\documentclass[11pt]{amsart}

\usepackage[english]{babel}
\usepackage{amsmath}
\usepackage{amsfonts}
\usepackage{amssymb}
\usepackage{amsthm}
\usepackage{color}
\definecolor{forestgreen}{rgb}{0.13, 0.55, 0.13}
\usepackage{caption}
\usepackage[english]{babel}
\usepackage{yfonts}
\usepackage[T1]{fontenc}
\usepackage[utf8x]{inputenc}
\usepackage{enumerate}
\usepackage{verbatim}
\usepackage{graphicx}
\usepackage{verbatim}
\usepackage{faktor}
\usepackage{xcolor}
\usepackage{xfrac}
\usepackage{tikz}
\usepackage[all]{xy}
\usepackage{hyperref}
\usepackage[backend = biber,citestyle=ieee-alphabetic, bibstyle=alphabetic ,maxnames=4,minnames=3,maxbibnames=99]{biblatex}
\usepackage{contour}
\usepackage{comment}
\usepackage{listings}
\usepackage{pythonhighlight}
\usepackage{mathabx}
\usepackage{graphicx}
\lstnewenvironment{Python}[1][]{\lstset{style=mypython, frame=none, breaklines=false #1}}{}
\title[Arc coordinates for maximal representations]{Arc coordinates for maximal representations}
\author{Marta Magnani}

\usepackage[margin=1.2in]{geometry}
\usepackage{bm}
\usepackage[percent]{overpic}
\usepackage{float}
\usepackage{enumerate}

\newcommand{\X}{\mathcal X}
\newcommand{\Y}{\mathcal Y}
\newcommand{\Yi}{\mathcal Y_{0,\infty}}
\renewcommand{\a}{\overline{\mathfrak{a}}^{+}}
\newcommand{\Sp}{{\rm Sp}}
\newcommand{\PSp}{{\rm PSp}}
\newcommand{\SL}{{\rm SL}}
\newcommand{\PSL}{{\rm PSL}}
\renewcommand{\O}{{\rm O}}

\newcommand{\PSO}{{\rm PSO}}
\newcommand{\PO}{{\rm PO}}

\newcommand{\GL}{{\rm GL}}
\newcommand{\Sym}{{\rm Sym}}
\newcommand{\Hom}{{\rm Hom}}
\newcommand{\Stab}{{\rm Stab}}
\newcommand{\Id}{{\rm Id}}
\newcommand{\R}{\mathbb R}
\renewcommand{\H}{\mathbb H}

\newcommand{\bpm}{\begin{pmatrix}}
\newcommand{\epm}{\end{pmatrix}}
\newcommand{\notpitchfork}{\mathrel{\ooalign{\hidewidth\reflectbox{$\pitchfork$}\hidewidth\cr$/$}}}

\theoremstyle{plain}
\newtheorem{thm}{Theorem}[section]
\newtheorem{lem}[thm]{Lemma}
\newtheorem{prop}[thm]{Proposition}
\newtheorem{cor}[thm]{Corollary}

\newtheorem*{introtheorem*}{Theorem}
\newtheorem*{introprop*}{Proposition}
\newtheorem*{introcor*}{Corollary}

\theoremstyle{definition}
\newtheorem*{introdefinition*}{Definition}
\newtheorem{definition}[thm]{Definition}
\newtheorem{remark}[thm]{Remark}
\newtheorem{example}[thm]{Example}
\newtheorem{notation}[thm]{Notation}

\begin{document}

\begin{abstract}
We generalize arc coordinates for maximal representations from a hyperbolic surface with boundary into $\PSp(4,\R)$, focusing on the case where the surface is a pair of pants. We introduce geometric parameters within the space of right-angled hexagons in the Siegel space $\X$. These parameters enable the visualization of a right-angled hexagon as a polygonal chain inside the hyperbolic plane $\H^{2}$. We explore the geometric properties of reflections in $\X$ and introduce the notion of maximal representation of the reflection group $W_{3}=\mathbb{Z}/2\mathbb{Z}*\mathbb{Z}/2\mathbb{Z}*\mathbb{Z}/2\mathbb{Z}$. We parametrize maximal representations from $W_{3}$ into $\PSp^{\pm}(4,\R)$, this induces a natural parametrization of a subset of maximal and Shilov hyperbolic representations into $\PSp(4,\R)$.  
\end{abstract}  

\maketitle

\section{Introduction}

\subsection{The space of maximal representations}

Given $\Sigma$ a closed oriented surface of negative Euler characteristic and fundamental group $\Gamma$, the Teichmüller space $\mathcal{T}(\Sigma)$ is the parameter space of marked hyperbolic structures on $\Sigma$. It is well known that, with the introduction of the holonomy map, one can associate to a point in $\mathcal{T}(\Sigma)$ a discrete and faithful representation $\rho: \Gamma \to \PSL(2,\R)$ so that the surface $\Sigma$ is realized by the quotient $\Sigma=\rho(\Gamma) \backslash \H^{2}$. This representation is well defined up to conjugation by an element in $\PSL(2,\R)$ so that the space $\mathcal{T}(\Sigma)$ can be identified with a connected component of the representation variety $\Hom(\Gamma,\PSL(2,\R))/\PSL(2,\R)$ which consists entirely of discrete and faithful representations \cite{goldman1980discontinuous}.

This phenomenon of the representation variety to admit components consisting only of injective homomorphisms with discrete image is still true if we substitute $\PSL(2,\R)$ with a semisimple real Lie group of higher rank $G$. In this sense higher rank Teichmüller space was developed as a generalization of classical Teichmüller space. More precisely given $G$ a semisimple real Lie group of higher rank, a \emph{higher Teichmüller space} is a subset of $\Hom(\Gamma,G)/G$ which is a union of connected components that consist entirely of discrete and faithful representations. To such a representation $\rho$ we can associate the quotient $\rho(\Gamma)\backslash \X$ where $\X$ is the symmetric space associated to $G$. The space $\X$ is a non-positively curved Riemannian symmetric manifold of higher rank, where \emph{rank} denotes the maximal dimension of an isometrically embedded flat inside $\X$. The quotient $\rho(\Gamma)\backslash \X$ is a locally symmetric space whose fundamental group is isomorphic to the fundamental group of $\Sigma$.

There are two well-known families of higher Teichmüller spaces: Hitchin components and maximal representations. In this paper we are interested in maximal representations. These are defined when $G$ is a Hermitian Lie group such as $\PSp(2n,\R)$ and are singled out by the maximal value of the Toledo number, which is a generalization of the Euler number. Burger Iozzi and Wienhard \cite{burger2005maximal}, \cite{burger2010surface}  studied the Toledo invariant for general Hermitian Lie groups and proved that maximal representations $\Hom^{\text{max}}(\Gamma,G)$ can be characterised as those representations admitting a monotone equivariant boundary map $\xi:S^{1}\to \widecheck{S}$, where $\widecheck{S}$ denotes the Shilov boundary and coincides with the set of real Lagrangians $\mathcal{L}(\R^{2n})$ when $G$ is $\PSp(2n,\R)$.

There is a related theory for surfaces with punctures or boundary components. The first thing to notice is that when $\partial \Sigma \neq \varnothing$ then $\Gamma$ is a free group and the whole representation variety is connected. Denote $\Sigma=\Sigma_{g,m}$ a surface of genus $g$ and $m$ boundary components with fundamental group $\Gamma_{g,m}$. Let
$$\Gamma_{g,m}= \langle a_{1},b_{1},...a_{g},b_{g},c_{1},...c_{m}| \ \prod_{i=1}^{g}[a_{i},b_{i}]\prod_{j=1}^{m}c_{j}=1 \rangle$$
be a presentation where the elements $c_{i}$ represent loops which are freely homotopic to the corresponding boundary components of $\partial \Sigma$ with positive orientation. A boundary condition might be imposed by considering 

\begin{equation}\label{boundcond}
\begin{aligned}
 \Hom^{\widecheck{S}}(\Gamma_{g,m},G) = \{ \rho \in \Hom(\Gamma_{g,m},G) \, | \, \rho(c_{i}) \text{ has at least} \\
 \text{one fixed point in } \widecheck{S}, \quad 1 \leq i \leq m \}
\end{aligned}
\end{equation}

\noindent where $\widecheck{S}$ denotes the Shilov boundary. In this case  $\Hom^{\text{max}}(\Gamma_{g,m},G) \subset \Hom^{\widecheck{S}}(\Gamma_{g,m},G)$ and in particular $\Hom^{\text{max}}(\Gamma_{g,m},G)$ is a union of connected components of the set $\Hom^{\widecheck{S}}(\Gamma_{g,m},G)$ \cite[Corollary 14]{burger2010surface}. In this paper we are interested in maximal representations inside (\ref{boundcond}) that satisfy a further condition: we will fix a union of conjugacy classes by imposing in (\ref{boundcond}) that every $\rho(c_{i})$ fixes exactly two points in $\widecheck{S}$ on which it acts expandingly and contractingly respectively. This is equivalent for the representation to be Anosov in the sense of \cite{guichardwienhard2012}.  We denote this space $\Hom^{\text{max,Shilov}}(\Gamma_{g,m},G)$.

\subsection{The results}

Coordinates on the space of maximal representations often arise as a generalization of well known coordinates on the classical Teichmüller space $\mathcal{T}(\Sigma)$. Analogues of Fenchel–Nielsen coordinates were developed by Strubel \cite{strubel2015fenchel}, whereas analogues of shear coordinates were developed by Alessandrini Guichard Rogozinnikov and Wienhard \cite{alessandrini2019noncommutative}. 

We want to generalize arc coordinates to the space of maximal representations. In classical Teichmüller theory arc coordinates were introduced by Harer \cite{Har86} and developed by Penner \cite{penner1987decorated} to decompose decorated Teichmüller space of punctured surface. This decomposition was generalized by \cite{ushijima1999canonical} \cite{penner2002} for surfaces with boundary. Similar coordinates were used in \cite{luo2007teichmuller} \cite{guo2009parameterizations} and in \cite{Mondello09}. 

We will consider the case where $\Sigma=\Sigma_{g,m}$ is a compact orientable smooth surface of genus $g$ and $m$ boundary components. We denote $\Gamma_{g,m}$ the fundamental group $\pi_{1}(\Sigma_{g,m})$, which is isomorphic to the free group $\mathbb{F}_{2g+m-1}$. An element $\gamma \in \Gamma_{g,m}$ is called \emph{peripheral} if it is represented by a loop that is freely homotopic into a boundary component of $\Sigma_{g,m}$. We can equip $\Sigma_{g,m}$ with a complete hyperbolic structure of finite volume with geodesic boundary. The universal covering $\widetilde{\Sigma}_{g,m}$ of $\Sigma_{g,m}$ is a closed subset of the hyperbolic plane $\H^{2}$ where boundary curves are geodesics.

\begin{figure}[!h]
   \centering
   \captionsetup{justification=centering,margin=2cm}
   \setlength{\unitlength}{0.1\textwidth}
   \begin{picture}(6,2.8)
     \put(0,0){\includegraphics[width=8cm,height=4.5cm]{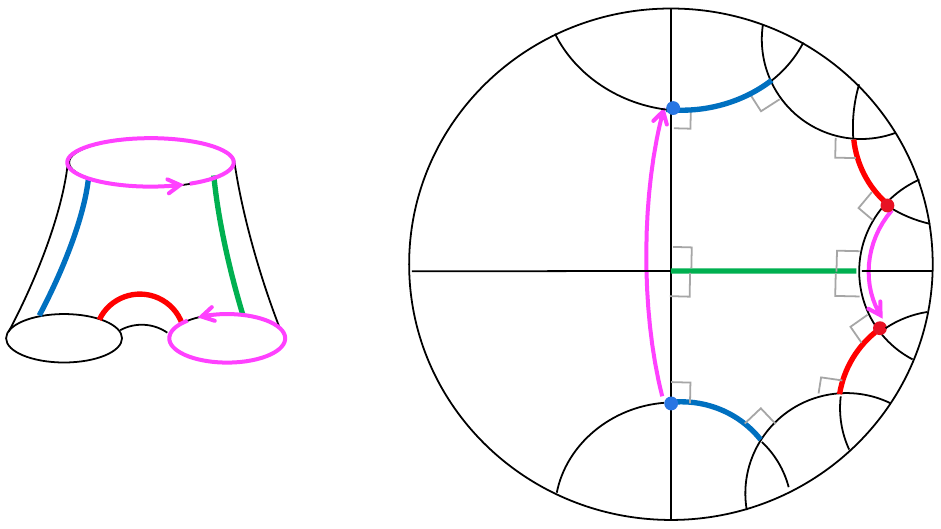}}
    \put(0.7,1.7){\textcolor{magenta}{$\alpha$}}
    \put(1,1.2){\textcolor{magenta}{$\beta$}}
    \put(3.1,1.1){\textcolor{magenta}{$\rho(\alpha$)}}
    \put(4.25,1.5){\textcolor{magenta}{$\rho(\beta$)}} 
    \put(2.5,2){\textcolor{gray}{$\H^{2}$}}
 \end{picture}   
\caption{The maximal representation $\rho: \Gamma_{0,3}\to \PSL(2,\R)$}
   \label{intropic1}
\end{figure}
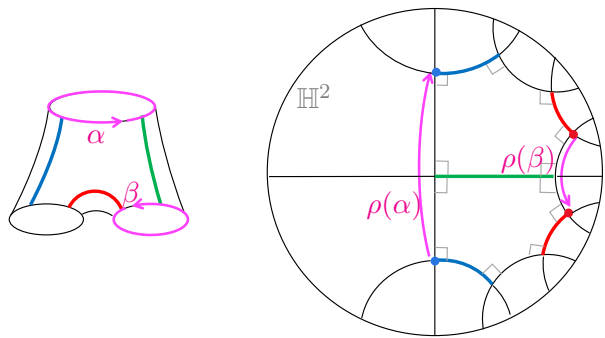

Arc coordinates are obtained by decomposing the surface in hexagons through the choice of a maximal collection $\{a_{1},...,a_{k}\}$ of pairwise disjoint arcs with starting and ending point on a boundary component which are essential and pairwise non-homotopic. For every hexagon in this decomposition there are exactly three alternating edges belonging to $\partial \Sigma_{g,m}$. We denote by $E$ the set of all edges and by $E_{bdry}$ the set of edges lying on a boundary component. For a fixed hyperbolic structure we can always realize the hexagon decomposition of $\Sigma_{g,m}$ in a way such that every edge is a geodesic and every arc is the unique geodesic which is orthogonal to the boundary at both endpoints. For each choice of $\{a_{1},...,a_{k}\}$ we get a parametrization of the Teichmüller space $\mathcal{T}(\Sigma_{g,m})$: once we fix the lengths $l(a_{1}),...,l(a_{k})$ there is a unique hyperbolic metric that makes $\Sigma_{g,m} \setminus \bigcup_{i} a_{i}$ a union of hyperbolic right-angled hexagons where each hexagon has exactly three alternating edges $a_{i_{1}},a_{i_{2}},a_{i_{3}}$ in $E \backslash E_{bdry}$ of length $l(a_{i_{1}}),l(a_{i_{2}}),l(a_{i_{3}})$ respectively, where $i_{1},i_{2},i_{3} \in \{1,...,k\}$. This is due to the
well known fact that given three real numbers $b,c,d>0$ there exists (up to
isometries) a unique right-angled hexagon in $\H^{2}$ with alternating sides of lengths
$b,c$ and $d$ (see for example \cite[Lemma 6.2.2]{martelli2016introduction}). A point in the Teichmüller space $\mathcal{T}(\Sigma_{g,m})$ is identified with a maximal representation $\rho: \Gamma_{g,m} \to \PSL(2,\R)$. Since we are considering surfaces with geodesic boundary, the image $\rho(\gamma)$ of every element $\gamma \in \Gamma_{g,m}$ is a hyperbolic isometry fixing exactly two points in $\partial \H^{2}$. The above discussion asserts that once we fix the lengths $l(a_{1}),...,l(a_{k})$ we can explicitly write (up to conjugation) the maximal representation $\rho$ such that $\Sigma_{g,m}=\rho(\Gamma_{g,m}) \backslash \H^{2}$. An example for the surface $\Sigma_{0,3}$ (pair of pants) is given in Figure \ref{intropic1}, where the fundamental group $\Gamma_{0,3}$ is isomorphic to the free group $\mathbb{F}_{2}$ generated by $\alpha$ and $\beta$.

More generally, given a maximal representation $\rho: \Gamma_{g,m} \to \PSp(2n,\R)$, the image $\rho(\gamma)$ of every non-peripheral element $\gamma \in \Gamma_{g,m}$ is Shilov hyperbolic (see \cite{strubel2015fenchel}). Equivalently, $\rho(\gamma)$ fixes two transverse Lagrangians $l_{\gamma}^{+}$ and $l_{\gamma}^{-}$ on which it acts expandingly and contractingly respectively. These Lagrangians are the images $\xi(\gamma^{+})$ and $\xi(\gamma^{-})$ where $\xi: S^{1} \to \mathcal{L}(\mathbb{R}^{2n})$ is the equivariant boundary map and $l^{\pm}_{\gamma}=\xi(\gamma^{\pm})$.  
We consider the set of maximal representations where the property of being Shilov hyperbolic is true also for peripheral elements.

\begin{introdefinition*} (see \ref{defmaxshil})
 A maximal representation  $\rho: \pi_{1}(\Sigma) \to \PSp(2n,\R)$ is called \emph{Shilov hyperbolic} if $\rho(\gamma)$ is Shilov hyperbolic for every $\gamma \in \pi_{1}(\Sigma)$. The set of maximal representations which are Shilov hyperbolic is denoted by $\text{Hom}^{\text{max,Shilov}} (\pi_{1}(\Sigma),\PSp(2n,\R))$. We define $\chi^{\text{max,Shilov}} (\pi_{1}(\Sigma),\PSp(2n,\R))$ as the quotient 
    $$\chi^{\text{max,Shilov}} (\pi_{1}(\Sigma),\PSp(2n,\R)):=\text{Hom}^{\text{max,Shilov}} (\pi_{1}(\Sigma),\PSp(2n,\R))/_{\PSp(2n,\R)}$$ 
    where $\PSp(2n,\R)$ is acting by conjugation: $\rho \sim \rho'$ if there exists $g \in \PSp(2n,\R)$ such that $\rho(\gamma)=g\rho'(\gamma)g^{-1}$ for all $\gamma \in \pi_{1}(\Sigma)$.
\end{introdefinition*}

 We focus on the case where $\Sigma=\Sigma_{0,3}$ and $n=2$: we consider the Siegel space $\X$ the symmetric space associated to $\Sp(4,\R)$ and we fix the Weyl chamber $\a$
 
\begin{equation*}
    \a = \{ (x_{1},x_{2}) \in \R^{2} | \ x_{1} \geq x_{2} \geq 0 \}
\end{equation*}

\noindent We denote by $\mathfrak{a}$ the set of regular vectors inside $\a$  
$$\mathfrak{a}=\{(x_{1},x_{2}) \in \R^{2}| \ x_{1}> x_{2}> 0 \}$$
and we further denote by $\mathfrak{d}$ the set 
$$\mathfrak{d}=\{(x_{1},x_{2}) \in \R^{2}| \ x_{1}=x_{2} \}$$
The first step is to introduce a parameter space for right-angled hexagons in $\X$. The subspaces of the Siegel space that play the role of geodesics in $\H^{2}$ are called $\R$-\emph{tubes}. In Section \ref{paramhex} we give the definition of a right-angled hexagon $H$ in $\X$, which is determined by a cyclic sequence of $\R$-tubes 
$$H=[\Y_{1},\Y_{2},\Y_{3},\Y_{4},\Y_{5},\Y_{6}]$$
where any two consecutive tubes are orthogonal. We further define the set of \emph{ordered right-angled hexagons} $\mathcal{H}$, this is given by the data $(H,\Y_{1})$ of a right-angled hexagon together with the choice of a tube $\Y_{1}$. We distinguish between generic (Definition \ref{defgenhexag}) and non-generic hexagons (Section \ref{secnongenhex}). A generic hexagon is parametrized by length parameters $\underline{b},\underline{c},\underline{d}$ inside $\mathfrak{a}$ and angle parameters $\alpha_{1},\alpha_{2}$ lying in $[0,2\pi)$ (Proposition \ref{parforhex}). In the non-generic case some length parameters lie in $\mathfrak{d}$ and some angle parameters vanish (Propositions \ref{parnongentype1}, \ref{parnongentype2} and \ref{hexparinH2}). This leads to a geometric visualization of a right-angled hexagon inside $\X$ in terms of a polygonal chain (Section \ref{secpolchain}). A parameter space which encloses both generic and non-generic hexagons is given by

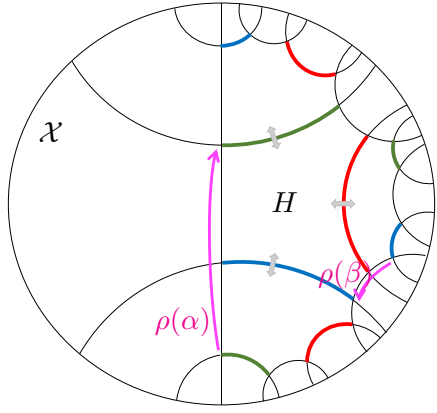
\begin{figure}[!h]
   \centering
\captionsetup{justification=centering,margin=2cm}
   \setlength{\unitlength}{0.1\textwidth}
   \begin{picture}(4.3,3.6)
\put(0,0){\includegraphics[width=5.8cm,height=5.5cm]{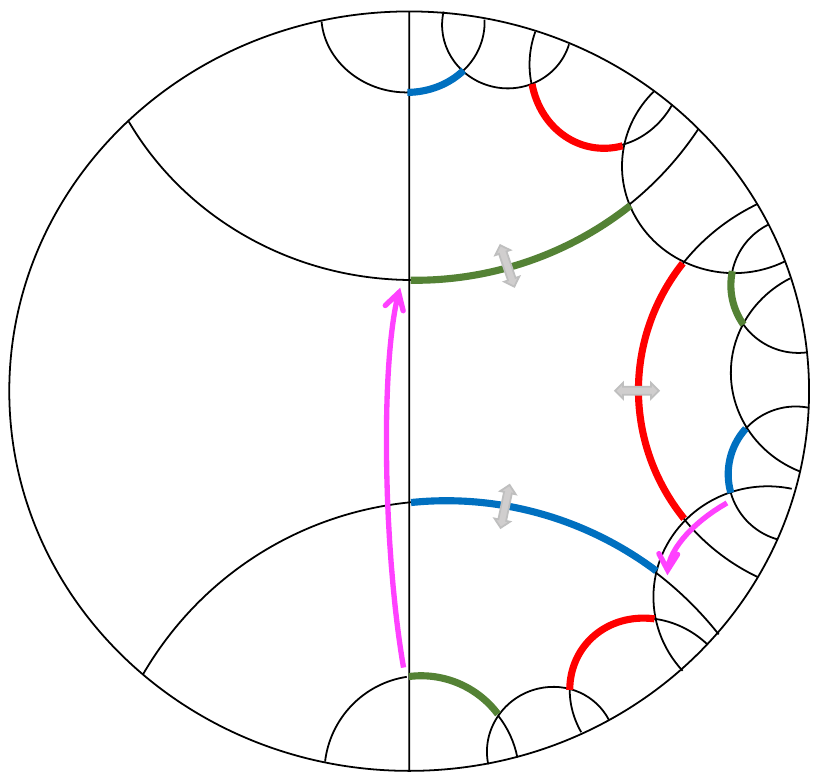}}
\put(0.3,2.3){$\X$}
\put(2.3,1.7){$H$}
\put(1.3,0.7){\textcolor{magenta}{$\rho(\alpha)$}}
\put(2.7,1.1){\textcolor{magenta}{$\rho(\beta)$}}  
   \end{picture}
\caption{The maximal representation $\rho:\Gamma_{0,3} \to \PSp(4,\R)$}
   \label{intropic3}
\end{figure}

\begin{introtheorem*} (see \ref{cpctarccoordthm})
The space $\mathcal{H}$ is parametrized up to isometry by $\mathcal{A}=\overline{\mathfrak{a}}^{3}\times [0,2\pi)^{2}/_{\sim}$.
\end{introtheorem*}

where $\overline{\mathfrak{a}}=\mathfrak{a}\cup\mathfrak{d}$. The equivalence relation collapses one of the angles to a point in the case where the hexagon degenerates to a non-generic one. These parameters were firstly introduced  with the aim of generalizing hexagon parameters in $\H^{2}$. This approach turned out to be very tricky and this is explained in detail in Section \ref{disconpar}. Geometric parameters for a maximal representation $\rho: \Gamma_{0,3} \to \PSp(4,\R)$ should be thought as the data of lengths and angles which uniquely determine two adjacent hexagons in $\X$ both having three alternating sides of length $\underline{b},\underline{c}$ and $\underline{d}$ respectively. The maximal representation is then determined by the image of the generators of the fundamental group generalizing the geometric construction of Figure \ref{intropic1}. If we want to extend our hexagon-parameters for two adjacent hexagons by only changing the angle parameters we can not guarantee that the constructed hexagons have the same alternating side-lengths. We will therefore construct two adjacent hexagons starting with one hexagon $H$ and obtaining the others by reflecting $H$ across a side (Figure \ref{intropic3}). 

\noindent This leads to the parametrization of a subset $$\chi^{\mathcal{S}}\subset \chi^{\text{max,Shilov}} (\Gamma_{0,3},\PSp(4,\R))$$
The idea is to see the fundamental group $\Gamma_{0,3}$ as a subgroup of the reflection group 
$$W_{3}=\mathbb{Z}/2\mathbb{Z}*\mathbb{Z}/2\mathbb{Z}*\mathbb{Z}/2\mathbb{Z}=\langle s_{1},s_{2},s_{3}| \ s_{1}^{2}=s_{2}^{2}=s_{3}^{2}=1 \rangle$$
through the following homomorphism $\phi$
$$\begin{aligned}
\phi: \Gamma_{0,3}& \to W_{3} \\
 \alpha& \mapsto s_{1}s_{2}\\
  \beta& \mapsto s_{2}s_{3}
\end{aligned}$$ 

\noindent We define the notion of maximal representation of the group $W_{3}$ into $\PSp^{\pm}(2n,\R)$, where $\PSp^{\pm}(2n,\R)$ denotes the union of symplectic and antisymplectic matrices (Definition \ref{defantisympl}).

\begin{figure}[!h]
   \centering
   \captionsetup{justification=centering,margin=2cm}
   \setlength{\unitlength}{0.1\textwidth}
   \begin{picture}(4.5,3)
     \put(0.8,0){\includegraphics[width=4.5cm,height=4.5cm]{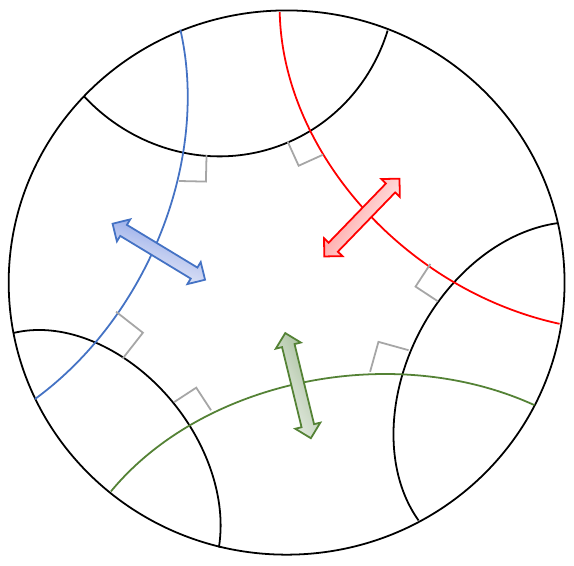}}
          \put(1.8,-0.16){$X_{2}$}
          \put(1.1,0.2){\textcolor{forestgreen}{$Q_{1}$}}
          \put(0.7,0.6){\textcolor{blue}{$P_{2}$}}
          \put(0.5,1.1){$X_{1}$}
          \put(1,2.5){$Z_{2}$}
          \put(1.5,2.8){\textcolor{blue}{$P_{1}$}}
          \put(2,2.95){\textcolor{red}{$R_{2}$}}
          \put(2.6,2.8){$Z_{1}$}
          \put(3.7,1.7){$Y_{2}$}
          \put(3.7,1.2){\textcolor{red}{$R_{1}$}}
          \put(3.6,0.9){\textcolor{forestgreen}{$Q_{2}$}}
          \put(3,0.1){$Y_{1}$}
          \put(0.9,1.8){\textcolor{blue}{$\rho(s_{1})$}}
          \put(2.2,0.4){\textcolor{forestgreen}{$\rho(s_{2})$}}
          \put(2.8,2){\textcolor{red}{$\rho(s_{3})$}}
 \end{picture} 
\caption{The reflections $\rho(s_{1}),\rho(s_{2}),\rho(s_{3})$ for $\rho: W_{3} \to \PSp^{\pm}(2n,\R)$ maximal}
   \label{intropic4}
\end{figure}
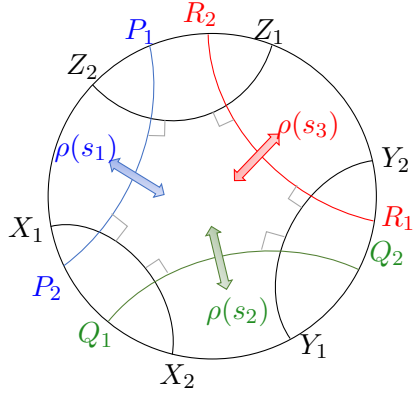

\begin{introdefinition*} (see \ref{defmaxW3repr})
 A representation $ \rho: W_{3} \to \PSp^{\pm}(2n,\R)$ is \emph{maximal} if there exists a maximal 6-tuple of Lagrangians $(P_{1},P_{2},Q_{1},Q_{2},R_{1},R_{2})$ such that $\rho(s_{1}),\rho(s_{2}),\rho(s_{3})$ are reflections of $\X$ fixing $(P_{1},P_{2}),(Q_{1},Q_{2}),(R_{1},R_{2})$ respectively and such that
\begin{equation*}
    \begin{cases}
        \rho(s_{1})(X_{1})=X_{2} \text{ and } \rho(s_{1})(Z_{1})=Z_{2}\\
\rho(s_{2})(X_{1})=X_{2} \text{ and } \rho(s_{2})(Y_{1})=Y_{2}\\
        \rho(s_{3})(Y_{1})=Y_{2} \text{ and } \rho(s_{3})(Z_{1})=Z_{2}
    \end{cases}
\end{equation*}
 where $X_{1},X_{2},Y_{1},Y_{2},Z_{1},Z_{2}$ are uniquely determined by (see Figure \ref{intropic4})
$$\Y_{P_{1},P_{2}} \perp \Y_{X_{1},X_{2}} \perp \Y_{Q_{1},Q_{2}} \perp \Y_{Y_{1},Y_{2}} \perp \Y_{R_{1},R_{2}} \perp \Y_{Z_{1},Z_{2}} \perp \Y_{P_{1},P_{2}}$$
\end{introdefinition*}

\begin{introtheorem*} (see \ref{paramforchimax})
The set $\chi^{\text{max}}(W_{3},\PSp^{\pm}(4,\R))=\Hom^{\text{max}}(W_{3},\PSp^{\pm}(4,\R))/\PSp(4,\R)$  is parametrized by $\mathcal{S}\subset \mathcal{A}\times\mathcal{K}^{3}$ where $\mathcal{A}$ is the parameter space of a right-angled hexagon and $\mathcal{K}$ is the set
$$ \mathcal{K}=\Big\{\bpm -K&0\\0&K\epm, \  K \in \PO(2), \ K^{2}=\Id \Big\}$$ 
\end{introtheorem*}

In Section \ref{sectiongeomintK} we give a geometrical interpretation of the set $\mathcal{K}$ in terms of the polygonal chain associated to a right-angled hexagon. We prove that the restriction to $\Gamma_{0,3}$ of such a maximal representation is maximal and Shilov hyperbolic.

\begin{introprop*} (see \ref{maximalityproof})
Fix $\widetilde{\rho} \in \Hom^{\text{max}} (W_{3},\PSp^{\pm}(4,\R))$. 
Then the representation $\rho:=\widetilde{\rho}|_{Im(\phi)}$ is inside $\Hom^{\text{max,Shilov}} (\Gamma_{0,3},\PSp(4,\R))$.
\end{introprop*}

This allows us to define $\chi^{\mathcal{S}}$ as the image $\chi^{\mathcal{S}}:=Im(f)$ where $f$ is the map
$$\begin{aligned}
f: \chi^{\text{max}}(W_{3},\PSp^{\pm}(4,\R))& \to \chi^{\text{max,Shilov}}(\Gamma_{0,3},\PSp(4,\R)) \\
 \bigl[\widetilde{\rho}\bigr]& \mapsto \bigl[\widetilde{\rho}|_{Im(\phi)}\bigr]
\end{aligned}
$$

\noindent The parametrization of $\chi^{\mathcal{S}}$ is obtained by imposing an equivalent relation on $\mathcal{S}$ which identifies the points that have same image under $f$.
\begin{introtheorem*} (see \ref{thesis})
 The set $\chi^{\mathcal{S}}$ is parametrized by $\mathcal{S}/_{\sim}$   
\end{introtheorem*}

In Corollary \ref{mapfnotinjnorsurj} we show that, contrary to the hyperbolic case (Proposition \ref{teichprop}), the map $f$ is not injective nor surjective. The set $\chi^{\mathcal{S}}$ is an 8-dimensional subspace inside the 10-dimensional space $\chi^{\text{max,Shilov}}(\Gamma_{0,3},\PSp(4,\R))$. The space $\chi^{\text{max,Shilov}}(\Gamma_{0,3},\PSp(4,\R))$ has 4 connected components (see \cite{alessandrini2019noncommutative}), and the paramater space for $\chi^{\mathcal{S}}$ gives 8 connected components. In Proposition \ref{propconncomp} we show that $\chi^{\mathcal{S}}$ hits all connected components of  $\chi^{\text{max,Shilov}}(\Gamma_{0,3},\PSp(4,\R))$. 

A motivation for this work is to study compactification of character varieties where similar arguments can be carried out with non-Archimedean Siegel spaces as in \cite{burger2017maximal}. We expect applications of this work in the study of the real spectrum compactification of maximal character varieties (see \cite{BIPP21}, \cite{burger2023real}) Of particular interest are rank two groups where \cite{burger2021weyl} and \cite{OuyTamb23} \cite{OuyangTambGoth} suggest a link with flat structures with angle multiple of $\frac{\pi}{2}$. Developing arc coordinates for those would be interesting.

\subsection{Organization of the paper}
Section \ref{chaptersiegel} recalls some fundamental properties of the Siegel space $\X$- the symmetric space associated to $\Sp(2n,\R)$. We fix notation and focus on the case of $\Sp(4,\R)$. In Section \ref{parforquin} we define the set of generic quintuples and give a parameter space for them (Proposition \ref{partrip}). These parameters will be very useful for the parametrization of right-angled hexagons. Section \ref{chapterhexagons} is dedicated to the study of hexagons. We define the set of ordered right-angled hexagons (Definition \ref{defordhex}) and distinguish between generic (Definition \ref{defgenhexag}) and non-generic hexagons (Section \ref{secnongenhex}). We introduce a parameter space for both cases (Proposition \ref{parforhex} for the generic case and Propositions \ref{parnongentype1}, \ref{parnongentype2} and \ref{hexparinH2} for the non-generic case). A parameter space which encloses both generic and non-generic hexagons is given in Theorem \ref{cpctarccoordthm}. These parameters will be called \emph{arc coordinates}. In Section \ref{disconpar} we show how arc coordinates arise from the idea of generalizing coordinates of a hexagon in $\H^{2}$ and explain the problems encountered in this approach. In Section \ref{refinsiegel} we discuss properties of reflections in the Siegel space $\X$. We define the \emph{reflection set associated to the side of a hexagon} (Definition \ref{reflsetside}) and give a geometric interpretation of it (Section \ref{sectiongeomintK}). In Section \ref{chapterparformax} we define the set $\chi^{\text{max,Shilov}} (\pi_{1}(\Sigma),\PSp(2n,\R))$ and we further define the notion of a maximal representation from the group $W_{3}=\mathbb{Z}/2\mathbb{Z}*\mathbb{Z}/2\mathbb{Z}*\mathbb{Z}/2\mathbb{Z}$ into $\PSp^{\pm}(2n,\R)$ (Definition \ref{defmaxW3repr}), providing a parameter space for the  $\PSp^{\pm}(4,\R)$-case (Theorem \ref{paramforchimax}). We further define the set $\chi^{\mathcal{S}}\subset \chi^{\text{max,Shilov}} (\Gamma_{0,3},\PSp(4,\R))$ (Definition \ref{defchiS}) to which we provide a parameter space in Theorem \ref{thesis}.

\subsection*{Acknowledgements} I would like to thank first and foremost my doctoral advisors Beatrice Pozzetti and Anna Wienhard for their constant support throughout my years as a PhD student. I would like to thank Eugen Rogozinnikov for all the helpful comments and suggestions in many parts of this article. Thanks also to Valentina Disarlo, Daniele Alessandrini, Gabriele Viaggi, Max Riestenberg and Mitul Islam for interesting discussions related to this work. I acknowledge funding by the DFG, through the RTG 2229 “Asymptotic Invariants and Limits of Groups and Spaces”- 281869850- and the Emmy Noether project 427903332 of B. Pozzetti.

\section{Preliminaries} \label{chaptersiegel}

\subsection{The Siegel space} \label{defsiegel}
The \emph{Siegel space} $\X$ is the symmetric space associated to the symplectic group $\Sp(2n,\R)$. Recall that the symplectic group $$\Sp(2n,\R)= \{ M \in \SL(2n,\R) \ | \ M^{T}J_{n}M=J_{n}   \}$$ is the subgroup of $\SL(2n,\R)$ preserving the symplectic form $\omega(\cdot, \cdot)$ represented, with respect to the standard basis, by the matrix \begin{equation*} 
J_{n}=
\bpm
0 & \Id_{n}\\
-\Id_{n} & 0 \epm
\end{equation*}
The group $\Sp(2n,\R)$ can also be described as the group of block matrices:

\begin{equation*}
\Sp(2n,\R)=\Big\{  \bpm
A & B\\
C & D 
\epm | \ A^{T}C,B^{T}D \text{ symmetric, and  } A^{T}D-C^{T}B=\Id_{n}\Big\} 
\end{equation*}

\noindent When $n=1$ the group $\Sp(2,\R)$ coincides with $\SL(2,\R)$. There are two models commonly used for the Siegel space: the \emph{upper-half space} and the \emph{Borel embedding} model.

The \emph{upper-half space model} is a generalization of the upper-half space model of the hyperbolic plane and is given by a specific set of symmetric matrices: $$\X= \{ X+iY,  \ X \in \Sym(n,\R), \ Y \in \Sym^{+}(n,\R) \}$$ where $\Sym(n,\R)$ denotes the set of $n$-dimensional symmetric matrices with coefficients in $\R$ and $\Sym^{+}(n,\R)$ is the subset of $\Sym(n,\R)$ given by positive definite matrices. The group $\Sp(2n,\R)$ acts on $\X$ by fractional linear transformations:
\begin{equation*}
\centering   
 \bpm
A & B\\
C & D 
\epm \cdot Z = (AZ+B)(CZ+D)^{-1}
\end{equation*}
The \emph{Borel embedding model} is given by $$\mathbb{X} = \{ l \in \mathcal{L} (\mathbb{C}^{2n}) | \ \  i\omega(\sigma(\cdot),\cdot)_{\mathbb{C}}\vert_{l\times l} \  \text{is positive definite} \}$$ where $\mathcal{L} (\mathbb{C}^{2n})$ is the set of Lagrangians and $\sigma : \mathbb{C}^{2n} \to \mathbb{C}^{2n} $ denotes complex conjugation. We consider the  affine chart $\iota: \Sym(n,\mathbb{C}) \to \mathcal{L} (\mathbb{C}^{2n})$ that associates to a symmetric matrix $Z$ the linear subspace of $\mathbb{C}^{2n}$ spanned by the columns of the matrix $\left(
\begin{array}{c}
Z\\
\Id_{n}\\
\end{array}
\right)_{2n \times n}$. This induces an $\Sp(2n,\R)$-equivariant identification $\iota: \X \mapsto \mathbb{X}$. The restriction of the affine chart $\iota$ to the subspace $\Sym(n,\R)$ provides a parametrization of the set of real Lagrangians that are transverse as linear subspaces to $ \langle e_{1},...,e_{n} \rangle $, which will be denoted by $l_{\infty}$ or just $\infty$.

\subsection{Boundary and Lagrangians} \label{Maslov}
The set of real Lagrangians $\mathcal{L}(\R^{2n})$ naturally arises as the unique closed $\Sp(2n,\R)$-orbit in the boundary of $\X$ in its Borel embedding and for
this reason $\mathcal{L}(\R^{2n})$ is the Shilov boundary of the bounded domain realization of $\X$ (see \cite{wienhard2004bounded}). Denote by $\mathcal{L}(\R^{2n})^{(k)}$ the set of $k$-tuples of pairwise transverse Lagrangians. It is easy to prove that the group $\Sp(2n,\R)$ acts transitively on $\mathcal{L}(\R^{2n})^{(2)}$. Moreover, it has $(n+1)$ orbits in $\mathcal{L}(\R^{2n})^{(3)}$, indexed by the Maslov index (see for example \cite{souriau2005construction}). The Maslov index is cyclically invariant, is invariant under the action of $\Sp(2n,\R)$ on $\mathcal{L}(\R^{2n})^{(3)}$ and the group $\Sp(2n,\R)$ acts transitively on the set of triples of pairwise transverse Lagrangians with the same Maslov index \cite{lion1980weil}. The value of the Maslov index is maximal on the orbit of $( \langle e_{1},...,e_{n} \rangle, \langle e_{n+1},...,e_{2n} \rangle, \langle e_{1}+e_{n+1},...,e_{n}+e_{2n} \rangle) =(l_{\infty},0,\Id)$.

\begin{definition} A triple  of pairwise transverse Lagrangians is called \emph{maximal} if it is in the $\Sp(2n,\R)$-orbit of $(l_{\infty},0,\Id)$. An m-tuple $(l_{1},...,l_{m})$ is \emph{maximal} if for every $i<j<k$ the triple $(l_{i},l_{j},l_{k})$ is maximal.
\end{definition}

Maximal triples are a generalization of positively oriented triples in the circle $S^{1}$ and they play a central role in the study of maximal representations. It is useful to have a concrete criterion to check when triples of Lagrangian are maximal. 

\begin{lem} (\cite{burger2017maximal}) \label{max}
The following hold:\\
\emph{(1)} Any cyclic permutation of a maximal triple is maximal;\\
\emph{(2)} The triple $(l_{\infty},X,Y)$ is maximal if and only if $Y-X$ is positive definite;\\
\emph{(3)} If $Z-X$ is positive definite, the triple $(X,Y,Z)$ is maximal if and only if $Z-Y$ and $Y-X$ are positive definite.
\end{lem}

From a given maximal $m$-tuple we can obtain a maximal $(m+k)$-tuple by adding a maximal $k$-tuple between two consecutive Lagrangians. 

\begin{lem} \label{maxlemma2}
Let $(P_{1},...,P_{m})$ be a maximal $m$-tuple. For $i \in \{ 1,...,m-1 \}$ and $k \geq 1$ let $(P_{i},Q_{1},...,Q_{k},P_{i+1})$ be maximal. Then the $(m+k)$-tuple $(P_{1},...,P_{i},Q_{1},...,Q_{k},P_{i+1},...,P_{m})$ is maximal.    
\end{lem}

\proof Up to isometry we reduce to the case where $P_{1}=0, P_{m}=l_{\infty}$, that is we consider the $(m+k)$-tuple $(0,P_{2},...P_{i},Q_{1},...,Q_{k},P_{i+1},...,l_{\infty})$ where $(0,P_{2},...,P_{m-1},l_{\infty})$ maximal and $(P_{i},Q_{1},...,Q_{k},P_{i+1})$ maximal. Using Lemma \ref{max} result follows immediately.
\endproof

\subsection{$\Sp(2n,\R)$-invariant distance} \label{invdist}

We introduce a $\Sp(2n,\R)$-invariant distance on the symmetric space $\X$. Fix a point $p$ in a maximal flat $F$ and a Weyl chamber $\a \subset T_{p}F$. This is a fundamental domain for the action of $\Sp(2n,\R)$ on the tangent bundle $T\X$. In our case we have 
\begin{equation*}
    \a = \{ (x_{1},...,x_{n}) \in \R^{n} | \ x_{1} \geq...\geq x_{n} \geq 0 \}
\end{equation*}

\noindent A vector in the Weyl chamber is \emph{regular} if all the inequalities are strict, which is equivalent to being contained in a unique flat. In order to define a vectorial $\Sp(2n,\R)$-invariant distance in $\X$ we need to recall from \cite{burger2017maximal} the definition of an endomorphism-valued cross-ratio. If $l_{1},l_{2} \in \mathcal{L}(\R^{2n})$ are transverse (denoted by $ l_{1} \pitchfork l_{2}$), we denote by $p_{l_{1}}^{||l_{2}}:\R^{2n} \to l_{1}$ the projection to $ l_{1}$ parallel to $ l_{2}$.

\begin{definition} \label{defcrossratio} For Lagrangians $l_{1},...,l_{4} \in \mathcal{L}(\mathbb{C}^{2n})$ such that
$ l_{1} \pitchfork l_{2}$ and $ l_{3} \pitchfork l_{4}$ the \emph{cross-ratio} $R(l_{1},l_{2},l_{3},l_{4})$ is given in the Borel embedding model by the endomorphism of $l_{1}$
\begin{equation*}
   R(l_{1},l_{2},l_{3},l_{4})= p_{l_{1}}^{||l_{2}} \circ p_{l_{4}}^{||l_{3}}|_{l_{1}}
\end{equation*}
 \end{definition}

In the upper half space model the explicit expression for the cross-ratio is given by (\cite[Lemma 4.2]{burger2017maximal}):
\begin{equation*} 
 R(X_{1},X_{2},X_{3},X_{4})=(X_{1}-X_{2})^{-1}(X_{4}-X_{2})(X_{4}-X_{3})^{-1}(X_{1}-X_{3})  
\end{equation*}

\noindent where $R$ is expressed with respect to the basis of $X_{1}$ given by the columns of the matrix $\left(
\begin{array}{c}
X_{1}\\
\Id_{n}\\
\end{array}
\right)$. 

\begin{lem} [\cite{burger2017maximal}] \label{crinfty}
    Assume $0,Z,X,l_{\infty}$ are pairwise transverse. Then
    $$R(0,Z,X,l_{\infty})=Z^{-1}X$$
\end{lem}

We can now define the vectorial distance $d^{\a}$. The fact that the cross-ratio can be used to describe the projection of a pair of points in $\X$ onto the Weyl chamber was proved by Siegel in \cite{siegel1943symplectic}.

\begin{definition} \label{defweyldist}
The \emph{vectorial distance} $d^{\a}$ is the projection onto the Weyl chamber $\a$:
\[
\begin{aligned}
 \X^{2}&\to \a\\
(X,Z)& \mapsto (\log(\lambda_{1}),...,\log(\lambda_{n}))
\end{aligned}
\]

where $\lambda_{i}=\frac{1+\sqrt{r_{i}}}{1-\sqrt{r_{i}}}$ and $1>r_{1} \geq ... \geq r_{n}  \geq 0$ are the eigenvalues of $R(X,\overline{Z},Z,\overline{X})$

\end{definition}

\noindent A pair of points in $\X$ can be mapped to any other pair of points if and only if their vectorial distance $d^{\a}$ is the same. For interesting properties about the distance $d^{\a}$ see \cite{parreau2010distance},\cite{kapovich2017dynamics}.

\begin{lem} [\cite{fanoni2020basmajian}] \label{dist}
Let $A$ and $B$ be positive definite symmetric matrices such that the difference $B-A$ is positive definite. Let $\mu_{1} \geq ... \geq \mu_{n}$ be the eigenvalues of $A^{-1}B$. Then
\begin{equation*}
    d^{\a}(iA,iB)=(\log\mu_{1},...,\log\mu_{n})
\end{equation*}

\end{lem}

\subsection{Copies of $\H^{2}$ inside the Siegel space $\X$} \label{H2insideX} 

\begin{definition}
    Let $\X$ be the symmetric space associated to $\Sp(2n,\R)$. A \emph{maximal polydisk} in $\X$ is the image of a totally geodesic and holomorphic embedding of the Cartesian product of $n$ copies of $\H^{2}$ into $\X$.
\end{definition}

We will focus on the symmetric space $\X$ associated to $\Sp(4,\R)$. In this case an example of a maximal polydisk is the image of the following map $\psi$:
\[
\begin{aligned}
\psi: \H^{2}\times\H^{2}&\to \X\\
(z_{1},z_{2})& \mapsto \bpm
z_{1} & 0\\
0 & z_{2} 
\epm
\end{aligned}
\]

\noindent We will refer to this polydisk as the \emph{model polydisk} since every other polydisk is translate of our model polydisk by an element in $\Sp(4,\R)$ (see \cite{wolf1972fine}). Let $(M_{1},M_{2})$ be an element of $\SL(2,\R)\times\SL(2,\R)$. Then $(M_{1},M_{2})$ acts on the model polydisk as following:
$$
(M_{1},M_{2}) \cdot \bpm
z_{1} & 0\\
0 & z_{2} 
\epm = \bpm
M_{1}(z_{1}) & 0\\
0 & M_{2}(z_{2}) 
\epm
$$
where $M(z)$ is the action on a point $z \in \H^{2}$ by Möbius transformation. Let $\Delta$ be the diagonal embedding given by
\begin{equation*}
\begin{aligned}
\Delta: \SL(2,\R)\times\SL(2,\R)&\to \Sp(4,\R)\\
\Big( \bpm
a_{1} & b_{1}\\
c_{1} & d_{1} 
\epm,\bpm
a_{2} & b_{2}\\
c_{2} & d_{2} 
\epm \Big)& \mapsto \bpm
a_{1} & 0 & b_{1} &0\\
0 & a_{2} & 0 & b_{2}\\
c_{1} & 0 & d_{1} &0\\
0 & c_{2} & 0 & d_{2}
\epm
\end{aligned}
\end{equation*}

\noindent then $\psi(M_{1}(z_{1}),M_{2}(z_{2}))=\Delta(M_{1},M_{2})\big(\psi(z_{1},z_{2})\big)$. In particular the set $\psi\big((z,z)\big)$ is a copy of $\H^{2}$ inside $\X$ and will be called the \emph{diagonal disc}.

\subsection{$\mathbb{R}$-tubes} \label{tubes}

The subspaces of the Siegel space that play the role of geodesics in $\H^{2}$ are called $\R$-\emph{tubes}. Let $\{a,b \}$ be an unordered pair of transverse Lagrangians.

\begin{definition} \label{defRtube} (\textbf{$\R$-\emph{tube}}) The $\R$-\emph{tube} associated to $\{ a,b \}$ is the set $$\Y_{ a,b }= \{ l\in \X | \  R(a,l,\sigma(l),b) = -\Id \}$$ 
\end{definition}

It can be proven (see \cite{burger2017maximal}) that $\Y_{ a,b }$ is a totally geodesic subspace of $\X$ of
the same real rank as $\X$ and that it is the parallel set of the Riemannian singular geodesics whose endpoints in the visual boundary of $\X$ are the Lagrangians $a$ and $b$.
The group $\Sp(2n,\R)$ acts transitively on $\mathcal{L}(\R^{2n})^{(2)}$ and for every $g \in \Sp(2n,\R)$ it holds $g \cdot \Y_{ a,b } = \Y_{ ga,gb }$. Up to the symplectic group action we can therefore reduce to a model $\R$-tube, the one with endpoints $0$ and $l_{\infty}$. In the upper-half space model this will be called the \emph{standard tube} and consists of matrices of the form 
\begin{equation*}
\Yi= \{ iY| Y \in \Sym^{+}(n,\R)  \}
\end{equation*}

\noindent Intersection patterns of $\R$-tubes in the Siegel space reflect the intersection patterns of geodesics in the hyperbolic plane. 

\begin{prop} [\cite{fanoni2020basmajian}]
 \label{prop}
If $(l_{1},l_{2},l_{3},l_{4})$ is maximal, the intersection $\Y_{l_{1},l_{3}} \cap \Y_{l_{2},l_{4}}$ consists of a single point and $\Y_{l_{1},l_{2}} \cap \Y_{l_{3},l_{4}}$ is empty.
\end{prop}

\begin{definition} Two $\R$-tubes $\Y_{a,b}$ and $\Y_{c,d}$ are \emph{orthogonal} if they are orthogonal as submanifolds of the symmetric space (where there is a well defined $\Sp(2n,\R)$-invariant scalar product).
\end{definition} 

\begin{remark} \label{ortocond}
The orthogonality relation can be expressed as a property of the cross-ratio of the boundary points:
if $(a,c,b,d)$ is maximal, the $\R$-tubes $\mathcal{Y}_{a,b}$ and $\Y_{c,d}$ are orthogonal if and only if $R(a,c,b,d)=2 \Id$ (see \cite[Definition 4.14]{burger2017maximal}).
\end{remark}

Denote by $((a,b)):= \{ l \in \mathcal{L}| \ (a,l,b) \text{ is maximal} \}$ and by 
\begin{equation*}
    p_{a,b}: \X \cup ((a,b)) \to \Y_{a,b}
\end{equation*}

\noindent the orthogonal projection. It will be useful to have concrete expressions for the orthogonal projection to $((a,b))$ when $(a,b)=(0,l_{\infty})$. 

\begin{lem} [\cite{fanoni2020basmajian}]  \label{orto}
For any $A\in \Sym^{+}(n,\R)$ the $\R$-tubes $\Y_{A,-A}$ and $\Yi$ are orthogonal and their unique intersection point is $iA$. In particular $p_{0,\infty}(A)=iA$.
\end{lem}

\begin{lem} [\cite{fanoni2020basmajian}]  \label{proj} 
If $(a,x,y,b) \in \mathcal{L}(\R^{2n})^{4}$ is a maximal 4-tuple and $p_{a,b}$ is the orthogonal projection onto $\Y_{a,b}$, the distance $$d^{\overline{\mathfrak{a}}^{+}}(p_{a,b}(x),p_{a,b}(y))=(\log\mu_{1},...,\log\mu_{n})$$ where $\mu_{i}$ are the eigenvalues of the cross-ratio $R(a,x,y,b)$.
\end{lem}

\subsection{Computing orthogonal tubes}

A crucial tool to construct right-angled hexagons in the Siegel space $\X$ is computing orthogonal $\R$-tubes. For this reason, this section lists concrete criteria to determine them. 

\begin{lem} \label{unicityoforth}
Let $(P_{1},P_{2},P_{3},P_{4})$ be a maximal $4$-tuple. Then there exists a unique tube $\Y_{P_{5},P_{6}}$ orthogonal to both $\Y_{P_{1},P_{4}}$ and $\Y_{P_{2},P_{3}}$.
\end{lem}

\proof Up to $\Sp(2n,\R)$-action we can consider $(P_{1},P_{2},P_{3},P_{4})=(0,\Id,P,\infty)$. By Lemma \ref{orto} the tubes orthogonal to $\Yi$ are of the form $Y_{-Q,Q}$ where $Q\in \Sym^{+}(n,\R)$. We want to find $Q$ such that the triple $(\Id,Q,P)$ is maximal and such that $\Y_{-Q,Q} \perp \Y_{\Id,P}$. By the orthogonality condition (see Remark \ref{ortocond}) this happens if and only if $R(P,-Q,\Id,Q)=2\Id$. Developing the left-hand side we obtain:
$$2(P+Q)^{-1}Q(Q-\Id)^{-1}(P-\Id)=2(P+Q)^{-1}\big((Q-\Id)Q^{-1}\big)^{-1}(P-\Id)=$$
$$=2\big((\Id-Q^{-1})(P+Q)\big)^{-1}(P-\Id)=2\Id$$
This simplifies to $P+Q-Q^{-1}P-\Id=P-\Id$. We obtain $Q^{2}=P$ which has unique solution $Q=\sqrt{P}$. In particular $P \in \Sym^{+}(n,\R)$ as $(0,\Id,P,\infty)$ is maximal (see Lemma \ref{max}).
\endproof

\begin{lem} \label{comput1}
Let $(0,\Id,P,\infty)$ be a maximal quadruple. Then $\Y_{-\sqrt{P},\sqrt{P}}$ is the unique $\R$-tube orthogonal to both $\Yi$ and $\Y_{\Id,P}$.
\end{lem}
\endproof

\proof Follows directly from the proof of Lemma \ref{unicityoforth}.
\endproof

\begin{lem} \label{comput2}
Let $(P_{1},P_{2},P_{3},P_{4})$ be a maximal quadruple. Then the unique tube orthogonal to both $\Y_{P_{1},P_{4}}$ and  $\Y_{P_{2},P_{3}}$ is $\Y_{Z_{1},Z_{2}}$ where
$$
Z_{1}=g^{-1}\big(-\sqrt{g P_{3}}\big), \  Z_{2}=g^{-1} \big(\sqrt{g P_{3}} \big)
$$
$$g=
\bpm A & 0\\
0 & A^{-T} \epm
\bpm \Id & (P_{1}-P_{4})^{-1}\\
0 & \Id \epm
\bpm 0 & \Id\\
-\Id & 0 \epm
\bpm \Id & -P_{4}\\
0 & \Id \epm$$
\\
and
$A=\sqrt{(P_{1}-P_{4})(P_{2}-P_{1})^{-1}(P_{2}-P_{4})}$.
\end{lem}

\proof The matrix $g\in \Sp(2n,\R)$ is an isometry such that $(gP_{1},gP_{2},gP_{4})= (0,\Id,\infty)$. Result follows from Lemma \ref{comput1}.
\endproof

\begin{lem} \label{quadrismax}
Let $(P_{1},P_{2},P_{3},P_{4},P_{5},P_{6})$ be a maximal 6-tuple and let $Q_{1},Q_{2},Q_{3},Q_{4}$ be such that
$\Y_{P_{1},P_{2}}\perp\Y_{Q_{1},Q_{2}}\perp\Y_{P_{3},P_{4}}\perp\Y_{Q_{3},Q_{4}}\perp\Y_{P_{5},P_{6}}$. Then the quadruple $(P_{3},Q_{2},Q_{3},P_{4})$ is maximal.
\end{lem}

\proof 
Let $g \in \Sp(2n,\R)$ be such that $(gQ_{1},gQ_{2})=(\infty,0)$ and $ (gP_{3},gP_{4})=(-\Id,\Id)$.
We obtain the tubes $g\cdot \Y_{P_{1},P_{2}}=\Y_{-M,M}$ and $g\cdot \Y_{P_{5},P_{6}}=\Y_{P,Q}$ for some $M,P,Q$ positive definite matrices. The tube  $g\cdot \Y_{Q_{3},Q_{4}}=\Y_{X,Y}$ is such that 
\begin{equation} \label{Xposdef}
  \Y_{-\Id,\Id}\perp\Y_{X,Y}\perp \Y_{P,Q}  
 \end{equation}
 where $P$ and $Q$ are positive definite matrices. By construction of orthogonal tubes we know $(-\Id,0,\Id)$ and $(X,\Id,Y)$ maximal (Lemma \ref{comput1}). It is not hard to show that the matrix $X$ needs to be positive definite for the condition (\ref{Xposdef}) to be satisfied. It follows $(-\Id,0,X,\Id)$ maximal and so is its preimage $(P_{3},Q_{2},Q_{3},P_{4})$. 
\endproof

We end this section by giving some concrete expressions to find two orthogonal tubes when one of them is of the form $\Y_{-P,P}$ for a positive definite matrix $P$. This configuration will turn out to be very useful when defining the parameter space of right-angled hexagons.
  
\begin{lem} \label{comput3}  
Let $(0,P_{1},P_{2},\infty)$ be a maximal quadruple. Then 
$$
\Y_{-P_{1},P_{1}} \perp \Y_{P_{1}P_{2}^{-1}P_{1},P_{2}} \text{ and } \ \Y_{-P_{2},P_{2}} \perp \Y_{P_{1},P_{2}P_{1}^{-1}P_{2}}
$$
\end{lem} 

\proof For the first case it is sufficient to find $X \in \Sym(n,\R)$ such that\\ 
 $R(X,P_{1},P_{2},-P_{1})=2\Id$ (see Remark \ref{ortocond}). Developing the left-hand side we obtain
$$
 (X-P_{1})^{-1}(-2P_{1})(-P_{1}-P_{2})^{-1}(X-P_{2})=2\Id
$$
which can be rewritten $
 P_{1}(P_{1}+P_{2})^{-1}(X-P_{2})=(X-P_{1})$. This simplifies to 
$$
(X-P_{2})=(P_{1}+P_{2})P_{1}^{-1}(X-P_{1})=(\Id+P_{2}P_{1}^{-1})(X-P_{1})=X-P_{1}+P_{2}P_{1}^{-1}X-P_{2}
$$
We obtain $
P_{2}P_{1}^{-1}X=P_{1}$. Result follows. The proof for the second case is the same. 
\endproof

\subsection{The symmetric spaces $\X_{\GL(n,\R)}$ and $\X_{\SL(n,\R)}$ } \label{sym}

Recall that the standard model for the symmetric space associated to $\GL(n,\R)$ is $$\X_{\GL(n,\R)}=\Sym^{+}(n,\R)$$ We endow $\X_{\GL(n,\R)}$ with the distance $d_{\GL}$ given by $$d_{\GL}(X,Y)=\sqrt{\sum_{i=1}^{n} (\log \lambda_{i})^{2}}$$ where $\lambda_{i}$ are the eigenvalues of $XY^{-1}$. With this choice of $d_{\GL}$ the natural identification $\X_{\GL(n,\R)}=\Yi$ is an isometry (where $\Yi$ is equipped with the induced Riemannian metric). Recall also that the symmetric space associated to $\SL(n,\R)$ is $$\X_{\SL(n,\R)}=\{ X \in \Sym^{+}(n,\R)| \ \det(X)=1 \}$$ Similarly, we endow $\X_{\SL(n,\R)}$ with the distance $d_{\SL}$ given by $$d_{\SL}(X,Y)=\sqrt{\sum_{i=1}^{n} (\log \lambda_{i})^{2}}$$ where $\lambda_{i}$ are the eigenvalues of $XY^{-1}$.
In particular, the symmetric space $\X_{\SL(2,\R)}$ can be identified with the hyperbolic upper-half plane $\H^{2}= \{ z=x+iy| \ x,y \in \R, \ y>0 \}$ via the following map:
\begin{equation*} 
\begin{aligned}
 h: \mathcal{X}_{\SL(2,\R)}&\to \H^{2}\\
B& \mapsto \big[ B \big] \cdot i
\end{aligned}
 \end{equation*}
 
\noindent where $B$ is an element of $\PSL(2,\R$) and acts on $\H^{2}$ via Möbius transformations. The inverse of $h$ is given by $h^{-1}(z)=\sqrt{AA^{T}}$ where $A \in \PSL(2,\R)$ and $A \cdot i =z$. 

\begin{remark} \label{scale}
If we endow $\H^{2}$ with the distance $d_{\H^{2}}$ relative to the standard metric $\frac{dx^{2}+dy^{2}}{y^{2}}$ on the upper-half plane then  $h$ is not an isometry and in general $(\X_{\SL(2,\R)},d_{\SL})$ and $(\H^{2},d_{\H^{2}})$ are not isometric. To have an isometry we have to scale $d_{\H^{2}}$ by a factor of $\frac{1}{\sqrt{2}}$. 
\end{remark}

\subsection{The geometry of the standard tube $\Yi$} \label{seccyl}

Let us consider the symmetric space $\X$ associated to $\Sp(4,\R)$.

\begin{lem} \label{HxR}
The tube $\Yi$ is isometrically identified with $\R \times \H^{2}$.
\end{lem}

\proof As seen in Section \ref{sym} there is a natural identification $\Yi = \X_{\GL(2,\R)}$, where $\X_{\GL(2,\R)}$ is the set of positive definite symmetric matrices. The map
\[
\begin{aligned}
f= \pi^{\R} \times \pi^{\H^{2}}: \X_{\GL(2,\R)}&\to \R\times\X_{\SL(2,\R)}\\
Q& \mapsto \Big( \frac{\log \det Q}{\sqrt{2}}, \frac{Q}{\sqrt{\det Q}} \Big)
\end{aligned}
\]

\noindent is a bijection, with inverse $f^{-1}\big( (r,B) \big)=\sqrt{e^{\sqrt{2}r}}B$. When $\R\times\X_{\SL(2,\R)}$ and $\X_{\GL(2,\R)}$ are considered as metric spaces (endowed with $d_{\R} \times d_{\SL}$ and  $d_{\GL}$ respectively), the map $f$ is an isometry (\cite{fanoni2020basmajian} Lemma 2.17). The identification $\Yi=\R \times \H^{2}$ follows from $\X_{\SL(2,\R)} = \H^{2}$ (see Section \ref{sym}). In particular all copies of $\H^{2}$ in $\Yi$ are canonically identified. Observe that to turn $\Yi=\R \times \H^{2}$ into an isometric identification we have to scale the metric $d_{\H^{2}}$ by a factor of $\frac{1}{\sqrt{2}}$ (Remark \ref{scale}).
\endproof

\begin{definition} \label{defhypcomp}
Given $A \in \Sym^{+}(2,\R)$ the \emph{hyperbolic component} of $A$ is the point $\pi^{\H^{2}}(p_{0,\infty}(A))$. Similarly the $\R$-component of $A$ is the point $\pi^{\R}(p_{0,\infty}(A))$ and will be called \emph{level} of $A$.
\end{definition}

Let $\H^{2}= \{x+iy| \ x,y \in \R, \ y>0 \}$, then for any fixed level in $\R \times \H^{2}$ the set of diagonal matrices coincides with the $y$-axis of $\H^{2}$, where the set  $ \big\{ \begin{pmatrix} 
\lambda_{1} & 0\\
0 & \lambda_{2}
\end{pmatrix}, \lambda_{1}>\lambda_{2}, \ \lambda_{1} \cdot \lambda_{2}=1 \big\}$ consists of points "above" $i \in \H^{2}$ in the vertical $y$-axis of the hyperbolic plane. 

\begin{remark} \label{meaningvector}
Let $iX,iY \in \Yi$ and $(d_{1},d_{2}) \in \a$ such that $d^{\a}(iX,iY)=(d_{1},d_{2})$. We can associate to $(d_{1},d_{2})$ a more meaningful vector $(\mathbf{r}_{(d_{1},d_{2})},\mathbf{h}_{(d_{1},d_{2})})$ (Figure \ref{cyl}). The vector $(\mathbf{r},\mathbf{h})$ based at $iX$ has first coordinate $\mathbf{r}$ equal to the difference between the levels of $iX$ and $iY$
$$\mathbf{r}=d^{\R}(\pi^{\R}(iX),\pi^{\R}(iY)) $$ and second coordinate $\mathbf{h}$ equal to the distance between the two points in $\H^{2}$
$$\mathbf{h}=d^{\H^{2}}(\pi^{\H^{2}}(iX),\pi^{\H^{2}}(iY))$$

\begin{figure}[!h]
   \centering
   \setlength{\unitlength}{0.1\textwidth}
   \begin{picture}(6,3)
     \put(0.8,0){\includegraphics[width=6.5cm,height=5.3cm]{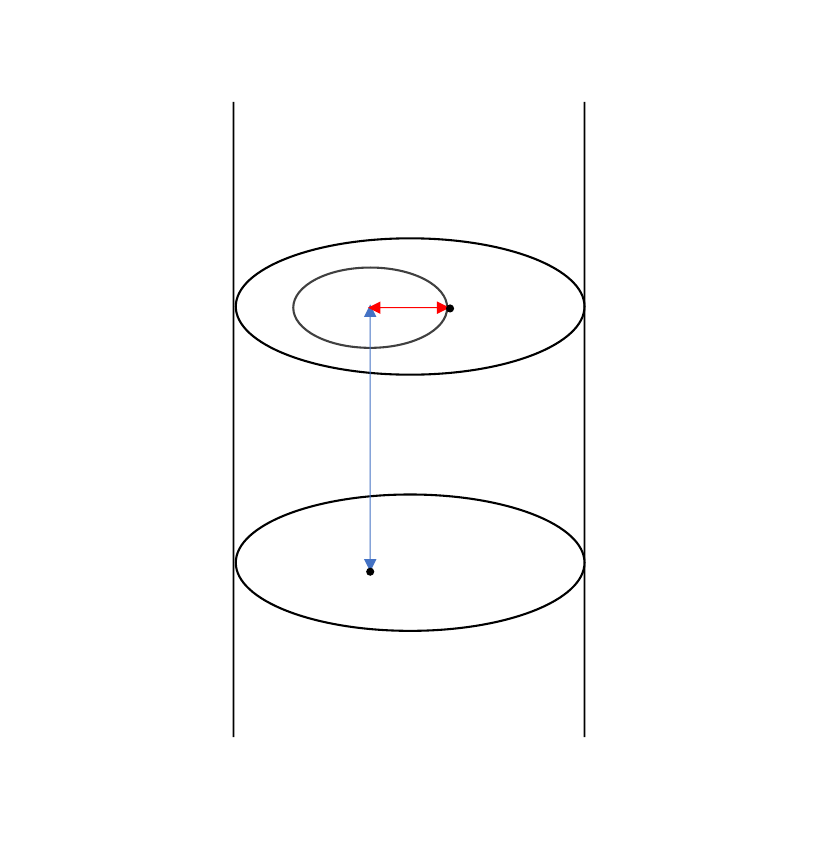}}
     \put(2.8,2.23){\textcolor{red}{$\mathbf{h}$}}
     \put(2.8,1.6){\textcolor{blue}{$\mathbf{r}$}}
     \put(2.3,1.1){$iX$}
      \put(3.2,2.2){$iY$}
      \put(2.8,0.4){$\Yi$}
   \end{picture}
   \caption{Geometric interpretation of $(\mathbf{r},\mathbf{h})$ for two points $iX,iY$ in $\Yi$}
   \label{cyl}
\end{figure}
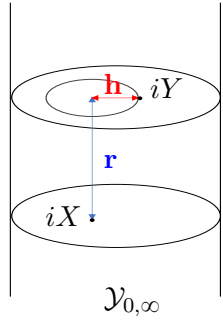

\noindent It is not hard to show that given $d^{\a}(iX,iY)=(d_{1},d_{2})$ then
$\mathbf{r}= \frac{d_{1}+d_{2}}{\sqrt{2}}$ and $\mathbf{h}=(d_{1}-d_{2})$. The vector $(\mathbf{r},\mathbf{h})$ also gives a geometric condition for the maximality of the triple $(l_{\infty},X,Y)$ i.e. for the matrix $Y-X$ to be positive definite (see Lemma \ref{max}). It holds (\cite[Corollary 2.21]{fanoni2020basmajian}):
$$
Y-X \text{ positive definite } \iff \mathbf{r}> \frac{1}{\sqrt{2}}\mathbf{h}
$$
\end{remark}

The following Lemma of linear algebra will play a crucial role in the definition of parameters for a hexagon inside $\X$.

\begin{lem} \label{diagonalizing_mat} 
For any $M \in \Sym^{+}(2,\R)$ with distinct eigenvalues there exist unique $S,Q$ in $\PSO(2)$ and $\PO(2)\backslash\PSO(2)$ respectively such that $$SMS^{T}=QMQ^{T}=\begin{pmatrix}
\lambda_{1} & 0\\
0 & \lambda_{2}
\end{pmatrix}, \text{where } \lambda_{1}>\lambda_{2}$$
\end{lem}

\proof Let $v_{1},v_{2}$ be orthonormal eigenvectors relative to the eigenvalues  $\lambda_{1}>\lambda_{2}>0$ respectively and let $L$ denote the orthogonal matrix $L=\Big( \begin{bmatrix}
v_{1}
\end{bmatrix} \begin{bmatrix}
v_{2}
\end{bmatrix} \Big)$. If $\det L =1$, it is a standard fact of linear algebra that $
S=\begin{pmatrix}
[ \  \ v_{1}^{T} \ \ ]\\
[ \  \ v_{2}^{T} \ \ ]
\end{pmatrix}=L^{T}
$ is the unique element of $\PSO(2)$ such that $SMS^{T}=\bpm
\lambda_{1} & 0\\
0 & \lambda_{2}
\epm, \ \lambda_{1}>\lambda_{2}$. Put $Q =  \begin{pmatrix}
-1 & 0\\
0 & 1
\end{pmatrix}S = \begin{pmatrix}
-[ \  \ v_{1}^{T} \ \ ]\\
 \ \ [ \  \ v_{2}^{T} \ \ ]
\end{pmatrix}$. Then $\det Q =-1$ and $Q$ is the desired matrix in $\PO(2)\backslash\PSO(2)$. If $\det L =-1$ then the two diagonalizing matrices are $Q=L^{T}$ and $S=\begin{pmatrix}
-1 & 0\\
0 & 1
\end{pmatrix}Q$.

\endproof

\subsection{Geometric interpretation of diagonalization matrix} \label{geomintdiag}
The group $\PSp(4,\R)$ acts on a point $iM \in \Yi$ via fractional linear transformations. In particular the isometries stabilizing the standard tube are of the form $\bpm A&0\\0&A^{-T}\epm$ where $A \in \GL(2,\R)$.  Recall that we identify $\Sym^{+}(2,\R)=\Yi= \R \times \H^{2}$ (see Sections \ref{sym} and \ref{seccyl}). In this identification the identity matrix is identified with the point $(0,i) \in \R\times\H^{2}$. Moreover, all copies of $\H^{2}$ in $\Yi$ are canonically identified.
For a matrix $S \in \PSO(2)$ we want to interpret the action 
\begin{equation} \label{diagaction}
\PSp(4,\R) \ni \begin{pmatrix}
S & 0\\
0 & S
\end{pmatrix} \cdot (iM)=
iSMS^{T}
\end{equation}
as a transformation which fixes the level of $M$ and rotates its hyperbolic component around $i \in \H^{2}$. If we consider the action on $\H^{2}$ through Möbius transformations we see that 
$$
\Stab_{\PSL(2,\R)}(i)= \PSO(2)
$$
For $\theta \in [0,\pi)$ the action (\ref{diagaction}) of a matrix $S=\begin{pmatrix}
\cos\theta & -\sin\theta\\
\sin\theta & \cos\theta
\end{pmatrix} \in \PSO(2)$ can be interpreted as a clockwise rotation  of angle $2\theta$ around $i \in \H^{2}$. For every $S \in \PSO(2)$ and every $\theta \in [0,\pi)$ there is a unique way to write $S$ as a rotation matrix of the form
$$S=\begin{pmatrix}
\cos\theta & -\sin\theta\\
\sin\theta & \cos\theta
\end{pmatrix} \sim \begin{pmatrix}
\cos(\pi+\theta) & -\sin(\pi+\theta)\\
\sin(\pi+\theta) & \cos(\pi+\theta)
\end{pmatrix}=-S$$

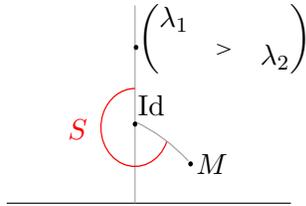
\begin{figure}[!h]
   \centering
   \captionsetup{justification=centering,margin=2cm}
   \setlength{\unitlength}{0.1\textwidth}
   \begin{picture}(3,1.9)
     \put(0,0){\includegraphics[width=4cm,height=3cm]{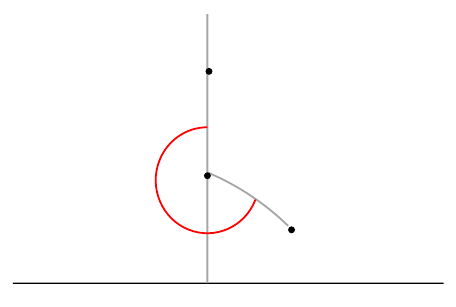}}
          \put(0.6,0.7){\textcolor{red}{$S$}}
          \put(1.2,0.9){$\Id$}
          \put(1.7,0.4){$M$}
          \put(1.2,1.5){$\begin{pmatrix}
\lambda_{1}&{}& {}\\
{}&^{>}& \lambda_{2}
\end{pmatrix}$}
\end{picture}
\caption{Geometric interpretation of diagonalization matrix $S$}
   \label{angle}
\end{figure}

\noindent Given M positive definite with distinct eigenvalues we interpret the unique $S \in \PSO(2)$ for which
$SMS^{T}=\bpm
\lambda_{1} & 0\\
0 & \lambda_{2}
\epm, \ \lambda_{1}>\lambda_{2}
$ as the angle formed by the semi-axis $\{(0,y)| \ y>1\}$ inside $\H^{2}$ and the geodesic segment connecting the hyperbolic components of $\Id$ and $M$ (Figure \ref{angle}).

In this paper we will  use both the matrix and the angle notation: angle parameters will be denoted by $S$ or $\alpha$ depending on the context, where
$$S=\bpm
\cos\frac{\alpha}{2} & -\sin\frac{\alpha}{2}\\
\sin\frac{\alpha}{2} & \cos\frac{\alpha}{2}
\epm,\alpha \in [0,2\pi)$$

\begin{remark} (\emph{Drawing angles "on the left"}) \label{left} The matrix $S=\begin{pmatrix}
\cos\theta & -\sin\theta\\
\sin\theta & \cos\theta
\end{pmatrix}$ acts on $M$ by clockwise rotation of center $\Id$ and angle $2\theta$ on the $\H^{2}$-component of the standard tube $\Yi$. For this reason to draw the angle parameters we will consider  the oriented geodesic going from $M$ to $\Id$ and draw the angle on the left of it.
\end{remark}

\begin{remark} \label{refl}
For $\Lambda=\bpm \lambda_{1} & 0\\ 0 & \lambda_{2}
\epm$ with $\lambda_{1} \neq \lambda_{2}$, the stabilizer of the quadruple $(0,\Id,\Lambda,\infty)$ is given by
\begin{equation*}
\text{Stab}_{\PSp(4,\R)}(0,\Id,\Lambda,\infty)=\Big\{ \Id, \bpm r & 0\\0 & r
\epm \Big\} \cong \mathbb{Z}/2\mathbb{Z}
\end{equation*}

\noindent where $r=\begin{pmatrix}
-1 & 0\\
0 & 1
\end{pmatrix}$. In the identification $\Yi=\R \times \H^{2}$ the action of the matrix $\bpm
r & 0\\
0 & r
\epm$ on $\Yi$ is a reflection across the $y$-axis of $\H^{2}$. Let $M=\begin{pmatrix}
m_{1} & m_{2}\\
m_{2} & m_{3}
\end{pmatrix}$ be positive definite and identify $M$ with $a+ib$ in $\H^{2}$ as in Section \ref{sym}. Then
$$rMr=\begin{pmatrix}
m_{1} & -m_{2}\\
-m_{2} & m_{3}
\end{pmatrix} \overset{(\S \ref{sym})}{=} -a+ib \in \H^{2}
$$
Put $M^{r}=rMr$. If $M$ is a point of angle $\alpha \in [0,2\pi)$ from $\bpm
\lambda_{1} & 0\\
0 & \lambda_{2}
\epm, \ \lambda_{1}>\lambda_{2}$, then $M^{\emph{r}}$ is a point of angle $(2\pi-\alpha)$ from $\bpm
\lambda_{1} & 0\\
0 & \lambda_{2}
\epm, \ \lambda_{1}>\lambda_{2}$ (Figure \ref{angle2}).

\begin{figure}[!h]
   \centering
   \captionsetup{justification=centering,margin=2cm}
   \setlength{\unitlength}{0.1\textwidth}
   \begin{picture}(5.3,1.7)
     \put(0.8,0){\includegraphics[width=5cm,height=3cm]{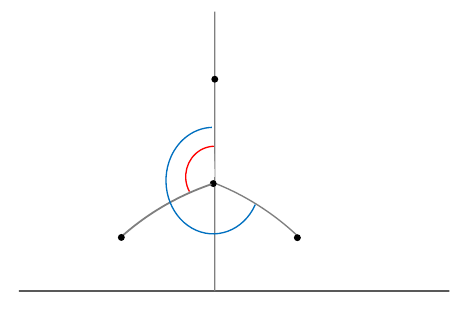}}
          \put(2.05,1.1){\textcolor{red}{$S_{1}$}}
          \put(1.7,1){\textcolor{blue}{$S_{2}$}}
          \put(2.35,0.8){$\Id$}
          \put(2.9,0.5){$M^{r}$}
          \put(1.4,0.5){$M$}
          \put(2.4,1.5){$\begin{pmatrix}
\lambda_{1}&{}& {}\\
{}&^{>}& \lambda_{2}
\end{pmatrix}$}
\end{picture}
\caption{The point $M^{r}$ is obtained by reflecting $M$ across the $y$-axis }
   \label{angle2}
\end{figure}
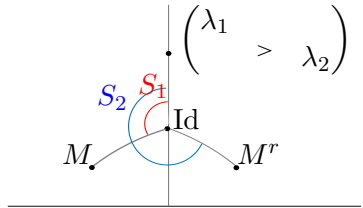

\noindent To see this using the angle interpretation of diagonalization matrices take the unique $S_{1} \in \PSO(2)$ diagonalizing $M$ as in Lemma \ref{diagonalizing_mat}. It holds:
\begin{equation} \label{clarity}
S_{1}MS_{1}^{T}=\bpm
\lambda_{1} & 0\\
0 & \lambda_{2}
\epm= \ 
S_{1}
rM^{r}r
S_{1}^{T}
\end{equation}
so that
$$
\begin{pmatrix}
\lambda_{1} & 0\\
0 & \lambda_{2}
\end{pmatrix}=r
\begin{pmatrix}
\lambda_{1} & 0\\
0 & \lambda_{2}
\end{pmatrix}
r\overset{(\ref{clarity})}{=}
\underbrace{rS_{1}r}_{S_{2}}M^{r}\underbrace{rS_{1}^{T}r}_{S_{2}^{T}}
$$
If 
$$S_{1}=\begin{pmatrix}
\cos\theta & -\sin\theta\\
\sin\theta & \cos\theta
\end{pmatrix}$$ then $$S_{2}=\begin{pmatrix}
\cos\theta & \sin\theta\\
-\sin\theta & \cos\theta
\end{pmatrix}=\begin{pmatrix}
\cos(2\pi-\theta) & -\sin(2\pi-\theta)\\
\sin(2\pi-\theta) & \cos(2\pi-\theta)
\end{pmatrix}$$
\end{remark}

\noindent More generally if $A,B$ are symmetric positive definite matrices such that $A^{-1}B$ has distinct eigenvalues $\lambda_{1}>\lambda_{2}$, then the stabilizer $\text{Stab}_{\PSp(4,\R)}(0,A,B,\infty)$ is isomorphic to $\mathbb{Z}/2\mathbb{Z}$ where the non-trivial element represents a reflection in the hyperbolic component of $\Yi$ across the geodesic going through $\pi^{\H^{2}}(iA)$ and $\pi^{\H^{2}}(iB)$.

\subsection{Orientation of the hyperbolic component of $\Yi$}    \label{secorienthypcomp}

\begin{prop} \label{orienthypcomp}
Consider the symmetric space $\X$ associated to $\Sp(4,\R)$. Choosing an orientation of  $\H^{2}$ inside $\Yi=\R \times \H^{2}$ is equivalent to choosing an orientation of  $\mathbb{P}(\infty) \simeq \mathbb{P}(0)$ where $\infty,0 \in \mathcal{L}(\R^{4})$.  
\end{prop}
\proof
Fix a basis $\mathcal{B}=\{e_{1},e_{2},e_{3},e_{4} \}$ of $\mathbb{R}^{4}$ and consider the two transverse Lagrangians \\ $0=\langle e_{3},e_{4} \rangle$ and $\infty=\langle e_{1},e_{2} \rangle$. The standard tube
$\Yi= \{ iY| \ Y \in \Sym^{+}(2,\R)  \}$ is isometrically identified with $\R \times \H^{2}$ (Lemma \ref{HxR}). The hyperbolic plane inside $\Yi$ is identified with the symmetric space associated to $\SL(2,\R)$, that is $\X_{\SL(2,\R)}=\{ X \in \Sym^{+}(2,\R)| \ \det(X)=1 \}$. All copies of $\H^{2}$ inside $\Yi$ are canonically identified and stabilized by the set of matrices
\begin{equation} \label{O(2)_set}
\Big\{ \bpm
R & 0\\
0 & R 
\epm, R \in \O(2)\Big\} 
\end{equation}
acting $\Sp(4,\R)$-equivariantly in the identification $ \X \overset{\iota}{\mapsto} \mathbb{X}$. Let $
\gamma(t)=
\begin{pmatrix}
ie^{t} & 0\\
0 & ie^{-t} 
\end{pmatrix}
$ be a geodesic ray
lying inside the hyperbolic component of $\Yi$. Then $\gamma(t)$ converges to the Lagrangian  $l_{+}=\langle e_{1},e_{4} \rangle$ when $t \to \infty$, and to the Lagrangian $l_{-}=\langle e_{2},e_{3} \rangle$ when $t \to -\infty$. To fix an orientation of $\H^{2}$ inside $\Yi$ it is sufficient to orient its visual boundary $\partial_{\infty} \H^{2}$. In Section \ref{geomintdiag} we have investigated the action of  orthogonal matrices on $\H^{2} \subset \Yi=\R \times \H^{2}$ and we have interpreted matrices as angles in the hyperbolic component. The set in (\ref{O(2)_set}) acts preserving the Lagrangians $0$ and $\infty$ respectively and the visual boundary of the hyperbolic component is realized by the $\O(2)$-orbit of the Lagrangian $l_{+}$ so that $\partial_{\infty} \H^{2}$ is given by $\{g \cdot l_{+}\}$ where $g$ is as in (\ref{O(2)_set}). Any such point $g \cdot l_{+}$ is a Lagrangian which intersects $\infty$ in one line. We get the following identification:
\begin{equation*} 
\begin{aligned}
 \partial_{\infty} \H^{2}& \to \mathbb{P}(\infty) \\
 l& \mapsto l \cap \infty
\end{aligned}
 \end{equation*}
\noindent To fix an orientation of $\H^{2}$ inside $\Yi$ it is therefore sufficient to orient $\mathbb{P}(\infty)$.
This set is canonically identified with $\mathbb{P}(0)$. To see this let $v^{\perp \omega} = \{ u \in \R^{4}| \ \omega(v,u)=0 \}$. 

\noindent In particular $v \in v^{\perp \omega}$ and dim($v^{\perp \omega}$)=3. Then $\mathbb{P}(0)$ and $\mathbb{P}(\infty)$ are identified through the map $[v] \mapsto [v^{\perp \omega} \cap \infty]$.
\endproof

\subsection{Isometries reflecting the hyperbolic component} \label{isomrefl}

\begin{prop}
For $A \in \GL(2,\R)$ let $f_{A}$ be an isometry stabilizing the standard tube $\Yi$:
 \begin{equation*} 
\begin{aligned}
f_{A}: \Yi&\to \Yi\\
iY& \mapsto iAYA^{T}
\end{aligned}
 \end{equation*}

Then $f_{A}$ is reversing the orientation of the hyperbolic component of $\Yi$ if and only if $\det A<0$.
\end{prop}

\proof 
Recall that $\Yi =\chi_{\GL_{2,\R}} = \Sym^{+}(2,\R)= \R \times \chi_{\SL_{2,\R}}= \R \times \H^{2} $. The isometry $f_{A}$ is linear in $Y$ and its differential is the map $X \mapsto AXA^{T}$ for any tangent vector $X=\bpm
x_{1} & x_{2} \\
x_{2}  & x_{3}  
\epm$, $X \cong \bpm
x_{1}\\
x_{2}\\  
x_{3}  
\epm$.
For $A=\bpm
a_{1} & a_{2} \\
a_{3}  & a_{4}  
\epm$ the tangent vector $AXA^{T}=df_{A}(X)$ can be rewritten as
$$
\bpm
a_{1}^{2}x_{1}+2a_{1}a_{2}x_{2}+a_{2}^{2}x_{3}\\
a_{1}a_{3}x_{1}+(a_{1}a_{4}+a_{2}a_{3})x_{2}+a_{2}a_{4}x_{3}\\  
a_{3}^{2}x_{1}+2a_{3}a_{4}x_{2}+a_{4}^{2}x_{3}  
\epm=\bpm
a_{1}^{2} & 2a_{1}a_{2} & a_{2}^{2} \\
a_{1}a_{3} & a_{1}a_{4}+a_{2}a_{3} & a_{2}a_{4}\\
a_{3}^{2} &  2a_{3}a_{4} & a_{4}^{2}
\epm \bpm
x_{1}\\
x_{2}\\  
x_{3}  
\epm
$$
where 
$$
\det \bpm
a_{1}^{2} & 2a_{1}a_{2} & a_{2}^{2} \\
a_{1}a_{3} & a_{1}a_{4}+a_{2}a_{3} & a_{2}a_{4}\\
a_{3}^{2} &  2a_{3}a_{4} & a_{4}^{2}
\epm= (\det A)^{3}
$$

\noindent The map $f_{A}$ is therefore reversing the orientation of the tube $\Yi$ if and only if $\det A <0$. To finish the proof we need to show that only the orientation of the hyperbolic component can be reversed, not the orientation of the $\R$-component. This holds as the action of a $f_{A}$ on the $\R$-component of  $iY \in \R \times \H^{2}$ is a translation: if $\pi^{\R}(iY)$ is the $\R$-component of $iY$ then the $\R$-component of $iAYA^{T}$ is given by $\pi^{\R}(iY)+\frac{2 \log |\det A|}{\sqrt{2}}$. 
\endproof

\begin{remark}
More generally given $l_{1},l_{2}$ in $\mathcal{L}(\R^{4})^{(2)}$ and $g$ an isometry stabilizing the tube $\Y_{l_{1},l_{2}}$, whether or not $g$ is reversing the orientation of the hyperbolic component of $\Y_{l_{1},l_{2}}$ is intrinsic and only depends on the sign of $\det A$.
\end{remark}

\begin{definition} \label{isomref} 
An isometry $g\in \PSp(4,\R)$ conjugate to  $\begin{pmatrix}
A & 0\\
0 & A^{-T}
\end{pmatrix}$ with $A \in \GL(2,\R)$ is called \emph{reflecting} (resp. \emph{non-reflecting} ) if $\det A<0$ (resp. $>0$). 
\end{definition}

\section{Parameters for quintuples} \label{parforquin}

\subsection{The sets $\mathcal{Q}^{gen}$ and $\mathcal{Q}^{st}$} In this section we introduce parameters for generic quintuples. 

\begin{definition} \label{genquadr}
Let $(P,X,Y,Q)$ be a maximal quadruple and let $\mu_{1},...,\mu_{n}$ be the eigenvalues of the cross-ratio $R(P,X,Y,Q)$. The quadruple $(P,X,Y,Q)$ is said to be \emph{generic} if for any $i \neq j$ it holds $\mu_{i} \neq \mu_{j}$.
\end{definition}

\begin{remark}
Recall that we denote by $p_{P,Q}$ the orthogonal projection on the tube $\Y_{P,Q}$. Let $\underline{b}=d^{\a}(p_{P,Q}(X),p_{P,Q}(Y))$ be the vector obtained by the orthogonal  projection of $X$ and $Y$ on the tube $\Y_{P,Q}$. From Lemma \ref{proj} it is easy to see that the quadruple $(P,X,Y,Q)$ is generic if and only if the vector $\underline{b}$ is a regular vector of the Weyl chamber.
\end{remark}

\begin{definition} \label{gentrip}
The set of \emph{generic quintuples} $\mathcal{Q}^{gen}$ is given by:
\begin{equation*}
\mathcal{Q}^{gen}:= \{(P,X,Y,Z,Q) \text{ maximal} |  \ (P,X,Y,Q) \text{ and }(P,Y,Z,Q) \text{ generic }\}
\end{equation*}
\end{definition}

\begin{remark}
Observe that the definition of generic quintuple strongly depends on the order of the quintuple: given $(P,X,Y,Z,Q)$ generic it is not necessarily true that a cyclic permutation of the quintuple is generic. 
\end{remark}

We will see in the next section how the parametrization of $\mathcal{Q}^{gen}$ is connected with the parametrization of right-angled hexagons of $\X$. Let us now consider the symmetric space associated to $\Sp(4,\R)$. 

\begin{definition} The set of \emph{standard quintuples} $ \mathcal{Q}^{st}\subset \mathcal{Q}^{gen}$ is given by
$$
\mathcal{Q}^{st}:= \{(0,X,\Id,\bpm
\lambda_{1} & 0\\
0 & \lambda_{2}
\epm,\infty) \in \mathcal{Q}^{gen}| \ \lambda_{1},\lambda_{2} \in \R, \ \lambda_{1}>\lambda_{2} \ \}
$$   
\end{definition}

\begin{remark} \label{iso}
Recall that for a diagonal matrix $\Lambda$ with different eigenvalues the stabilizer of $(0,\Id,\Lambda,\infty)$ is isomorphic to $\mathbb{Z}/2\mathbb{Z}$ (see Remark \ref{refl}). For any $(P,X,Y,Z,Q) \in \mathcal{Q}^{gen}$ we can always find a $g \in \PSp(4,\R)$ such that $g(P,Y,Z,Q)=(0,\Id,\bpm
\lambda_{1} & 0\\
0 & \lambda_{2}
\epm,\infty)$ where $\lambda_{1} > \lambda_{2}$. It is therefore clear that $
\mathcal{Q}^{gen}/_{\PSp(4,\R)} \cong \mathcal{Q}^{st}/_{\mathbb{Z}/2\mathbb{Z}}$.
\end{remark}

We denote by $\mathfrak{a}$ the set of regular vectors inside $\a$  
$$\mathfrak{a}=\{(x_{1},x_{2}) \in \R^{2}| \ x_{1}> x_{2}> 0 \}$$

\begin{prop} \label{partrip} 
The set $\mathcal{Q}^{gen}/_{\PSp(4,\R)}$ is parametrized by $\mathfrak{a}^{2} \times \PSO(2) / _{\sim}$ where for $S$ in $\PSO(2)$  and $r=\bpm -1 & 0\\0 & 1\epm$ the equivalent relation is $
S \sim rSr$.

The parametrization is given by 
$$
\Big((c_{1},c_{2}),(d_{1},d_{2}),[S]\Big) \mapsto  \Big[ \Big( 0,S^{T}\begin{pmatrix}
\frac{1}{e^{c_{2}}} & 0\\
0 & \frac{1}{e^{c_{1}}}
\end{pmatrix}S,\ \Id, \ \begin{pmatrix}
e^{d_{1}} & 0\\
0 & e^{d_{2}}
\end{pmatrix}, \infty \Big)\Big] \in \mathcal{Q}^{st}/ _{\mathbb{Z}/2\mathbb{Z}}
$$

\noindent with inverse  
$$
[(P,X,Y,Z,Q)] \mapsto  \Big( \ (c_{1},c_{2})=d^{\mathfrak{a}^{+}}\big(p_{P,Q}(X), p_{P,Q}(Y)\big),(d_{1},d_{2})=d^{\mathfrak{a}^{+}}\big(p_{P,Q}(Y),p_{P,Q}(Z)\big),[S] \ \Big)
$$

\noindent where 
$$
S (g X) S^{T}= \bpm
\frac{1}{e^{c_{2}}} & 0\\
0 & \frac{1}{e^{c_{1}}}
\epm, \ \frac{1}{e^{c_{2}}} > \frac{1}{e^{c_{1}}}
$$
and $g$ is a map in $\PSp(4,\R)$ such that $g(P,X,Y,Z,Q) \in \mathcal{Q}^{st}$.

\noindent The parameter space can be rewritten as $\mathfrak{a}^{2} \times [0,2\pi) / _{\sim}$ where for $\alpha \in [0,2\pi)$ the corresponding $\PSO(2)$-parameter is $S=\bpm
\cos\frac{\alpha}{2} & -\sin\frac{\alpha}{2}\\
\sin\frac{\alpha}{2} & \cos\frac{\alpha}{2}
\epm$ and the equivalence relation is $\alpha \sim (2 \pi- \alpha)$.

\end{prop} 

\proof
We first show how to find parameters $(\underline{c},\underline{d},[S])$ for a given quintuple $[(P,X,Y,Z,Q)]$ in $\mathcal{Q}^{gen} / _{\PSp(4,\R)}$. We want to use the fact: 
$$
    \mathcal{Q}^{gen}/_{\PSp(4,\R)} \cong \mathcal{Q}^{st}/_{\mathbb{Z}/2\mathbb{Z}}
    $$ 
    
\begin{figure}[!h]
   \centering
   \captionsetup{justification=centering,margin=2cm}
   \setlength{\unitlength}{0.1\textwidth}
   \begin{picture}(4.6,2.5)
     \put(0.8,0){\includegraphics[width=4cm,height=4cm]{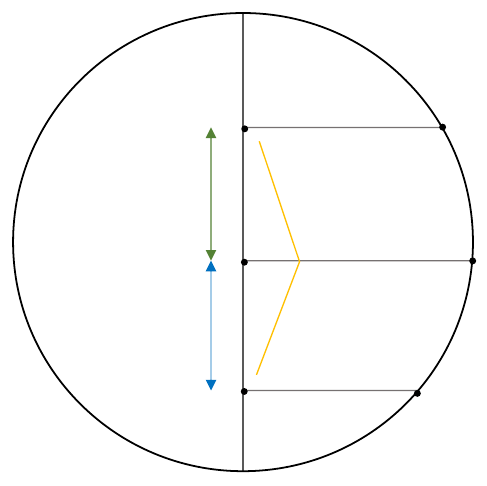}}
          \put(2,-0.1){$0$}
          \put(2,2.55){$\infty$}
          \put(3.4,1.1){$\Id= gY$}
          \put(3.1,0.4){$gX$}
          \put(3.2,1.9){$\bpm
e^{d_{1}} & 0\\
0 & e^{d_{2}}
\epm=gZ$}
          \put(1.7,1.5){\textcolor{forestgreen}{$\underline{d}$}}
          \put(1.7,0.9){\textcolor{blue}{$\underline{c}$}}
          \put(2.4,1.3){\textcolor{orange}{$S$}}
   \end{picture}
\caption{The isometry $g$ sends the quintuple $(P,X,Y,Z,Q)$ to a standard one}
   \label{gentripl}
\end{figure}
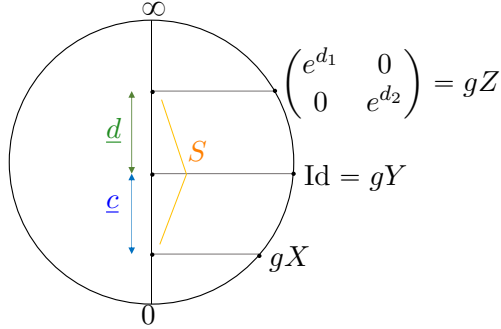

\noindent Let $(P,X,Y,Z,Q)\in \mathcal{Q}^{gen}$. Up to isometry  we can consider $P=0$ and $Q=\infty$. Put
$$\underline{c}=(c_{1},c_{2})=d^{\mathfrak{a}^{+}}\big(p_{0,\infty}(X), p_{0,\infty}(Y)\big) \in \mathfrak{a}$$
$$\underline{d}=(d_{1},d_{2})=d^{\mathfrak{a}^{+}}\big(p_{0,\infty}(Y),p_{0,\infty}(Z)\big) \in \mathfrak{a}$$
\noindent Recall that the isometries stabilizing the standard tube are of the form $\bpm A&0\\0&A^{-T}\epm$ where $A \in \GL(2,\R)$. Let $g \in \Stab(\Yi)$ be such that
$gY=\Id$ and $gZ=\bpm e^{d_{1}} & 0\\0 & e^{d_{2}}\epm$ (Figure \ref{gentripl}).

The first equality forces $A=\mathcal{O} \sqrt{Y^{-1}}$ where $\mathcal{O}\in \PO(2)$. The second equality forces $\mathcal{O}=P,Q$ where $P$ and $Q$ are the unique matrices in $\PSO(2)$ and $\PO(2) \backslash \PSO(2)$ respectively (see Lemma \ref{diagonalizing_mat}) such that  $$P(\sqrt{Y^{-1}}Z\sqrt{Y^{-1}})P^{T}=Q(\sqrt{Y^{-1}}Z\sqrt{Y^{-1}})Q^{T}=\begin{pmatrix}
e^{d_{1}} & 0\\
0 & e^{d_{2}}
\end{pmatrix}$$

\noindent Accordingly, the two only possibilities for $g$ are:
$$
g_{1}=\begin{pmatrix}
P\sqrt{Y^{-1}} & 0\\
0 & P\sqrt{Y}
\end{pmatrix} \text{ and } g_{2}=\begin{pmatrix}
Q\sqrt{Y^{-1}} & 0\\
0 & Q\sqrt{Y}
\end{pmatrix}
$$
In particular we know $Q=rP$ where $r=\bpm -1 & 0\\0 & 1\epm$, so that $g_{2}X=r g_{1}Xr$. 

\noindent Let $S$ be the unique matrix in $\PSO(2)$ such that $Sg_{1}XS^{T}=\bpm
\lambda_{1}&0\\0&\lambda_{2}\epm$ for some $\lambda_{1}>\lambda_{2}$. Then $d^{\a}(g_{1}X,\Id)=(c_{1},c_{2})$ so that $\lambda_{1}=\frac{1}{e^{c_{2}}}$ and $\lambda_{2}=\frac{1}{e^{c_{1}}}$. 

\noindent On the other hand $S'=rSr$ is the unique matrix in $\PSO(2)$ such that $S'g_{2}XS'^{T}=\bpm
\frac{1}{e^{c_{2}}} & 0\\
0 & \frac{1}{e^{c_{1}}}
\epm$. The point $g_{2}X=rg_{1}Xr$ is the image of $g_{1}X$ under a reflection on the hyperbolic component (see Remark \ref{refl}) and the two quintuples $$
(0,g_{1}X,\Id,\bpm e^{d_{1}} & 0\\0 &e^{d_{2}}\epm ,\infty ) \ \ \text{ and } \ \ ( 0,g_{2}X,\Id,\bpm e^{d_{1}} & 0\\0 & e^{d_{2}}\epm ,\infty)
$$
are equivalent in $\mathcal{Q}^{st}/_{\mathbb{Z}/2\mathbb{Z}}$. The third parameter $[S] \in \PSO(2)/_{\sim}$ is given by the diagonalization matrix and has the geometric interpretation of an angle: given $S \in \PSO(2)$ we write it as a matrix $
S=\begin{pmatrix}
\cos\frac{\alpha}{2} & -\sin\frac{\alpha}{2}\\
\sin\frac{\alpha}{2} & \cos\frac{\alpha}{2}
\end{pmatrix},\alpha \in [0,2\pi)
$ and the equivalence relation is the identification of angle $\alpha$ with angle $(2\pi-\alpha)$ (see Figure \ref{angle2}). We obtain parameters
$(\underline{c},\underline{d},[\alpha])$ in $\mathfrak{a}^{2}\times [0,2 \pi)/_{\sim}$. The parameter $S$ or $\alpha$ will be called the \emph{angle parameter} of the generic quintuple and provides information about the angle between the hyperbolic components of $X$ and $Z$. We will draw the angle on the left as explained in Remark \ref{left}.

\noindent For the inverse map, to any element of the parameter space $\mathfrak{a}^{2} \times \PSO(2) / _{\sim}$ we can associate a unique quintuple inside $\mathcal{Q}^{st}/_{\mathbb{Z}/2\mathbb{Z}}$. To $\big( (c_{1},c_{2}),(d_{1},d_{2}),[S] \big)$ in $\mathfrak{a}^{2} \times \PSO(2) / _{\sim}$ we associate the standard quintuple 
$$
(P,X,Y,Z,Q)=\Big(0, S^{T}\bpm
\frac{1}{e^{c_{2}}} & 0\\
0 & \frac{1}{e^{c_{1}}}
\epm S,\ \Id, \ \bpm e^{d_{1}} & 0\\0 & e^{d_{2}}\epm, \infty \Big)
$$
Then $X$ is a matrix such that $d^{\mathfrak{a}^{+}}(iX,i\Id)=(c_{1},c_{2})$ and $SXS^{T}=\begin{pmatrix}
\frac{1}{e^{c_{2}}} & 0\\
0 & \frac{1}{e^{c_{1}}}
\end{pmatrix}$. For any $S' \sim S$ we obtain an equivalent quintuple $(P,X',Y,Z,Q)$ inside $\mathcal{Q}^{st}/_{\mathbb{Z}/2\mathbb{Z}}$ where
$$
(P,X',Y,Z,Q)=(P,rXr,Y,Z,Q)=\begin{pmatrix}
r & 0\\
0 & r
\end{pmatrix} \cdot (P,X,Y,Z,Q)
$$

\endproof

\begin{cor}
 The set $\mathcal{Q}^{gen}/_{\PSp(4,\R)}$ is parametrized by $
\mathfrak{a}^{2}\times [0,\pi]$.   
\end{cor}

\proof  
The equivalence relation of Proposition \ref{partrip} is given by $\alpha \sim 2 \pi - \alpha$. We can always choose $\alpha \in [0,\pi]$ as representative of the equivalence class.
\endproof

\noindent To conclude this section we state two technical lemmas that will be useful later.

\begin{lem} \label{teclem1}
Let $p=\big( (c_{1},c_{2}),(d_{1},d_{2}),[S] \big) \in \mathfrak{a}^{2} \times \PSO(2) / _{\sim}$ and let $X,Y$ be positive definite such that $d^{\mathfrak{a}^{+}}(iX,iY)=(c_{1},c_{2})$. Then the unique $Z$ such that $(0,X,Y,Z,\infty)$ corresponds to $p$ in the parametrization of Proposition \ref{partrip} is given by
$$
Z=\sqrt{Y}R^{T}S\begin{pmatrix}
e^{d_{1}} & 0\\
0 & e^{d_{2}}
\end{pmatrix}S^{T}R\sqrt{Y}
$$
where  $R$ is the unique matrix in $\PSO(2)$ such that $
R(\sqrt{Y}^{-1}X\sqrt{Y}^{-1})R^{T}=\begin{pmatrix}
\frac{1}{e^{c_{2}}} & 0\\
0 & \frac{1}{e^{c_{1}}}
\end{pmatrix}$.
\end{lem}

\proof
It is easy to check that for such a $Z$ it holds
$d^{\mathfrak{a}^{+}}(iY,iZ)=(d_{1},d_{2})$. By Proposition \ref{partrip} we know that  $SgXS^{T}=\bpm \lambda_{1}&0\\0&\lambda_{2}\epm$ where $\lambda_{1}>\lambda_{2}$ and where $g$ is such that $g(0,X,Y,Z,\infty)$ is a standard quintuple i.e.
$gY=\Id$ and $gZ=\bpm e^{d_{1}}&0\\0&e^{d_{2}}\epm$. Then $g=\bpm
A & 0\\0 & A^{-T}\epm$ where $A=S^{T}R\sqrt{Y}^{-1}$ and $gX=S^{T}R\sqrt{Y}^{-1}X\sqrt{Y}^{-1}SR^{T}$. It holds
$$
S (gX) S^{T}=R(\sqrt{Y}^{-1}X\sqrt{Y}^{-1})R^{T}=\bpm\frac{1}{e^{c_{2}}} & 0\\0 & \frac{1}{e^{c_{1}}}\epm
$$
To finish the proof we need to check that for $S' \sim S$ we obtain the same point in $\mathcal{Q}^{gen} / _{\PSp(4,\R)}$. Take $S'=rSr$, where $r=\bpm-1 & 0\\0 & 1\epm$. Then $
Z'=\sqrt{Y}R^{T}S'\begin{pmatrix}
e^{d_{1}} & 0\\
0 & e^{d_{2}}
\end{pmatrix}S'^{T}R\sqrt{Y}$.
Consider
$$
h=\begin{pmatrix}
\sqrt{Y}R^{T}S'rS^{T}R\sqrt{Y}^{-1} & 0\\
0 & \sqrt{Y}^{-1}R^{T}S'rS^{T}R\sqrt{Y}
\end{pmatrix}
$$
Then $h(0,X,Y,Z,\infty)=(0,X,Y,Z',\infty)$ so that $[(0,X,Y,Z,\infty)]=[(0,X,Y,Z',\infty)]$ in $\mathcal{Q}^{gen}/ _{\PSp(4,\R)}$. Geometrically the map $h$ is a reflection in the $\H^{2}$-component across the geodesic passing through $\pi^{\H^{2}}(iX)$ and $\pi^{\H^{2}}(iY)$, denoted for simplicity by $X$ and $Y$ respectively. This is shown in Figure \ref{hmap} below.
\endproof

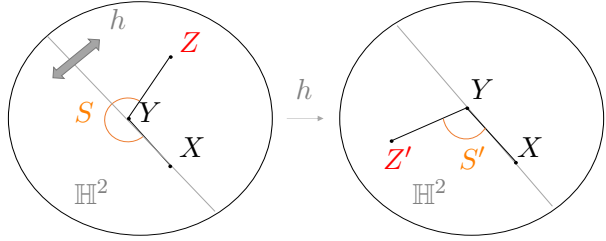
\begin{figure}[!h]
   \centering
   \captionsetup{justification=centering,margin=2cm}
   \setlength{\unitlength}{0.1\textwidth}
   \begin{picture}(3.5,2)
     \put(-1,0){\includegraphics[width=8cm,height=3.2cm]{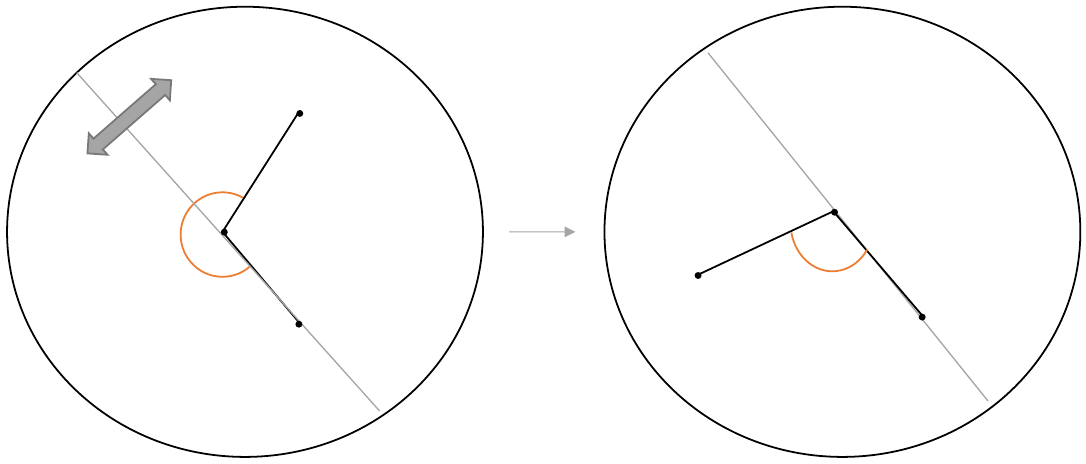}}
          \put(0.15,1){$Y$}
          \put(-0.1,1.8){\textcolor{gray}{$h$}}
          \put(0.5,0.7){$X$}
          \put(1.5,1.2){\textcolor{gray}{$h$}}
          \put(0.5,1.6){\textcolor{red}{$Z$}}
          \put(-0.4,1){\textcolor{orange}{$S$}}
          \put(2.9,0.6){\textcolor{orange}{$S'$}}
          \put(2.25,0.6){\textcolor{red}{$Z'$}}
          \put(3.4,0.7){$X$}
          \put(3,1.2){$Y$}
          \put(-0.4,0.3){\textcolor{gray}{$\H^{2}$}}
          \put(2.5,0.3){\textcolor{gray}{$\H^{2}$}}
          \end{picture}  
\caption{The map $h$ in the Poincaré disk model of $\H^{2}$}
   \label{hmap}
\end{figure}

\begin{lem} \label{teclem2}
Let $p=\big( (c_{1},c_{2}),(d_{1},d_{2}),[S] \big)$ be inside $\mathfrak{a}^{2} \times \PSO(2) / _{\sim}$ and let $Y,Z$ be positive definite such that $d^{\mathfrak{a}^{+}}(iY,iZ)=(d_{1},d_{2})$. Then the unique $X$ such that $(0,X,Y,Z,\infty)$ corresponds to $p$ in the parametrization of Proposition \ref{partrip} is given by
$$
X=\sqrt{Y}P^{T}S^{T}\begin{pmatrix}
\frac{1}{e^{c_{2}}} & 0\\
0 & \frac{1}{e^{c_{1}}}
\end{pmatrix}SP\sqrt{Y}
$$
where  $P$ is the unique matrix in $\PSO(2)$ such that $
P(\sqrt{Y}^{-1}Z\sqrt{Y}^{-1})P^{T}= \begin{pmatrix}
e^{d_{1}} & 0\\
0 & e^{d_{2}}
\end{pmatrix}
$.
\end{lem}

\proof
It is easy to check that for such $X$ it holds
$d^{\mathfrak{a}^{+}}(iX,iY)=(c_{1},c_{2})$. Take $
g=\begin{pmatrix}
P\sqrt{Y}^{-1} & 0\\
0 & P\sqrt{Y}
\end{pmatrix}
$ then $gY = \Id$, $gZ =\begin{pmatrix}
e^{d_{1}} & 0\\
0 & e^{d_{2}}
\end{pmatrix}$ and $gX=P\sqrt{Y}^{-1}X\sqrt{Y}^{-1}P^{T}=S^{T}\begin{pmatrix}
\frac{1}{e^{c_{2}}} & 0\\
0 & \frac{1}{e^{c_{1}}}
\end{pmatrix}S$, so that 
$$
S (g\cdot X) S^{T}=\bpm\frac{1}{e^{c_{2}}} & 0\\0 & \frac{1}{e^{c_{1}}}\epm
$$
To finish the proof we need to check that for $S' \sim S$ we obtain the same point in $\mathcal{Q}^{gen} / _{\PSp(4,\R)}$. Take $S'=rSr$ where $r=\bpm-1 & 0\\0 & 1\epm$. Then 
$$
X'=\sqrt{Y}P^{T}S'^{T}\begin{pmatrix}
\frac{1}{e^{c_{2}}} & 0\\
0 & \frac{1}{e^{c_{1}}}
\end{pmatrix}S'P\sqrt{Y}=\sqrt{Y}P^{T}rS^{T}\begin{pmatrix}
\frac{1}{e^{c_{2}}} & 0\\
0 & \frac{1}{e^{c_{1}}}
\end{pmatrix}SrP\sqrt{Y}
$$
Consider
$$
h=\begin{pmatrix}
\sqrt{Y}P^{T}rP\sqrt{Y}^{-1} & 0\\
0 & \sqrt{Y}^{-1}P^{T}rP\sqrt{Y}
\end{pmatrix}
$$
Then $h(0,X,Y,Z,\infty)=(0,X',Y,Z,\infty)$
so that $[(0,X,Y,Z,\infty)]=[(0,X',Y,Z,\infty)]$ in $\mathcal{Q}^{gen}/_{\PSp(4,R)}$. Geometrically the map $h$ can be seen as a reflection in the $\H^{2}$-component across the geodesic passing through $Y$ and $Z$ (similar to Figure \ref{hmap}).
\endproof

\section{Parameters for right-angled hexagons} \label{chapterhexagons}

In this section we define ordered right-angled hexagons in the Siegel space $\X$.  We distinguish between generic and non-generic hexagons and introduce a parameter space for both. A parameter space which encloses both cases is given in Theorem \ref{cpctarccoordthm}.

\subsection{Definition of hexagon, the sets $\mathcal{H},\mathcal{H}^{gen}$ and $\mathcal{H}^{st}$} \label{paramhex}

\begin{definition} \label{defhex}
A \emph{right-angled hexagon} in $\X$ is a cyclic sequence of six $\R$-tubes \\$H=[\Y_{1},\Y_{2},\Y_{3},\Y_{4},\Y_{5},\Y_{6}]$ where any two consecutive tubes are orthogonal and such that
$$\Y_{1}=\Y_{P_{1},P_{2}}, \ \Y_{2}= \Y_{Q_{1},Q_{2}}, \ \Y_{3}=  \Y_{P_{3},P_{4}} , \ \Y_{4}= \Y_{Q_{3},Q_{4}} , \ \Y_{5}= \Y_{P_{5},P_{6}},\ \Y_{6}= \Y_{Q_{5},Q_{6}}
$$
 for a maximal 12-tuple $(P_{1},Q_{6},Q_{1},P_{2},P_{3},Q_{2},Q_{3},P_{4},P_{5},Q_{4},Q_{5},P_{6})$.
\end{definition}

\begin{figure}[!h]
   \centering
   \captionsetup{justification=centering,margin=2cm}
   \setlength{\unitlength}{0.1\textwidth}
   \begin{picture}(2.6,2.7)
     \put(0,0){\includegraphics[width=4cm,height=4cm]{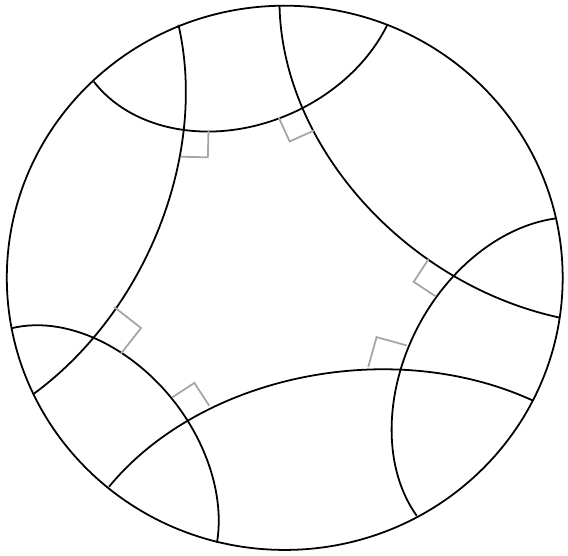}}
          \put(0.8,-0.1){$Q_{2}$}
          \put(0.3,0.1){$P_{3}$}
          \put(-0.2,0.7){$P_{2}$}
          \put(-0.3,1.1){$Q_{1}$}
          \put(0.2,2.3){$Q_{6}$}
          \put(0.7,2.6){$P_{1}$}
          \put(1.2,2.6){$P_{6}$}
          \put(1.8,2.5){$Q_{5}$}
          \put(2.6,1.5){$Q_{4}$}
          \put(2.6,1.1){$P_{5}$}
          \put(2.4,0.6){$P_{4}$}
          \put(1.8,0){$Q_{3}$}
          \put(1.1,1.3){$H$}
          \put(0.4,1.4){$\Y_{1}$}
          \put(0.45,0.65){$\Y_{2}$}
          \put(1.2,0.6){$\Y_{3}$}
          \put(1.9,1){$\Y_{4}$}
          \put(1.7,1.7){$\Y_{5}$}
          \put(1,2.1){$\Y_{6}$}
 \end{picture}
\caption{The maximal 12-tuple determining the right-angled hexagon $H=[\Y_{1},\Y_{2},\Y_{3},\Y_{4},\Y_{5},\Y_{6}]$}
   \label{max12tuple}
\end{figure}

The maximal 12-tuple determining a right hexagon $H$ in $\X$ is illustrated in Figure \ref{max12tuple}.

\begin{definition} \label{defstabhex}
 Let $H=[\Y_{1},\Y_{2},\Y_{3},\Y_{4},\Y_{5},\Y_{6}]$ be a right-angled hexagon in $\X$. We define the \emph{stabilizer of $H$} and denote it by $\Stab(H)$ the stabilizer
    $$\Stab(H)=\big\{g \in \PSp(2n,\R)| \ g\cdot \Y_{i}=\Y_{i}, \ i\in \{1,...,6\} \big\}$$   
\end{definition}

\begin{definition} \label{defordhex}
The set $\mathcal{H}$ of \emph{ordered right-angled hexagons} in $\X$ is defined by 
$$\mathcal{H}:=\{(H,\Y_{1})| \  H=[\Y_{1},\Y_{2},\Y_{3},\Y_{4},\Y_{5},\Y_{6}] \text{ right-angled hexagon }\}$$
\end{definition}

We want to be able to determine a point $(H,\Y_{1})$ inside $ \mathcal{H}$ by giving the data of an ordered maximal 6-tuple. There are many ways to do this, as explained in the following lemma.

\begin{lem} \label{6tupledethex}
Let $H=[\Y_{1},\Y_{2},\Y_{3},\Y_{4},\Y_{5},\Y_{6}]$ be a right-angled hexagon with associated maximal 12-tuple $(P_{1},Q_{6},Q_{1},P_{2},P_{3},Q_{2},Q_{3},P_{4},P_{5},Q_{4},Q_{5},P_{6})$.
Then $(H,\Y_{1}) \in \mathcal{H}$ is uniquely determined by the following ordered maximal 6-tuples:
\begin{gather} 
(P_{1},P_{2},P_{3},P_{4},P_{5},P_{6}) \label{6tuple1}
\\ (Q_{1},Q_{2},Q_{3},Q_{4},Q_{5},Q_{6})\label{6tuple2}
\\(P_{2},Q_{2},P_{4},P_{5},Q_{5},P_{1})\label{used6tuple}
\end{gather}
\end{lem}

\proof
Given the maximal 6-tuple $(P_{1},P_{2},P_{3},P_{4},P_{5},P_{6})$ we use Lemma \ref{comput2} to uniquely determine $Q_{1},Q_{2},Q_{3},Q_{4},Q_{5},Q_{6}$ such that
\begin{equation} \label{orthogonaltubes}
\Y_{P_{1},P_{2}}\perp \Y_{Q_{1},Q_{2}}\perp \Y_{P_{3},P_{4}}\perp\Y_{Q_{3},Q_{4}}\perp\Y_{P_{5},P_{6}}\perp\Y_{Q_{6},Q_{5}}\perp\Y_{P_{1},P_{2}}
\end{equation}
The hexagon $H$ is determined by the $\R$-tubes in (\ref{orthogonaltubes}) and we put $\Y_{1}=\Y_{P_{1},P_{2}}$. The quadruples  $(P_{1},Q_{6},Q_{1},P_{2}),(P_{3},Q_{2},Q_{3},P_{4}),(P_{5},Q_{4},Q_{5},P_{6})$ are maximal by Lemma \ref{quadrismax} and we obtain  a maximal 12-tuple $(P_{1},Q_{6},Q_{1},P_{2},P_{3},Q_{2},Q_{3},P_{4},P_{5},Q_{4},Q_{5},P_{6})$ by Lemma \ref{maxlemma2}. The proof for the 6-tuple in (\ref{6tuple2}) is similar and we put again $\Y_{1}=\Y_{P_{1},P_{2}}$ where the 6-tuple $(P_{1},P_{2},P_{3},P_{4},P_{5},P_{6})$ is uniquely determined by  the orthogonality conditions in (\ref{orthogonaltubes}). Given $(P_{2},Q_{2},P_{4},P_{5},Q_{5},P_{1})$ maximal 
we construct the hexagon $H$ as following: let $g \in \Sp(2n,\R)$ such that $(gP_{1},gP_{2})=(\infty,0)$. Let us denote
$$gQ_{2}=A, \ gP_{4}=B, \ gP_{5}=C, \ gQ_{5}=D $$
We use Lemma \ref{comput2} and Lemma \ref{comput3} to uniquely determine the right-angled hexagon $H_{0,A,B,C,D,\infty}$ as shown in in Figure \ref{hexag} below. The maximality of the 12-tuple at the boundary is again guaranteed by Lemma \ref{quadrismax}. We put $H=g^{-1}(H_{0,A,B,C,D,\infty})$ and $\Y_{1}=\Y_{P_{1},P_{2}}$.
\endproof

\begin{figure}[!h]
   \centering
   \captionsetup{justification=centering,margin=2cm}
   \setlength{\unitlength}{0.1\textwidth}
   \begin{picture}(2.6,2.7)
     \put(0,0){\includegraphics[width=4.2cm,height=4cm]{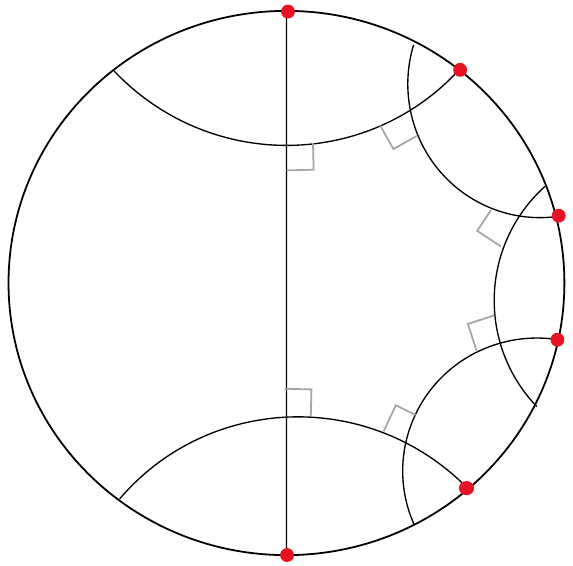}}
          \put(1.3,-0.16){\textcolor{red}{$0$}}
          \put(1.25,2.6){\textcolor{red}{$\infty$}}
          \put(2.7,1){\textcolor{red}{$B$}}
          \put(2.5,0.6){$Z_{1}$}
          \put(2.7,1.5){\textcolor{red}{$C$}}
          \put(2.6,1.8){$Z_{2}$}
          \put(1.9,2.5){$DC^{-1}D$}
          \put(2.2,2.3){\textcolor{red}{$D$}}
          \put(0.2,2.3){$-D$}
          \put(1.8,-0.05){$AB^{-1}A$}
          \put(2.2,0.2){\textcolor{red}{$A$}}
          \put(0.2,0.1){$-A$}
   \end{picture}
   
\caption{The right-angled hexagon $H_{0,A,B,C,D,\infty}$ is uniquely determined by the maximal 6-tuple $(0,A,B,C,D,\infty)$}
   \label{hexag}
\end{figure}
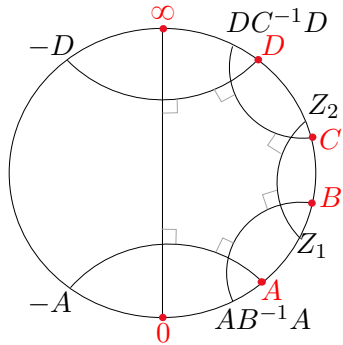

\begin{notation}
In this paper we will use a maximal 6-tuple as in (\ref{used6tuple}) to uniquely determine a right-angled hexagon $(H,\Y_{1})$ inside $\mathcal{H}$. In order to simplify the notation we will write 
$$H=(P,A,B,C,D,Q)$$
to refer to a hexagon $(H,\Y_{P,Q})$ where $H$ is uniquely determined by $(P,A,B,C,D,Q)$ as in (\ref{used6tuple}). The choice of the tube $\Y_{P,Q}$ is encoded in the order of the maximal 6-tuple.
When $P=0$ and $Q=\infty$ then $A,B,C,D$ are positive definite and we obtain a hexagon $(H, \Yi)$ where $H$ is shown in Figure \ref{hexag}. In particular the maximal 12-tuple associated to $H$ is given by $$H=(\infty,-D,-A,0,AB^{-1}A,A,Z_{1},B,C,Z_{2},D,DC^{-1}D)$$
where $Z_{1},Z_{2}$ are uniquely defined by requiring $\Y_{AB^{-1}A,B} \perp \Y_{Z_{1},Z_{2}} \perp \Y_{C,DC^{-1}D}$.
\end{notation}

\begin{definition} \label{defgenhexag}
The set of \emph{generic hexagons} $\mathcal{H}^{gen} \subset \mathcal{H}$ is given by ordered 6-tuples of the form
$$
\mathcal{H}^{gen}:= \{(P,A,B,C,D,Q) \text{ maximal}| \ (P,A,B,Q),(P,B,C,Q),(P,C,D,Q) \text{ generic} \}
$$
\end{definition}

\begin{figure}[!h]
   \centering
   \captionsetup{justification=centering,margin=2cm}
   \setlength{\unitlength}{0.1\textwidth}
   \begin{picture}(2.6,2.5)
     \put(0,0){\includegraphics[width=4cm,height=4cm]{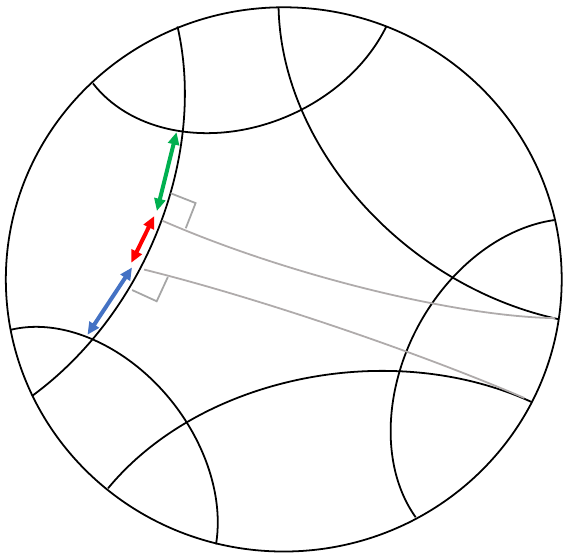}}
     \put(0.7,2.6){$Q$}
    \put(-0.1,0.6){$P$}
    \put(0.8,-0.1){$A$}
    \put(2.4,0.6){$B$}
    \put(2.6,1.1){$C$}   
    \put(1.8,2.5){$D$}
    \put(0.45,1.5){\textcolor{red}{$\underline{c}$}}
     \put(0.3,1.15){\textcolor{blue}{$\underline{b}$}}
      \put(0.6,1.8){\textcolor{forestgreen}{$\underline{d}$}}
 \end{picture}
   
\caption{The hexagon $(H,\Y_{1})$ is generic if and only if $\underline{b},\underline{c}$ and $\underline{d}$ are regular}
   \label{casesproof}
\end{figure}
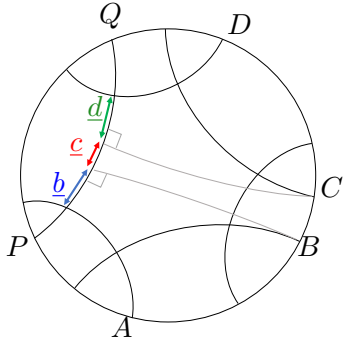

 \begin{remark} \label{equivdefgenhex}
Let $H=(P,A,B,C,D,Q)$ be a generic hexagon. Project $A,B,C$ and $D$ orthogonally on $\Y_{P,Q}$ and denote by  $\underline{b},\underline{c},\underline{d}$ the vectors $d^{\a}(p_{P,Q}(A),p_{P,Q}(B))$, $d^{\a}(p_{P,Q}(B),p_{P,Q}(C))$ and $d^{\a}(p_{P,Q}(C),p_{P,Q}(D))$ respectively (Figure \ref{casesproof}). It is easy to see that the hexagon $(H,\Y_{P,Q})$ is generic if and only if the vectors $\underline{b},\underline{c},\underline{d}$ are regular.
\end{remark}

Let us now focus on the symmetric space associated to $\Sp(4,\R)$.

\begin{definition} \label{standhex}
The set of \emph{standard hexagons} $\mathcal{H}^{st} \subset \mathcal{H}^{gen}$ is given by
   $$
\mathcal{H}^{st}:= \{(0,A,\Id,\bpm \lambda_{1}&0\\0&\lambda_{2}\epm,D,\infty) \in \mathcal{H}^{gen}| \ \lambda_{1},\lambda_{2} \in \R, \  \lambda_{1}>\lambda_{2} \}
$$ 
\end{definition}

\begin{remark} \label{iso2}
Similarly to what we have seen for quintuples for any $(H,\Y_{1})$ in $\mathcal{H}^{gen}$ we can always find an isometry $g \in \PSp(4,\R)$ such that $(gH,g\Y_{1}) \in \mathcal{H}^{st}$. For a diagonal matrix $\Lambda$ with different eigenvalues the stabilizer of $(0,\Id,\Lambda,\infty)$ is isomorphic to $\mathbb{Z}/2\mathbb{Z}$ (Remark \ref{refl}). It holds $\mathcal{H}^{gen}/_{\PSp(4,\R)} \cong \mathcal{H}^{st}/_{\mathbb{Z}/2\mathbb{Z}}$. 
\end{remark}

\subsection{Non-generic hexagons: the set $\mathcal{H}^{nongen}$} \label{secnongenhex}
We define three different types of non-generic hexagons depending on how many quadruples in $H=(P,A,B,C,D,Q)$ fail to be generic: a non-generic hexagon of type $k$ is a hexagon where $k$ quadruples are non-generic. 

\begin{definition} \label{nongenquadr} (\textbf{Non-generic quadruple}) Let $(P,X,Y,Q)$ be a maximal quadruple and let $(\mu_{1},\,\mu_{2})$ be the eigenvalues of the cross-ratio $R(P,X,Y,Q)$. The quadruple $(P,X,Y,Q)$ is said to be \emph{non-generic} if $\mu_{1}=\mu_{2}$.
\end{definition}

\begin{definition} The set $\mathcal{H}^{nongen}_{type1}$ is given by 
$$\mathcal{H}^{nongen}_{type1}:=\mathcal{H}^{nongen}_{type1.1}\cup\mathcal{H}^{nongen}_{type1.2}\cup\mathcal{H}^{nongen}_{type1.3}$$
where
$$\mathcal{H}^{nongen}_{type1.1}:= \{(P,A,B,C,D,Q) \text{ maximal}| \ (P,A,B,Q) \text{ non-generic},(P,B,C,Q),(P,C,D,Q) \text{ generic} \}$$
$$\mathcal{H}^{nongen}_{type1.2}:= \{(P,A,B,C,D,Q) \text{ maximal}| \  (P,B,C,Q) \text{ non-generic},(P,A,B,Q),(P,C,D,Q) \text{ generic} \}$$
$$\mathcal{H}^{nongen}_{type1.3}:= \{(P,A,B,C,D,Q) \text{ maximal}| \  (P,C,D,Q)\text{ non-generic},(P,A,B,Q) ,(P,B,C,Q)\text{ generic} \}$$
 \end{definition}

 \begin{definition} The set $\mathcal{H}^{nongen}_{type2}$ is given by 
$$\mathcal{H}^{nongen}_{type2}:=\mathcal{H}^{nongen}_{type2.1}\cup\mathcal{H}^{nongen}_{type2.2}\cup\mathcal{H}^{nongen}_{type2.3}$$
where
$$\mathcal{H}^{nongen}_{type2.1}:= \{(P,A,B,C,D,Q) \text{ maximal}| \ (P,A,B,Q),(P,B,C,Q) \text{ non-generic},(P,C,D,Q) \text{ generic} \}$$
$$\mathcal{H}^{nongen}_{type2.2}:= \{(P,A,B,C,D,Q) \text{ maximal}| \  (P,A,B,Q),(P,C,D,Q) \text{ non-generic},(P,B,C,Q) \text{ generic} \}$$
$$\mathcal{H}^{nongen}_{type2.3}:= \{(P,A,B,C,D,Q) \text{ maximal}| \ (P,B,C,Q), (P,C,D,Q)\text{ non-generic},(P,A,B,Q) \text{ generic} \}$$
 \end{definition}

\begin{definition} The set $\mathcal{H}^{nongen}_{type3}$ is given by
$$\mathcal{H}^{nongen}_{type3}:= \{(P,A,B,C,D,Q) \text{ maximal}| \ (P,A,B,Q),(P,B,C,Q),(P,C,D,Q) \text{ non-generic} \}$$
\end{definition}

\begin{prop}
$$\mathcal{H}=\mathcal{H}^{gen}\cup\mathcal{H}^{nongen}_{type1}\cup\mathcal{H}^{nongen}_{type2}\cup\mathcal{H}^{nongen}_{type3}$$
\end{prop}

\proof The inclusion $\mathcal{H}^{gen}\cup\mathcal{H}^{nongen}_{type1}\cup\mathcal{H}^{nongen}_{type2}\cup\mathcal{H}^{nongen}_{type3} \subset \mathcal{H}$ is trivial.\\
Let $(H,\Y_{1}) \in \mathcal{H}$. By Lemma \ref{6tupledethex} we can uniquely determine $(H,\Y_{1})$ from a maximal 6-tuple $(H,\Y_{1})=(P,A,B,C,D,Q)$ where $\Y_{1}=\Y_{P,Q}$. Let $\underline{b},\underline{c},\underline{d}$ be the three vectors in Figure \ref{casesproof}. The hexagon $(H,\Y_{1})$ is generic if and only if the vectors $\underline{b},\underline{c},\underline{d}$ are regular and is non-generic if one of them is inside $\mathfrak{d}$. For each vector $\underline{b},\underline{c}$ and $\underline{d}$ we see if it is generic or not and we list all possible configurations. The hexagon $(H,\Y_{1})$ must be contained in one of this exhaustive list.
We obtain $2^{3}=8$ possible configurations, one in $\mathcal{H}^{gen}$, three in $\mathcal{H}^{nongen}_{type1}$, three in $\mathcal{H}^{nongen}_{type2}$ and one in $\mathcal{H}^{nongen}_{type3}$.
\endproof

\subsection{Arc coordinates for generic hexagons} \label{secgen}

In this section we parametrize generic hexagons in $\X$ up to isometry. We will concentrate on the case where $\X$ is the symmetric space associated to $\Sp(4,\R)$. Recall that $\mathfrak{a}$ denotes the set of regular vectors inside $\a$.

\begin{prop} \label{parforhex}
The set $\mathcal{H}^{gen}/_{\PSp(4,\R)}$ is parametrized by $\mathfrak{a}^{3} \times  \PSO(2)^{2} / _{\sim}$ where for $(S_{1},S_{2}) \in \PSO(2)^{2}$ and $r=\bpm
-1 & 0\\
0 & 1
\epm$ the equivalent relation is $
(S_{1},S_{2}) \sim (r S_{1}r, r S_{2}r)$.

\noindent The parametrization is given by
$$
(\underline{b},\underline{c},\underline{d},[S_{1},S_{2}]) \mapsto [(0,A,\Id,C,D,\infty)]\in \mathcal{H}^{st}/_{\mathbb{Z}/2\mathbb{Z}}
$$
where $\underline{b}=(b_{1},b_{2}), \ \underline{c}=(c_{1},c_{2}), \ \underline{d}=(d_{1},d_{2})$ and
$$
A= S_{1}^{T}\bpm
\frac{1}{e^{b_{2}}} & 0\\
0 & \frac{1}{e^{b_{1}}}
\epm S_{1} \\
$$
$$
C= \bpm
e^{c_{1}} & 0\\
0 & e^{c_{2}}
\epm \\
$$
$$
D=\bpm
0 & \sqrt{e^{c_{1}}}\\
-\sqrt{e^{c_{2}}} & 0
\epm S_{2}\bpm e^{d_{1}} & 0\\
0 & e^{d_{2}}
\epm S_{2}^{T} \bpm
0 & -\sqrt{e^{c_{2}}}\\
\sqrt{e^{c_{1}}} & 0
\epm
$$

\noindent The parameter space can be rewritten as $\mathfrak{a}^{3} \times [0,2\pi)^{2} / _{\sim}$ where for $(\alpha_{1},\alpha_{2}) \in [0,2\pi)^{2}$ the corresponding $\PSO(2)^{2}$-parameter is $(S_{1},S_{2})$ with $S_{i}=\bpm
\cos\frac{\alpha_{i}}{2} & -\sin\frac{\alpha_{i}}{2}\\
\sin\frac{\alpha_{i}}{2} & \cos\frac{\alpha_{i}}{2}
\epm$. The equivalence relation is $(\alpha_{1},\alpha_{1})\sim (2\pi- \alpha_{1},2\pi- \alpha_{2})$.

\end{prop}

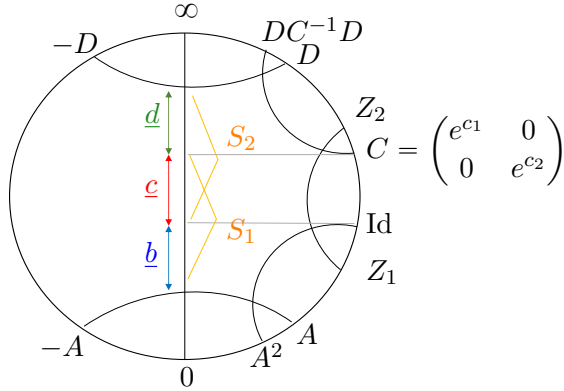
\begin{figure}[!h]
   \centering
   \captionsetup{justification=centering,margin=2cm}
   \setlength{\unitlength}{0.1\textwidth}
   \begin{picture}(5,2.95)
     \put(0.8,0){\includegraphics[width=4.8cm,height=4.5cm]{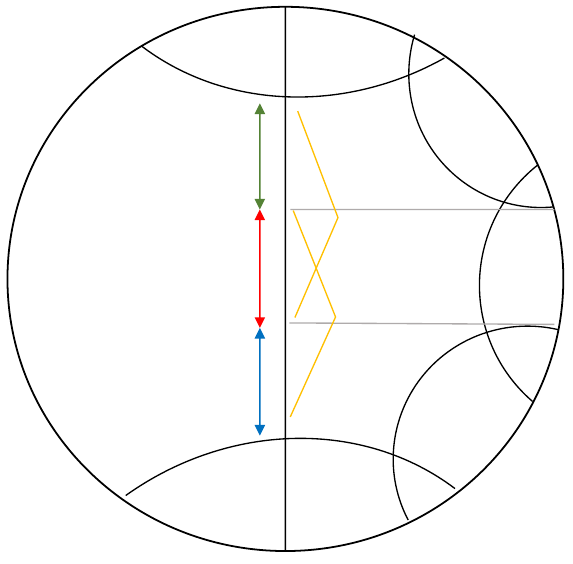}}
          \put(2.3,-0.16){$0$}
          \put(2.25,3){$\infty$}
          \put(3.9,1.15){$\Id$}
          \put(3.9,0.75){$Z_{1}$}
          \put(3.9,1.8){$C=\bpm e^{c_{1}} & 0\\0 & e^{c_{2}}\epm$}
          \put(3.8,2.15){$Z_{2}$}
          \put(3,2.8){$DC^{-1}D$}
          \put(3.3,2.6){$D$}
          \put(1.2,2.7){$-D$}
          \put(2.9,0){$A^{2}$}
          \put(3.3,0.2){$A$}
          \put(1.1,0.1){$-A$}
          \put(2.7,1.1){\textcolor{orange}{$S_{1}$}}
          \put(2.7,1.9){\textcolor{orange}{$S_{2}$}}
          \put(2,0.9){\textcolor{blue}{$\underline{b}$}}
          \put(2,1.5){\textcolor{red}{$\underline{c}$}}
          \put(2,2.1){\textcolor{forestgreen}{$\underline{d}$}}
          \end{picture}
\caption{The standard right-angled hexagon $(H^{st},\Yi)$ with parameters $(\underline{b},\underline{c},\underline{d},S_{1},S_{2})$}
   \label{Parhex}
\end{figure}

\proof 
We first show how to find parameters $\big((b_{1},b_{2}),(c_{1},c_{2}),(d_{1},d_{2}),[S_{1},S_{2}]\big)$ for a given $(H,\Y_{1})$ inside $\mathcal{H}^{gen}$.
Let $(H,\Y_{1})=(P,A,B,C,D,Q)$, $\Y_{1}=\Y_{P,Q}$. Up to isometry we can consider $P=0$ and $Q=\infty$. As $(H,\Yi) \in \mathcal{H}^{gen}$ the quintuples $(0,A,B,C,\infty)$ and $(0,B,C,D,\infty)$ both belong to $\mathcal{Q}^{gen}$.  We use Proposition \ref{partrip} to find parameters $\underline{b},\underline{c},\underline{d}$ in $\mathfrak{a}^{3}$: 
\begin{gather*}
(b_{1},b_{2})=d^{\mathfrak{a}^{+}}(iA,iB)\\  (c_{1},c_{2})=d^{\mathfrak{a}^{+}}(iB,iC)\\  
(d_{1},d_{2})=d^{\mathfrak{a}^{+}}(iC,iD)
\end{gather*}

\begin{figure}[!h]
   \centering
   \captionsetup{justification=centering,margin=2cm}
   \setlength{\unitlength}{0.1\textwidth}
   \begin{picture}(3.5,2.3)
     \put(-1,0){\includegraphics[width=8.2cm,height=3.5cm]{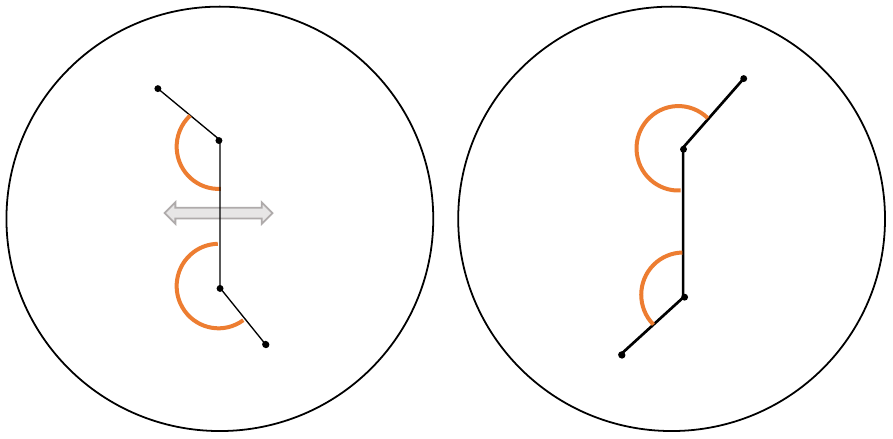}}
\put(-0.2,0.5){\textcolor{orange}{$S_{1}$}}
\put(-0.25,1.5){\textcolor{orange}{$S_{2}$}}
\put(0.35,0.8){$\Id$}
\put(0.3,1.6){$C$}
\put(0.55,0.3){$A$}  
\put(-0.3,1.8){$D$}
\put(-0.8,1.1){\textcolor{gray}{$\H^{2}$}}
\put(2.5,0.7){\textcolor{orange}{$S'_{1}$}}
\put(2.5,1.5){\textcolor{orange}{$S'_{2}$}}
\put(3.1,0.7){$\Id$}
\put(3.1,1.4){$C$}
\put(2.4,0.3){$A$}  
\put(3.2,1.9){$D$}
\put(1.9,1.1){\textcolor{gray}{$\H^{2}$}}
 \end{picture}
\caption{Visualization of the equivalence relation $[S_{1},S_{2}]=[S'_{1},S'_{2}]$ in the Poincaré disk model of $\H^{2}$}
   \label{anglesH2}
\end{figure}
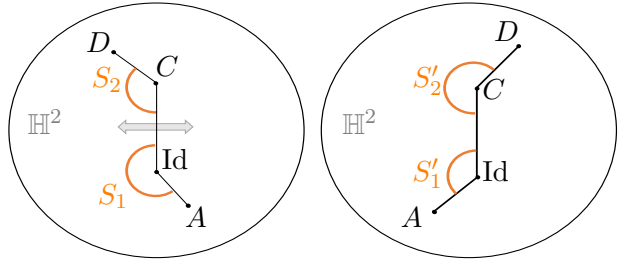

\noindent Let $g \in \Stab(\Yi)$ be such that $gB=\Id$ and $gC= \bpm e^{c_{1}}&0\\0&e^{c_{2}}\epm$. We obtain exactly two possibilities $g_{1},g_{2}$ sending $(H,\Yi)$ to the standard hexagons 
$$(H_{1},\Yi)=(0,g_{1}A,\Id,C,g_{1}D,\infty) \text{ and } (H_{2},\Yi)=(0,g_{2}A,\Id,C,g_{2}D,\infty)$$
It holds $g_{2}A= rg_{1}Ar$, $g_{2}D= rg_{1}Dr$, where $r=\bpm-1 & 0\\0&1\epm$. The points $rg_{1}Ar$ and $rg_{1}Dr$ are the image under a reflection in the hyperbolic component of $\Yi$. Let $(\underline{b},\underline{c},[S_{1}])$ and $(\underline{c},\underline{d},[S_{2}])$ be the parameters associated to the quintuples $(0,g_{1}A,\Id,C,\infty)$ and  $(0,\Id,C,g_{1}D,\infty)$ respectively. We associate to $(H,\Yi)$ the point $(\underline{b},\underline{c},\underline{d},[S_{1},S_{2}])$, where $S_{i}\sim rS_{i}r$ for $i\in \{1,2\}$. Parameters $S_{1},S_{2}$ are diagonalization matrices and have geometric interpretation of an angle. For $\alpha_{i} \in [0,2\pi)$ it holds $S_{i}=\bpm
\cos\frac{\alpha_{i}}{2} & -\sin\frac{\alpha_{i}}{2}\\
\sin\frac{\alpha_{i}}{2} & \cos\frac{\alpha_{i}}{2}
\epm, \ i \in {1,2}$. The equivalence relation $S\sim rSr$ is the identification of angle $\alpha_{i}$ with angle $(2\pi-\alpha_{i})$. See Figures \ref{Parhex} and \ref{anglesH2} for a visualization of the parameters.

\noindent Given $(\underline{b},\underline{c},\underline{d},[S_{1},S_{2}])$ in $\mathfrak{a}^{3} \times \PSO(2)^{2} / _{\sim}$,  we construct a standard hexagon $(H^{st},\Yi)$ in the following way: up to $\PSp(4,\R)$-action we can consider $B=\Id$ and $C$ diagonal. It is sufficient to determine $A,C$ and $D$ for $(H^{st},\Yi)$ to be uniquely determined. The equality $(c_{1},c_{2})=d^{\mathfrak{a}^{+}}(i\Id,iC)$
forces $C=\bpm
e^{c_{1}} & 0\\
0 & e^{c_{2}}
\epm$.
We use Lemma \ref{teclem1} and Lemma \ref{teclem2} to uniquely determine $A$ and $D$ respectively.
We use Lemma \ref{comput2} and Lemma \ref{comput3} to compute the corresponding orthogonal tubes. To finish the proof, we need to check that for $(S_{1}',S_{2}') \sim (S_{1},S_{2})$ we obtain an equivalent hexagon $(H'^{st},\Yi)$ inside $\mathcal{H}^{st} / _{\mathbb{Z}/2\mathbb{Z}}$. For $i=1,2$ let $S_{i}'=rS_{i}r$. Then
$$
A'= r S_{1}^{T}r\bpm
\frac{1}{e^{b_{2}}} & 0\\
0 & \frac{1}{e^{b_{1}}}
\epm r S_{1} r=r A r 
$$

$$
C'=C
$$
$$
D'=\bpm
0 & \sqrt{e^{c_{1}}}\\
-\sqrt{e^{c_{2}}} & 0
\epm r S_{2} r \bpm e^{d_{1}} & 0\\
0 & e^{d_{2}}
\epm r S_{2}^{T} r \bpm
0 & -\sqrt{e^{c_{2}}}\\
\sqrt{e^{c_{1}}} & 0
\epm=$$
$$=\bpm
0 & \sqrt{e^{c_{1}}}\\
\sqrt{e^{c_{2}}} & 0
\epm S_{2} \bpm e^{d_{1}} & 0\\
0 & e^{d_{2}}
\epm S_{2}^{T} \bpm
0 & \sqrt{e^{c_{2}}}\\
\sqrt{e^{c_{1}}} & 0
\epm = r D r 
$$

\noindent so that $(gH^{st},g\Yi)=(H'^{st},\Yi)$ where 
$g =  \bpm
r & 0\\
0 & r 
\epm$.
\endproof

\begin{cor}
 The set $\mathcal{H}^{gen}/_{\PSp(4,\R)}$ is parametrized by $\mathfrak{a}^{3} \times [0,\pi] \times [0,2\pi)$.  
\end{cor}

\proof  
The equivalence relation of Proposition \ref{parforhex} is given by $(\alpha_{1},\alpha_{2}) \sim (2\pi - \alpha_{1},2\pi - \alpha_{2})$. We choose $\alpha_{1} \in [0,\pi]$ as representative of the equivalence class.
\endproof

\begin{definition} \label{defarccoords} The parameters of Proposition \ref{parforhex} will be called \emph{arc coordinates} for a generic right-angled hexagon $(H,\Y_{1})$. The vectors $\underline{b},\underline{c},\underline{d}$ will be called \emph{length parameters} and $\alpha_{1},\alpha_{2}$ will be called \emph{angle parameters}.
\end{definition}

\begin{remark}
The term \emph{arc coordinates} introduced in Definition \ref{defarccoords} could be misleading as we also use it for the parametrization of classical Teichmüller space and for its generalization in the case of maximal representations. Nevertheless, we have decided to keep this name also for the parameters of a hexagon as they are crucial for the construction of parameters for maximal representations and will appear in their parameter space (Theorem \ref{thesis}).
\end{remark}

\subsection{Polygonal chain associated to a right-angled hexagon} \label{secpolchain}

In this section we define the polygonal chain associated to an ordered right-angled hexagon and show how this is related to length and angle parameters. For the purposes of this paper we will define the polygonal chain of $(H,\Y_{1})$ in the case where $\Y_{1}=\Yi$. Recall that given $A \in \Sym^{+}(2,\R)$, the hyperbolic component of $A$ is the point $\pi^{\H^{2}}(p_{0,\infty}(A))$. Recall also that for two points $iA,iB$ with $d^{\a}(iA,iB)=(d_{1},d_{2})$ the hyperbolic distance $d^{\H^{2}}(\pi^{\H^{2}}(iA),\pi^{\H^{2}}(iB))=\mathbf{h}$ is given by
\begin{equation} \label{maph}
\mathbf{h}=\mathbf{h}(\underline{d})=d_{1}-d_{2}
\end{equation}

\begin{definition} \label{defpolchain}
Let $\X$ be the symmetric space associated to $\Sp(4,\R)$ and let $(H,\Yi) \in \mathcal{H}$ be an ordered right-angled hexagon in $\X$. Let $$H=(0,A,B,C,D,\infty)$$ where $(0,A,B,C,D,\infty)$ is a maximal 6-tuple.
The \emph{polygonal chain associated to} $(H,\Yi)$ is the connected series of geodesic segments with vertices given by the ordered sequence of points (possibly coinciding)
$$\Big(\pi^{\H^{2}}(iA),\pi^{\H^{2}}(iB),\pi^{\H^{2}}(iC) ,\pi^{\H^{2}}(iD) \Big)$$
The \emph{segments of the polygonal chain} are the oriented geodesic segments (possibly collapsing to one point):
$$\overrightarrow{\pi^{\H^{2}}(iA)\pi^{\H^{2}}(iB)}, \ \overrightarrow{\pi^{\H^{2}}(iB)\pi^{\H^{2}}(iC)}, \ \overrightarrow{\pi^{\H^{2}}(iC)\pi^{\H^{2}}(iD)}$$
The \emph{angles of the polygonal chain} are the angles formed by two consecutive segments (measured on the left-hand side of the oriented segments).
\end{definition}

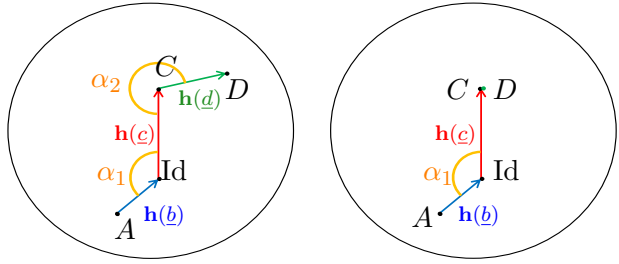
\begin{figure}[!h]
   \centering\captionsetup{justification=centering,margin=2cm}
   \setlength{\unitlength}{0.1\textwidth}
   \begin{picture}(3.5,2.3)
     \put(-1,0){\includegraphics[width=8.2cm,height=3.5cm]{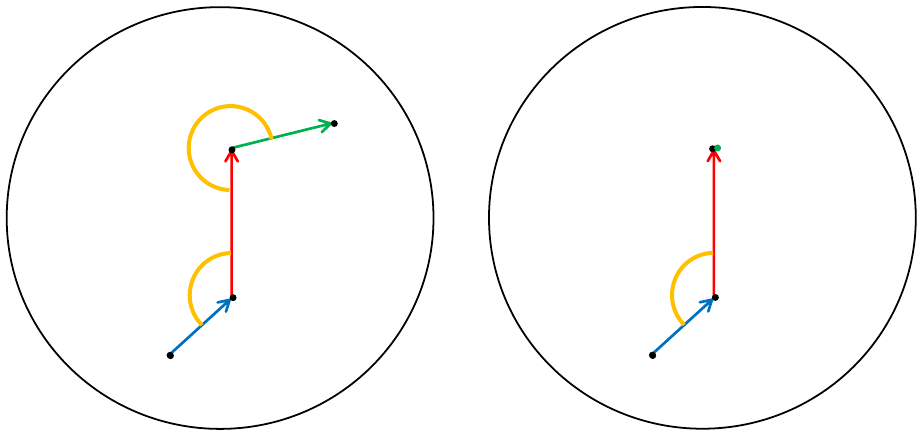}}
\put(0.2,0.4){\textcolor{blue}{$_{\mathbf{h}(\underline{b})}$}}
\put(-0.05,1.1){\textcolor{red}{$_{\mathbf{h}(\underline{c})}$}}
\put(0.5,1.4){\textcolor{forestgreen}{$_{\mathbf{h}(\underline{d})}$}}
\put(0.35,0.7){$\Id$}
\put(-0.05,0.2){$A$}
\put(0.3,1.6){$C$}
\put(0.9,1.4){$D$}
\put(2.9,0.4){\textcolor{blue}{$_{\mathbf{h}(\underline{b})}$}}
\put(2.7,1.1){\textcolor{red}{$_{\mathbf{h}(\underline{c})}$}}
\put(3.2,0.7){$\Id$}
\put(2.5,0.3){$A$}
\put(2.8,1.4){$C$}
\put(3.2,1.4){$D$}
\put(-0.2,0.7){\textcolor{orange}{$\alpha_{1}$}}
\put(-0.25,1.5){\textcolor{orange}{$\alpha_{2}$}}
\put(2.6,0.7){\textcolor{orange}{$\alpha_{1}$}}
 \end{picture}
   
\caption{Polygonal chains of a generic hexagon and of a non-generic hexagon of type 1.3}
   \label{polchain}
\end{figure}

\noindent If $(H,\Yi) \in \mathcal{H}^{gen}$ the segments of the polygonal chain of $(H,\Yi)$ have hyperbolic length given by $\mathbf{h}(\underline{b}),\mathbf{h}(\underline{c})$ and $\mathbf{h}(\underline{d})$ respectively where $\mathbf{h}$ is the map in (\ref{maph}) and $\underline{b},\underline{c},\underline{d}$ are length parameters of the arc coordinates. Up to an isometry $g \in \Sp(4,\R)$ we can consider the case where $B=\Id$ and $C$ is diagonal. Observe that to a generic hexagon $(H,\Yi)$ we can associate exactly two polygonal chains up to isometry, and these are drawn in Figure \ref{anglesH2}.
If the hexagon is non-generic some segments contract to a point. The hyperbolic length of the segment is zero as the corresponding length parameter is inside $\mathfrak{d}=\{(x_{1},x_{2})| \ x_{1}=x_{2})\}$. This will be made more clear in the next section. The polygonal chain of both a generic and a non-generic hexagon is illustrated in the Poincaré disc model in Figure \ref{polchain}. For simplicity for any $X \in \Sym^{+}(2,\R)$ we have denoted the point $\pi^{\H^{2}}(iX)$ as $X$.

\subsection{Arc coordinates for non-generic hexagons} \label{secnongen}

We denote by $\mathfrak{d}$ the following set
$$\mathfrak{d}=\{(x_{1},x_{2}) \in \R^{2}| \ x_{1}=x_{2} \}$$
The following Proposition arise as a natural generalisation of Proposition \ref{parforhex}. 

\begin{prop} \label{parnongentype1}
Non-generic hexagons of type 1 are parametrised up to isometry by
 $$\mathcal{H}^{nongen}_{type1.1}/_{\PSp(4,\R)} \ \cong \  \mathfrak{d} \times \mathfrak{a}^{2} \times [0,2\pi)/_{\sim}$$
 $$\mathcal{H}^{nongen}_{type1.2}/_{\PSp(4,\R)} \ \cong \ \mathfrak{a} \times \mathfrak{d} \times \mathfrak{a} \times [0,2\pi)/_{\sim}$$
 $$\mathcal{H}^{nongen}_{type1.3}/_{\PSp(4,\R)}  \ \cong \ \mathfrak{a}^{2} \times\mathfrak{d}  \times [0,2\pi)/_{\sim}$$

where for $\alpha \in [0,2\pi)$ the equivalence relation is given by $\alpha \sim 2 \pi- \alpha$.
 \end{prop}

\begin{figure}[!h]
   \centering\captionsetup{justification=centering,margin=2cm}
   \setlength{\unitlength}{0.1\textwidth}
   \begin{picture}(6.3,1.6)
     \put(0,0){\includegraphics[width=10.3cm,height=2.8cm]{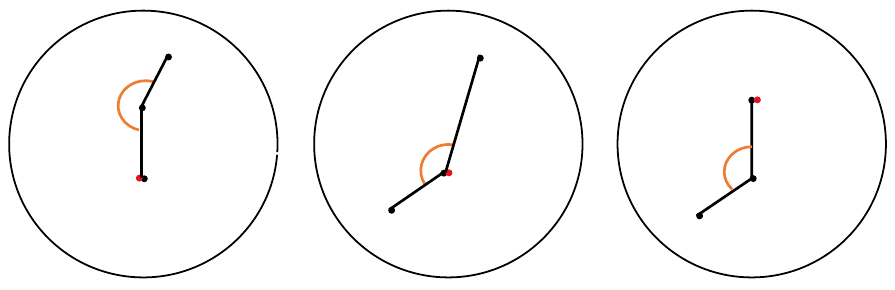}}
\put(1.1,0.5){$\Id$}
\put(0.8,0.5){\textcolor{red}{$A$}}
\put(1.1,1){$C$}
\put(1.3,1.5){$D$}
\put(3.6,1.5){$D$}
\put(3.1,0.5){$\Id$}
\put(3.4,0.5){\textcolor{red}{$C$}}
\put(2.6,0.5){$A$}
\put(5.7,1.2){\textcolor{red}{$D$}}
\put(5.3,1.2){$C$}
\put(5.6,0.6){$\Id$}
\put(5.1,0.2){$A$}
 \end{picture}
   
\caption{Polygonal chains of non-generic hexagons of type 1.1, 1.2 and 1.3 respectively}
   \label{nongen1}
\end{figure}
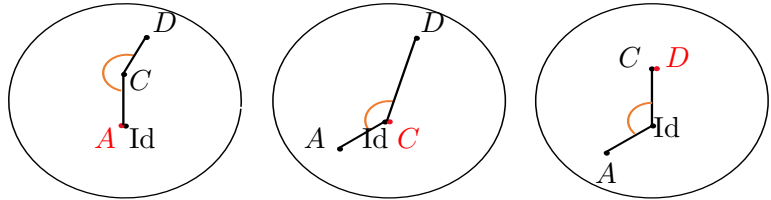

\proof Let $H=(P,A,B,C,D,Q) \in \mathcal{H}^{nongen}_{type1.1}$. Up to isometry we can assume $P=0$, $Q=\infty$, $B=\Id$ and $C$ diagonal. As $(0,A,\Id,\infty)$ non-generic we can not define an angle parameter between the hyperbolic components of $A$ and $\Id$ and the parameter $d^{\a}(iA,i\Id)$ is inside $\mathfrak{d}$. Geometrically this means that the two points coincide in the $\H^{2}$-component of $\Yi$ (Figure \ref{nongen1}). The quintuple $(0,\Id,C,D,\infty)$ is generic and we use Proposition \ref{partrip} to determine the angle parameter. Up to reflection on the hyperbolic component we can always choose the angle parameter to lie inside $[0,\pi] \cong \PSO(2)/_{\sim}$.

\noindent Conversely, given $\big( (b,b),(c_{1},c_{2}),(d_{1},d_{2}),[S] \big)  \in \mathfrak{d} \times \mathfrak{a}^{2} \times \PSO(2)/_{\sim}$ we construct the hexagon $H_{type1.1}^{nongen}=(0,A,\Id,C,D,\infty)$ where
$$
A= \bpm
\frac{1}{e^{b}} & 0\\
0 & \frac{1}{e^{b}}
\epm  \\
$$
$$
C= \bpm
e^{c_{1}} & 0\\
0 & e^{c_{2}}
\epm \\
$$
$$
D=\bpm
0 & \sqrt{e^{c_{1}}}\\
-\sqrt{e^{c_{2}}} & 0
\epm S \bpm e^{d_{1}} & 0\\
0 & e^{d_{2}}
\epm S^{T} \bpm
0 & -\sqrt{e^{c_{2}}}\\
\sqrt{e^{c_{1}}} & 0
\epm
$$

\noindent The proofs for type 1.2 and 1.3 are similar.
\endproof

\begin{prop} \label{parnongentype2}
Non-generic hexagons of type 2 are parametrized up to isometry by
 $$\mathcal{H}^{nongen}_{type2.1}/_{\PSp(4,\R)} \  \cong \ \mathfrak{d}^{2} \times \mathfrak{a}$$
 $$\mathcal{H}^{nongen}_{type2.2}/_{\PSp(4,\R)} \  \cong \  \mathfrak{d} \times \mathfrak{a} \times \mathfrak{d} $$
 $$\mathcal{H}^{nongen}_{type2.3}/_{\PSp(4,\R)} \  \cong \ \mathfrak{a} \times\mathfrak{d}^{2}  $$
\end{prop}

\begin{figure}[!h]
   \centering\captionsetup{justification=centering,margin=2cm}
   \setlength{\unitlength}{0.1\textwidth}
   \begin{picture}(6.3,1.6)
     \put(0,0){\includegraphics[width=10.3cm,height=2.8cm]{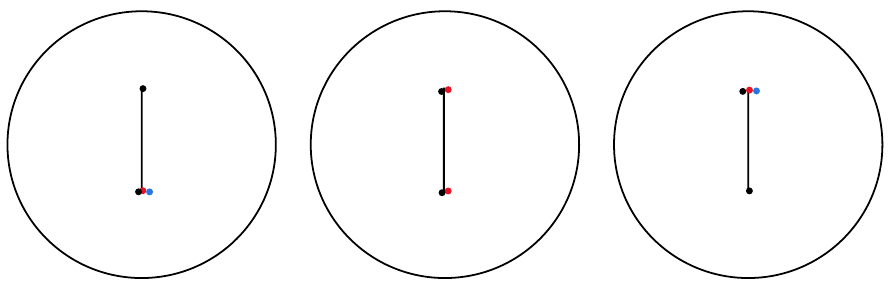}}
\put(0.8,0.3){A}
\put(1,0.3){\textcolor{red}{$\Id$}}
\put(1.2,0.3){\textcolor{blue}{$C$}}
\put(1,1.3){$D$}
\put(3.1,0.3){$A$}
\put(3.4,0.3){\textcolor{red}{$\Id$}}
\put(3.2,1.3){$C$}
\put(3.4,1.3){\textcolor{red}{$D$}}
\put(5.6,0.6){$\Id$}
 \end{picture}
\caption{Polygonal chains of three non-generic hexagons of type 2.1, 2.2 and 2.3 respectively}
   \label{nongen2}
\end{figure}
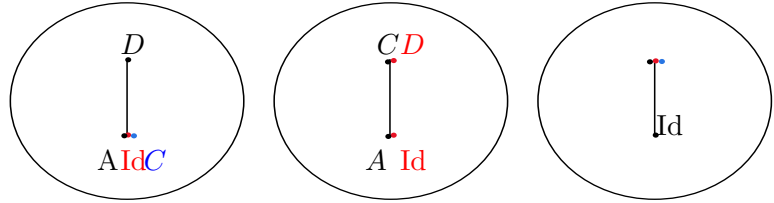

\proof This is similar to the proof of Proposition \ref{parnongentype1}. Since two quintuples are non-generic, we do not have any angle in the parameter space. Up to isometry we can move the polygonal chains of the hexagons in a configuration shown in Figure \ref{nongen2}. The vector parameters are the same of Proposition \ref{parforhex}, where two length are not regular and lie inside $\mathfrak{d}$.
\endproof

\begin{prop} \label{hexparinH2}
Non-generic hexagons of type 3 are parametrized up to isometry by
$$\mathcal{H}^{nongen}_{type3} \cong \mathfrak{d}^{3} \cong \R^{3}_{>0}$$
\end{prop}

\proof 
Let $(H,\Y_{1})=(P,A,B,C,D,Q)$ be inside $ \mathcal{H}^{nongen}_{type3}$ and up to isometry let us consider $P=0$, $Q=\infty$ and $B=\Id$. The matrices $A,C$ and $D$ are all multiples of the identity matrix. Equivalently, the points $A,\Id,C,D$ all coincide in the $\H^{2}$-component of $\Yi$ and there is no angle parameter. The length parameters are of the form $(b,b)$, $(c,c)$ and $(d,d)$. This case corresponds to hexagon-parameters in the hyperbolic case: we obtain the 3-dimensional space of right-angled hexagons of $\H^{2}$.
\endproof

\subsection{Arc coordinates for $\mathcal{H}$}

In this section we present arc coordinates in a more compact way. We introduce a parameter space for $\mathcal{H}$ which encloses both the generic and the non-generic case. Again we will focus on the case where $\X$ is the symmetric space associated to $\Sp(4,\R)$. We introduce the symbol $\overline{\mathfrak{a}}$ to denote the union $\overline{\mathfrak{a}}=\mathfrak{a}\cup \mathfrak{d}$ that is the set
$$\overline{\mathfrak{a}}=\{(x_{1},x_{2})| \ x_{1}\geq x_{2}>0\}$$

\begin{definition}
The space of \emph{decorated arc coordinates} $\mathcal{A}^{dec}_{(H,\Y_{1})}$ is given by 
$$ \mathcal{A}^{dec}_{(H,\Y_{1})}:= \overline{\mathfrak{a}}^{3} \times [0,2\pi)^{2} $$
We further define $\mathcal{A}_{(H,\Y_{1})}$ as the set
$$ \mathcal{A}_{(H,\Y_{1})}:=\mathcal{A}^{dec}_{(H,\Y_{1})} /_{\sim}$$
where the equivalence relation is given by
$(\underline{b},\underline{c},\underline{d},\alpha_{1},\alpha_{2}) \sim (\underline{b},\underline{c},\underline{d},2\pi-\alpha_{1},2\pi-\alpha_{2})$.
\end{definition}

\begin{remark}
It is straightforward to see that if $\underline{b},\underline{c},\underline{d} \in \mathfrak{a}^{3}$ then $\mathcal{A}_{(H,\Y_{1})}$ is the space of arc coordinates for generic hexagons described in Proposition \ref{parforhex}.
\end{remark}

\begin{thm}  \label{cpctarccoordthm}
Let $\X$ be the symmetric space associated to $\Sp(4,\R)$. The space $\mathcal{H}$ of ordered right-angled hexagons in $\X$ is parametrized up to isometry by $\mathcal{A}:=\mathcal{A}_{(H,\Y_{1})} / _{\sim}$
where for $(\underline{b},\underline{c},\underline{d},\alpha_{1},\alpha_{2}) \in \mathcal{A}_{(H,\Y_{1})}$ we have the following equivalent relation $\sim$ :
\begin{enumerate}
\item[(i)] If $\underline{b} \in \mathfrak{d}$:  $(\underline{b},\underline{c},\underline{d},\alpha_{1},\alpha_{2})\sim (\underline{b},\underline{c},\underline{d},\overline{\alpha}_{1},\alpha_{2}) \ \forall \overline{\alpha}_{1}$
\item[(ii)] If $\underline{c} \in \mathfrak{d}$:  $(\underline{b},\underline{c},\underline{d},\alpha_{1},\alpha_{2})\sim (\underline{b},\underline{c},\underline{d},\overline{\alpha}_{1},\overline{\alpha}_{2})$ 
for $\overline{\alpha}_{1},\overline{\alpha}_{2}$ such that   $\overline{\alpha}_{1}+\overline{\alpha}_{2}=\alpha_{1}+\alpha_{2}$

\item[(iii)]
If $\underline{d} \in \mathfrak{d}$:  $(\underline{b},\underline{c},\underline{d},\alpha_{1},\alpha_{2})\sim (\underline{b},\underline{c},\underline{d},\alpha_{1},\overline{\alpha}_{2}) \ \forall \overline{\alpha}_{2} $ 
\end{enumerate}
\end{thm}

\proof Let $(\underline{b},\underline{c},\underline{d},[\alpha_{1},\alpha_{2}]) \in \mathcal{A}$. We construct a right-angled hexagon $(H,\Yi)=(0,A,\Id,C,D,\infty)$ where $(0,A,\Id,C,D,\infty)$ is a maximal 6-tuple and $C$ is diagonal. We do it in the following way: we look at the length parameters $(\underline{b},\underline{c},\underline{d})$ which uniquely determine the genericity type of the hexagon and then we use one of Propositions \ref{parnongentype1}, \ref{parnongentype2} and \ref{hexparinH2} to construct $(H,\Yi)$. In the case of non-generic hexagons some of the angle parameters vanish and this is translated in the equivalent relations of $\mathcal{A}$ by collapsing the angle parameter in one point. More precisely:
\begin{itemize}

\item[0.] If $\underline{b},\underline{c},\underline{d} \in \mathfrak{a}^{3}$ we construct a generic hexagon with arc coordinates $(\underline{b},\underline{c},\underline{d},[\alpha_{1},\alpha_{2}])$ using Proposition \ref{parforhex}: $(H,\Yi)=(\underline{b},\underline{c},\underline{d},[\alpha_{1},\alpha_{2}])$.

\item[1.1] If $\underline{b} \in \mathfrak{d}, \ \underline{c},\underline{d} \in \mathfrak{a}^{2}$ the angle parameter $\alpha_{1}$ is collapsed into a point and we use $\alpha_{2}$ to construct a non-generic hexagon of type 1.1 using Proposition \ref{parnongentype1}: \\$(H,\Yi)=(\underline{b},\underline{c},\underline{d},[\alpha_{2}])$. 

\item[1.2] If $\underline{c} \in \mathfrak{d}, \ \underline{b},\underline{d} \in \mathfrak{a}^{2}$ we use Proposition \ref{parnongentype1} to construct a non-generic hexagon of type 1.2 where 
$(H,\Yi)=(\underline{b},\underline{c},\underline{d},[\alpha_{1}+\alpha_{2}-\pi])$. The reason why we choose to translate the angle by $\pi$ is the following: in the procedure of constructing a hexagon $(H,\Yi)$ from a point $(\underline{b},\underline{c},\underline{d},[\alpha_{1},\alpha_{2}])$ inside $\mathcal{A}$ we know that angle parameters $\alpha_{1},\alpha_{2}$ have a geometric interpretation realised in the polygonal chain of $(H,\Yi)$. If the hexagon is non-generic of type 1 then we only need one angle parameter to construct $(H,\Yi)$. In this construction procedure, when moving continuously from a point $(\underline{b},\underline{c},\underline{d},[\alpha_{1},\alpha_{2}])$ with $\underline{c} \in \mathfrak{a}$ to a point $(\underline{b},\underline{c},\underline{d},[\alpha_{1},\alpha_{2}])$ with $\underline{c} \in \mathfrak{d}$ we want the constructed hexagons to be close. To do this we need to construct $(H,\Yi)$ using the angle parameter $[\alpha_{1}+\alpha_{2}-\pi]$. This is illustrated in Figure \ref{translpi} below.

 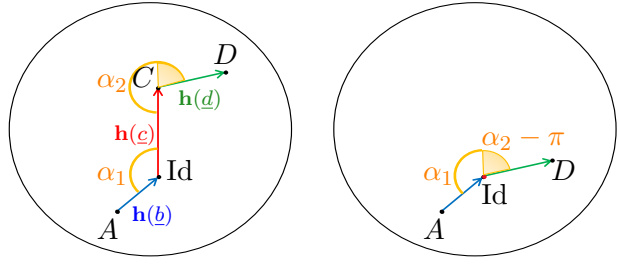
\begin{figure}[!h]
   \centering\captionsetup{justification=centering,margin=2cm}
   \setlength{\unitlength}{0.1\textwidth}
   \begin{picture}(3.5,2.3)
     \put(-1,0){\includegraphics[width=8.2cm,height=3.5cm]{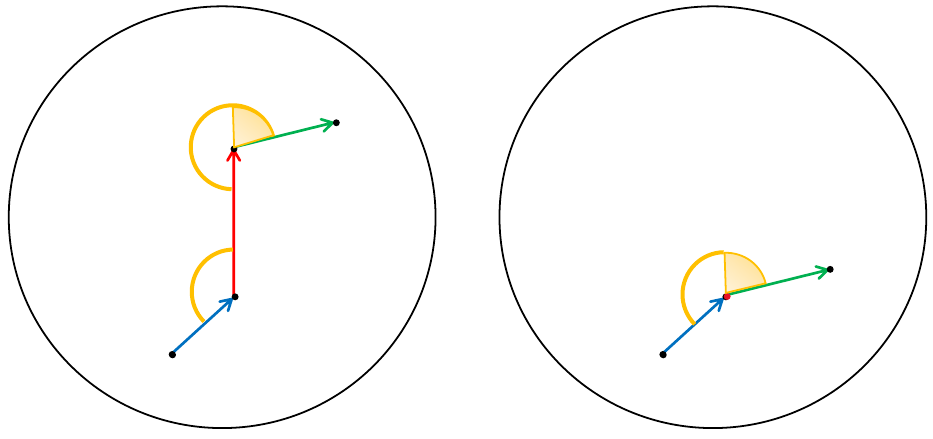}}
\put(0.1,0.4){\textcolor{blue}{$_{\mathbf{h}(\underline{b})}$}}
\put(-0.05,1.1){\textcolor{red}{$_{\mathbf{h}(\underline{c})}$}}
\put(0.5,1.4){\textcolor{forestgreen}{$_{\mathbf{h}(\underline{d})}$}}
\put(0.4,0.7){$\Id$}
\put(-0.2,0.2){$A$}
\put(0.1,1.5){$C$}
\put(0.8,1.7){$D$}
\put(3.1,0.5){$\Id$}
\put(2.6,0.2){$A$}
\put(-0.2,0.7){\textcolor{orange}{$\alpha_{1}$}}
\put(-0.2,1.5){\textcolor{orange}{$\alpha_{2}$}}
\put(3.7,0.7){$D$}
\put(2.6,0.7){\textcolor{orange}{$\alpha_{1}$}}
\put(3.1,1){\textcolor{orange}{$\alpha_{2}-\pi$}}
 \end{picture}
   
\caption{Construction of $(H,\Yi)$ of type 1.2 when $\underline{c} \to \mathfrak{d}$}
   \label{translpi}
\end{figure}

\item[1.3] If $\underline{d} \in \mathfrak{d},\  \underline{c},\underline{d} \in \mathfrak{a}^{2}$ then the angle parameter $\alpha_{2}$ is collapsed into a point and we use $\alpha_{1}$ to construct a non-generic hexagon of type 1.1 using Proposition \ref{parnongentype1}: \\$(H,\Yi)=(\underline{b},\underline{c},\underline{d},[\alpha_{1}])$.
\end{itemize}
\noindent If two length parameters are inside $\mathfrak{d}$ then two of the three equivalence relations of $\mathcal{A}$ are satisfied. In this case both angle parameters are collapsed into a point. 

\noindent For example if \emph{(i)} $\underline{b}\in \mathfrak{d}$ and \emph{(ii)} $\underline{c}\in \mathfrak{d}$ then $(\underline{b},\underline{c},\underline{d},\alpha_{1},\alpha_{2}) \sim (\underline{b},\underline{c},\underline{d},\overline{\alpha}_{1},\overline{\alpha}_{2}) \ \forall \overline{\alpha}_{1},\forall\overline{\alpha}_{2}$ as 
$$(\underline{b},\underline{c},\underline{d},\alpha_{1},\alpha_{2})\overset{(ii)}{\sim}(\underline{b},\underline{c},\underline{d},\alpha_{1}+\alpha_{2}-\overline{\alpha}_{2},\overline{\alpha}_{2}) \ \forall \overline{\alpha}_{2}\overset{(i)}{\sim}(\underline{b},\underline{c},\underline{d},\overline{\alpha}_{1},\overline{\alpha}_{2}) \ \forall \overline{\alpha}_{1}$$
We construct $(H,\Yi)$ in the following way:
\begin{itemize}
\item[2.1] If $\underline{b},\underline{c} \in \mathfrak{d}^{2}, \ \underline{d} \in \mathfrak{a}$ we use Proposition \ref{parnongentype2} to construct a non-generic hexagon of type 2.1 where $(H,\Yi)=(\underline{b},\underline{c},\underline{d})$.

\item[2.2] If $\underline{b},\underline{d} \in \mathfrak{d}^{2}, \ \underline{c} \in \mathfrak{a}$ we use Proposition \ref{parnongentype2} to construct a non-generic hexagon of type 2.2 where $(H,\Yi)=(\underline{b},\underline{c},\underline{d})$.

\item[2.3] If $\underline{c},\underline{d} \in \mathfrak{d}^{2}, \ \underline{b} \in \mathfrak{a}$ we use Proposition \ref{parnongentype2} to construct a non-generic hexagon of type 2.3 where $(H,\Yi)=(\underline{b},\underline{c},\underline{d})$.

\item[3.] If $\underline{b},\underline{c},\underline{d} \in \mathfrak{d}^{3}$ we use Proposition \ref{hexparinH2} to construct a non-generic hexagon of type 3 where $(H,\Yi)=(\underline{b},\underline{c},\underline{d})=(b,c,d) \in \R^{3}_{+}$.
\end{itemize}

\noindent It is clear that any equivalent point $(\underline{b},\underline{c},\underline{d},\overline{\alpha}_{1},\overline{\alpha}_{2})\sim(\underline{b},\underline{c},\underline{d},\alpha_{1},\alpha_{2})$ in $\mathcal{A}$ induces an isometric hexagon $(H',\Yi)$ in $\mathcal{H}$. 

\noindent Conversely, let $(H,\Y_{1})$ be a hexagon in $\mathcal{H}$ and let us write $H=(P,A,B,C,D,Q)$. Up to isometry we can consider $P=0$, $B=\Id$, $Q=\infty$ and $C$ diagonal, so that $\Y_{1}=\Yi$. We put $\underline{b}=d^{\mathfrak{a}^{+}}(iA,i\Id)$, $\underline{c}=d^{\mathfrak{a}^{+}}(i\Id,iC)$ and $\underline{d}=d^{\mathfrak{a}^{+}}(iC,iD)$. Again we use Propositions \ref{parnongentype1} and \ref{parnongentype2} to determine arc coordinates. More precisely:

\begin{itemize}
\item[0.] If $\underline{b},\underline{c},\underline{d} \in \mathfrak{a}^{3}$ we associate to $(H,\Yi)$ the point $(\underline{b},\underline{c},\underline{d},[\alpha_{1},\alpha_{2}])$ using Proposition \ref{parforhex}.

\item[1.1] If $\underline{b} \in \mathfrak{d}, \ \underline{c},\underline{d} \in \mathfrak{a}^{2}$ then for $(\underline{b},\underline{c},\underline{d},\alpha_{1},\alpha_{2}) \in \mathcal{A}_{H,\Y_{1}}$ it holds $ (\underline{b},\underline{c},\underline{d},\alpha_{1},\alpha_{2})\sim (\underline{b},\underline{c},\underline{d},\overline{\alpha}_{1},\alpha_{2})$ for all $ \overline{\alpha}_{1}$. We compute the point $(\underline{b},\underline{c},\underline{d},[\alpha_{2}])$ using Proposition \ref{parnongentype1} and we associate to $(H,\Yi)$ the point $(\underline{b},\underline{c},\underline{d},[\bullet,\alpha_{2}])$.

\item[1.2] If $\underline{c} \in \mathfrak{d}, \ \underline{b},\underline{d} \in \mathfrak{a}^{2}$ then for $(\underline{b},\underline{c},\underline{d},\alpha_{1},\alpha_{2}) \in \mathcal{A}_{H,\Y_{1}}$ it holds $(\underline{b},\underline{c},\underline{d},\alpha_{1},\alpha_{2})\sim (\underline{b},\underline{c},\underline{d},\overline{\alpha}_{1},\overline{\alpha}_{2})$ for $\overline{\alpha}_{1},\overline{\alpha}_{2}$ such that   $\overline{\alpha}_{1}+\overline{\alpha}_{2}=\alpha_{1}+\alpha_{2}$. 
We compute the point $(\underline{b},\underline{c},\underline{d},[\alpha])$ using Proposition \ref{parnongentype1} and we associate to $(H,\Yi)$ the point $(\underline{b},\underline{c},\underline{d},[\frac{\alpha}{2},\frac{\alpha}{2}+\pi])$.

\item[1.3] If $\underline{d} \in \mathfrak{d},\  \underline{c},\underline{d} \in \mathfrak{a}^{2}$ then for $(\underline{b},\underline{c},\underline{d},\alpha_{1},\alpha_{2}) \in \mathcal{A}_{H,\Y_{1}}$ it holds $ (\underline{b},\underline{c},\underline{d},\alpha_{1},\alpha_{2})\sim (\underline{b},\underline{c},\underline{d},\alpha_{1},\overline{\alpha}_{2})$ for all $ \overline{\alpha}_{2}$. We compute the point $(\underline{b},\underline{c},\underline{d},[\alpha_{1}])$ using Proposition \ref{parnongentype1} and we associate to $(H,\Yi)$ the point $(\underline{b},\underline{c},\underline{d},[\alpha_{1},\bullet])$.
\end{itemize}
For the cases 2.1, 2.2, 2.3 and 3 all the angle parameters vanish and we associate to $(H,\Yi)$ the point $(\underline{b},\underline{c},\underline{d},[\bullet,\bullet])$.
\endproof

\begin{notation}
Given $(H,\Y_{1})$ inside $\mathcal{H}$ its arc coordinates $(\underline{b},\underline{c},\underline{d},[\alpha_{1},\alpha_{2}])$ will be denoted $\mathcal{A}(H,\Y_{1})$. The vectors $\underline{b},\underline{c},\underline{d}$ are \emph{length parameters} and $\alpha_{1},\alpha_{2}$ are \emph{angle parameters}.
\end{notation}

\subsection{Hexagons inside a maximal polydisk} \label{hexinmaxpol}

Let $\X$ be the symmetric space associated to $\Sp(4,\R)$. We have seen in Section \ref{H2insideX} how to embed $\H^{2} \times \H^{2}$ inside $\X$. A right-angled hexagon $H=[\Y_{1},\Y_{2},\Y_{3},\Y_{4},\Y_{5},\Y_{6}]$ is contained in a maximal polydisk if there exists an isometry $g$ such that $gH$ is contained in the model polydisk. In particular a hexagon $H$ is contained in the model polydisk of $\X$ if and only if all tubes have diagonal matrices as endpoints. 
 
\begin{prop} \label{prophexinmaxpol}
The subspace $\mathcal{D}\subset \mathcal{A}$  
$$
\mathcal{D}=\Big\{ (\underline{b},\underline{c},\underline{d},[\alpha_{1},\alpha_{2}]) \in \mathcal{A}|\ [\alpha_{1},\alpha_{2}] \in \{[0,0],[0,\pi],[\pi,0],[\pi,\pi]\}  \Big\} \subset \mathcal{A}
$$
corresponds to right-angled hexagons inside a maximal polydisk in $\X$.  
\end{prop}

\proof
In the case where $\underline{b},\underline{c},\underline{d} \in \mathfrak{a}^{3}$ the point $p=(\underline{b},\underline{c},\underline{d},[\alpha_{1},\alpha_{2}])$ corresponds to a generic hexagon. Using Proposition \ref{parforhex} we know 
$p=[(H,\Yi)]$ where $(H,\Yi)=(0,A,\Id,C,D,\infty)$ with $C=\bpm e^{c_{1}}&0\\0&e^{c_{2}}\epm$ and
$$
A=\bpm \frac{1}{e^{b_{2}}}&0\\0&\frac{1}{e^{b_{1}}}\epm, \ D=\bpm e^{c_{1}+d_{2}}&0\\0&e^{c_{2}+d_{1}}\epm \text{ for } \alpha_{1}=\alpha_{2}=0
$$
$$
A=\bpm \frac{1}{e^{b_{1}}}&0\\0&\frac{1}{e^{b_{2}}}\epm, \ D=\bpm e^{c_{1}+d_{1}}&0\\0&e^{c_{2}+d_{2}}\epm \text{ for } \alpha_{1}=\alpha_{2}=\pi
$$
$$
A=\bpm \frac{1}{e^{b_{2}}}&0\\0&\frac{1}{e^{b_{1}}}\epm, \ D=\bpm e^{c_{1}+d_{1}}&0\\0&e^{c_{2}+d_{2}}\epm \text{ for } \alpha_{1}=0, \ \alpha_{2}=\pi
$$
$$
A=\bpm \frac{1}{e^{b_{1}}}&0\\0&\frac{1}{e^{b_{2}}}\epm, \ D=\bpm e^{c_{1}+d_{2}}&0\\0&e^{c_{2}+d_{1}}\epm \text{ for } \alpha_{1}=\pi, \ \alpha_{2}=0
$$
All four cases correspond to hexagons consisting of tubes that have diagonal matrices as endpoints. This is consistent with the geometrical meaning of the angle parameter described in Section \ref{geomintdiag}. A similar argument can be used for the case of non-generic hexagons of type 1.  All non-generic right-angled hexagons of type 2 and 3 are contained in a maximal  polydisk in $\X$ and in these cases for all $(\underline{b},\underline{c},\underline{d},[\alpha_{1},\alpha_{2}]) \in \mathcal{D}$ it holds $(\underline{b},\underline{c},\underline{d},\alpha_{1},\alpha_{2}) \sim (\underline{b},\underline{c},\underline{d},\overline{\alpha}_{1},\overline{\alpha}_{2}) $ for all $\overline{\alpha}_{1},\overline{\alpha}_{2}$.

\noindent Conversely, if $[(H,\Y_{1})] \in \mathcal{H}/_{\PSp(4,\R)}$ is contained in a maximal polydisk then we can move it into the model polydisk through an isometry. It is easy to see that in this case the point $p \in \mathcal{A}$ corresponding to $[(H,\Y_{1})]$ must be inside $\mathcal{D}$.
\endproof

\begin{definition} \label{defmathcalDH2}
    We define $\mathcal{D}_{\H^{2}}:=\big\{ (\underline{b},\underline{c},\underline{d},[\alpha_{1},\alpha_{2}]) \in \mathcal{A}|\ \underline{b},\underline{c},\underline{d} \in   \mathfrak{d}^{3}  \big\}$.
\end{definition}

Recall \ref{defstabhex} for the definition of the stabilizer of a right-angled hexagon.
 
\begin{prop} \label{stabofhex}
Let $\X$ be the symmetric space associated to $\Sp(4,\R)$ and let $H=[\Y_{1},\Y_{2},\Y_{3},\Y_{4},\Y_{5},\Y_{6}]$ be a right-angled hexagon in $\X$. It holds
    \begin{itemize}
        \item[(i)] If $H$ is contained in a copy of $\H^{2}$ inside $\X$ then $\Stab(H)\cong \PO(2)$
        \item[(ii)] If $H$ is contained in a maximal polydisk then $\Stab(H)\cong \mathbb{Z}/2\mathbb{Z}$
        \item[(iii)] If $H$ is not contained in any maximal polydisk then $\Stab(H)=\{id\}$
    \end{itemize}
\end{prop}

\proof Up to isometry we can consider $H=(0,A,\Id,C,D,\infty)$ where $C$ is diagonal. 
\begin{itemize}
    \item[\emph{(i)}] If $H$ is contained in the diagonal disc then the matrices $A,C,D$ are all multiples of $\Id$ and so are all endpoints of the $\R$-tubes of $H$. It is clear that the stabilizer is $\PO(2)$.  
    \item[\emph{(ii)}] The matrices $A,C,D$ are all diagonal and so are all endpoints of the $\R$-tubes of $H$. The stabilizer is given by the identity together with $\bpm r&0\\0&r\epm$ where $r=\bpm -1&0\\0&1\epm$.
    \item[\emph{(iii)}] This is clear.
\end{itemize}
\endproof

\section{Discussion about the parameters} \label{disconpar}

In $\H^{2}$ a right-angled hexagon is uniquely determined by the length of three alternating sides. On the other hand, all length parameters of a right-angled hexagon $[(H,\Y_{1})]$ in $\X$ are lying in the tube $\Y_{1}$. In this section we will recall the proof of the parametrization in $\H^{2}$ of right-angled hexagons as done in \cite[Lemma 6.2.2]{martelli2016introduction} suitably adapted to the upper-half space model. We discuss the problems that arise when generalizing these hexagon-parameters in the symmetric space $\X$ associated to $\Sp(4,\R)$. 

\subsection{The $\H^{2}$-case} We recall two standard results in hyperbolic geometry. 

\begin{lem} \label{bijofcrossratios}
Let $\gamma_{a,b},\gamma_{c,d},\gamma_{e,f}$ be three infinite geodesics in $\H^{2}$ with endpoints $\{a,b\},\{c,d\}$ and $\{e,f\}$ respectively. Suppose
$\gamma_{a,b} \perp \gamma_{c,d} \perp \gamma_{e,f}$. Then there exists a bijective map
\begin{equation*} 
\begin{aligned}
T: \R^{+} &\to \R^{+} \\
R(c,b,e,d)& \mapsto R(a,b,e,f)
\end{aligned}
 \end{equation*}
given by $T(x)=\frac{(x+1)^{2}}{4x}$
\end{lem}

\begin{prop} \label{bijmapH2}
Let $\gamma_{x,z} \perp \gamma_{y,w}$ be two orthogonal geodesics in $\H^{2}$ with endpoints $\{x,z \}$ and $\{y,w \}$ respectively and let $P$ be their intersection point. Then there exists a bijective map $f=f(\gamma_{x,z},\gamma_{y,w},P)$ defined as follows: for $d>0$ let $P'$ be one of the two points in $\gamma_{y,w}$ at distance $d$ from $P$. Let $\gamma_{P'}$ be the geodesic through $P'$ orthogonal to $\gamma_{y,w}$ and denote by $b$ one of its endpoints. We can define
\begin{equation*} 
\begin{aligned}
f: \R^{+} &\to \R^{+} \\
 d^{\H^{2}}(P,P')& \mapsto d^{\H^{2}}(P,p_{x,z}(b))
\end{aligned}
 \end{equation*}
where $p_{x,z}$ denotes the orthogonal projection on the geodesic $\gamma_{x,z}$. The map $f$ is given by $f(d)=\log\Big(\frac{e^{d}+1}{e^{d}-1}\Big)$. This expression does not depend on the choice of the points $P',b$.
\end{prop}

\begin{figure}[!h]
   \centering
   \captionsetup{justification=centering,margin=2cm}
   \setlength{\unitlength}{0.1\textwidth}
   \begin{picture}(6,2)
     \put(0,0){\includegraphics[width=9cm,height=3cm]{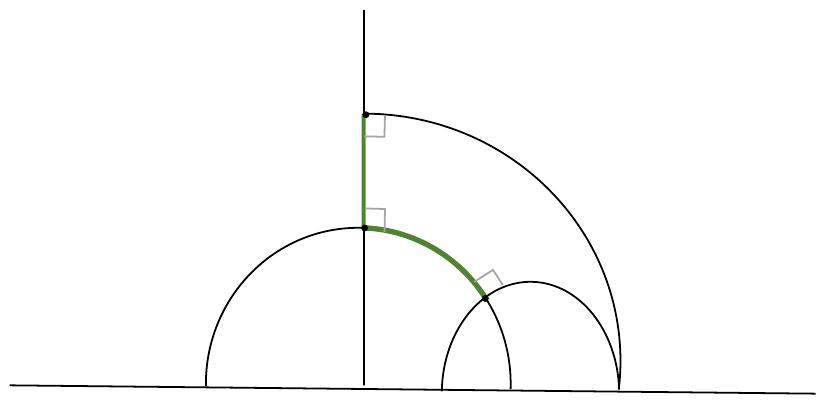}}
    \put(1.2,-0.2){$x$}
    \put(2.4,-0.2){$y$}
    \put(3.5,-0.2){$z$}
    \put(3.3,0.7){$p_{x,z}(b)$}
    \put(4.3,-0.2){$b$}
    \put(2.6,1.1){\textcolor{forestgreen}{$d$}}
    \put(2.7,0.5){\textcolor{forestgreen}{$f(d)$}}
    \put(2.3,0.9){$P$}
    \put(2.3,1.3){$P'$}
    \put(2.4,2){$w$}
    \put(0.9,1){\textcolor{gray}{$\H^{2}$}}    \end{picture}
\caption{The map $f$}
   \label{premartelli}
\end{figure}
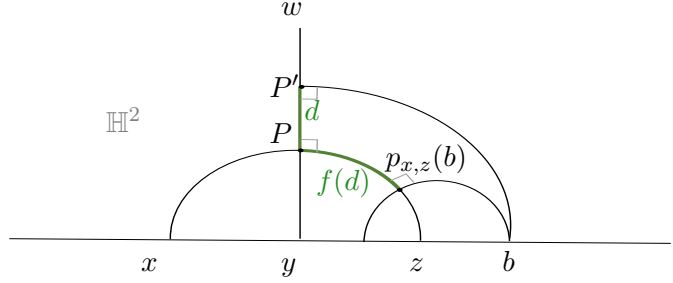

\begin{lem} \label{Martelli}(\cite[Lemma 6.2.2]{martelli2016introduction})
    Given three real numbers $b,c,d>0$ there exists (up to isometries) a unique hyperbolic right-angled hexagon with three alternate sides of length $b,c$ and $d$ respectively.
\end{lem}

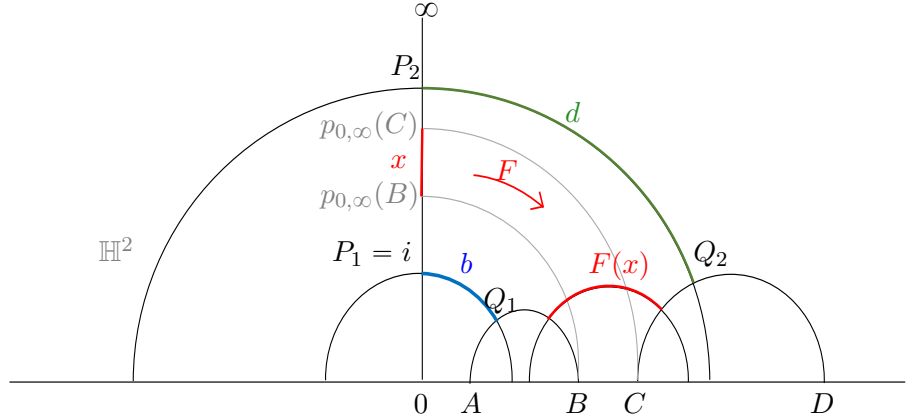
\begin{figure}[!h]
   \centering
   \captionsetup{justification=centering,margin=2cm}
   \setlength{\unitlength}{0.1\textwidth}
   \begin{picture}(7.4,3)
     \put(0,0){\includegraphics[width=12cm,height=5cm]{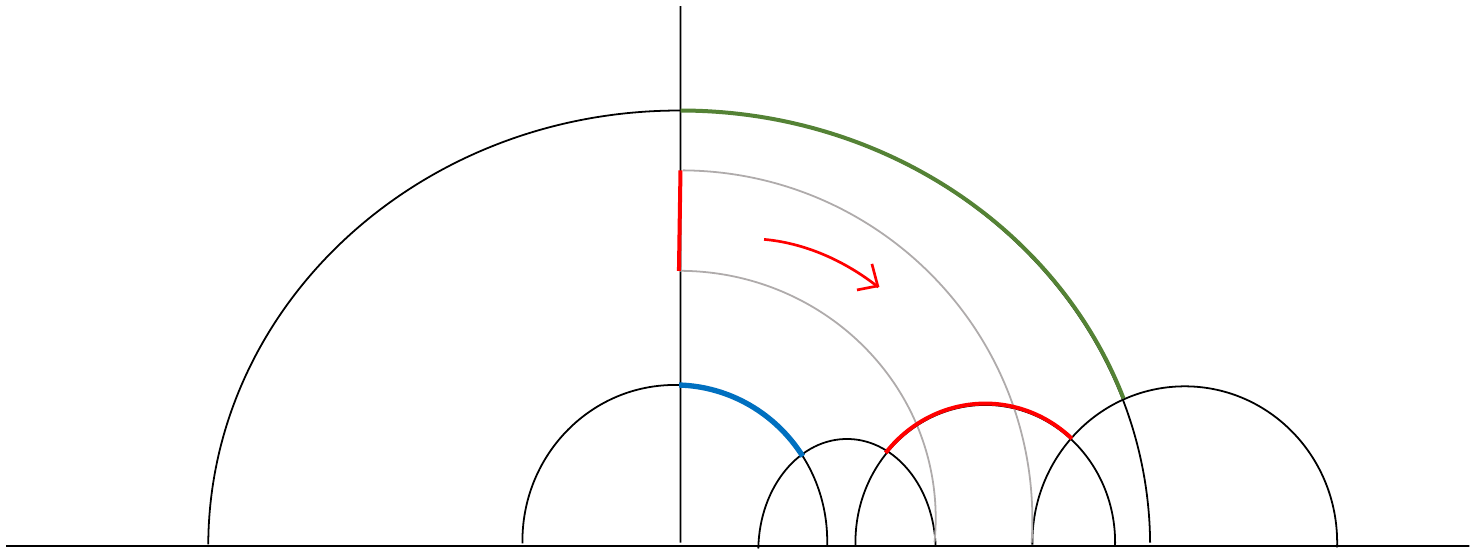}}
    \put(6.9,-0.2){$D$}
    \put(5.3,-0.2){$C$}
    \put(5,1){\textcolor{red}{$F(x)$}}
    \put(4.8,-0.2){$B$}
    \put(3.9,-0.2){$A$}
    \put(3.9,1){\textcolor{blue}{$b$}}
    \put(4.1,0.7){$Q_{1}$}
    \put(3.3,1.9){\textcolor{red}{$x$}}
    \put(2.8,1.1){$P_{1}=i$}
    \put(4.2,1.8){\textcolor{red}{$F$}}
    \put(4.8,2.3){\textcolor{forestgreen}{$d$}}
    \put(2.7,2.2){\textcolor{gray}{$p_{0,\infty}(C)$}}
    \put(3.3,2.7){$P_{2}$}
    \put(5.9,1.1){$Q_{2}$}
    \put(2.7,1.6){\textcolor{gray}{$p_{0,\infty}(B)$}}
    \put(3.5,-0.2){$0$}
    \put(3.5,3.2){$\infty$}
    \put(0.8,1.1){\textcolor{gray}{$\H^{2}$}}\end{picture}
\caption{Construction of a right-angled hexagon in the $\H^{2}$-case}
   \label{martelli}
\end{figure}

\proof Let $b,d>0$. The construction of the hexagon goes as follows: take a geodesic $\gamma$ with two arbitrary points $P_{1},P_{2}$ in it. Without loss of generality we can assume $\gamma$ to be the vertical geodesic $\gamma_{0,\infty}=\{iy| \ y>0\}$ and $P_{1}=i$ (Figure \ref{martelli}). Draw the perpendicular from $P_{1}=i$ and from $P_{2}$. At distances $b$ and $d$ we find two points $Q_{1}$ and $Q_{2}$ and we draw again two perpendiculars $\gamma_{A,B}$ and $\gamma_{C,D}$, with some points at infinity $A,B$ and $C,D$ respectively. Draw the unique perpendiculars to $\gamma$ pointing to $B$ and $C$: they intersect $\gamma$ in two points $p_{0,\infty}(B)$ and $p_{0,\infty}(C)$. Note that the lengths $d^{\H^{2}}(P_{1},p_{0,\infty}(B))$ and $d^{\H^{2}}(p_{0,\infty}(C),P_{2})$ have some fixed length depending only on $b$ and $d$ through a bijective map given explicitly in Proposition \ref{bijmapH2} (this is the map $f=f^{-1}$). We can vary the parameter $x=d^{\H^{2}}(p_{0,\infty}(B),p_{0,\infty}(C))$, the geodesics $\gamma_{A,B}$ and $\gamma_{C,D}$ are
ultra-parallel and there is a unique segment orthogonal to both which has length $F(x)$. The function $F: (0,+ \infty) \to (0,+ \infty)$ is continuous, strictly monotonic, and with $\lim_{x \to \infty}= \infty$: therefore there is precisely one $c$ such that $F(x)=c$.
\endproof

\begin{remark} \label{mapF}
By Proposition \ref{bijmapH2} we know that there is a bijection between the length $b,d$ of Figure \ref{martelli} and the segments $d^{\H^{2}}(P_{1},p_{0,\infty}(B)),d^{\H^{2}}(p_{0,\infty}(C),P_{2})$ respectively. We also know that $F(x)=c$ for a bijective map $F$. We can therefore think as the lengths $b,c,d$ determining the hexagon as all lying on the vertical geodesic $\gamma_{0,\infty}$. To explicitly write the map $F$ we consider the case where $A=e^{-2b}$, $B=1$, $C=e^{c}$, $D=e^{c+2d}$, $P_{1}=ie^{-b}$ and $P_{2}=ie^{c+d}$. For  $b,d$ fixed we obtain
$$
T \circ F(c)=\frac{(e^{c+2d}-1)(1-e^{2b+c})}{e^{c}(1-e^{2b})(e^{2d}-1)}=y
$$
where $T$ is the bijective map of Lemma \ref{bijofcrossratios} and $F(c)=T^{-1}(y)$.
\end{remark}

\subsection{The $\X$-case}

In the Siegel space $\X$ the analogue of geodesics are $\R$-tubes and length-parameters take value in the Weyl chamber. 

\begin{lem} \label{bijofcrossratios2}
Let $\Y_{A,B},\Y_{C,D},\Y_{E,F}$ be three $\R$-tubes inside $\X$ such that
$\Y_{A,B} \perp \Y_{C,D} \perp \Y_{E,F}$.
Let us denote by $(x_{1} \geq... \geq x_{n})$ the eigenvalues of $R(C,B,E,D)$ and by $(y_{1} \geq... \geq y_{n})$ the eigenvalues of $R(A,B,E,F)$. Then there exist a bijective map $T(x_{1},...x_{n})=(y_{1},...y_{n})$ where  $y_{i}=\frac{(x_{i}+1)^{2}}{4x_{i}}$ for $i \in \{ 1,...n \}$.

\end{lem}

\proof Up to isometry we can consider $\Y_{C,D}=\Yi$, $\Y_{A,B}=\Y_{-\Id,\Id}$ and $\Y_{E,F}=\Y_{-X,X}$ where $(0,\Id,X,\infty)$ is maximal and $X$ is diagonal with diagonal entries $(x_{1},...,x_{n})$. We obtain $R(C,B,E,D)=R(0,\Id,X,\infty)=(x_{1},...,x_{n})$ and
$R(A,B,E,F)=R(-\Id,\Id,X,-X)$. The matrix $R(-\Id,\Id,X,-X)$ is diagonal with entries $\big(\frac{(x_{1}+1)^{2}}{4x_{1}},...,\frac{(x_{n}+1)^{2}}{4x_{n}}\big)$.
\endproof

\begin{prop} \label{bijmap1}
Let $\Y_{X,Z} \perp \Y_{Y,W}$ be two orthogonal $\R$-tubes in $\X$  and let $P$ be their intersection point. Then there exists a bijective map $f=f(\Y_{X,Z},\Y_{Y,W},P)$ defined as follows: for $(d_{1},...d_{n}) \in \a$ let $P'$ be a point in $\Y_{Y,W}$ at distance $(d_{1},...d_{n})$ from $P$. Let $\Y_{P'}$ be the tube through $P'$ orthogonal to $\Y_{Y,W}$ and denote by $B$ one of its endpoints. Then we can define
\begin{equation*} 
\begin{aligned}
f: \overline{\mathfrak{a}} &\to \overline{\mathfrak{a}} \\
 d^{\a}(P,P')& \mapsto d^{\a}(P,p_{X,Z}(B))
\end{aligned}
 \end{equation*}
where $p_{X,Z}$ denotes the orthogonal projection on the tube $\Y_{X,Z}$. The map $f$ is given by $f(d_{1},...d_{n})=(\log\Big(\frac{e^{d_{n}}+1}{e^{d_{n}}-1}\Big),...,\log\Big(\frac{e^{d_{1}}+1}{e^{d_{1}}-1}\Big) )$. This expression does not depend on the choice of the points $P',B$. In particular the image of a regular point inside $\a$ is a regular point.
\end{prop}

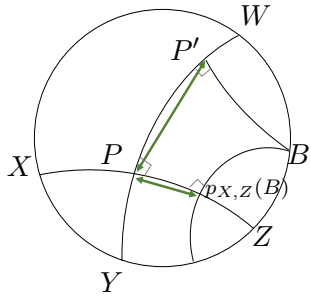
\begin{figure}[!h]
   \centering
   \captionsetup{justification=centering,margin=2cm}
   \setlength{\unitlength}{0.1\textwidth}
   \begin{picture}(2.5,2.3)
     \put(0,0){\includegraphics[width=3.5cm,height=3.5cm]{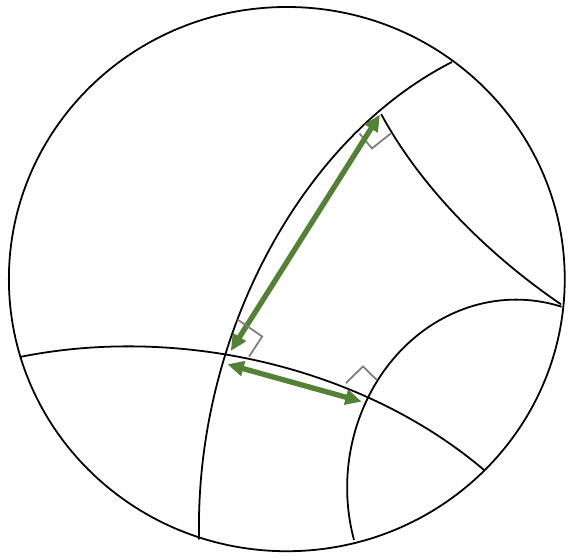}}
          \put(0.6,-0.2){$Y$}
          \put(1.8,2.1){$W$}
          \put(2.2,0.9){$B$}
          \put(1.9,0.2){$Z$}
          \put(-0.2,0.8){$X$}
          \put(1.2,1.8){$P'$}
          \put(0.6,0.9){$P$}
          \put(1.5,0.7){$_{p_{X,Z}(B)}$}
   \end{picture}   
\caption{There is a bijection between the green vectors}
\label{bij1}
\end{figure}

\proof Without loss of generality we can assume $\Y_{Y,W}=\Yi$ and $\Y_{X,Z}=\Y_{-\Id,\Id}$. The intersection point $\Yi \cap \Y_{-\Id,\Id}$ is the point $P=i\Id$. Let $P'=iB$  where $B$ is a diagonal matrix with entries $(e^{d_{1}},...,e^{d_{n}})$. By Lemma \ref{proj} and Lemma \ref{comput3} we know that $d^{\a}(P,p_{-\Id,\Id}(B))$ is given by $(\log \mu_{1},...,\log \mu_{1})$, where $\mu_{1}>...>\mu_{n}$ are the eigenvalues of $R(-\Id,0,B^{-1},\Id)$. Calculations give
$R(-\Id,0,B^{-1},\Id)=(\Id-B^{-1})^{-1}(\Id+B^{-1})$. The matrix  $R(-\Id,0,B^{-1},\Id)$ is diagonal with entries  $(\Big(\frac{e^{d_{1}}+1}{e^{d_{1}}-1}\Big),...,\Big(\frac{e^{d_{n}}+1}{e^{d_{n}}-1}\Big))$. It follows $f(d_{1},...,d_{n})=(\log\Big(\frac{e^{d_{n}}+1}{e^{d_{n}}-1}\Big),...,\log\Big(\frac{e^{d_{1}}+1}{e^{d_{1}}-1}\Big) )
$. Observe that we need to invert the order of $d_{1},...,d_{n}$ since the function $h(x)=\log\Big(\frac{e^{x}+1}{e^{x}-1}\Big)$ is decreasing for $x>0$. From the expression of $f$ it is clear that regular points of $\a$ are sent to regular points. It is trivial to show that the expression of $f$ does not depend on the choice of the points $P',B$.
\endproof

\begin{cor} \label{genofquadr}
Let $(H,\Yi) \in \mathcal{H}^{gen}$ be a generic hexagon $(H,\Yi)=(0,A,\Id,C,D,\infty)$. Then the quadruples
$(-A,0,A^{2},A),(-D,0,C,D)$ are generic.    
\end{cor}

\proof The quadruple $(0,A,\Id,\infty)$ is generic and the parameter $\underline{b}=(b_{1},b_{2})$ of Figure \ref{nobij} lies inside the set of regular vectors $\mathfrak{a}$. To show genericity of $(-A,0,A^{2},A)$ we need to prove that the cross-ratio $R(-A,0,A^{2},A)$ has distinct eigenvalues. By Lemma \ref{proj} we know that taking the logarithm of these ordered eigenvalues gives the distance $d^{\a}(p_{-A,A}(0),p_{-A,A}(A^{2}))$. By Proposition \ref{bijmap1} we know that this vector is the image under the bijection $f$ of the vector $\underline{b}$ and that $f$ is sending regular points to regular points. It follows $d^{\a}(p_{-A,A}(0),p_{-A,A}(A^{2})) \in \mathfrak{a}$ and so $(-A,0,A^{2},A)$ is generic. 
\endproof

\subsection{Changing side of the hexagon} 
The following Proposition relates length parameters of arc coordinates when we change side of the ordered right-angled hexagon.

\begin{prop}
Let $H=[\Y_{1},\Y_{2},\Y_{3},\Y_{4},\Y_{5},\Y_{6}]$ be a right-angled hexagon. Let $(\underline{b_{1}},\underline{c_{1}},\underline{d_{1}})$, $(\underline{b_{2}},\underline{c_{2}},\underline{d_{2}})$, $(\underline{b_{3}},\underline{c_{3}},\underline{d_{3}})$ denote length parameters of $\mathcal{A}(H,\Y_{1}), \mathcal{A}(H,\Y_{3})$ and $\mathcal{A}(H,\Y_{5})$ respectively. Then $\underline{b_{1}}=\underline{d_{2}}$, $\underline{b_{2}}=\underline{d_{3}}$ and $\underline{b_{3}}=\underline{d_{1}}$.
\end{prop}

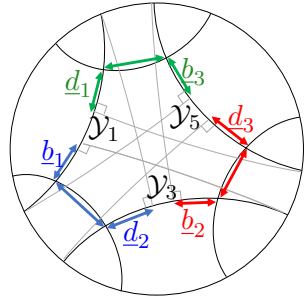
\begin{figure}[!h]
   \centering
   \captionsetup{justification=centering,margin=2cm}
   \setlength{\unitlength}{0.1\textwidth}
   \begin{picture}(2.6,2.5)
     \put(0,0){\includegraphics[width=4cm,height=4cm]{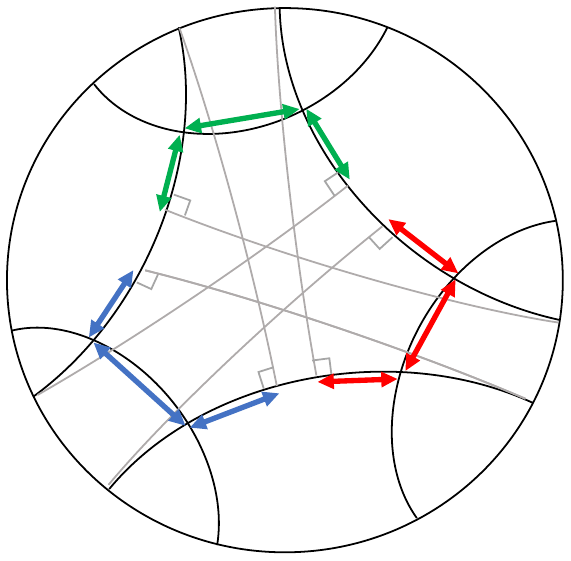}}
        \put(0.3,1.2){\textcolor{blue}{$\underline{b}_{1}$}}
        \put(0.5,1.8){\textcolor{forestgreen}{$\underline{d}_{1}$}}
        \put(1.5,0.6){\textcolor{red}{$\underline{b}_{2}$}}
        \put(1,0.5){\textcolor{blue}{$\underline{d}_{2}$}}
        \put(1.9,1.5){\textcolor{red}{$\underline{d}_{3}$}}
        \put(1.5,1.9){\textcolor{forestgreen}{$\underline{b}_{3}$}}
        \put(0.7,1.4){$\Y_{1}$}
          \put(1.2,0.9){$\Y_{3}$}
          \put(1.4,1.5){$\Y_{5}$}
 \end{picture}
   \caption{There is a bijection between vectors of the same colour}
   \label{star}
\end{figure}

\proof Let us prove $\underline{b_{1}}=\underline{d_{2}}$. Let (see the configuration of Figure \ref{max12tuple})
$$\Y_{1}=\Y_{P_{1},P_{2}}, \ \Y_{2}= \Y_{Q_{1},Q_{2}}, \ \Y_{3}=  \Y_{P_{3},P_{4}} , \ \Y_{4}= \Y_{Q_{3},Q_{4}} , \ \Y_{5}= \Y_{P_{5},P_{6}},\ \Y_{6}= \Y_{Q_{5},Q_{6}}
$$
By definition of $\mathcal{A}(H,\Y_{1})$ we know 
$\underline{b_{1}}=d^{\a}(p_{P_{1},P_{2}}(Q_{2}),p_{P_{1},P_{2}}(P_{4}))$ where \\$(H,\Y_{1})=(P_{2},Q_{2},P_{4},P_{5},Q_{5},P_{1})$, $\Y_{1}=\Y_{P_{1},P_{2}}$. Similarly by definition of $\mathcal{A}(H,\Y_{3})$ we know $\underline{d_{2}}=d^{\a}(p_{P_{3},P_{4}}(P_{1}),p_{P_{3},P_{4}}(Q_{1}))$ where $(H,\Y_{3})=(P_{4},Q_{4},P_{6},P_{1},Q_{1},P_{3})$, $\Y_{3}=\Y_{P_{3},P_{4}}$. Let $\underline{v}=d^{\a}(p_{Q_{1},Q_{2}}(P_{2}),p_{Q_{1},Q_{2}}(P_{3}))$
and let $f:\overline{\mathfrak{a}} \to \overline{\mathfrak{a}}$ be the map $f(a_{1},a_{2})=(\log\Big(\frac{e^{a_{2}}+1}{e^{a_{2}}-1}\Big),\log\Big(\frac{e^{a_{1}}+1}{e^{a_{1}}-1}\Big) )$. By Proposition \ref{bijmap1} it holds $\underline{b_{1}}=f^{-1}(\underline{v})=\underline{d_{2}}$. The proof for $\underline{b_{2}}=\underline{d_{3}}$ and $\underline{b_{3}}=\underline{d_{1}}$ is similar.
\endproof
 
\noindent An analogue Proposition which relates the length-parameters $\underline{c}_{i}$ is trickier. 
More generally when we try to generalize the map $F$ described in Lemma \ref{Martelli} we can not guarantee bijectivity. This will be explained in the next section.

\subsection{Constraints in generalizing hexagon parameters of $\H^{2}$} \label{constrains}

\begin{figure}[!h]
   \centering
   \captionsetup{justification=centering,margin=2cm}
   \setlength{\unitlength}{0.1\textwidth}
   \begin{picture}(5,3)
     \put(0.8,0){\includegraphics[width=4.8cm,height=4.5cm]{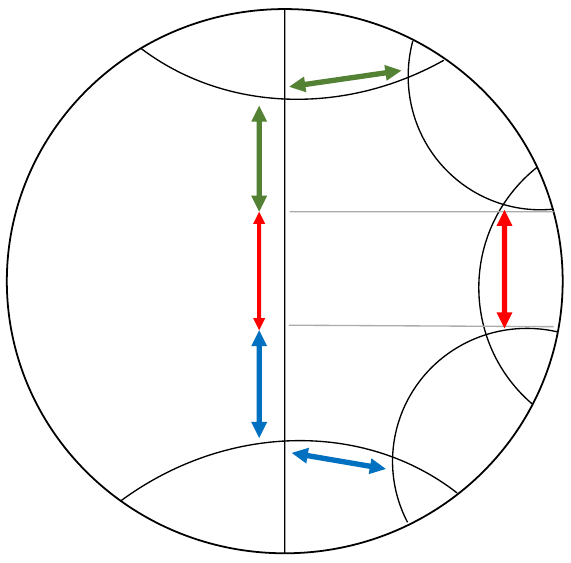}}
         \put(2.3,-0.16){$0$}
          \put(2.25,3){$\infty$}
          \put(3.9,1.15){$\Id$}
          \put(3.9,0.75){$Z_{1}$}
          \put(3.3,0.9){$P$}
          \put(3.9,1.8){$C=\bpm e^{c_{1}} & 0\\0 & e^{c_{2}}\epm$}
          \put(3.8,2.15){$Z_{2}$}
          \put(3.4,1.9){$Q$}
          \put(3,2.8){$DC^{-1}D$}
          \put(3.3,2.6){$D$}
          \put(1.2,2.7){$-D$}
          \put(2.9,0){$A^{2}$}
          \put(3.3,0.2){$A$}
          \put(1.1,0.1){$-A$}
           \put(2.7,1.5){$	\nleftrightarrow$}
          \put(2,0.9){\textcolor{blue}{$\underline{b}$}}
          \put(2,1.5){\textcolor{red}{$\underline{c}$}}
          \put(2,2.1){\textcolor{forestgreen}{$\underline{d}$}}
          \end{picture}
\caption{There is no bijective map between the red vectors}
   \label{nobij}
\end{figure}

In the previous section we have shown that given a right-angled hexagon $(H,\Y_{1})$ we can find a bijection between length parameters $\underline{b},\underline{d}$ of $\mathcal{A}(H,\Y_{1})$ and the vectorial length of two alternating sides (see Figure \ref{nobij}). This is analogue to the hyperbolic case. It is natural to ask whether for a right-angled hexagon in $\X$ there exists a bijective map also between the vector-parameter $\underline{c}$ of Figure \ref{nobij} and the missing alternating side of the hexagon. When the hexagon is non-generic of type 3 this is trivially true and corresponds to the immersion of hyperbolic hexagons inside $\X$. In this section we show that this is not the case for a general right-angled hexagon $H$ inside $\X$. More precisely, let $(H,\Yi)$ be a right-angled hexagon inside $\mathcal{H}$. We can determine $(H,\Yi)$ by the following maximal 12-tuple:
$$(H,\Yi)=(\infty,-D,-A,0,A^{2},A,Z_{1},\Id,C=\bpm 
e^{c_{1}} & 0\\
0 & e^{c_{2}}
\epm,Z_{2},D,DC^{-1}D)$$ 

\noindent Let $P,Q$ be the intersection points $P=\Y_{A^{2},\Id} \cap \Y_{Z_{1},Z_{2}}$, $Q=\Y_{Z_{1},Z_{2}} \cap \Y_{C,DC^{-1}C}$ and let $F$ be the following map
\begin{equation} \label{lemmaleficomap}
\begin{aligned}
F: \overline{\mathfrak{a}} &\to \overline{\mathfrak{a}} \\
 d^{\a}(i\Id,iC)=(c_{1},c_{2})& \mapsto d^{\a}(P,Q)
\end{aligned}
 \end{equation}

\noindent One can ask if the map $F$ is bijective. In this section we show that this is not the case and we provide a counterexample in the case where $(H,\Yi)$ is contained in a maximal polydisk.

\begin{remark} \label{mapFlemmamal}
The image of the map $F$ in (\ref{lemmaleficomap}) is the distance $d^{\a}(P,Q)=(\log\mu_{1},\log\mu_{2})$ where $\mu_{1}>\mu_{2}$ are the eigenvalues of $R(Z_{1},\Id,C,Z_{2})$. Asking for the existence of such a map $F$ is equivalent to ask for the existence of a map $T \circ F$  
    \begin{equation*} 
\begin{aligned} 
T \circ F: \overline{\mathfrak{a}} &\to \overline{\mathfrak{a}} \\
(c_{1},c_{2})& \mapsto (\log\lambda_{1},\log\lambda_{2})
\end{aligned}
 \end{equation*}
where $\lambda_{1} \geq \lambda_{2}$ are the eigenvalues of $R(A^{2},\Id,C,DC^{-1}D)$ and $T$ is the bijective map of Lemma \ref{bijofcrossratios2} (composed with the logarithm map). By abuse of notation we will write this map as $F$ and we can express the cross-ratio with respect to arc coordinates. This is made more precise in the following definition.
\end{remark}

\begin{definition} \label{defmaleficmap}
Let $\underline{b},\underline{d} \in \overline{\mathfrak{a}}$ and $\alpha_{1},\alpha_{2} \in [0,2\pi)$. We will call the \emph{malefic map} $F_{\underline{b}, \underline{d}, \alpha_{1},\alpha_{2}}$ the map defined as following:
\begin{equation*} 
\begin{aligned} 
F_{\underline{b}, \underline{d}, \alpha_{1},\alpha_{2}}: \overline{\mathfrak{a}} &\to \overline{\mathfrak{a}} \\
(c_{1},c_{2})& \mapsto (\log\lambda_{1},\log\lambda_{2})
\end{aligned}
 \end{equation*}
where $\lambda_{1} \geq \lambda_{2}$ are the eigenvalues of the cross-ratio $R(A^{2},\Id,C,DC^{-1}D)$ where \\
$(H,\Yi)=(0,A,\Id,C,D,\infty)$ is the right-angled hexagon with arc coordinates $\mathcal{A}(H,\Yi)$ equal to $(\underline{b},\underline{c},\underline{d},[\alpha_{1},\alpha_{2}])$.
\end{definition}

\begin{example} \label{malmapzpi}
The malefic map $F_{\underline{b}, \underline{d}, \alpha_{1},\alpha_{2}}$ clearly depends on the choice of the parameters $\underline{b}, \underline{d}, \alpha_{1},\alpha_{2}$. It is not hard to show that for $(\alpha_{1},\alpha_{2})=(0,0)$ and $(\alpha_{1},\alpha_{2})=(\pi,\pi)$ respectively we obtain 
\begin{equation*}    
 F_{\underline{b}, \underline{d}, 0,0}(c_{1},c_{2})=\Big(\log \frac{(e^{c_{1}+2d_{2}}-1)(1-e^{2b_{2}+c_{1}})}{e^{c_{1}}(1-e^{2b_{2}})(e^{2d_{2}}-1)},\log \frac{(e^{c_{2}+2d_{1}}-1)(1-e^{2b_{1}+c_{2}})}{e^{c_{2}}(1-e^{2b_{1}})(e^{2d_{1}}-1)} \Big) 
\end{equation*}
\begin{equation*}
F_{\underline{b},\underline{d},\pi,\pi}(c_{1},c_{2})= \Big(\log \frac{(e^{c_{1}+2d_{1}}-1)(1-e^{2b_{1}+c_{1}})}{e^{c_{1}}(1-e^{2b_{1}})(e^{2d_{1}}-1)},\log \frac{(e^{c_{2}+2d_{2}}-1)(1-e^{2b_{2}+c_{2}})}{e^{c_{2}}(1-e^{2b_{2}})(e^{2d_{2}}-1)} \Big)   
\end{equation*} 

\noindent where $\underline{b}=(b_{1},b_{2})$ and $\underline{d}=(d_{1},d_{2})$. Observe that in both cases the constructed hexagon lies inside a maximal polydisk (see Proposition \ref{prophexinmaxpol}).
\end{example}

\begin{lem}
Let $\underline{b},\underline{d} \in \overline{\mathfrak{a}}$, $\alpha_{1},\alpha_{2} \in [0,2\pi)$ and let $F_{\underline{b}, \underline{d}, \alpha_{1},\alpha_{2}}$ be the malefic map. It holds 
$$
F_{\underline{b}, \underline{d}, \alpha_{1},\alpha_{2}}(c_{1},c_{2})=F_{\underline{b}, \underline{d}, 2\pi-\alpha_{1},2\pi-\alpha_{2}}(c_{1},c_{2})
$$
\end{lem}

\proof This is straightforward by Proposition \ref{parforhex} in the generic case and more generally by Theorem \ref{cpctarccoordthm}: for angle parameters $(\alpha_{1},\alpha_{2})$ and $(2\pi-\alpha_{1},2\pi-\alpha_{2})$ we obtain two isometric hexagons.
\endproof

We can extend the malefic map $F_{\underline{b}, \underline{d}, \alpha_{1},\alpha_{2}}$ of Definition \ref{defmaleficmap} to the set $\big(\R_{\geq 0}\times\R_{\geq 0}\big) \backslash \{(0,0)\}$, that is we allow the case where $(c_{1},c_{2})$ is such that $c_{1}<c_{2}$ or $c_{i}=0$ for $i \in \{1,2\}$.
The image $F_{\underline{b}, \underline{d}, \alpha_{1},\alpha_{2}}(c_{1},c_{2})$ for a point $(c_{1},c_{2}) \in \overline{\mathfrak{a}}$ is obtained by computing the cross-ratio $R(A^{2},\Id,C,DC^{-1}D)$. In Theorem \ref{cpctarccoordthm} we have provided an explicit way to construct a hexagon $(H,\Yi) \in \mathcal{H}$ from arc coordinates $(\underline{b},\underline{c},\underline{d},[\alpha_{1},\alpha_{2}])$. More precisely we have shown how to compute positive definite symmetric matrices $A,C,D$ where $(H,\Yi)=(0,A,\Id,C,D,\infty)$. The explicit formulas appear in Proposition \ref{parforhex} for the generic case and are suitably adapted to the non-generic case in Proposition \ref{parnongentype1}, \ref{parnongentype2} and \ref{hexparinH2}. We extend these formulas to the case where $(c_{1},c_{2})$ is such that $c_{1}<c_{2}$ or $c_{i}=0$ for $i \in \{1,2\}$.

\begin{prop} \label{lemmalnottrue}
Let $\underline{b},\underline{d} \in \overline{\mathfrak{a}}$, $\alpha_{1},\alpha_{2} \in [0,2\pi)$ and let $\widetilde{F}$ denote the malefic map extended to $\big(\R_{\geq 0}\times\R_{\geq 0}\big) \backslash \{(0,0)\}$.   
It holds
$$
\widetilde{F}_{\underline{b}, \underline{d},\alpha_{1},\alpha_{2}}(c_{1},c_{2})=\widetilde{F}_{\underline{b}, \underline{d},\pi-\alpha_{1},\pi-\alpha_{2}}(c_{2},c_{1})
$$
Furthermore, if $(c_{1},c_{2})$ is a point lying on one of the semi-axis of
$\big(\R_{\geq 0}\times\R_{\geq 0}\big) \backslash \{(0,0)\}$ then $\widetilde{F}_{\underline{b}, \underline{d},\alpha_{1},\alpha_{2}}(c_{1},c_{2})$ is also lying on a semi-axis.
\end{prop}

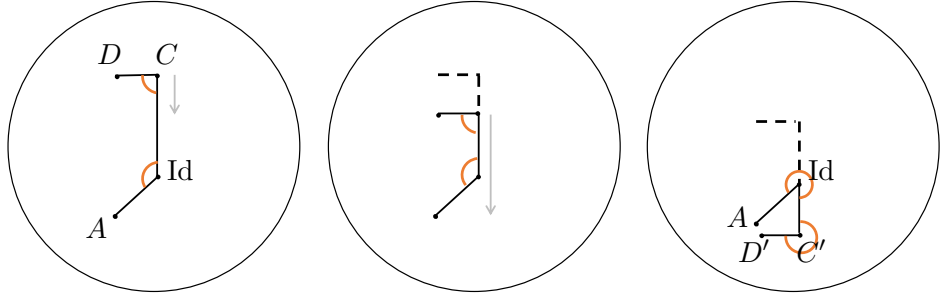
\begin{figure}[!h]
   \centering\captionsetup{justification=centering,margin=2cm}
   \setlength{\unitlength}{0.1\textwidth}
   \begin{picture}(7.5,2.3)
     \put(0,0){\includegraphics[width=12.5cm,height=4cm]{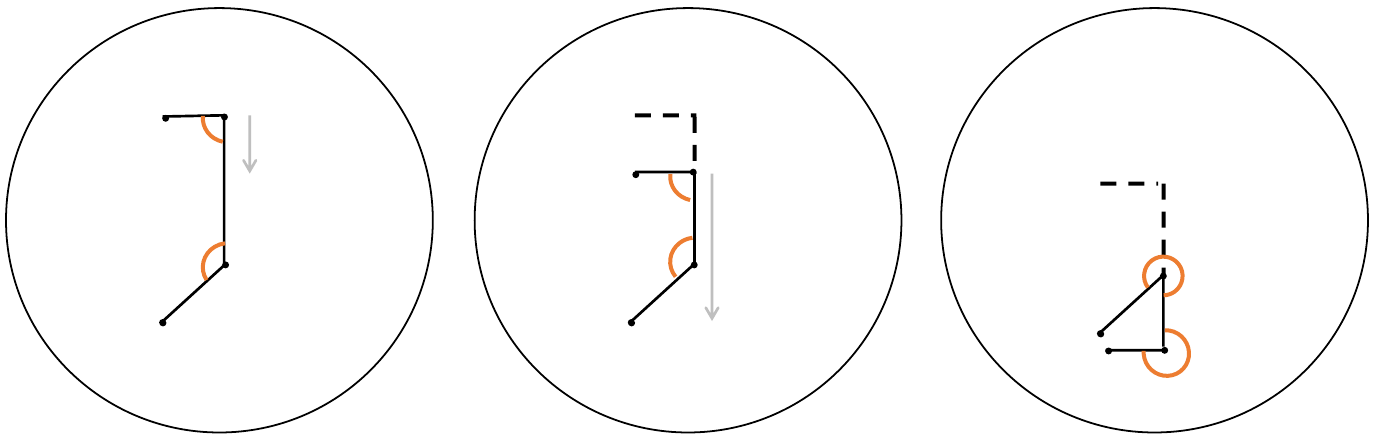}}
\put(1.4,1){$\Id$}
\put(0.7,0.5){$A$}
\put(1.3,2){$C$}
\put(0.8,2){$D$}
\put(6.8,0.3){$C'$}
\put(6.9,1){$\Id$}
\put(6.2,0.6){$A$}
\put(6.3,0.3){$D'$}
 \end{picture}
\caption{Continuous transformation of the polygonal chain when going from $\mathfrak{a}$ to $\mathfrak{a}^{-}$}
   \label{lemmal}
\end{figure}

\proof Let us understand the geometrical meaning of 
$\widetilde{F}_{\underline{b}, \underline{d},\alpha_{1},\alpha_{2}}(c_{1},c_{2})$ for a point $(c_{1},c_{2})$ with $c_{1}<c_{2}$. If in the parametrization of Proposition \ref{parforhex} we consider the set $\mathfrak{a}^{-}=\{ 0<x_{1}<x_{2}\}$ instead of $\mathfrak{a}=\{ x_{1}>x_{2}>0\}$ we are choosing to diagonalize the matrix $C$ with an increasing order of the eigenvalues. In the geometric interpretation of angle parameters of Section \ref{geomintdiag} the angle $\alpha$ denotes the angle from the semi-axis $\{(0,y)| \ y>1\} \in \H^{2}$. By picking the set $\mathfrak{a}^{-}$ we are considering the angle $\alpha+\pi$ when $\alpha \in [0,\pi)$ and the angle $\alpha-\pi$ when $\alpha \in [\pi,2\pi)$. From the equivalent relations of the angle parameters we know $\alpha+\pi \sim 2\pi-(\alpha+\pi)=\pi-\alpha$ and $\alpha-\pi \sim 2\pi-(\alpha-\pi)=\pi-\alpha$, so that $\widetilde{F}_{\underline{b}, \underline{d},\alpha_{1},\alpha_{2}}(c_{1},c_{2})=\widetilde{F}_{\underline{b}, \underline{d},\pi-\alpha_{1},\pi-\alpha_{2}}(c_{2},c_{1})
$. We should think at the extended map $\widetilde{F}_{\underline{b}, \underline{d},\alpha_{1},\alpha_{2}}$ as a way to construct right-angled hexagons in a continuous way by moving the point $C$. The polygonal chain of the hexagon is transformed as shown in Figure \ref{lemmal}.

\vspace{0.5cm}

\noindent Let us show that $\widetilde{F}_{\underline{b}, \underline{d},\alpha_{1},\alpha_{2}}$ preserves semi-axes. Let $(c_{1},c_{2})$ be such that $c_{1}=0$. This means $C=\bpm
e^{0} & 0\\
0 & e^{c_{2}}
\epm=\begin{pmatrix}
1 & 0\\
0 & \lambda
\end{pmatrix}$ for a $\lambda>0$, so that $C$ and $\Id$ are not transverse. Furthermore, there exists a $g \in \Sp(4,\R)$ such that
$g(A^{2},\Id,C,DC^{-1}D)=(0,\Id,M,\infty)
$ where $M$ is positive definite and such that $\Id$ and $M$ are not transverse. This means
$M=\begin{pmatrix}
1 & 0\\
0 & \mu
\end{pmatrix}$, $\mu>0
$ and we know $R(0,\Id,M,\infty)=M$ so that $\widetilde{F}_{\underline{b}, \underline{d},\alpha_{1},\alpha_{2}}(0,c_{2})=(\log(1),\log(\mu))=(0,y)$ for some $y>0$.
\endproof

\begin{cor}
Let $\underline{b}, \underline{d} \in \mathfrak{a}^{2}$ and let $F_{0}, F_{\pi}$ be the maps $F_{0}=F_{\underline{b},\underline{d},0,0}$, $F_{\pi}=F_{\underline{b},\underline{d},\pi,\pi}$
where $F_{\underline{b}, \underline{d},\alpha_{1},\alpha_{2}}$ denotes the malefic map. Then $F_{0}$ is not surjective and $F_{\pi}$ is not injective.
\end{cor}

\begin{figure}[!h]
\centering
\begin{tikzpicture} 
    \fill[green] (0,0) -- (45:2) arc (45:0:2) -- (0,0);
    \fill[yellow] (0,0) -- (45:2) arc (45:90:2) -- cycle;
    \draw[red] (0,0) -- (2,2);
    \draw[dashed] (0,0) -- (2,0);
    \draw[dashed] (0,0) -- (0,2);
\begin{scope}[shift={(3.5,0)}]
\fill[green] (0,0) -- (37:2) arc (37:0:2) -- (0,0);
    \fill[yellow] (0,0) -- (37:2) arc (37:90:2) -- cycle;
    \draw (0,0) -- (2,2);
    \draw[red] (0,0) -- (2,3/2);
    \draw[dashed] (0,0) -- (2,0);
    \draw[dashed] (0,0) -- (0,2);
     \end{scope}
    \draw[->, thick] (2.5,1) -- (3,1) node[midway, above]{$\widetilde{F}_{0}$};
\begin{scope}[shift={(0,-3.5)}]    
\fill[yellow] (0,0) -- (45:2) arc (45:0:2) -- (0,0);
    \fill[green] (0,0) -- (45:2) arc (45:90:2) -- cycle;
    \draw[red] (0,0) -- (2,2);
    \draw[dashed] (0,0) -- (2,0);
    \draw[dashed] (0,0) -- (0,2);
\end{scope}
\draw[->, thick] (0.5,-0.5) -- (0.5,-1) node[midway, right]{$c_{1}\leftrightarrow c_{2}$};  
\draw[->, thick] (2.5,-2.5) -- (4,-0.6) node[midway, right]{$\widetilde{F}_{\pi}$};
\end{tikzpicture}
\caption{$\widetilde{F}_{0}(c_{1},c_{2})=\widetilde{F}_{\pi}(c_{2},c_{1})$}
\label{tikz}
\end{figure}
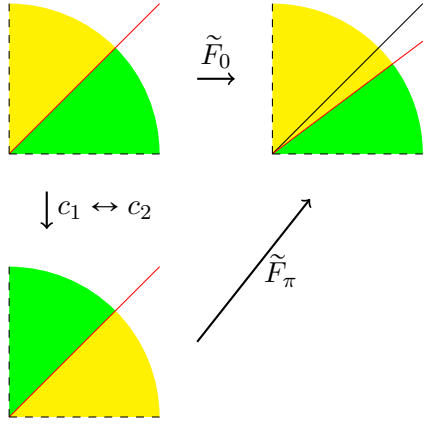

\proof
Consider the extended malefic maps $\widetilde{F}_{0}=\widetilde{F}_{\underline{b},\underline{d},0,0}$ and  $\widetilde{F}_{\pi}=\widetilde{F}_{\underline{b},\underline{d},\pi,\pi}$. The map $\widetilde{F}_{\underline{b}, \underline{d},\alpha_{1},\alpha_{2}}$ is continuous and from example \ref{malmapzpi} it is easy to see that $\widetilde{F}_{0}(x,x) \neq (X,X)$ (Figure \ref{tikz}). 

\noindent We deduce that when restricting to $\mathfrak{a}=\{x_{1}>x_{2}>0\}$ (i.e. considering the malefic map $F$) the map $F_{0}$ is not surjective and the map $F_{\pi}$ is not injective. This is illustrated in Figures \ref{plot1} and \ref{plot2} below. The program to generate these figures can be found in the github repository \url{https://github.com/martamagnani/Arc-coord/blob/main/Lemma_is_false.py}.  
\endproof

\begin{figure}[!h]
\centering
   \setlength{\unitlength}{0.1\textwidth}
   \begin{picture} (7,2.5)
      \put(0,0){\includegraphics[width=10cm, height=3.5cm]{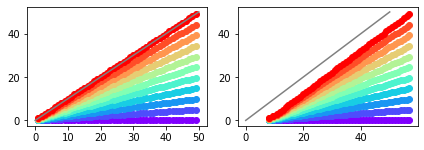}}
  \end{picture}
   \caption{The right-hand side shows the image of the map $F_{0}$ (not surjective) when $\underline{b}=(40,0.01)$ and $\underline{d}=(35,0.01)$}
   \label{plot1}
\end{figure}

\begin{figure}[!h]
\centering
   \setlength{\unitlength}{0.1\textwidth}
   \begin{picture} (7,2.5)
      \put(0,0){\includegraphics[width=10cm, height=3.5cm]{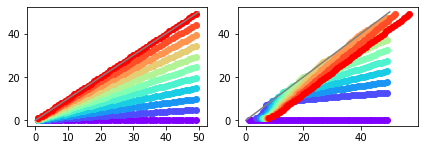}}
  \end{picture}
   \caption{The right-hand side shows the image of the map $F_{\pi}$ (not injective) when $\underline{b}=(40,0.01)$ and $\underline{d}=(35,0.01)$}
   \label{plot2}
\end{figure}

\newpage \begin{remark} \label{genwelldef} In Corollary \ref{genofquadr} we have seen how the genericity of the hexagon \\$(H,\Yi)=(0,A,\Id,C,D,\infty)$ induces the genericity of the quadruples $(-A,0,A^{2},A)$ and $(-D,0,C,D)$ respectively (see Figure \ref{nobij}). On the other hand, the quadruple $(Z_{1},\Id,C,Z_{2})$ is not necessarily generic. The extended map $\widetilde{F}$ of Proposition \ref{lemmalnottrue} does not preserve the diagonal $\mathfrak{d}=\{x_{1}=x_{2}\}$. By continuity of $\widetilde{F}$ we deduce that in general $(c_{1},c_{2})$ regular does not imply $F(c_{1},c_{2})$ regular. The parameters of Proposition \ref{parforhex} strongly depend on the order of the 6-tuple defining the hexagon or equivalently on the choice of a tube $\Yi$.
\end{remark}

\section{Reflections in the Siegel space} \label{refinsiegel}

In this section we study reflections in the Siegel space. We first recall properties of reflections in the hyperbolic plane $\H^{2}$ and we then generalize the results for the Siegel space $\X$. We define the notion of reflection set associated to the side of a hexagon, which will be used in the next section to parametrize maximal representations. 

\subsection{Reflections in $\H^{2}$} \label{secreflH2}

Let $\H^{2}$ denote the upper-half space model of the hyperbolic plane $$\H^{2}=\{x+iy| \ x,y \in \R, \ y>0\}$$
A reflection in $\H^{2}$ can be defined as a non-trivial isometry fixing an infinite geodesic $\gamma \in \H^{2}$. We propose an equivalent definition of reflection that will be generalized to define reflections in the Siegel space. 

\begin{definition}
Let $\SL^{-}(2,\R)$ be the set $\SL^{-}(2,\R):=\{M \in \GL(2,\R)| \det M=-1 \}$. The union $\SL(2,\R)\cup \SL^{-}(2,\R)$ forms a group that we denote $\SL^{\pm}(2,\R)$.
\end{definition}

\begin{definition}
A reflection in $\H^{2}$ is an involution of $\SL^{-}(2,\R)$.    
\end{definition}

We will denote by $\PSL^{\pm}(2,\R)$ the group $\SL^{\pm}(2,\R)/_{\{ \pm \Id\}}$. The action of $\PSL^{\pm}(2,\R)$ on $\H^{2}$ by Möbius transformations is not well defined for $R \in \SL^{-}(2,\R)$ as the resulting point may not lie in $\H^{2}$. To define an action of $\PSL^{\pm}(2,\R)$ on $\H^{2}$ let us denote by $\H^{2}_{\pm}$ the extended hyperbolic plane 
$$
\H^{2}_{\pm}=\{ x \pm iy| \ x,y \in \R, \ y>0\}
$$
so that $\H^{2}=\H^{2}_{\pm}/_{\sim}$ where $x+iy \sim x-iy$. A matrix $R \in \SL^{-}(2,\R)$ acts on $\H^{2}$ through Möbius transformations in the following way 
$$
R\cdot z := [R\cdot z] \in \H^{2}_{\pm}/_{\sim}
$$

\begin{lem}
All reflections of $\H^{2}$ are conjugated by an element of $\SL(2,\R)$.    
\end{lem}

\proof The proof is given in the more general case in Lemma \ref{reflconj} where it is shown for the Siegel space $\X$ and the group $\Sp(2n,\R)$. The proof for $\H^{2}$ is the case $n=1$.
\endproof

\begin{definition}
We call the \emph{standard reflection} in $\H^{2}$ the map $r=\bpm-1&0\\0&1\epm$.
\end{definition}

\begin{prop} \label{propreflH2_1}
Let $R$ be a reflection in $\H^{2}$. Then $R$ fixes exactly two boundary points $p,q \in \partial \H^{2}$. Moreover, $R$ fixes the infinite geodesic $\gamma_{p,q}$ that has $p,q$ as endpoints.    
\end{prop}

\proof The proof is given in the general case in Proposition \ref{reflfixp}.
\endproof

\begin{prop} \label{propreflH2_2}
Given $p,q \in \partial \H^{2}$, there is a unique reflection $R$ fixing both $p$ and $q$. The map $R$ is an isometry sending any boundary point $x$ to the unique boundary point $R(x)$ such that $\gamma_{p,q} \perp \gamma_{x,R(x)}$
\end{prop}

\proof The proof is given in the general case in Proposition \ref{lemnumberofrefl}.
\endproof

\begin{prop} \label{reflinH2}
Let $(q_{1},q_{2},q_{3},q_{4})$ be a positive quadruple in $\partial \H^{2}$ and let $R$ be the reflection fixing two boundary points $p_{1},p_{2} \in \partial \H^{2}$. If $(p_{2},q_{1},q_{2},q_{3},q_{4},p_{1}) \text{ is positive (possibly } p_{2}=q_{1} \text{ or } p_{1}=q_{4})$,
then $(p_{1},R(q_{4}),R(q_{3}),R(q_{2}),R(q_{1}),p_{2})$ is positive.
\end{prop}

\proof The proof is given in the general case in Proposition \ref{propreflinsiegel}
\endproof

\subsection{Reflections in $\X$}

\begin{definition} \label{defantisympl}
Let $\omega(\cdot,\cdot)$ be the symplectic form represented, with respect to the standard basis, by the matrix
$J_{n}=
\bpm
0 & \Id_{n}\\
-\Id_{n} & 0 
\epm
$. A matrix $M \in \GL(2n,\R)$ is  \emph{antisymplectic} if
$M^{T}J_{n}M=-J_{n}$. The set of antisymplectic matrices will be denoted by $\Sp^{-}(2n,\R)$. More precisely $\Sp^{-}(2n,\R)$ is the set 
\begin{equation*}
\Sp^{-}(2n,\R)=\Big\{  \bpm
A & B\\
C & D 
\epm | \ A^{T}C,B^{T}D \text{ symmetric, and  } A^{T}D-C^{T}B=-\Id_{n}\Big\} 
\end{equation*} 

\noindent The union of symplectic and antisymplectic matrices forms a group that will be denoted by $\Sp^{\pm}(2n,\R)$.
\end{definition}

\begin{definition} \label{reflection}
    A \emph{reflection} $R$ in $\X$ is an antisymplectic involution of $\X$.
\end{definition}

We will denote by $\PSp^{\pm}(2n,\R)$ the group $\Sp^{\pm}(2n,\R)/_{ \{\pm \Id \}}$. The action of $\PSp^{\pm}(2n,\R)$ on $\X$ by fractional linear transformations is not well defined for $R \in \Sp^{-}(2n,\R)$ as the resulting point may not lie in $\X$. To define an action of $\PSp^{\pm}(2n,\R)$ on $\X$ let us denote by $\X^{\pm}$ the extended Siegel space:
$$
\X^{\pm}= \{ X \pm iY | \ X \in \Sym(n,\R), Y \in \Sym^{+}(n,\R) \}
$$
Then $\X=\X^{\pm}/_{\sim}$ where $X+iY \sim X-iY$. For $R$ antisymplectic and $Z \in \X$ we define the action
$$
R\cdot Z := [R\cdot Z] \in \X^{\pm}/_{\sim}
$$

\begin{lem} \label{reflconj}
All reflections of $\X$ are conjugated by an element of $\Sp(2n,\R)$.
\end{lem}

\proof Let $R$ be a reflection of $\X$. Since $R$ is an involution we know that its eigenvalues are given by the set $\{ \pm 1 \}$. Recall that we denote by $\mathcal{L}(\R^{2n})^{(k)}$ the set of $k$-tuples of real pairwise transverse Lagrangians. As $\Sp(2n,\R)$ acts transitively on pairs of transverse Lagrangians, it is sufficient to show that the $R$-eigenspaces $E_{1},E_{-1} $ are inside $\mathcal{L}(\R^{2n})^{(2)}$.
For $u,v \in E_{1}$ it holds $\omega(u,v)=\omega(R(u),R(v))=-\omega(u,v)
$ where the first equality holds since $u,v \in E_{1}$ and the second one since $R$ is antisymplectic. It follows that $\omega(u,v)=0$ for any $u,v \in E_{1}$, that is $E_{1}$ is a Lagrangian subspace.
Similarly one can show that  $E_{-1}$ is also a Lagrangian subspace.
Since a real matrix with real eigenvalues has real eigenvectors, we conclude that $E_{1},E_{-1} \in \mathcal{L}(\R^{2n})^{(2)}$.
\endproof

\begin{definition}
We call the \emph{standard reflection} the map $R_{st}=\bpm -\Id&0\\0&\Id\epm$.
\end{definition}

\begin{lem}
Let $R$ be a reflection inside $\X$. Then  for any $X,Z$ in $\X$ it holds $$d^{\a}(R(X),R(Z))=d^{\a}(X,Z)$$
\end{lem}

\proof This follows immediately from the definition of $d^{\a}$ (Definition \ref{defweyldist}) and from the properties of the cross-ratio (Definition \ref{defcrossratio}).
\endproof

\begin{prop} \label{reflfixp}
Let $\X$ be the symmetric space associated to $\Sp(4,\R)$ and let $R$ be a reflection of $\X$. Then the set
$$Fix_{\mathcal{L}(\R^{4})}(R)=\{l \in \mathcal{L}(\R^{4}) | \ R(l)=l \}$$
is given by the $R$-eigenspaces $E_{1},E_{-1}$ together with an $S^{1}$-isomorphic family $\mathcal{F}$ of pairwise transverse Lagrangians each of which is not transverse to neither $E_{1}$ nor $E_{-1}$. Moreover, $R$ fixes the tube $\Y_{E_{1},E_{-1}}$ and fixes a flat inside any $Y_{l_{i},l_{j}}$ where $l_{i},l_{j} \in \mathcal{F}$.
\end{prop}

\proof Since any $R$ is conjugated to the standard reflection through an element of $\Sp(4,\R)$, let us prove the proposition for 
$R_{st}=\bpm -\Id&0\\0&\Id\epm$. Let $(e_{1},e_{2},e_{3},e_{4})$ denote the standard basis of $\R^{4}$. We have
$E_{1}=\langle e_{3},e_{4}\rangle$ and $E_{-1}=\langle e_{1},e_{2}\rangle$
where $E_{1},E_{-1} \in \mathcal{L}(\R^{4})^{(2)}$. For any $u\in \mathbb{P}(E_{1})$ there exists a unique $v\in \mathbb{P}(E_{-1})$ such that $\omega(u,v)=0$. For any $w_{1},w_{2} \in \langle u,v \rangle$ it holds $\omega(w_{1},w_{2})=0$ so that $l=\langle u,v \rangle$ is in $\mathcal{L}(\R^{4})$ and $R_{st}(l)=l$. We obtain the set $\mathcal{F} \subset Fix_{\mathcal{L}(\R^{4})}(R_{st})$:
$$\mathcal{F}= \big\{ l=\langle u,v \rangle\big\} \cong  \mathbb{P}(E_{1}) \cong \mathbb{P}(E_{-1})\cong S^{1} $$
\noindent where $u \in \mathbb{P}(E_{1})$ and $v$ is the unique element of $\mathbb{P}(E_{-1})$ such that $\omega(u,v)=0$.

\noindent We want to show that for any $l \in \mathcal{F}$ it holds $E_{1} \notpitchfork l \notpitchfork E_{-1}$. Fix $\alpha \in \R$ and consider $u \in \mathbb{P}(E_{1})$, $u \neq e_{4}$ to be the vector $u=e_{3}+\alpha e_{4}$. Then the corresponding $v$ in $ \mathbb{P}(E_{-1}), \ v \neq e_{1}$ such that $\omega(u,v)=0$ is given by $v=-\alpha e_{1}+e_{2}$. 
Let $l=\langle u,v \rangle=\langle e_{3}+\alpha e_{4}, -\alpha e_{1}+e_{2} \rangle \in \mathcal{F} \backslash \{ \langle e_{4},e_{1}\rangle \}$, then $l$ intersects $E_{1}$ in the line $\langle e_{3}+\alpha e_{4} \rangle \subset E_{1}$  and intersects $E_{-1}$ in the line $\langle -\alpha e_{1}+e_{2}  \rangle \subset E_{-1}$. We are left with the case $l=\langle e_{4},e_{1} \rangle$ which is clearly not transverse to $E_{1}$ nor $E_{-1}$. We have showed $E_{1} \notpitchfork l \notpitchfork E_{-1}$ for every $l \in \mathcal{F}$, we are left to show that for any $l_{1},l_{2}$ inside $\mathcal{F}$  it holds $l_{1} \pitchfork l_{2}$. Similarly as before let $l_{1},l_{2} \in \mathcal{F} \backslash \{ \langle e_{4},e_{1}\rangle \}$ where
$l_{1}=\langle u_{1},v_{1} \rangle=\langle e_{3}+\alpha e_{4}, -\alpha e_{1}+e_{2} \rangle$ and 
$l_{2}=\langle u_{2},v_{2} \rangle=\langle e_{3}+\beta e_{4}, -\beta e_{1}+e_{2} \rangle$
for $\alpha,\beta \in \R, \ \beta \neq \alpha$. It is easy to see that $l_{1} \pitchfork l_{2}$. It is also trivial to show that transversality holds in the case $l_{1} =\langle e_{4},e_{1}\rangle$.\\
The reflection $R_{st}$ is fixing the tube $\Y_{E_{1},E_{-1}}$: to see this recall that the affine chart $\iota$ in Section \ref{defsiegel} identifies $\mathcal{L}(\mathbb{C}^{4})$ with $\Sym(2,\mathbb{C})$. In this chart the Lagrangian $E_{-1}=\langle e_{1},e_{2} \rangle$ is the point at infinity in the Shilov boundary $\mathcal{L}(\R^{4})$ of $\X$ and the expression for the tube $\Y_{E_{1},E_{-1}}$ is given by the standard tube $\Yi=\{ iY| \ Y \in \Sym^{+}(2,\R)\}$. For any $iY \in \Yi$ it holds 
$R_{st}(iY)=-iY=iY$ in $\X^{\pm}/_{\sim}=\X$. Observe that with respect to the tube $\Y_{E_{1},E_{-1}}$ the reflection $R_{st}$ is the analogue of a reflection in $\H^{2}$: it is sending any boundary point $X \in \Sym(2,\R)$ (transverse to both $E_{1}$ and $E_{-1}$) to the unique $R_{st}(X)=-X$ such that $\Yi \perp \Y_{X,R_{st}(X)}$.
Let us now consider $l_{1}, l_{2} \in \mathcal{F}$ where $l_{1}=\langle e_{1}, e_{4} \rangle$ and $l_{2}=\langle e_{2}, e_{3} \rangle$. Let us change the standard basis $\mathcal{B}=(e_{1},e_{2},e_{3},e_{4})$ with the basis $\mathcal{B}'=(e_{3},e_{2},e_{1},e_{4})$. Writing vectors of $\mathbb{C}^{4}$ in this new basis means considering the chart $T \circ \iota :\Sym(2,\mathbb{C}) \to \mathcal{L}(\mathbb{C}^{4})$ where $T(\mathcal{B})=\mathcal{B}'$. In particular in this chart the tube $\Y_{l_{1},l_{2}}$ has the standard form
$\Y_{l_{1},l_{2}}= \{ iY| \ Y \in \Sym^{+}(2,\R)\}$ and the reflection $R_{st}$ written in basis $\mathcal{B}'$ is given by $\widetilde{R}=R_{\mathcal{B}'}=\bpm -r&0\\0&r\epm$ where $r=\bpm -1&0\\0&1\epm$. It holds $\widetilde{R}(iY)=-iY^{r}=iY^{r}$ in $ \X^{\pm}/_{\sim}=\X$
where by $Y^{r}$ we denote the point in the $\H^{2}$-component of the standard tube which is obtained by reflecting $Y$ across the standard vertical geodesic of the hyperbolic plane (see Section \ref{refl}). The reflection $\widetilde{R}$ fixes the flat $\mathbb{D}$ where
$$\mathbb{D}= i \bpm d_{1}&0\\0&d_{2}\epm \cong \R \times \bpm \lambda&0\\0&\frac{1}{\lambda}\epm \subset \Yi $$
The reflection $\widetilde{R}$ is reflecting across a geodesic $\gamma$ in the $\H^{2}$-component of the tube $\Y_{l_{1},l_{2}}$ and is therefore fixing the flat $\R \times \gamma$ inside the tube.  Since $\Sp(4,\R)$ acts transitively on the space of transverse Lagrangians we deduce that the same holds for any $l_{1},l_{2} \in \mathcal{F}$.
\endproof

\begin{cor}
There is no maximal triple in the $S^{1}$-family $\mathcal{F}$ of Proposition \ref{reflfixp}.
\end{cor} 
\proof
Let $(e_{1},e_{2},e_{3},e_{4})$ be the standard basis in $\R^{4}$ and as usual let us denote by $l_{\infty},0$ and $\Id$ the Lagrangians
$l_{\infty}=\langle e_{1},e_{2}\rangle$, $0=\langle e_{3},e_{4} \rangle$ and $\Id=\langle e_{1}+e_{3},e_{2}+e_{4} \rangle$ respectively. Since any reflection is conjugated to the standard one, we prove the result for the standard reflection $R_{st}=\bpm -\Id&0\\0&\Id  \epm$. In the proof of Proposition \ref{reflfixp}) we have seen that for $R_{st}$ we have $E_{1}=0$ and $E_{-1}=l_{\infty}$. Each Lagrangian of $\mathcal{F}$ intersects $l_{\infty}$ and $0$ in one line and $\mathcal{F} \cong \mathbb{P}(E_{1}) \cong \mathbb{P}(E_{-1})$. Let $l_{1},l_{2},l_{3}$ be three points in $\mathcal{F}$. Up to $\GL(2,\R) \cong \Stab(E_{1},E_{-1})$-action we can choose the three vectors of $\mathbb{P}(E_{-1})$ to be $e_{1},e_{2}$ and $e_{1}+e_{2}$ respectively ($\GL(2,\R)$ acts three-transitively on the lines of $\R^{2}$) and we obtain
$l_{1}=\langle e_{1},e_{4} \rangle$, $ l_{2}=\langle e_{2},e_{3} \rangle$ and $ l_{3}=\langle e_{1}+e_{2},e_{3}-e_{4} \rangle$. Let $g \in \Sp(4,\R)$ be such that $g(l_{1},l_{2})=(l_{\infty},0)$. Then
$$g=\bpm A&0\\0&A^{-T}\epm \circ \bpm 1&0&0&0\\0&0&0&-1\\0&0&1&0\\0&1&0&0 \epm$$

\noindent and let us choose for simplicity $A=\Id$. Then $gl_{3}= \langle e_{2}+e_{3}, e_{1}+e_{4} \rangle$ which corresponds to the matrix $M=\bpm 0&1\\1&0\epm$ in the identification of Section \ref{defsiegel}. The triple $(l_{\infty},0,M)$ is not maximal as its Maslov index is zero (see Section \ref{Maslov} for the definition of Maslov index).
\endproof

In Proposition \ref{reflfixp} we have seen that for a given reflection $R \in \PSp(4,\R)^{-}$ there is a different geometrical behaviour when considering what $R$ is doing with respect to the tube $\Y_{E_{1},E_{-1}}$ or to the tube $\Y_{l_{1},l_{2}}$, where $l_{1},l_{2}$ are two arbitrary points inside $\mathcal{F}$.

\begin{definition}  \label{reflexotic}
The reflection $R_{ex}=\bpm -r&0\\0&r\epm$ where $r=\bpm -1&0\\0&1\epm$ will be called the \emph{exotic reflection}.
\end{definition}

\subsection{Reflection set associated to the side of a hexagon} \label{reflsetsection}

\begin{prop} \label{lemnumberofrefl}
Let $\X$ be the symmetric space associated to $\Sp(4,\R)$. Let $(P,X,Y,Q)$ be a maximal quadruple in $\X$ and let $g$ be an isometry such that $g(P,X,Y,Q)=(0,\Id,Y',\infty)$ for $Y'$ diagonal. Let $\mathcal{R}(P,X,Y,Q) \subset \PSp^{\pm}(4,\R)$ be the set defined by
$$\mathcal{R}(P,X,Y,Q):=\{R \text{ reflection } | \ R(P)=P, \ R(Q)=Q \text{ and } \Y_{X,R(X)} \perp \Y_{P,Q} \perp \Y_{Y,R(Y)} \}$$
It holds:
\begin{itemize}
\item[(i)] If $(P,X,Y,Q)$ is generic then $\mathcal{R}(P,X,Y,Q)= \{g^{-1}R_{st}g,g^{-1}R_{ex}g  \}$
\item[(ii)] If $(P,X,Y,Q)$ is non-generic then $\mathcal{R}(P,X,Y,Q)=g^{-1}\mathcal{K}g$ where
\end{itemize}

$$ \mathcal{K}=\Big\{\bpm -K&0\\0&K\epm, \  K \in \PO(2), \ K^{2}=\Id \Big\}$$
\end{prop}

 \begin{figure}[!h]
   \centering
   \captionsetup{justification=centering,margin=2cm}
   \setlength{\unitlength}{0.1\textwidth}
   \begin{picture}(3,3)
     \put(0,0){\includegraphics[width=4cm,height=4cm]{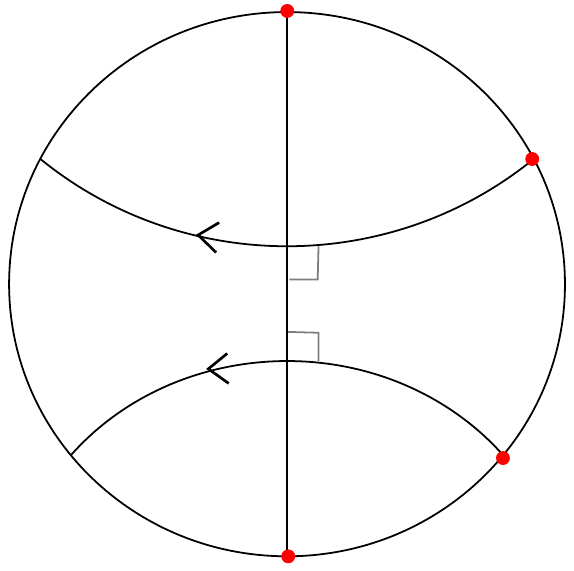}}
          \put(0.8,-0.2){$P=R(P)$}
          \put(0.8,2.6){$Q=R(Q)$}
          \put(2.5,1.8){$Y$}
\put(-0.45,1.8){$R(Y)$}
\put(2.4,0.4){$X$}
\put(-0.4,0.4){$R(X)$}
\end{picture} 
\vspace{3mm}
\caption{The set $\mathcal{R}(P,X,Y,Q)$ depends on the genericity of the quadruple $(P,X,Y,Q)$}
\label{howmanyrefl}
\end{figure}
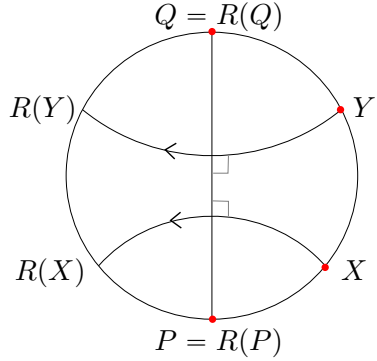
 
\proof Let $g$ be an isometry such that $g(P,X,Y,Q)=(0,\Id,Y',\infty)$ for $Y'$ diagonal. As 
$\mathcal{R}(P,X,Y,Q)=g^{-1}\mathcal{R}(0,\Id,Y',\infty)g$, it is sufficient to show the proposition for $(0,\Id,Y',\infty)$. Let $R$ be a reflection such that
\begin{equation*} 
\begin{cases}
R(0)=0 \text{ and } R(\infty)=\infty\\
\Y_{\Id,R(\Id)} \perp \Y_{0,\infty} \perp \Y_{Y',R(Y')}
 \end{cases}  
\end{equation*}
We obtain two possibilities for $\mathcal{R}(0,\Id,Y',\infty)$:\\
\emph{(i)} If $(0,\Id,Y',\infty)$ is generic then $Y'$ is of the form $Y'=\bpm y_{1}&0\\0&y_{2}\epm$ where $y_{1}\neq y_{2}$. Then $\mathcal{R}(0,\Id,Y',\infty)=\{R_{st},R_{ex}\}$. The reflection $R_{st}$ fixes the tube $\Yi$ (as $E_{1}^{R_{st}}=0, E_{-1}^{R_{st}}=l_{\infty}$) whereas $R_{ex}$ fixes a flat inside $\Yi$ (see Proposition \ref{reflfixp}).\\
\emph{(ii)} If $(0,\Id,Y',\infty)$ is non-generic then $Y'=\bpm y&0\\0&y\epm=y\Id$. Let  $R=\bpm -K&0\\0&K\epm$ for a $K \in \PO(2)$. Then $R$ is antisymplectic and it holds $R(\Id)=-\Id$ so that $\Y_{0,\infty} \perp \Y_{\Id,R(Y')}$. Similarly
$R(Y')=-Y'$ so that $\Y_{0,\infty} \perp \Y_{Y',R(Y')}$. We further need $R^{2}=\Id$ for $R$ to be an involution which is satisfied exactly when $K^{2}=\Id$. 
\endproof

\begin{prop} \label{propreflinsiegel}
Let $(l_{1},l_{2},l_{3},l_{4})$ be a maximal quadruple.  Let $R$ be a reflection inside $\mathcal{R}(P,X,Y,Q)$ where $(P,X,Y,Q)$ is a maximal quadruple. Suppose $(X,l_{1},l_{2},l_{3},l_{4},Y)$ maximal (possibly $X=l_{1}$ or $Y=l_{4}$), then $(Q,R(Y),R(l_{4}),R(l_{3}),R(l_{2}),R(l_{1}),R(X),P)$ is maximal.
\end{prop}

\proof Let $g$ be an isometry such that $g(P,X,Y,Q)=(0,\Id,Y',\infty)$ for $Y'$ diagonal. We want to show that the image 
\begin{equation} \label{toshowmax}
\big(R(\infty),R(Y'),R(l_{4}),R(l_{3}),R(l_{2}),R(l_{1}),R(\Id),R(0)\big)
\end{equation}
is maximal for $(\Id,l_{1},l_{2},l_{3},l_{4},Y')$ maximal and $R \in \mathcal{R}(0,\Id,Y',\infty)$. 
By Proposition \ref{lemnumberofrefl}:
\begin{itemize}
\item[(i)] $\mathcal{R}(0,\Id,Y',\infty)= \{R_{st}R_{ex}\}$ if $Y'=\bpm y_{1}&0\\0&y_{2}\epm, y_{1}\neq y_{2}$
\item[(ii)] $\mathcal{R}(0,\Id,Y',\infty)= \mathcal{K}$ if $Y'=\bpm y&0\\0&y\epm$
\end{itemize}
Observe that $\{R_{st}R_{ex}\}\subset \mathcal{K}$. Using Lemma \ref{max} it is not hard to show that (\ref{toshowmax}) is maximal.
\endproof

\begin{definition}\label{reflsetside} 
Let $H=[\Y_{1},\Y_{2},\Y_{3},\Y_{4},\Y_{5},\Y_{6}]$ be a right-angled hexagon. The \emph{reflection set $\mathcal{R}_{\Y_{k}}^{\Y_{k-1},\Y_{k+1}}\subset \PSp^{\pm}(4,\R)$ associated to} $\Y_{k}$ is the set of reflections which are fixing the endpoints of $\Y_{k}$ and are switching the endpoints of $\Y_{k-1}$ and $\Y_{k+1}$ respectively.  
\end{definition}

\begin{cor} \label{reflsetcoro}
Let $H=[\Y_{1},\Y_{2},\Y_{3},\Y_{4},\Y_{5},\Y_{6}]$ be a right-angled hexagon.\\ Let $\Y_{k-1}=\Y_{P_{1},P_{2}}$, $\Y_{k}=\Y_{Q_{1},Q_{2}}$, $\Y_{k+1}=\Y_{P_{4},P_{5}}$ and let $g$ be an isometry such that \\$g(Q_{1},P_{2},P_{3},Q_{2})=(0,\Id,Y',\infty)$ for $Y'$ diagonal.
Then
\begin{equation*}
\mathcal{R}_{\Y_{k}}^{\Y_{k-1},\Y_{k+1}}=\begin{cases}
\{g^{-1}R_{st}g,g^{-1}R_{ex}g\}, \ if \   (Q_{1},P_{2},P_{3},Q_{2}) \  generic\\
g^{-1}\mathcal{K}g \ if \   (Q_{1},P_{2},P_{3},Q_{2}) \  non \ generic
\end{cases}    
\end{equation*}
\end{cor}

\proof
Follows directly from Proposition \ref{lemnumberofrefl}.
\endproof

\noindent We can rewrite Corollary \ref{reflsetcoro} in terms of arc coordinates.

\begin{cor} \label{reflsetwitharcoords}
Let $H=[\Y_{1},\Y_{2},\Y_{3},\Y_{4},\Y_{5},\Y_{6}]$ be a right-angled hexagon where 
$$\Y_{1}=\Y_{P_{1},P_{2}}, \ \Y_{2}= \Y_{Q_{1},Q_{2}}, \ \Y_{3}=  \Y_{P_{3},P_{4}} , \ \Y_{4}= \Y_{Q_{3},Q_{4}} , \ \Y_{5}= \Y_{P_{5},P_{6}},\ \Y_{6}= \Y_{Q_{5},Q_{6}}
$$
Suppose $(H,\Y_{1})$ has arc coordinates $\mathcal{A}(H,\Y_{1})=(\underline{b},\underline{c},\underline{d},[\alpha_{1},\alpha_{2}])$. Let $F_{\underline{b}, \underline{d}, \alpha_{1},\alpha_{2}}$ be the malefic map of Definition \ref{defmaleficmap} and let $g_{1},g_{2},g_{3}$ be isometries such that
$$g_{1}(Q_{1},P_{2},P_{3},Q_{2})=(0,\Id,Y_{1},\infty)$$
$$g_{2}(Q_{5},P_{6},P_{1},Q_{6})=(0,\Id,Y_{1},\infty)$$
$$g_{3}(Q_{3},P_{4},P_{5},Q_{4})=(0,\Id,Y_{3},\infty)$$
where $Y_{1},Y_{2},Y_{3}$ are diagonal matrices.
It holds
\begin{equation} \label{firstreflset}
\mathcal{R}_{\Y_{2}}^{\Y_{1},\Y_{3}}=\begin{cases}
\{g_{1}^{-1}R_{st}g_{1},g_{1}^{-1}R_{ex}g_{1}\}, \ if \ \underline{b} \in \mathfrak{a}\\
g_{1}^{-1}\mathcal{K}g_{1}   \ if \ \underline{b} \in \mathfrak{d}
\end{cases}    
\end{equation}

\begin{equation} \label{secondreflset}
\mathcal{R}_{\Y_{6}}^{\Y_{5},\Y_{1}}=\begin{cases}
\{g_{2}^{-1}R_{st}g_{2},g_{2}^{-1}R_{ex}g_{2}\}, \ if \ \underline{d} \in \mathfrak{a}\\
g_{2}^{-1}\mathcal{K}g_{2}    \ if \ \underline{d} \in \mathfrak{d}
\end{cases}    
\end{equation}

\begin{equation}  \label{thirdreflset}
\mathcal{R}_{\Y_{4}}^{\Y_{3},\Y_{5}}=\begin{cases}
\{g_{3}^{-1}R_{st}g_{3},g_{3}^{-1}R_{ex}g_{3}\}, \ if \ F_{\underline{b}, \underline{d}, \alpha_{1},\alpha_{2}}(\underline{c}) \in \mathfrak{a}\\
g_{3}^{-1}\mathcal{K}g_{3}    \ if \ F_{\underline{b}, \underline{d}, \alpha_{1},\alpha_{2}}(\underline{c}) \in \mathfrak{d}
\end{cases}    
\end{equation}

\end{cor}

\begin{figure}[!h]
   \centering
   \captionsetup{justification=centering,margin=2cm}
   \setlength{\unitlength}{0.1\textwidth}
   \begin{picture}(3.5,3.5)
     \put(0.0,0){\includegraphics[width=5cm,height=4.7cm]{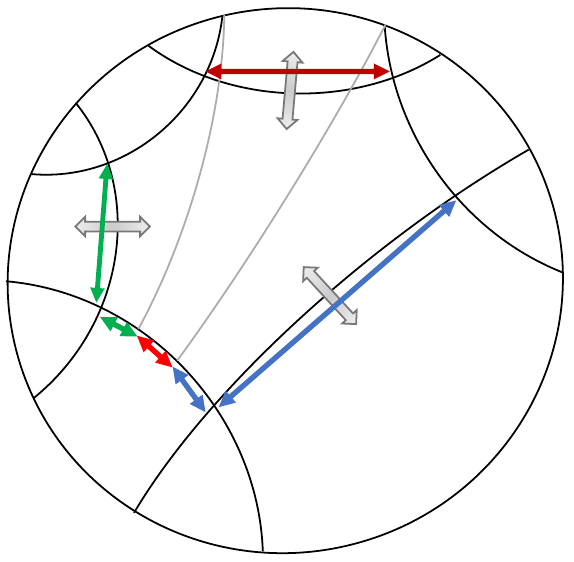}}
     \put(1.4,-0.26){$P_{2}$}
          \put(0.7,0){$Q_{1}$}
          \put(-0.1,0.8){$Q_{6}$}
          \put(-0.2,1.4){$P_{1}$}
          \put(0.1,2.5){$Q_{5}$}
          \put(0.7,2.9){$Q_{4}$}
          \put(1.1,3){$P_{5}$}
          \put(-0.1,2.1){$P_{6}$}
          \put(2.1,3){$P_{4}$}
          \put(2.4,2.8){$Q_{3}$}
          \put(3.1,2.2){$Q_{2}$}
          \put(3.2,1.6){$P_{3}$}
    \put(2.2,1.5){$\Y_{2}$}
    \put(1.4,0.35){$\Y_{1}$}
    \put(2.4,2.3){$\Y_{3}$}
    \put(1.8,2.6){$\Y_{4}$}
    \put(0.7,2.4){$\Y_{5}$}
    \put(0.35,1.55){$\Y_{6}$}
    \put(0.95,0.75){\textcolor{blue}{$\underline{b}$}}
    \put(0.75,1){\textcolor{red}{$\underline{c}$}}
    \put(0.5,1){\textcolor{forestgreen}{$\underline{d}$}}
    \put(1.8,1.1){\textcolor{gray}{$R_{1}$}}
    \put(0.9,1.8){\textcolor{gray}{$R_{3}$}}
    \put(1.4,2.2){\textcolor{gray}{$R_{2}$}}
 \end{picture}
  \vspace{0.3cm} 
\caption{$R_{1},R_{2},R_{3}$ are in $\mathcal{R}_{\Y_{2}}^{\Y_{1},\Y_{3}},\mathcal{R}_{\Y_{4}}^{\Y_{3},\Y_{5}},\mathcal{R}_{\Y_{6}}^{\Y_{5},\Y_{1}}$ respectively}
   \label{reflsetwarc}
\end{figure}

\proof Let us prove (\ref{firstreflset}). By Corollary \ref{reflsetcoro} we know 
\begin{equation*}
\mathcal{R}_{\Y_{2}}^{\Y_{1},\Y_{3}}=\begin{cases}
\{g_{1}^{-1}R_{st}g_{1},g_{1}^{-1}R_{ex}g_{1}\}, \ if \   (Q_{1},P_{2},P_{3},Q_{2}) \  generic\\
g_{1}^{-1}\mathcal{K}g_{1} \ if \   (Q_{1},P_{2},P_{3},Q_{2}) \  non \ generic
\end{cases}    
\end{equation*}
The quadruple $(Q_{1},P_{2},P_{3},Q_{2})$ is generic if the matrix given by the cross-ratio $Cr(Q_{1},P_{2},P_{3},Q_{2})$ has distinct eigenvalues $\mu_{1},\mu_{2}$ and non-generic if $\mu_{1}=\mu_{2}$. By Lemma \ref{proj} it holds $(\mu_{1},\mu_{2})=d^{\a}(p_{Q_{1},Q_{2}}(P_{2}),p_{Q_{1},Q_{2}}(P_{3}))$. The vector $(\mu_{1},\mu_{2})$ is the image $f(\underline{b})$ where $f$ is the bijective map of Proposition \ref{bijmap1}. These vectors are drawn in blue in Figure \ref{reflsetwarc}. In particular $f$ preserves regular vectors. Equality (\ref{firstreflset}) follows and the proof for (\ref{secondreflset}) is similar. For (\ref{thirdreflset}) we need to write $F_{\underline{b}, \underline{d}, \alpha_{1},\alpha_{2}}(\underline{c})$ instead of $\underline{c}$ as there is no bijective map as in the cases (\ref{firstreflset}) and (\ref{secondreflset}). This is explained in Section \ref{constrains}.
\endproof

\subsection{Geometrical interpretation of the set $\mathcal{K}$} \label{sectiongeomintK}

In this section we give a geometrical interpretation to the reflection set associated to the side of a hexagon (and in particular to the set $\mathcal{K}$) in terms of the polygonal chain associated to the hexagon. For simplicity we will consider an ordered hexagon of the form $(H,\Yi)$ and will study the associated polygonal chain defined in \ref{defpolchain}. 

\begin{definition}
Two right-angled hexagons $H_{1},H_{2}$ are said to be \emph{adjacent at $\Y_{1}$} if\\ 
$H_{1}=[\Y_{1},\Y_{2},\Y_{3},\Y_{4},\Y_{5},\Y_{6}]$ and $H_{2}=[\Y_{1},\Y_{6},\Y_{7},\Y_{8},\Y_{9},\Y_{2}]$. Two such adjacent hexagons will be denoted $H_{1} \#_{\Y_{1}} H_{2}$.
\end{definition}

When two hexagons are adjacent they share one vertex of the correspondent polygonal chains and we can look at the "attached" polygonal chains. 

\begin{definition} Let $(H,\Yi)$ and $(H',\Yi)$ be two adjacent right-angled hexagons where
$$(H,\Yi)=(0,A,B,C,D,\infty)$$
$$(H',\Yi)=(0,A',B',C',D',\infty), \ D'=A$$

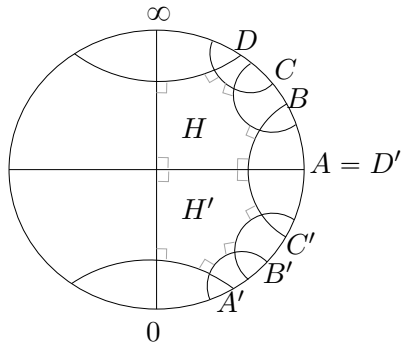
\begin{figure}[!h]
   \centering
   \captionsetup{justification=centering,margin=2cm}
   \setlength{\unitlength}{0.1\textwidth}
   \begin{picture}(3,2.6)
     \put(0.0,0){\includegraphics[width=4cm,height=3.8cm]{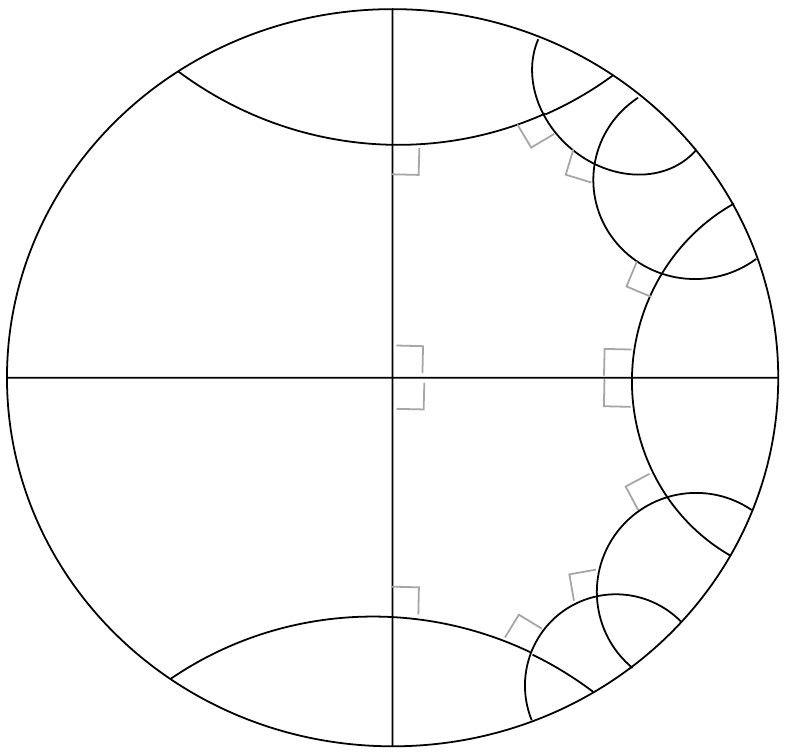}}
          \put(1.95,2.25){$D$}
          \put(2.3,2){$C$}
          \put(2.4,1.75){$B$}
          \put(2.6,1.2){$A=D'$}
          \put(2.4,0.5){$C'$}
          \put(2.2,0.25){$B'$}
          \put(1.8,0){$A'$}
          \put(1.2,-0.25){$0$}
          \put(1.2,2.5){$\infty$}
          \put(1.5,0.8){$H'$}
          \put(1.5,1.5){$H$}
 \end{picture}
  \vspace{0.2cm} 
\caption{The two adjacent hexagons $H$ and $H'$}
   \label{attachang1}
\end{figure}

\noindent We can look at the ordered sequence of points 
$$\Big(\pi^{\H^{2}}(iA'),\pi^{\H^{2}}(iB'),\pi^{\H^{2}}(iC'),\pi^{\H^{2}}(iA),\pi^{\H^{2}}(iB),\pi^{\H^{2}}(iC) ,\pi^{\H^{2}}(iD) \Big)$$
obtained by the union of vertices of the polygonal chains associated to $(H,\Yi)$ and $(H',\Yi)$ respectively. This induces an orientation on the segments forming the polygonal chains. The \emph{attachment angle} $\beta$ between these two polygonal chains is the angle (measured on the left) formed by the two (non-vanishing) segments attached at the point $\pi^{\H^{2}}(iD')=\pi^{\H^{2}}(iA)$.  \end{definition}

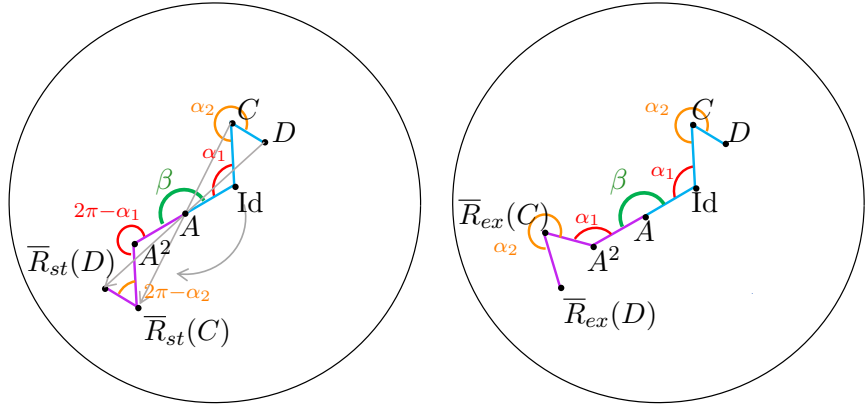
\begin{figure}[!h]
   \centering\captionsetup{justification=centering,margin=2cm}
   \setlength{\unitlength}{0.1\textwidth}
   \begin{picture}(7.5,3.5)
     \put(0,0){\includegraphics[width=11.5cm,height=5.5cm]{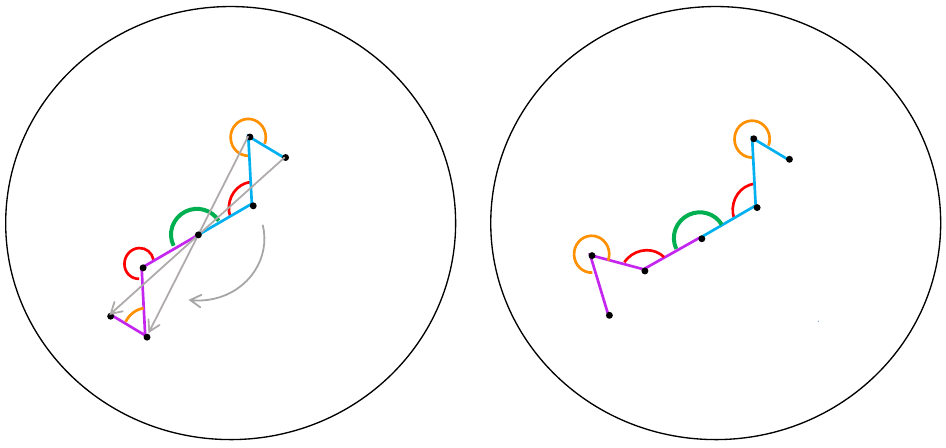}}
\put(0.6,1.7){\textcolor{red}{$_{2\pi-\alpha_{1}}$}}
\put(1.15,1.25){$A^{2}$}
\put(5,1.2){$A^{2}$}
\put(5.55,2.05){\textcolor{red}{$_{\alpha_{1}}$}}
\put(4.2,1.4){\textcolor{orange}{$_{\alpha_{2}}$}}
\put(3.9,1.6){$\overline{R}_{ex}(C)$}
\put(1.6,2.6){\textcolor{orange}{$_{\alpha_{2}}$}}
\put(5.5,2.6){\textcolor{orange}{$_{\alpha_{2}}$}}
\put(2,2.5){$C$}
\put(5.9,2.5){$C$}
\put(1.5,1.5){$A$}
\put(5.4,1.45){$A$}
\put(1.2,0.6){$\overline{R}_{st}(C)$}
\put(4.8,0.75){$\overline{R}_{ex}(D)$}
\put(2,1.7){$\Id$}
\put(5.9,1.7){$\Id$}
\put(1.7,2.2){\textcolor{red}{$_{\alpha_{1}}$}}
\put(4.9,1.65){\textcolor{red}{$_{\alpha_{1}}$}}
\put(1.3,1.9){\textcolor{forestgreen}{$\beta$}}
\put(5.2,1.9){\textcolor{forestgreen}{$\beta$}}
\put(2.3,2.3){$D$}
\put(6.2,2.3){$D$}
\put(1.2,1){\textcolor{orange}{$_{2\pi-\alpha_{2}}$}}
\put(0.2,1.2){$\overline{R}_{st}(D)$}

 \end{picture}
   
\caption{Polygonal chains of $(\overline{R}_{st}(H),\Yi)$ and $(\overline{R}_{ex}(H),\Yi)$ obtained from the polygonal chain of $(H,\Yi)$}
   \label{geomintS3}
\end{figure}

For a visualization of the attachment angle see for example Figure \ref{geomintS3}. We want to study the case where $H'$ is obtained by reflecting $H$ along a side.

\begin{definition} \label{adjsymmhex}
Let $H_{1}=[\Y_{1},\Y_{2},\Y_{3},\Y_{4},\Y_{5},\Y_{6}]$ and $H_{2}=[\Y_{1},\Y_{6},\Y_{7},\Y_{8},\Y_{9},\Y_{2}]$ be two hexagons adjacent at $\Y_{1}$. The hexagons $H_{1}\#_{\Y_{1}}H_{2}$ are said to be \emph{symmetric} if $H_{2}=R(H_{1})$ for a reflection $R \in \mathcal{R}_{\Y_{1}}^{\Y_{6},\Y_{2}}$.
\end{definition}

Our aim is to give a geometric interpretation of the set $\mathcal{K}$.  More precisely we show how the choice of $R \in \mathcal{K}$ is equivalent to choosing an attachment angle $\beta$ between the polygonal chains of $H$ and $\overline{R}(H)$, where $H\#\overline{R}(H)$ are adjacent symmetric and $\overline{R}$ is conjugate to $R$. We start by stating a proposition which will be useful later.

\begin{prop} \label{colinearity}
    Let $(0,P,Q,\infty)$ be a maximal quadruple and consider the orthogonal tubes $\Y_{-Q,Q} \perp \Y_{P,QP^{-1}Q}$. Suppose $(-Q,0,P,Q)$ generic. Then the hyperbolic components of $P,Q$ and $QP^{-1}Q$ lie on the same geodesic in $\H^{2}$ and $\pi^{\H^{2}}(p_{0,\infty}(Q))$ is the middle point of the three. If $(-Q,0,P,Q)$ is non-generic then the hyperbolic components coincide in $\H^{2}$.
\end{prop}

\proof By Proposition \ref{bijmap1} we know that $(-Q,0,P,Q)$ is generic if and only if the quadruple $(0,Q,QP^{-1}Q,\infty)$ is generic. Up to $\Sp(4,\R)$-action we can consider $Q=\Id$ and $ QP^{-1}Q=Y$ where $Y=\bpm y_{1}&0\\0&y_{2}\epm$, $y_{1}>y_{2}$. The hyperbolic component of $\Id$ is $i \in \H^{2}$ in the identification of Section \ref{seccyl}.  It is trivial to show that the hyperbolic components of $Y$ and $Y^{-1}$ lie on the same geodesic in $\H^{2}$ (the $y$-axis), where the point $i$ is in the middle. Since isometries preserve geodesics, the same is true more generally for tubes $\Y_{-Q,Q} \perp \Y_{P,PQP^{-1}Q}$. The non-generic case is trivial.
\endproof

We also recall a standard fact of linear algebra.

\begin{lem} \label{Kdecomp}
Let $K \in \PO(2)$ such that $K^{2}=\Id$. Then 
$$K=\Id \ \ \text{ or } \ \ K=\bpm-1&0\\0&1\epm \bpm \cos\theta&-\sin\theta\\\sin\theta&\cos\theta\epm$$
for a unique $\theta \in [0,\pi)$.
\end{lem}

\noindent Let us denote by $\beta$ the following map 

\begin{equation} \label{beta}
    \begin{aligned}
        \beta: \mathcal{K}&\to[0,2\pi)\\
        \bpm-K&0\\0&K\epm&\mapsto \begin{cases}
            \pi+2\theta\\
            \pi \ \ \text{ if } K=\Id
        \end{cases}
    \end{aligned}
\end{equation}
where 
$$K=\bpm-1&0\\0&1\epm \bpm \cos\theta&-\sin\theta\\\sin\theta&\cos\theta\epm, \ \theta \in [0,\pi)$$

\begin{prop} \label{propattachangle}
Let $(H,\Yi) \in \mathcal{H}$ where $H=(0,A,\Id,C,D,\infty)$ and has arc coordinates $(\underline{b},\underline{c},\underline{d},[\alpha_{1},\alpha_{2}])$. Suppose the angles of the polygonal chain associated to $(H,\Yi)$ are $\alpha_{1},\alpha_{2}$ (possibly only $\alpha$ or no angle).
Consider a reflection $\overline{R}$ inside $\mathcal{R}_{\Y_{-A,A}}^{\Yi,\Y_{A^{2},\Id}}$:
$$\overline{R}=g^{-1}Rg \ \in \mathcal{R}_{\Y_{-A,A}}^{\Yi,\Y_{A^{2},\Id}}$$
for $R=\bpm-K&0\\0&K\epm\in \mathcal{K}$ and $g$ an isometry such that $g(-A,0,A^{2},A)=(0,\Id,Y,\infty)$ with $Y$ diagonal.\\
Then the attachment angle between the polygonal chains of $(H,\Yi)$ and $(\overline{R}(H),\Yi)$ is given by $\beta(K)$ where $\beta$ is the map in (\ref{beta}). Moreover the polygonal chain associated to $(\overline{R}(H),\Yi)$ has
\begin{itemize}
    \item[(i)] segments of lengths $h(\underline{d}),h(\underline{c}),h(\underline{b})$ where $h$ is the map $h(d_{1},d_{2})=d_{1}-d_{2}$
\item[(ii)] angles (if there):
$\begin{cases}
\alpha_{2},\alpha_{1}  \ \ (or \ \alpha) \\
2\pi-\alpha_{2},2\pi-\alpha_{1} (or \ 2\pi-\alpha) \  \text{ if } K=\Id
\end{cases}$
\end{itemize}
\end{prop}

\begin{figure}[!h]
   \centering
   \captionsetup{justification=centering,margin=2cm}
   \setlength{\unitlength}{0.1\textwidth}
   \begin{picture}(3,2.8)
     \put(0.0,0){\includegraphics[width=4.5cm,height=4.2cm]{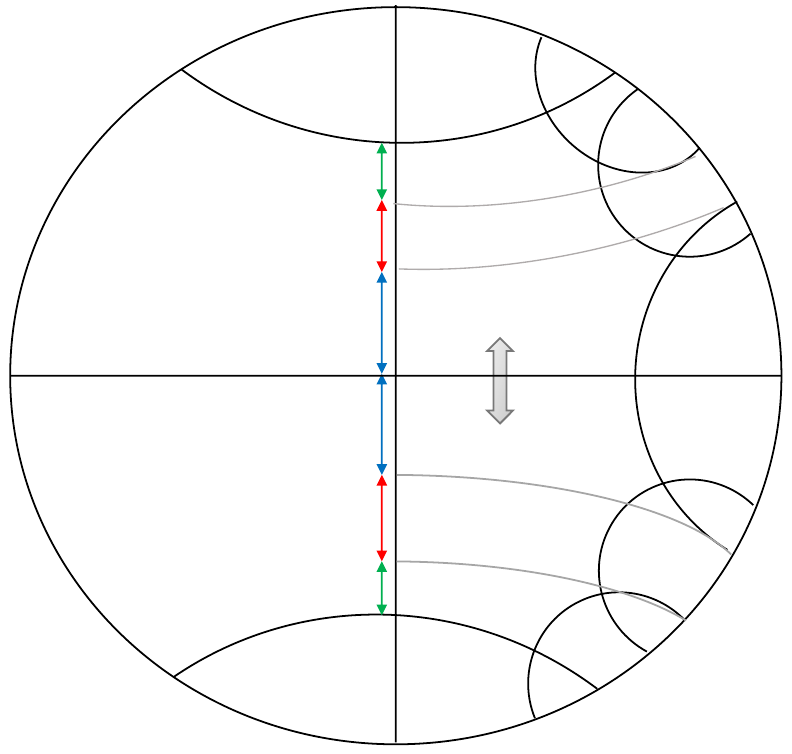}}
          \put(2.25,2.5){$D$}
          \put(2.6,2.2){$C$}
          \put(0.3,2.5){$-D$}
         \put(2,2.7){$DC^{-1}D$}
          \put(2.7,2){$\Id$}
          \put(2.9,1.3){$A=\overline{R}(A)$}
          \put(-0.5,1.3){$-A$}
          \put(2.8,0.6){$A^{2}=\overline{R}(\Id)$}
          \put(2.6,0.2){$\overline{R}(C)$}
          \put(2.1,-0.05){$\overline{R}(D)$}
          \put(1.4,-0.2){$0$}
          \put(1.35,2.7){$\infty$}
          \put(1.6,0.7){$\overline{R}(H)$}
           \put(1.9,1){\textcolor{gray}{$\overline{R}$}}
           \put(1.1,0.85){\textcolor{red}{$\underline{c}$}}
           \put(1.1,0.6){\textcolor{forestgreen}{$\underline{d}$}}
           \put(1.1,1.1){\textcolor{blue}{$\underline{b}$}}
           \put(1.1,1.75){\textcolor{red}{$\underline{c}$}}
           \put(1.1,1.5){\textcolor{blue}{$\underline{b}$}}
           \put(1.1,2){\textcolor{forestgreen}{$\underline{d}$}}
          \put(1.7,1.8){$H$}
 \end{picture}
  \vspace{0.2cm} 
\caption{The adjacent symmetric hexagons $H\#_{\Y_{-A,A}}\overline{R}(H)$}
   \label{geomintS1}
\end{figure}

\proof Let us first consider the case where $(H,\Yi)$ is generic. The two adjacent symmetric hexagons are illustrated in Figure \ref{geomintS1}. By Proposition \ref{propreflinsiegel} the 6-tuple \\$(0, \overline{R}(D), \overline{R}(C), \overline{R}(\Id)=A^{2}, \overline{R}(A)=A, \infty)$ is maximal. This 6-tuple determines the ordered sequence of vertices in the polygonal chain associated to $(\overline{R}(H),\Yi)$.

\noindent By  Corollary \ref{reflsetwitharcoords} we know 
$$\mathcal{R}_{\Y_{-A,A}}^{\Yi,\Y_{A^{2},\Id}}=\begin{cases}
\{g^{-1}R_{st}g,g^{-1}R_{ex}g\}, \ if \ \underline{b} \in \mathfrak{a}\\
g^{-1}\mathcal{K}g   \ if \ \underline{b} \in \mathfrak{d}
\end{cases} $$
where $g(-A,0,A^{2},A)=(0,\Id,Y,\infty)$ for $Y$ diagonal and $\{R_{st},R_{ex}\} \subset \mathcal{K}$. As $(H,\Yi)$ generic we know $\underline{b} \in \mathfrak{a}$ so that
$\overline{R}=g^{-1}R_{st}g$ or $\overline{R}=g^{-1}R_{ex}g$. By Proposition \ref{colinearity} we know that the attachment angle is $\beta(K)=\pi$ as the hyperbolic components of the points $A^{2}$, $A$ and $\Id$ are colinear. Put $\overline{R}_{st}=g^{-1}R_{st}g$. Computations give $\overline{R}_{st}=\bpm 0&A\\A^{-1}&0\epm$ and one can immediately see that $\overline{R}_{st} \in \mathcal{R}_{\Y_{-A,A}}^{\Yi,\Y_{A^{2},\Id}}$ as
$$\overline{R}_{st}(-A)=-A, \ \ \overline{R}_{st}(A)=A, \ \ \overline{R}_{st}(0)=\infty, \ \ \overline{R}_{st}(A^{2})=\Id$$
It is straightforward to see  that the segments of the polygonal chain associated to $(\overline{R}_{st}(H),\Yi)$ have length $h(\underline{d}),h(\underline{c}),h(\underline{b})$ respectively. The eigenspaces $E_{\pm 1}^{\overline{R}_{st}}$ of $\overline{R}_{st}$ are given by
$$E_{1}^{\overline{R}_{st}}=A, \ \ E_{-1}^{\overline{R}_{st}}=-A$$
By Proposition \ref{reflfixp} we know that $\overline{R}_{st}$ is fixing the tube $\Y_{-A,A}$ and sending any transverse $X \pitchfork A$ to the unique $\overline{R}_{st}(X)$ such that $\Y_{-A,A} \perp \Y_{X,\overline{R}_{st}(X)}$. By Proposition \ref{colinearity} the hyperbolic components of $\overline{R}_{st}(X)$ and $X$ lie therefore on the same geodesic inside $\H^{2}$. The polygonal chain associated to $(\overline{R}_{st}(H),\Yi)$ is obtained by rotating the polygonal chain of $H$ of an angle $\pi$ around $A$. This is illustrated on the left side of Figure \ref{geomintS3} where the polygonal chain of $H$ is drawn in blue and the polygonal chain of $\overline{R}_{st}(H)$ is drawn in purple. We see that the angles of the polygonal chain associated to $(\overline{R}_{st}(H),\Yi)$ are given by $2\pi-\alpha_{2}$ and $2\pi-\alpha_{1}$.

Put $\overline{R}_{ex}=g^{-1}R_{ex}g$. Instead of explicitly computing $\overline{R}_{ex}$ observe that if we denote by $f$ the map $f \in \Stab_{\PSp(4,\R)}(0,A^{2},A,\infty)$, $ \ f \neq \Id$ then the map $f \circ \overline{R}_{st}$ satisfies
$$f \circ \overline{R}_{st}(-A)=-A, \ \ f \circ \overline{R}_{st}(A)=A, \ \ f \circ \overline{R}_{st}=\infty, \ \ f \circ \overline{R}_{st}(A^{2})=\Id$$
so that $\overline{R}_{ex}=f \circ \overline{R}_{st}$. The geometric interpretation of $f$ is the reflection across the geodesic going through the hyperbolic components of $A$ and $A^{2}$ respectively and this geodesic also goes through the hyperbolic component of $\Id$ (Proposition \ref{colinearity}). The angles of the polygonal chain associated to $(\overline{R}_{ex}(H),\Yi)$ are therefore given by $\alpha_{2},\alpha_{1}$ and this is illustrated on the right-hand side of Figure \ref{geomintS3}.

\begin{figure}[!h]
   \centering\captionsetup{justification=centering,margin=2cm}
   \setlength{\unitlength}{0.1\textwidth}
   \begin{picture}(8,2.5)
     \put(0,0){\includegraphics[width=12cm,height=3.8cm]{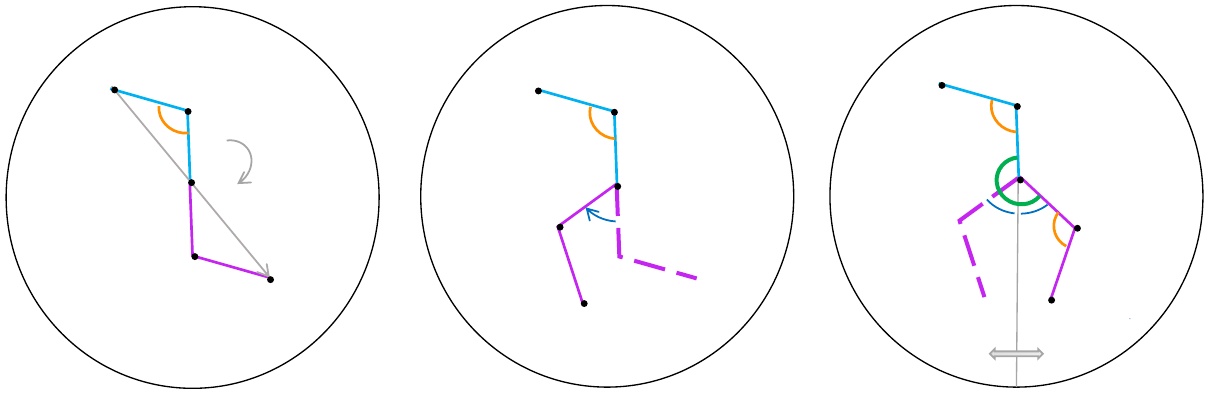}}
\put(1.3,1.7){\textcolor{orange}{$\alpha$}}
\put(1.7,1.5){\textcolor{gray}{$M$}}
\put(6.66,1.7){\textcolor{orange}{$\alpha$}}
\put(6.2,0.3){\textcolor{gray}{$r$}}
\put(6.9,0.9){\textcolor{orange}{$\alpha$}}
\put(6.25,1.){\textcolor{blue}{$_{2\theta}$}}
\put(6.55,1.){\textcolor{blue}{$_{2\theta}$}}
\put(3.7,0.9){\textcolor{blue}{$_{2\theta}$}}
\put(3.95,1.1){\textcolor{blue}{$S$}}
\put(6.2,1.3){\textcolor{forestgreen}{$\beta$}}
 \end{picture}
   
\caption{Geometrical interpretation of $\overline{R}=rSM$}
   \label{geomintS4}
\end{figure}

If $(H,\Yi)$ is non-generic of type 1.1 the length parameters $\underline{b},\underline{c},\underline{d} $ lie inside $\mathfrak{d}\times\mathfrak{a}^{2}$. The polygonal chain associated to $(H,\Yi)$ has only one angle $\alpha$ as illustrated in Figure \ref{nongen1}. By Corollary \ref{reflsetwitharcoords} we know $\overline{R}=g^{-1}\mathcal{K}g \ \in \mathcal{R}_{\Y_{-A,A}}^{\Yi,\Y_{A^{2},\Id}}$
where $g(-A,0,A^{2},A)=(0,\Id,Y,\infty)$ for $Y$ diagonal. As $\underline{b} \in \mathfrak{d}$ we know
$A=a\cdot \Id$.
Computations give 
$$\mathcal{R}_{\Y_{-A,A}}^{\Yi,\Y_{A^{2},\Id}}=\Big\{\bpm 0&aK\\a^{-1}K&0\epm, \ K \in \PO(2), \ K^{2}=\Id\Big\}$$
Given $\overline{R} \in \mathcal{R}_{\Y_{-A,A}}^{\Yi,\Y_{A^{2},\Id}}$ let us decompose $\overline{R}$ as following: \\let $r=\bpm-1&0\\0&1\epm$ and $S=\bpm \cos\theta&-\sin\theta\\\sin\theta&\cos\theta\epm$ for $\theta \in [0,\pi)$. We write
$$\overline{R}=\underbrace{\bpm r&0\\0&r\epm}_{r} \underbrace{\bpm S&0\\0&S\epm}_{S} \underbrace{\bpm 0&a\Id\\a^{-1}\Id&0\epm}_{M}$$
where $\bpm-1&0\\0&1\epm \bpm \cos\theta&-\sin\theta\\\sin\theta&\cos\theta\epm \neq R_{st}$ is the decomposition of Lemma \ref{Kdecomp}. 
The geometrical interpretation of this decomposition is illustrated in Figure \ref{geomintS4}. The map $M$ is analogue to the rotation of Figure \ref{geomintS3} (left). The map $S$ is a rotation of angle $2\theta$ around $i$ on the hyperbolic component of $\Yi$ (see Section \ref{geomintdiag}) and the map $r$ is a reflection across the vertical axis. We obtain an attachment angle $\beta(K)=\pi+2\theta$ and polygonal chain angle $\alpha$. When $\overline{R}=R_{st}$ we only have $\overline{R}=\bpm 0&a\Id\\a^{-1}\Id&0\epm$ and we are neither rotating nor reflecting. 

\noindent The proof for the other cases where $(H,\Yi)$ is non-generic are similar.
\endproof

\begin{remark} \label{reflsamehex} 
If $H$ is contained in a maximal polydisk it can happen that $\mathcal{R}_{\Y_{-A,A}}^{\Yi,\Y_{A^{2},\Id}}$ contains two different reflections $\overline{R},\overline{R'}$ for which $\overline{R}(H)=\overline{R'}(H)$. Let us denote for simplicity $p=(\underline{b},\underline{c},\underline{d},[\alpha_{1},\alpha_{2}])$ the arc coordinates associated to $(H,\Yi)$. By the geometrical interpretation of Proposition \ref{propattachangle} it is not hard to show that the case $\overline{R}(H)=\overline{R'}(H)$ happens exactly for 
\begin{equation} \label{bothinmaxpol}
    \overline{R}, \overline{R'}\in \{g^{-1}R_{st}g, \ g^{-1}R_{ex}g\} \  \ \text{ if } p \in \mathcal{D}\backslash \mathcal{D}_{\H^{2}}
\end{equation}
and for any
\begin{equation} \label{bothindiagdisc}
    \overline{R},\overline{R'} \in g^{-1}\mathcal{K}g \ \ \text{ if } p \in \mathcal{D}_{\H^{2}}
\end{equation}
The two hexagons $H$ and $\overline{R}(H)=\overline{R'}(H)$ lie both inside the model polydisk if $\overline{R},\overline{R'}$ are as in (\ref{bothinmaxpol}) and all the points of the two polygonal chains lie on the vertical geodesic of $\H^{2}$. In (\ref{bothindiagdisc}), the two hexagons $H$ and $\overline{R}(H)=\overline{R'}(H)$ lie both inside the diagonal disc and all the points of the polygonal chains coincide with $\pi^{\H^{2}}(i\Id)$.
\end{remark}

\begin{remark} 
In Proposition \ref{propattachangle} we have shown how to draw the polygonal associated to $(\overline{R}_{2}(H),\Yi)$ where $\overline{R}_{2} \in \mathcal{R}_{\Y_{-A,A}}^{\Yi,\Y_{A^{2},\Id}}$. We can state a similar result for $\overline{R}_{1} \in \mathcal{R}_{\Y_{-D,D}}^{\Yi,\Y_{C,DC^{-1}D}}$. 
We can then draw the three attached polygonal chains $\overline{R}_{2}(H)\#H\#\overline{R}_{1}(H)$. Figure \ref{combpolchains} illustrates all possible polygonal chains of $\overline{R}_{2}(H)\#H\#\overline{R}_{1}(H)$ in the case that $(H,\Yi)$ is a generic right-angled hexagon. The polygonal chains are drawn up to isometry, this means that we consider two polygonal chains to be equivalent if there exists an isometry  $g \in \PSp(4,\R) $ sending all the vertices of one to the vertices of the other.

\begin{figure}[!h]
   \centering\captionsetup{justification=centering,margin=2cm}
   \setlength{\unitlength}{0.1\textwidth}
   \begin{picture}(8.5,2.5)
     \put(0,0){\includegraphics[width=13.5cm,height=3.5cm]{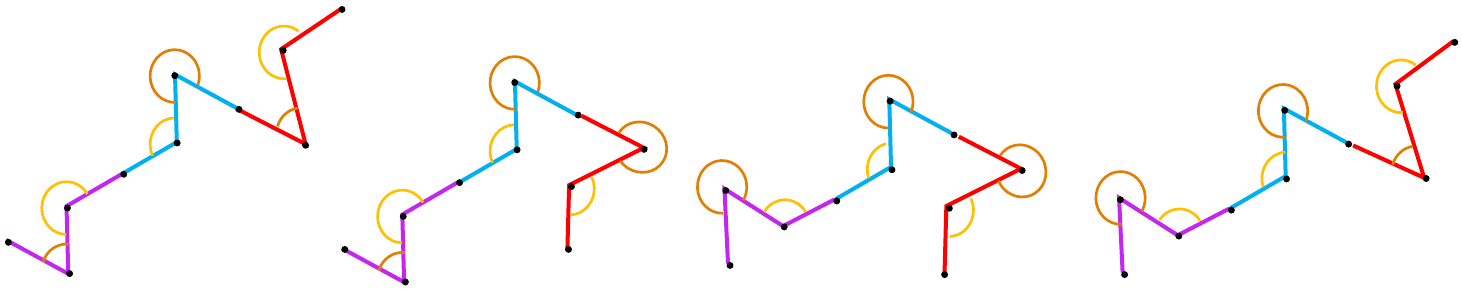}}
\put(1.1,1){\textcolor{cyan}{$H$}}
\put(1.8,1.8){\textcolor{red}{$\overline{R}_{1}(H)$}}
\put(0.6,0.4){\textcolor{violet}{$\overline{R}_{2}(H)$}}
\end{picture}
\caption{All possible polygonal chains (up to isometry) for $\overline{R}_{2}(H)\#H\#\overline{R}_{1}(H)$ when $(H,\Yi)$ is generic}
   \label{combpolchains}
\end{figure}

\end{remark}

\section{Parameters for maximal representations} \label{chapterparformax}

In this section we use arc coordinates of right-angled hexagons to parametrize maximal representations. We start by discussing geometric properties of Shilov hyperbolic isometries in $\PSp(2n,\R)$. We introduce the notion of a maximal representation from the reflection group
$\mathbb{Z}/2\mathbb{Z}*\mathbb{Z}/2\mathbb{Z}*\mathbb{Z}/2\mathbb{Z}$ into $\PSp^{\pm}(2n,\R)$. This will lead to the parametrization of a subset of maximal representations arising from a very geometric picture (Theorem \ref{thesis}).

\subsection{Shilov hyperbolic isometries}

\begin{definition} \label{defshilhyp}
   An element $g \in \PSp(2n,\R)$ is called \emph{Shilov hyperbolic} if it is conjugate to $\begin{pmatrix}
A & 0\\
0 & A^{-T} 
\end{pmatrix}$ for a matrix $A \in \GL(n,\R)$ with complex eigenvalues with modulus greater than one.
\end{definition}

An element $g\in \PSp(2n,\R)$ is Shilov hyperbolic if and only if $g$ fixes two transverse Lagrangians $l_{g}^{+},l_{g}^{-}$ on which it acts expandingly and contractingly respectively. In the following example we classify Shilov hyperbolic elements in $\PSp(4,\R)$.

\begin{example} \label{exampleshil} Let $(e_{1},e_{2},e_{3},e_{4})$ be the standard basis of $\R^{4}$. Recall that we denote $0, l_{\infty}$ the Lagrangians $0=\langle e_{3},e_{4} \rangle$ and $l_{\infty}=\langle e_{1},e_{2} \rangle$ respectively, and that the standard tube $\Yi= \{ iY| \ Y \in \Sym^{+}(2,\R)  \}$
is isometrically identified with $\R \times \H^{2}$ (see Lemma \ref{HxR}) through the map 
\[
\begin{aligned}
\pi^{\R} \times \pi^{\H^{2}}: \Yi&\to \R\times\Sym^{+}(2,\R)\\
iY& \mapsto \Big( \frac{\log \det Y}{\sqrt{2}}, \frac{Y}{\sqrt{\det Y}} \Big)
\end{aligned}
\]
Moreover, in the proof of Proposition \ref{orienthypcomp} we have shown how the visual boundary of the hyperbolic component of $\Yi$ can be realized as the $\O(2)$-orbit of the Lagrangian $l=\langle e_{1},e_{4}\rangle$. 
An isometry $g_{A}=\bpm
A & 0\\
0 & A^{-T} 
\epm$ in $\PSp(4,\R)$ stabilizes the standard tube $\Yi$. The action of $g_{A}$ on the $\R$-component of $\Yi$ is given by 
$$\pi^{\R}(g_{A}(iY))=\frac{\log \det^{2}(A)\det Y}{\sqrt{2}}=\frac{\log \det^{2}(A) }{\sqrt{2}}+ \pi^{\R}(iY)$$
We study the action of $g_{A}$ on the hyperbolic component of the standard tube $\Yi$ and give a geometric interpretation for $g_{A}$ to be Shilov hyperbolic. We have the following possibilities:
\begin{itemize}
\item  $A$ has one eigenvalues $\lambda \in \R$ and is conjugate to a matrix
$$  \bpm \bpm \lambda & 0\\ 0 & \lambda \epm & 0\\
0 & \bpm \frac{1}{\lambda} & 0\\ 0 & \frac{1}{\lambda} \epm \epm
\sim g_{A}$$
\noindent The isometry $g_{A}$ acts on $\H^{2}$ as the identity map and $\pi^{\R}(g_{A}(iY))-\pi^{\R}(iY) = \frac{4\log |\lambda|}{\sqrt{2}}$. It is clear that $g_{A}$ is Shilov hyperbolic if and only if $|\lambda| \neq 1$ that is $g_{A} \neq id$.
\item  $A$ has eigenvalues $\lambda > \mu \in \R$ and is conjugate to a matrix
$$  \bpm \bpm \lambda & 0\\ 0 & \mu \epm & 0\\
0 & \bpm \frac{1}{\lambda} & 0\\ 0 & \frac{1}{\mu} \epm \epm
\sim g_{A}$$

\noindent The isometry $g_{A}$ acts on the hyperbolic component of $\Yi$ as an hyperbolic isometry, it fixes exactly two points in the boundary of $\H^{2}$. The axis is the infinite geodesic on the $\H^{2}$-component having as endpoints the two fixed points and has translation length given by $\log|\frac{\lambda}{\mu}|$. Observe that if $\lambda$ and $\mu$ have different signs then $\det A <0$ and $g_{A}$ is reversing the orientation of the $\H^{2}$-component (see Definition \ref{isomref}). On the $\R$-component we have $\pi^{\R}(g_{A}(iY))-\pi^{\R}(iY) = \frac{2\log |\lambda\mu|}{\sqrt{2}}$. We should think at the isometry $g_{A}$ being Shilov hyperbolic if we can find a point in the tube that moves vertically (on the $\R$-component) more than horizontally (on the $\H^{2}$-component). This will be explained in detail in Lemma \ref{checkshil}.

\item  $A$ has one eigenvalues $\lambda\in \R$ and is conjugate to a matrix
$$  \bpm \bpm \lambda & 1\\ 0 & \lambda \epm & 0\\
0 & \bpm \frac{1}{\lambda} & 0\\ -\frac{1}{\lambda^{2}} & \frac{1}{\lambda} \epm \epm
\sim g_{A}$$
The map $g_{A}$ acts on the hyperbolic component of $\Yi$ as a parabolic isometry, it fixes exactly one point in the boundary of $\H^{2}$. We have $\pi^{\R}(g_{A}(iY))-\pi^{\R}(iY) = \frac{4\log |\lambda|}{\sqrt{2}}$. The isometry $g_{A}$ is Shilov hyperbolic if and only if $|\lambda| \neq 1$. 
\item  $A$ has two complex eigenvalues $\lambda e ^{i\theta},\lambda e ^{-i\theta}, \theta \neq 2k \pi $ and is conjugate to a matrix
$$  \bpm \bpm \lambda \cos\theta & -\lambda \sin\theta\\ \lambda \sin\theta & \lambda \cos\theta\epm & 0\\
0 & \bpm \frac{1}{\lambda} \cos\theta & -\frac{1}{\lambda} \sin\theta\\ \frac{1}{\lambda} \sin\theta & \frac{1}{\lambda} \cos\theta\epm \epm
\sim g_{A}$$
The map $g_{A}$ acts on the hyperbolic component of $\Yi$ as an elliptic isometry, it fixes exactly one point inside $\H^{2}$ and the angle of rotation is given by $\theta=-2\log |\lambda|$. The isometry $g_{A}$ is Shilov hyperbolic if and only if $|\lambda| \neq 1$. 
\end{itemize}

The geometrical interpretation of the action of $g_{A}$ on $\Yi$ in the hyperbolic, parabolic and elliptic case is illustrated in Figure \ref{g_Aaction}.

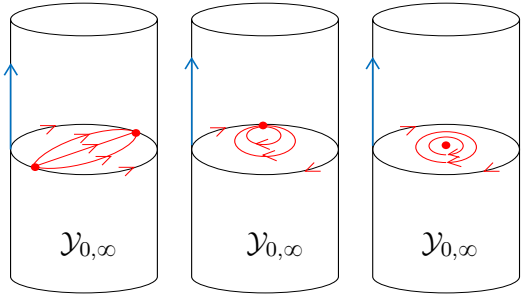
\begin{figure}[!h]
   \centering
   \setlength{\unitlength}{0.1\textwidth}
   \begin{picture}(4.5,2.7)
     \put(00,0){\includegraphics[width=7cm,height=4cm]{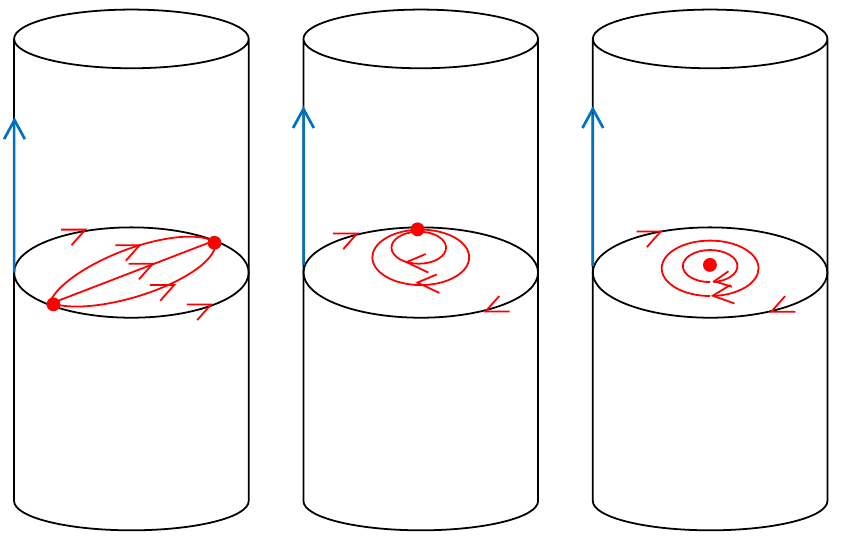}}
      \put(0.5,0.4){$\Yi$}
      \put(2.1,0.4){$\Yi$}
      \put(3.6,0.4){$\Yi$}
   \end{picture}
   \caption{The action of $g_{A}$ on $\Yi=\R \times \H^{2}$ in the hyperbolic (with $\det A>0$), parabolic and elliptic case}
   \label{g_Aaction}
\end{figure}

\end{example}

\begin{lem} \label{checkshil}
Let $g$ be an element of $\PSp(4,\R)$ fixing two Lagrangians $l_{1},l_{2}$ in $ \mathcal{L}(\R^{4})^{(2)}$ i.e. $g$ is conjugated to $g_{A}=\bpm
A & 0\\
0 & A^{-T} 
\epm$ for a matrix $A \in \GL(2,\R)$. Denote by $|\lambda| \geq |\mu|$ the modulus of the eigenvalues of $A$. Then
\begin{itemize}
    \item[(i)]  There exists $l\in \mathcal{L}(\R^{4})$ such that $(l_{1},l,l_{2})$ maximal and $(l_{1},l,g(l),l_{2})$ maximal if and only if $|\mu|>1$ (that is $g$ is Shilov hyperbolic).
\item[(ii)] $(l_{1},l,g(l),l_{2})$ maximal for all $l$ such that $(l_{1},l,l_{2})$ maximal if and only if $A=\lambda \Id$ for $\lambda \in \R, \  \lambda>1$.
\end{itemize}
\end{lem}

\proof
It is sufficient to prove the lemma for $g=g_{A}$, that is $l_{1}=0$ and $l_{2}=\infty$.
\begin{itemize}
    \item[\emph{(i)}]
We want to show that there exists $Y$ such that $(0,Y,\infty)$ maximal and $(0,Y,g_{A}Y,\infty)$ maximal if and only if $|\mu|>1$. Let us write $Y>0$ for a matrix $Y$ which is positive definite. By Lemma \ref{max} we know that $(0,Y,\infty)$ is maximal if and only if $Y>0$. Suppose there exists $Y >0$ such that $(0,Y,g_{A}Y,\infty)$ maximal, that is \\$g_{A}Y-Y$=$AYA^{T}-Y>0$. Recall that a matrix $M$ is positive definite if and only if $NMN^{T}>0$ for every invertible matrix $N$. In particular for $N=\sqrt{Y^{-1}}$ we obtain
\begin{equation} \label{cond}
\sqrt{Y^{-1}}(AYA^{T}-Y)\sqrt{Y^{-1}}=(\sqrt{Y^{-1}}A\sqrt{Y})(\sqrt{Y}A^{T}\sqrt{Y^{-1}})-\Id>0
\end{equation}
The matrix $\sqrt{Y}A^{T}\sqrt{Y^{-1}} \in \GL(2,\R)$ has the same eigenvalues of $A$. Let $v$ be the orthonormal eigenvector associated to $\mu$. Then $\sqrt{Y}A^{T}\sqrt{Y^{-1}}v=\mu v$ and we deduce $ \sqrt{Y}A^{T}\sqrt{Y^{-1}}\overline{v}=\overline{\mu}\overline{v}$. It follows from (\ref{cond}):
$$
v^{T}(\sqrt{Y^{-1}}A\sqrt{Y})(\sqrt{Y}A^{T}\sqrt{Y^{-1}})\overline{v}-v^{T}\overline{v}=|\mu|^{2}-1>0
$$
Suppose $|\mu|>1$. We want to find $Y>0$ such that $(0,Y,g_{A}Y,\infty)$ maximal that is we want to find $Y>0$ such that $g_{A}Y-Y>0$.  In Remark \ref{meaningvector} we have given an equivalent condition for $g_{A}Y-Y$ to be positive definite: let $\mathbf{r}=d^{\R}(\pi^{\R}(iY),\pi^{\R}(ig_{A}Y))$ and $ \mathbf{h}=d^{\H^{2}}(\pi^{\H^{2}}(iY),\pi^{\H^{2}}(ig_{A}Y))$, then 
 \begin{equation} \label{rhineq}
     g_{A}Y-Y>0 \iff \mathbf{r}>\frac{1}{\sqrt{2}}\mathbf{h}
 \end{equation}

As $g_{A}$ Shilov hyperbolic (this is the assumption $|\mu|<1$) we know that $g_{A}$ acts as a translation of distance $\mathbf{r}$ on the $\R$-component of $\Yi$ and as an isometry on the $\H^{2}$-component which can be hyperbolic parabolic or elliptic (see Example \ref{exampleshil}). Observe that for a fixed $A$ the distance $\mathbf{r}$ only depends on the eigenvalues of $A$ and not on the point $Y$, whereas the distance $\mathbf{h}$ depends on $Y$ and decreases the more $Y$ is close to the axis (if the isometry is hyperbolic) or to a fixed point (if the isometry is parabolic or elliptic). We can always find an open neighbour of the axis (or of a fixed point) such that the condition $\mathbf{r}>\frac{1}{\sqrt{2}}\mathbf{h}$ is satisfied.
\item[\emph{(ii)}] We want to show that $(0,Y,g_{A}Y,\infty)$ is maximal $\forall Y$ such that $(0,Y,\infty)$ maximal if and only if $A=\lambda \Id$ for $\lambda \in \R, \  \lambda>1$. This is clear after the discussion in \emph{(i)}: recall that $\lambda \Id$ for $\lambda \in \R$ acts on the $\H^{2}$-component of the tube as the identity map, so if $A=\lambda \Id$ the inequality in (\ref{rhineq}) is clearly satisfied. Conversely, suppose that $(0,Y,g_{A}Y,\infty)$ is maximal for all $Y>0$. Equivalently, for any $Y>0$ the inequality in (\ref{rhineq}) is satisfied, where the distance $\mathbf{r}$ is a fixed length depending only on the eigenvalues of $A$. This implies that the action of $g_{A}$ on the hyperbolic component of the tube is the identity map i.e. $g_{A}=\lambda \Id, \ \lambda \in \R$. If this were not the case, we could always find an $Y$ which does not satisfy the inequality in (\ref{rhineq}) by stepping away from the axis (in the hyperbolic case) or the fixed points (in the parabolic or elliptic case) and moving towards the boundary of $\H^{2}$. 
\end{itemize}
\endproof

\subsection{Maximal representations} \label{maxreprsec}

Let $\Sigma$ be an oriented surface with negative Euler characteristic and boundary $\partial \Sigma$. Fix a finite area hyperbolization on $\Sigma$ inducing an action of the fundamental group $\pi_{1}(\Sigma)$ on $S^{1}=\partial \H^{2}$. An element $\gamma \in \pi_{1}(\Sigma)$ is called \emph{peripheral} if it is freely homotopic to a boundary component. Maximal representations are representations that maximize the Toledo invariant, an invariant defined using bounded cohomology (see \cite{toledo1989representations}, \cite{burger2010surface}). It is a deep result from Burger Iozzi and Wienhard (\cite[Theorem 8]{burger2010surface} ) that maximal representations can
be equivalently characterized as representations admitting a well-behaved boundary map.

\begin{definition} \label{defmax}
    A representation $\rho: \pi_{1}(\Sigma) \to \PSp(2n,\R)$ is \emph{maximal} if there exists a $\rho$-equivariant map $\xi: S^{1} \to \mathcal{L}(\R^{2n})$ which is monotone (i.e. the image of any positively oriented triple in the circle is a maximal triple) and right continuous.
\end{definition}

Given a maximal representation $\rho: \pi_{1}(\Sigma) \to \PSp(2n,\R)$, the image $\rho(\gamma)$ of every non-peripheral element $\gamma \in \pi_{1}(\Sigma)$ is Shilov hyperbolic (see \cite{strubel2015fenchel}). Equivalently, $\rho(\gamma)$ fixes two transverse Lagrangians $l_{\gamma}^{+}$ and $l_{\gamma}^{-}$ on which it acts expandingly and contractingly respectively. These Lagrangians are the images $\xi(\gamma^{+})$ and $\xi(\gamma^{-})$ where $\xi: S^{1} \to \mathcal{L}(\mathbb{R}^{2n})$ is the equivariant boundary map and $l^{\pm}_{\gamma}=\xi(\gamma^{\pm})$.  
We want to parametrize the set of maximal representations where the property of being Shilov hyperbolic is true also for peripheral elements. This is equivalent to the requirement that the representations are Anosov in the sense of \cite{guichardwienhard2012}.

\begin{definition} \label{defmaxshil}
    A maximal representation  $\rho: \pi_{1}(\Sigma) \to \PSp(2n,\R)$ will be called \emph{Shilov hyperbolic} if $\rho(\gamma)$ is Shilov hyperbolic for every $\gamma \in \pi_{1}(\Sigma)$. The set of maximal representations which are Shilov hyperbolic will be denoted by $\text{Hom}^{\text{max,Shilov}} (\pi_{1}(\Sigma),\PSp(2n,\R))$. We define $\chi^{\text{max,Shilov}} (\pi_{1}(\Sigma),\PSp(2n,\R))$ the quotient 
    $$\chi^{\text{max,Shilov}} (\pi_{1}(\Sigma),\PSp(2n,\R)):=\text{Hom}^{\text{max,Shilov}} (\pi_{1}(\Sigma),\PSp(2n,\R))/_{\PSp(2n,\R)}$$ 
    where $\PSp(2n,\R)$ is acting by conjugation: $\rho \sim \rho'$ if there exists $g \in \PSp(2n,\R)$ such that $\rho(\gamma)=g\rho'(\gamma)g^{-1}$ for all $\gamma \in \pi_{1}(\Sigma)$.
\end{definition}

\begin{notation}
For the rest of the paper we wll denote by $W_{3}$ the reflection group \\$\mathbb{Z}/2\mathbb{Z} * \mathbb{Z}/2\mathbb{Z}* \mathbb{Z}/2\mathbb{Z}= \langle s_{1},s_{2},s_{3}| \ s_{1}^{2}=s_{2}^{2}=s_{3}^{2}=1 \rangle$.
\end{notation}

\begin{figure}[!h]
   \centering
   \captionsetup{justification=centering,margin=2cm}
   \setlength{\unitlength}{0.1\textwidth}
   \begin{picture}(4.5,3)
     \put(0.8,0){\includegraphics[width=4.5cm,height=4.5cm]{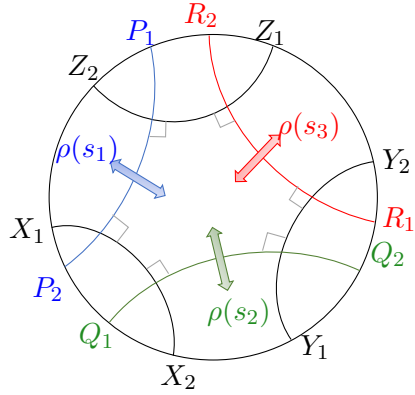}}
          \put(1.8,-0.16){$X_{2}$}
          \put(1.1,0.2){\textcolor{forestgreen}{$Q_{1}$}}
          \put(0.7,0.6){\textcolor{blue}{$P_{2}$}}
          \put(0.5,1.1){$X_{1}$}
          \put(1,2.5){$Z_{2}$}
          \put(1.5,2.8){\textcolor{blue}{$P_{1}$}}
          \put(2,2.95){\textcolor{red}{$R_{2}$}}
          \put(2.6,2.8){$Z_{1}$}
          \put(3.7,1.7){$Y_{2}$}
          \put(3.7,1.2){\textcolor{red}{$R_{1}$}}
          \put(3.6,0.9){\textcolor{forestgreen}{$Q_{2}$}}
          \put(3,0.1){$Y_{1}$}
          \put(0.9,1.8){\textcolor{blue}{$\rho(s_{1})$}}
          \put(2.2,0.4){\textcolor{forestgreen}{$\rho(s_{2})$}}
          \put(2.8,2){\textcolor{red}{$\rho(s_{3})$}}
 \end{picture} 
\caption{The reflections $\rho(s_{1}),\rho(s_{2}),\rho(s_{3})$ for $\rho: W_{3} \to \PSp^{\pm}(2n,\R)$ maximal}
   \label{maxrefl}
\end{figure}

\begin{definition} \label{defmaxW3repr}
    A representation $ \rho: W_{3} \to \PSp^{\pm}(2n,\R)$ is \emph{maximal} if there exists a maximal 6-tuple of Lagrangians $(P_{1},P_{2},Q_{1},Q_{2},R_{1},R_{2})$ such that $\rho(s_{1}),\rho(s_{2}),\rho(s_{3})$ are reflections of $\X$ fixing $(P_{1},P_{2}),(Q_{1},Q_{2}),(R_{1},R_{2})$ respectively and such that
\begin{equation*}
    \begin{cases}
        \rho(s_{1})(X_{1})=X_{2} \text{ and } \rho(s_{1})(Z_{1})=Z_{2}\\
\rho(s_{2})(X_{1})=X_{2} \text{ and } \rho(s_{2})(Y_{1})=Y_{2}\\
        \rho(s_{3})(Y_{1})=Y_{2} \text{ and } \rho(s_{3})(Z_{1})=Z_{2}
    \end{cases}
\end{equation*}
 where $X_{1},X_{2},Y_{1},Y_{2},Z_{1},Z_{2}$ are uniquely determined by
 \begin{equation} \label{perpcondmaxreflW3}
\Y_{P_{1},P_{2}} \perp \Y_{X_{1},X_{2}} \perp \Y_{Q_{1},Q_{2}} \perp \Y_{Y_{1},Y_{2}} \perp \Y_{R_{1},R_{2}} \perp \Y_{Z_{1},Z_{2}} \perp \Y_{P_{1},P_{2}}     
 \end{equation}
The space of maximal representations will be denoted by $\text{Hom}^{\text{max}} (W_{3},\PSp^{\pm}(2n,\R))$. We further define 
    $$\chi^{\text{max}} (W_{3},\PSp^{\pm}(2n,\R)):=\text{Hom}^{\text{max}} (W_{3},\PSp^{\pm}(2n,\R))/_{\PSp(2n,\R)}$$
\end{definition}

\begin{remark} \label{possforrho} 
Given the maximal 6-tuple $(P_{1},P_{2},Q_{1},Q_{2},R_{1},R_{2})$ the set of reflections $\rho(s_{i}), i=1,2,3$ for which $ \rho: W_{3} \to \PSp(4,\R)^{\pm}$ is maximal as in Definition \ref{defmaxW3repr} are given by the sets $\mathcal{R}(P_{1},Z_{2},X_{1},P_{2}),\mathcal{R}(Q_{1},X_{2},Y_{1},Q_{2})$ and $\mathcal{R}(R_{1},Y_{2},Z_{1},R_{2})$ respectively.
\end{remark}

\begin{lem} \label{Shilovproof}
Let $W_{3}= \langle s_{1},s_{2},s_{3}| \ s_{1}^{2}=s_{2}^{2}=s_{3}^{2}=1 \rangle$ and let $ \rho: W_{3} \to \PSp(4,\R)^{\pm}$ be maximal. Then the composition $\rho(s_{i}s_{j})=\rho(s_{i})\rho(s_{j})$ is a Shilov hyperbolic element of $\PSp(4,\R)$ for any $i \neq j$ where $i,j \in \{1,2,3 \}$.
\end{lem}

\proof The product of any two reflections is an element of $\PSp(2n,\R)$: for two reflections $\rho(s_{i})=R_{i}, \rho(s_{j})=R_{j}$ it holds
$$(R_{i}R_{j})^{T}J(R_{i}R_{j})=R_{j}^{T}R_{i}^{T}JR_{i}R_{j}=R_{j}^{T}(-J)R_{j}=J$$
Let $\rho: W_{3} \to \PSp(4,\R)^{\pm}$ be maximal, we want to show that $\rho(s_{i})\rho(s_{j})$ is Shilov hyperbolic. Without loss of generality let us assume $i=1$, $j=2$. By definition of maximality it is clear that $\rho(s_{1})\rho(s_{2})$ fixes $X_{1}$ and $X_{2}$, where $(P_{1},X_{1},P_{2},Q_{1},X_{2},Q_{2})$ is a maximal 6-tuple and $\Y_{P_{1},P_{2}} \perp \Y_{X_{1},X_{2}} \perp \Y_{Q_{1},Q_{2}}$. Up to isometry let us consider $(P_{1},P_{2},Q_{1},Q_{2},R_{1},R_{2})$ to be $(0,A,\Id,C,D,\infty)$ where $A,C,D$ are positive definite and $C$ is diagonal. The map $\rho(s_{1})\rho(s_{2})$ is inside $\PSp(4,\R)$ and fixes $0$ and $\infty$. This map is Shilov hyperbolic if and only if there exists a positive definite $Y$ such that $(0,Y,\rho(s_{1})\rho(s_{2})Y,\infty)$ is maximal (see Lemma \ref{checkshil}). Let $Y=A$. Then $\rho(s_{1})\rho(s_{2})(A)=\rho(s_{1})(A)$. We want to show that $(0,A,\rho(s_{1})A,\infty)$ is maximal. We know $(0,A,C)$ maximal and $\rho(s_{1}) \in \mathcal{R}(-D,0,C,D)$ (Remark \ref{possforrho}). Result follows by Proposition \ref{propreflinsiegel} and \ref{maxlemma2}.
\endproof

\begin{lem} \label{maxlemma1}
Let $ \rho: W_{3} \to \PSp^{\pm}(4,\R)$ be maximal and let $(X_{1},X_{2},Y_{1},Y_{2},Z_{1},Z_{2})$ be a maximal 6-tuple as in Definition \ref{defmaxW3repr}.
Then for $l_{1},l_{2},l_{3},l_{4} \in \mathcal{L}(\R^{4})$ it holds:
\begin{itemize}
\item[(i)] If $(X_{2},l_{1},l_{2},l_{3},l_{4},Z_{1})$ is maximal then $ (Z_{2},\rho(s_{1})l_{4},\rho(s_{1})l_{3},\rho(s_{1})l_{2},\rho(s_{1})l_{1},X_{1})$ is maximal
\item[(ii)] If $(Y_{2},l_{1},l_{2},l_{3},l_{4},X_{1})$  is maximal then $(X_{2},\rho(s_{2})l_{4},\rho(s_{2})l_{3},\rho(s_{2})l_{2},\rho(s_{2})l_{1},Y_{1})$ is maximal
\item[(iii)] If $(Z_{2},l_{1},l_{2},l_{3},l_{4},Y_{1})$ is maximal then $ (Y_{2},\rho(s_{3})l_{4},\rho(s_{3})l_{3},\rho(s_{3})l_{2},\rho(s_{3})l_{1},Z_{1})$ is maximal
\end{itemize}
\end{lem}

\proof Follows directly from Proposition \ref{propreflinsiegel}. 
\endproof

\subsection{Arc coordinates in classical Teichmüller} \label{arcoordinclassicteich}

Given a hyperbolic surface with boundary, arc coordinates provide a parametrization of the Teichmüller space. 

Given $\Sigma=\Sigma_{g,m}$ a compact orientable smooth surface of genus $g$ and $m$ boundary components, we can equip $\Sigma_{g,m}$ with a complete hyperbolic structure of finite volume with geodesic boundary. The universal covering $\widetilde{\Sigma}_{g,m}$ of $\Sigma_{g,m}$ is a closed subset of the hyperbolic plane $\H^{2}$ where boundary curves are geodesics. Let $\{ a_{1},...,a_{k} \}$ be a maximal collection of pairwise disjoint arcs in $\Sigma_{g,m}$ with starting and ending point on a boundary component which are essential and pairwise non-homotopic. The connected components of $\Sigma_{g,m} \setminus \bigcup_{i} a_{i}$ are given by a union of hexagons. Every arc will be called an \emph{edge} of the hexagon decomposition. For every hexagon there are exactly three alternating edges belonging to one boundary component of $\Sigma_{g,m}$. We denote by $E$ the set of all edges, $E_{bdry}$ the set of edges lying on a boundary component and by $\mathcal{H}$ the set of all hexagons of the decomposition. It can be shown that for such a collection $\{ a_{1},...,a_{k} \}$ it holds 
$$k= \# E \backslash E_{bdry}= 3 |\chi(\Sigma_{g,m})|= 3(2g-2+m)$$ 
and that the number of hexagons is given by $2 |\chi(\Sigma_{g,m})|= 2(2g-2+m)$. For a fixed hyperbolic structure we can always realize the hexagon decomposition of $\Sigma_{g,m}$ in a way such that every edge is a geodesic and every arc $a_{i} \in \{a_{1},...,a_{k} \}$ is the unique geodesic which is orthogonal to the boundary at both endpoints. We fix an orientation on the boundary components such that the surface lies to the right of the boundary.
For each choice of $\{ a_{1},...,a_{k} \}$ we get a parametrization of the Teichmüller space $\mathcal{T}(\Sigma_{g,m})$: once we fix the lengths $l(a_{1}),...,l(a_{k})$ there is a unique hyperbolic metric that makes $\Sigma_{g,m} \setminus \bigcup_{i} a_{i}$ a union of hyperbolic right-angled hexagons where each hexagon has exactly three alternating edges $a_{i_{1}},a_{i_{2}},a_{i_{3}}$ in $E \backslash E_{bdry}$ of length $l(a_{i_{1}}),l(a_{i_{2}}),l(a_{i_{3}})$ respectively, where $i_{1},i_{2},i_{3} \in \{1,...,k\}$. This is due to the well known fact that given three real numbers $b,c,d>0$ there exists (up to isometries) a unique right-angled hexagon in $\H^{2}$ with  alternating sides of lengths $b,c$ and $d$. 

Let us denote by $\Gamma_{g,m}$ the fundamental group $\pi_{1}(\Sigma_{g,m})$. It is well known that $\Gamma_{g,m}$ is isomorphic to the free group $\mathbb{F}_{2g+m-1}$ and that one can define the Teichmüller space $\mathcal{T}(\Sigma_{g,m})$ as the set of conjugacy classes of discrete and faithful representations $\rho: \Gamma_{g,m} \to  \PSL(2,\R)$.

In Definition \ref{defmaxshil} we have defined the space $\chi^{\text{max,Shilov}} (\Gamma_{g,m},\PSp(2n,\R))$. When $n=1$ the group $\PSp(2,\R)$ coincides with $\PSL(2,\R)$ and a Shilov hyperbolic element in $\PSL(2,\R)$ is conjugated to a matrix of the type $\bpm \lambda&0\\0&\lambda^{-1}\epm$ where $|\lambda|>1$. A representation $\rho$ in $ \Hom(\Gamma_{g,m},\PSL(2,\R))$ is discrete and faithful if and only if $\rho$ is maximal. The surface $\Sigma_{g,m}$ is then realized by the quotient $\Sigma_{g,m}= _{\rho(\Gamma_{g,m})} \backslash \H^{2}$ where  $\rho(\Gamma_{g,m})$ acts freely and properly discontinuously on $\H^{2}$. Restricting to Shilov hyperbolic representations is equivalent to require that the surface $\Sigma_{g,m}$ does not have cusps as boundary components. The above discussion asserts that once we fix the lengths $l(\alpha_{1}),...,l(\alpha_{k})$ we can explicitly write the representation $\rho:\Gamma_{g,m}\to \PSL(2,\R)$ in $\Hom^{\text{max},\text{Shilov}}(\Gamma_{g,m},\PSL(2,\R))$ such that $\Sigma_{g,m}=\rho(\Gamma_{g,m})\backslash \H^{2}$.

\begin{example} \label{arccoordbasic1}
Let $\Gamma_{0,3}$ be the fundamental group $\pi_{1}(\Sigma_{0,3})$, which is isomorphic to $\mathbb{F}_{2}$. We denote $\Gamma_{0,3}=\langle \alpha, \beta \rangle$ and consider $a_{1},a_{2},a_{3}$ three arcs as in Figure \ref{intropic1} which decompose $\Sigma_{0,3}$ in two hexagons. Once we fix the lengths $l(a_{1}),l(a_{2}),l(a_{3})$ we can uniquely draw two adjacent isometric hexagons in $\H^{2}$ up to isometry and we can reconstruct the generators  $\rho(\alpha),\rho(\beta)$ of the maximal representation which "closes up" the pair of pants. These are two hyperbolic isometries inside $\PSL(2,\R)$. 
\end{example}

\subsection{The group $\Gamma_{0,3}$ as a subgroup of $W_{3}$} \label{teichspacepop}
The fundamental group $\Gamma_{0,3}$ is isomorphic to the free group $\mathbb{F}_{2}$. On the other hand, the group $W_{3}= \langle s_{1},s_{2},s_{3}| \ s_{1}^{2}=s_{2}^{2}=s_{3}^{2}=1 \rangle$
has a normal subgroup $\Gamma= \langle s_{1}s_{2}, s_{2}s_{3} \rangle$ isomorphic to the free group $\mathbb{F}_{2}$. 
This allows us to see $\Gamma_{0,3}$ as a subgroup of $W_{3}$.

\begin{prop} \label{teichprop}
Let $\Gamma_{0,3}=\pi_{1}(\Sigma_{0,3})=\langle \alpha,\beta \rangle$ and $W_{3}= \langle s_{1},s_{2},s_{3}| \ s_{1}^{2}=s_{2}^{2}=s_{3}^{2}=1 \rangle$. Fix $\widetilde{\rho} \in \Hom^{\text{max}} (W_{3},\PSL^{\pm}(2,\R))$. Denote by  $\phi$ the following homomorphism
$$\begin{aligned}
\phi: \Gamma_{0,3}& \to W_{3} \\
 \alpha& \mapsto s_{1}s_{2}\\
  \beta& \mapsto s_{2}s_{3}
\end{aligned}$$  
It holds:
\begin{itemize}
    \item[(i)]The representation $\rho:=\widetilde{\rho}|_{Im(\phi)}$ is inside $\Hom^{\text{max,Shilov}} (\Gamma_{0,3},\PSL(2,\R))$.

\item [(ii)] For any $\rho \in \Hom^{\text{max,Shilov}} (\Gamma_{0,3},\PSL(2,\R))$ there exists a unique $\widetilde{\rho}\in \Hom^{\text{max}} (W_{3},\PSL^{\pm}(2,\R))$ such that $\rho=\widetilde{\rho}\circ \phi$
\item [(iii)]
The following map $f$ is a homeomorphism:
    $$\begin{aligned}
f: \chi^{\text{max}}(W_{3},\PSL^{\pm}(2,\R))& \to \chi^{\text{max,Shilov}}(\Gamma_{0,3},\PSL(2,\R)) \\
 \bigl[\widetilde{\rho}\bigr]& \mapsto \bigl[\widetilde{\rho}|_{Im(\phi)}\bigr]
\end{aligned}
$$
\end{itemize}
\end{prop}

\proof 
\begin{itemize}
    \item[\emph{(i)}] This will be proven for $\PSp(4,\R)$ in Proposition \ref{maximalityproof}. The proof for $\PSL(2,\R)$ is similar.
    \item[\emph{(ii)}] Let $\rho \in \Hom^{\text{max,Shilov}} (\Gamma_{0,3},\PSL(2,\R))$. Denote by $\{(x_{1},x_{2}), (y_{1},y_{2}), (z_{1},z_{2})\}  \subset \partial \H^{2}$ the fixed points of $\rho(\alpha), \rho(\beta)$ and $\rho(\beta^{-1}\alpha^{-1})$ respectively. Choose an orientation of the boundary $\partial \H^{2}$ such that $(x_{1},x_{2},y_{1},y_{2},z_{1},z_{2})$ is positive. In Section \ref{secreflH2} we have defined 
a reflection in $\H^{2}$ as an involution of $\SL^{-}(2,\R)$. Reflections in $\H^{2}$ fix an infinite geodesic $\gamma$ and are uniquely determined by the endpoints of $\gamma$ at the boundary of $\H^{2}$. For $p,q \in \partial \H^{2}$ let $\gamma_{p,q}$ denote the infinite geodesic having $p,q$ as endpoints.

\begin{notation}
For $p,q \in \partial \H^{2}$ we denote $R_{p,q}$ the unique non trivial reflection fixing the infinite geodesic $\gamma_{p,q}$. In other words $R_{p,q}$ is the unique non trivial isometry such that
$\gamma_{p,q} \perp \gamma_{x,R_{p,q}(x)}$
for any $x \in \partial \H^{2}$.
\end{notation}

\noindent Let $(p_{1},p_{2},q_{1},q_{2},r_{1},r_{2})$ be the positive 6-tuple inside $\partial \H^{2}$ uniquely determined by 
$$\gamma_{x_{1},x_{2}} \perp \gamma_{p_{1},p_{2}} \perp \gamma_{y_{1},y_{2}} \perp \gamma_{q_{1},q_{2}} \perp \gamma_{z_{1},z_{2}} \perp \gamma_{r_{1},r_{2}} \perp \gamma_{x_{1},x_{2}}$$ 
Define $\widetilde {\rho}: W_{3} \to \PSL^{\pm}(2,\R)$ such that (Figure \ref{Teichpop})
$$\widetilde {\rho}(s_{1})=R_{r_{1},r_{2}}, \ \ \widetilde {\rho}(s_{2})=R_{p_{1},p_{2}}, \ \ \widetilde {\rho}(s_{3})=R_{q_{1},q_{2}}$$

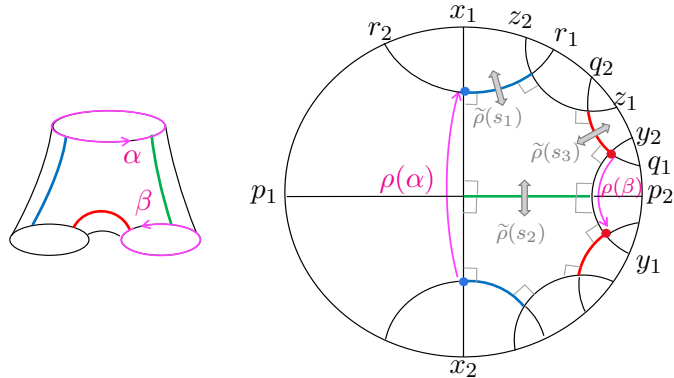
\begin{figure}[!h]
   \centering
   \captionsetup{justification=centering,margin=2cm}
   \setlength{\unitlength}{0.1\textwidth}
   \begin{picture}(5.5,2.9)
     \put(0,0){\includegraphics[width=8.5cm,height=4.5cm]{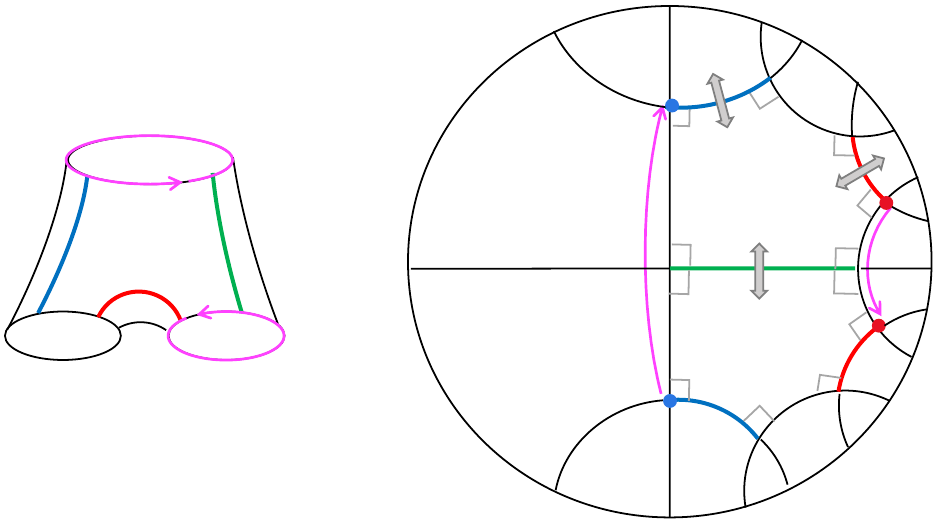}}
    \put(1,1.7){\textcolor{magenta}{$\alpha$}}
    \put(1.1,1.3){\textcolor{magenta}{$\beta$}}
    \put(3.2,1.5){\textcolor{magenta}{$\rho(\alpha$)}}
    \put(5.1,1.5){\textcolor{magenta}{$_{\rho(\beta)}$}}
    \put(3.8,-0.1){$x_{2}$}
    \put(5.4,0.8){$y_{1}$}
    \put(3.8,2.95){$x_{1}$}
    \put(2.1,1.4){$p_{1}$}
    \put(5.5,1.4){$p_{2}$}
    \put(5.4,1.9){$y_{2}$}
    \put(5.5,1.65){$q_{1}$}
    \put(5,2.5){$q_{2}$}
    \put(5.2,2.2){$z_{1}$}
    \put(4.7,2.75){$r_{1}$}
    \put(4.3,2.9){$z_{2}$}
    \put(3.1,2.8){$r_{2}$}
    \put(4,2.1){\textcolor{gray}{$_{\widetilde {\rho}(s_{1})}$}}
    \put(4.2,1.1){\textcolor{gray}{$_{\widetilde {\rho}(s_{2})}$}}
    \put(4.5,1.8){\textcolor{gray}{$_{\widetilde {\rho}(s_{3})}$}}
\end{picture}   
\caption{The maximal representation $\rho$ as a restriction of $\widetilde{\rho}$}
   \label{Teichpop}
\end{figure}

\noindent Then $\widetilde{\rho}$ is maximal. 
Moreover it is easy to show that (see for example \cite{martelli2016introduction} Proposition 6.2.1) $\rho(\alpha)= R_{q_{1},q_{2}} \circ R_{p_{1},p_{2}}$ and 
$\rho(\beta)= R_{r_{1},r_{2}}\circ R_{q_{1},q_{2}}$. It follows $$\widetilde{\rho}(s_{1}s_{2})=\widetilde{\rho}(s_{1})\widetilde{\rho}(s_{2})=\rho(\alpha)\text{ and } \widetilde{\rho}(s_{2}s_{3})=\widetilde{\rho}(s_{2})\widetilde{\rho}(s_{3})=\rho(\beta)$$  
so that $\widetilde{\rho}\circ \phi (\gamma)=\rho(\gamma)$ for all $\gamma \in \Gamma_{0,3}$. It is clear that $\widetilde{\rho}$ is the unique maximal representation such that $\widetilde{\rho}\circ \phi =\rho$.
\item[\emph{(iii)}] 
This follows directly from \emph{(ii)}. In particular recall that as $\Gamma_{0,3}$ is free, the set of representations $\Hom(\Gamma_{0,3},\PSL(2,\R))$ can be identified with $\PSL(2,\R)^{2}$ and we can carry over its topology.  
\end{itemize}
\endproof

\subsection{The set $\chi^{\mathcal{S}}$}

In Definition \ref{defmaxW3repr} we have defined the  set of maximal representations $\Hom^{\text{max}} (W_{3},\PSp^{\pm}(4,\R))$ and we know that we can see the fundamental group $\Gamma_{0,3}$ as a subgroup of $W_{3}$. In this section we define the set $\chi^{\mathcal{S}} \subset\chi^{\text{max,Shilov}} (\Gamma_{0,3},\PSp(4,\R))$. We start by giving an analogue of Proposition \ref{teichprop}\emph{(i)} in the case of $\PSp(4,\R)$, that is we show that the restriction to $\Gamma_{0,3}$ of a maximal representation as in Definition \ref{defmaxW3repr} is a maximal and Shilov hyperbolic representation as in Definition \ref{defmaxshil}. 

\begin{prop} \label{maximalityproof}
Let $\Gamma_{0,3}=\pi_{1}(\Sigma_{0,3})=\langle \alpha,\beta \rangle$ and $W_{3}= \langle s_{1},s_{2},s_{3}| \ s_{1}^{2}=s_{2}^{2}=s_{3}^{2}=1 \rangle$. Fix $\widetilde{\rho} \in \Hom^{\text{max}} (W_{3},\PSp^{\pm}(4,\R))$. Denote by  $\phi$ the following homomorphism
$$\begin{aligned}
\phi: \Gamma_{0,3}& \to W_{3} \\
 \alpha& \mapsto s_{1}s_{2}\\
  \beta& \mapsto s_{2}s_{3}
\end{aligned}$$ 
Then the representation $\rho:=\widetilde{\rho}|_{Im(\phi)}$ is inside $\Hom^{\text{max,Shilov}} (\Gamma_{0,3},\PSp(4,\R))$.

\end{prop}

\proof It is not hard to show that $\widetilde{\rho}\circ \phi (\gamma)$ is an element of $\PSp(4,\R)$ for any $\gamma \in \Gamma_{0,3}$. By abuse of notation let us denote the subgroup $\phi(\Gamma_{0,3})\trianglelefteq W_{3}$ just as $\Gamma_{0,3}$.
Given $\widetilde{\rho}$ in $\Hom^{\text{max}} (W_{3},\PSp^{\pm}(4,\R))$ we want to prove that $\rho=\widetilde{\rho}\circ \phi : \Gamma_{0,3} \to \PSp(4,\R)$ is maximal. As $\widetilde{\rho}:W_{3} \to \PSp^{\pm}(4,\R)$ is maximal we know that there exists a maximal 6-tuple $(P_{1},P_{2},Q_{1},Q_{2},R_{1},R_{2})$ satisfying the conditions of Definition \ref{defmaxW3repr}. In particular $(P_{1},P_{2})$, $(Q_{1},Q_{2})$ and $(R_{1},R_{2})$ are the fixpoints of $\widetilde{\rho}(s_{1}),\widetilde{\rho}(s_{2})$ and $\widetilde{\rho}(s_{3})$ respectively and  let us denote $(X_{1},X_{2},Y_{1},Y_{2},Z_{1},Z_{2})$ the Lagrangians satisfying (\ref{perpcondmaxreflW3}). We obtain the representation $\rho$ where 
\begin{equation*}
\begin{cases}
    \rho(\alpha)=\widetilde{\rho}\circ \phi(s_{1}s_{2}) \text{ fixes } X_{1},X_{2}\\
    \rho(\beta)=\widetilde{\rho}\circ \phi(s_{2}s_{3}) \text{ fixes } Y_{1},Y_{2}\\
    \rho(\beta^{-1}\alpha^{-1})=\widetilde{\rho}\circ \phi(s_{3}s_{1}) \text{ fixes } Z_{1},Z_{2}
\end{cases}
 \end{equation*}
To prove maximality of $\rho$ we use \cite[Corollary 3.4.2.]{strubel2015fenchel}, which gives a sufficient condition for a representation to be maximal by computing the Maslov index of elements in the Shilov boundary fixed by $\rho(C_{i})$, where $C_{i}$ are the peripheral elements of $\Gamma_{0,3}$. Finally, we know by Lemma \ref{Shilovproof} that $\widetilde{\rho}(s_{i},s_{j})$ is Shilov hyperbolic for any $i \neq j$. Any maximal representation sends non peripheral elements to Shilov hyperbolic elements in $\PSp(4,\R)$ (\cite{strubel2015fenchel}).
\\
\\
Additionally, we provide an explicit construction of the boundary map $\xi: S^{1} \to \mathcal{L}(\R^{4})$ in the spirit of \cite{fock2006moduli}. Fix $\rho_{0}$ a hyperbolization of $\Sigma_{0,3}$. 
Denote $\{(x_{1},x_{2}), (y_{1},y_{2}), (z_{1},z_{2})\}$ the fixed points of $\rho_{0}(\alpha), \rho_{0}(\beta)$ and $\rho_{0}(\beta^{-1}\alpha^{-1})$ respectively. Choose an orientation of the boundary $\partial \H^{2}$ such that $(x_{1},x_{2},y_{1},y_{2},z_{1},z_{2})$ is positive. By Proposition \ref{teichprop}\emph{(ii)} we know that there is a unique way to extend the action of $\Gamma_{0,3}$ on $\H^{2}$ (and on its boundary) to the group $W_{3}$ such that $\rho_{0}=\widetilde{\rho}_{0}\circ\phi$ where $\widetilde{\rho}_{0}:W_{3}\to \PSL^{\pm}(2,\R)$ is maximal. For simplicity given $s \in W_{3}$ and $p \in \partial \H^{2}$ we will denote the action $\widetilde{\rho}_{0}(s)\cdot p$ simply as $s \cdot p$.

\noindent As $\widetilde{\rho}:W_{3} \to \PSp^{\pm}(4,\R)$ is maximal there exists a maximal 6-tuple $(P_{1},P_{2},Q_{1},Q_{2},R_{1},R_{2})$ satisfying the conditions of Definition \ref{defmaxW3repr}. Let us denote $(X_{1},X_{2},Y_{1},Y_{2},Z_{1},Z_{2})$ the Lagrangians satisfying (\ref{perpcondmaxreflW3}). We define the following sets
$$H^{\H^{2}}:=\{x_{1},x_{2},y_{1},y_{2},z_{1},z_{2} \}$$
$$H^{\mathcal{L}}:=\{X_{1},X_{2},Y_{1},Y_{2},Z_{1},Z_{2} \}$$
$$\mathcal{O}^{\H^{2}}_{n}:= \bigcup_{|s| \leq n}s \cdot H^{\H^{2}}, \  s \in W_{3}$$
$$\mathcal{O}^{\mathcal{L}}_{n}:= \bigcup_{|s| \leq n} \widetilde{\rho}(s) \cdot H^{\mathcal{L}}$$
Define $\xi_{n}: \mathcal{O}^{\H^{2}}_{n} \to  \mathcal{O}^{\mathcal{L}}_{n}$ as
\begin{equation*}
    \begin{cases}
        \Big(\xi_{n}(x_{1}),\xi_{n}(x_{2}),\xi_{n}(y_{1}),\xi_{n}(y_{2}),\xi_{n}(z_{1}),\xi_{n}(z_{2})\Big)=(X_{1},X_{2},Y_{1},Y_{2},Z_{1},Z_{2})\\
        \xi_{n}(s \cdot p)= \widetilde{\rho}(s)\xi_{n}(p) \text{  for  } s \in W_{3}, |s| \leq n, \ p \in H^{\H^{2}}
    \end{cases}
\end{equation*}

\begin{figure}[!h]
   \centering
\captionsetup{justification=centering,margin=2cm}
   \setlength{\unitlength}{0.1\textwidth}
   \begin{picture}(5.7,6)
\put(0,0){\includegraphics[width=9.3cm,height=9cm]{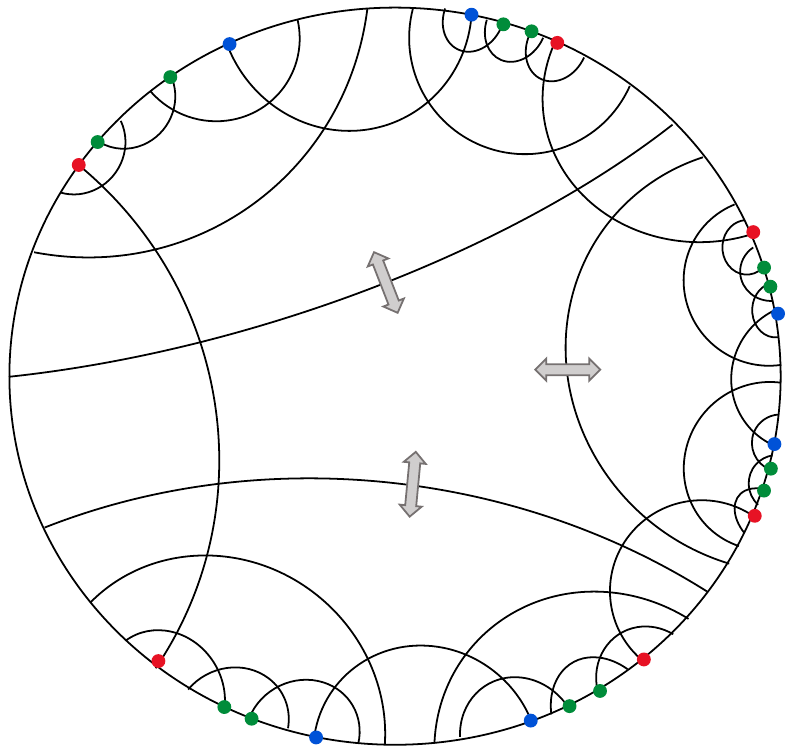}}
\put(2.3,-0.1){\textcolor{blue}{$s_{2}z_{2}$}}
\put(4,0){\textcolor{blue}{$s_{2}z_{1}$}}
\put(4.35,0.15){\textcolor{forestgreen}{$s_{2}s_{3}x_{1}$}}
\put(4.7,0.4){\textcolor{forestgreen}{$s_{2}s_{3}x_{2}$}}
\put(5,0.6){\textcolor{red}{$y_{1}$}}
\put(5.8,1.7){\textcolor{red}{$y_{2}$}}
\put(5.9,1.9){\textcolor{forestgreen}{$s_{3}s_{2}z_{1}$}}
\put(6,2.1){\textcolor{forestgreen}{$s_{3}s_{2}z_{2}$}}
\put(5.2,2.25){$_{s_{3}s_{2}H}$}
\put(6,2.4){\textcolor{blue}{$s_{3}x_{2}$}}
\put(4.8,2.8){$s_{3}H$}
\put(5.2,3.5){$_{s_{3}s_{1}H}$}
\put(6,3.3){\textcolor{blue}{$s_{3}x_{1}$}}
\put(6,3.6){\textcolor{forestgreen}{$s_{3}s_{1}y_{2}$}}
\put(5.9,3.8){\textcolor{forestgreen}{$s_{3}s_{1}y_{1}$}}
\put(5.9,4.1){\textcolor{red}{$z_{1}$}}
\put(4.4,5.45){\textcolor{red}{$z_{2}$}}
\put(3.4,5){$s_{1}s_{3}H$}
\put(4.1,5.65){\textcolor{forestgreen}{$s_{1}s_{3}x_{2}$}}
\put(3.8,5.85){\textcolor{forestgreen}{$s_{1}s_{3}x_{1}$}}
\put(3.1,5.8){\textcolor{blue}{$s_{1}y_{2}$}}
\put(1.3,5.6){\textcolor{blue}{$s_{1}y_{1}$}}
\put(0.6,5.2){\textcolor{forestgreen}{$s_{1}s_{2}z_{1}$}}
\put(0,4.8){\textcolor{forestgreen}{$s_{1}s_{2}z_{2}$}}
\put(1.2,4.5){$s_{1}s_{2}H$}
\put(0.3,4.5){\textcolor{red}{$x_{1}$}}
\put(0.9,0.6){\textcolor{red}{$x_{2}$}}
\put(1.6,0){\textcolor{forestgreen}{$s_{2}s_{1}y_{2}$}}
\put(1,0.3){\textcolor{forestgreen}{$s_{2}s_{1}y_{1}$}}
\put(2.7,1.3){$s_{2}H$}
\put(3.8,0.8){$s_{2}s_{3}H$}
\put(1.7,1){$s_{2}s_{1}H$}
\put(2.8,2.8){$H$}
\put(3.1,3.5){\textcolor{gray}{$s_{1}$}}
\put(2.8,2.2){\textcolor{gray}{$s_{2}$}}
\put(3.8,2.8){\textcolor{gray}{$s_{3}$}}
\put(2.5,4.3){$s_{1}H$}
   \end{picture}
\caption{Configuration of $\mathcal{O}_{2}^{\H^{2}}$}
   \label{maxproof}
\end{figure}
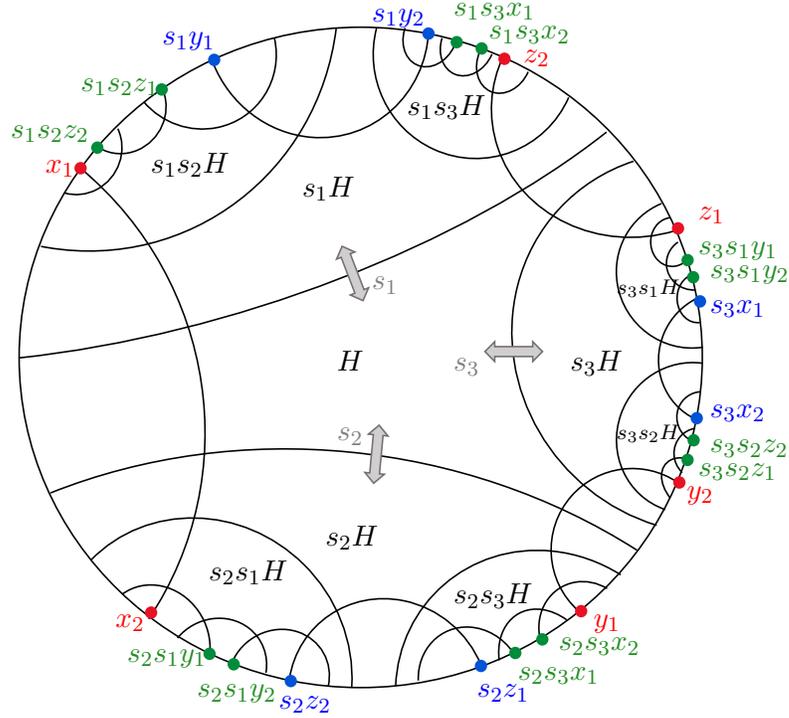

\noindent We show that the map $\xi_{n}$ is monotone by induction on $n$.\\
\underline{$\boldsymbol{n=0}$}: From the definition of $\xi$ it is clear that the map $\xi_{0}:H^{\H^{2}}\to H^{\mathcal{L}}$ is monotone.\\
\underline{$\boldsymbol{n=1}$}: We obtain the map $\xi_{1}:\mathcal{O}^{\H^{2}}_{1}\to \mathcal{O}^{\mathcal{L}}_{1}$ where $$\mathcal{O}_{1}^{\H^{2}}=H^{\H^{2}} \cup \{ s_{1}H^{\H^{2}},s_{2}H^{\H^{2}},s_{3}H^{\H^{2}} \}$$
$$\mathcal{O}_{1}^{\mathcal{L}}= H^{\mathcal{L}} \cup \{  \widetilde{\rho}(s_{1})H^{\mathcal{L}}, \widetilde{\rho}(s_{2})H^{\mathcal{L}}, \widetilde{\rho}(s_{3})H^{\mathcal{L}} \}$$
The set $\mathcal{O}_{1}^{\H^{2}}$ is given by $H^{\H^{2}}$ together with other six points, two for every $s_{i}H^{\H^{2}}, i\in \{1,2,3 \}$. For $s_{1}H^{\H^{2}}$ we only add the two points \{$s_{1}y_{1}, s_{1}y_{2}$\} as $s_{1}x_{1}=x_{2}$, $s_{1}x_{2}=x_{1}$ and $s_{1}z_{1}=z_{2}$, $ s_{1}z_{2}=z_{1}$. The same holds for $s_{2}H^{\H^{2}}$ and $s_{3}H^{\H^{2}}$. This is illustrated in Figure \ref{maxproof}.

\noindent The set $\mathcal{O}_{1}^{\H^{2}}$ is therefore formed by 12 points. The order on $\mathcal{O}_{1}^{\H^{2}}$ is given by the orientation of $\partial \H^{2}$. To explicitly write $\mathcal{O}_{1}^{\H^{2}}$ as a positive 12-tuple we use Proposition \ref{reflinH2} to show that the quadruples $(z_{2},s_{1}y_{2},s_{1}y_{1}x_{1})$, $ (x_{2},s_{2}z_{2},s_{2}z_{1},y_{1})$ and $(y_{2},s_{3}x_{2},s_{3}x_{1},z_{1})$ are positive. We obtain the following positive 12-tuple:
$$\mathcal{O}_{1}^{\H^{2}}=(z_{2},s_{1}y_{2},s_{1}y_{1},x_{1},x_{2},s_{2}z_{2},s_{2}z_{1},y_{1},y_{2},s_{3}x_{2},s_{3}x_{1},z_{1})$$
Similarly, the set $\mathcal{O}_{1}^{\mathcal{L}}$ consists of 12 Lagrangians: it is given by $H^{\mathcal{L}}$ together with six Lagrangians, two for every $ \widetilde{\rho}(s_{i})H^{\mathcal{L}}$. For $ \widetilde{\rho}(s_{1})H^{\mathcal{L}}$ we only add the Lagrangians $\{ \widetilde{\rho}(s_{1})Y_{1}, \widetilde{\rho}(s_{1})Y_{2}\}$. By maximality of $ \widetilde{\rho}$ we know (\ref{defmaxW3repr}) 
$$ \widetilde{\rho}(s_{1})X_{1}=X_{2}, \ \widetilde{\rho}(s_{1})X_{2}=X_{1}  \text{ and }  \widetilde{\rho}(s_{1})Z_{1}=Z_{2}, \  \widetilde{\rho}(s_{1})Z_{2}=Z_{1}$$
and the same is true for $ \widetilde{\rho}(s_{2})H^{\mathcal{L}}$ and $ \widetilde{\rho}(s_{3})H^{\mathcal{L}}$. To prove monotonic behaviour of $\xi_{1}$ we need to show maximality of the 12-tuple given by
$$\mathcal{O}_{1}^{\mathcal{L}}=(Z_{2}, \widetilde{\rho}(s_{1})Y_{2}, \widetilde{\rho}(s_{1})Y_{1},X_{1},X_{2}, \widetilde{\rho}(s_{2})Z_{2}, \widetilde{\rho}(s_{2})Z_{1},Y_{1},Y_{2}, \widetilde{\rho}(s_{3})X_{2}, \widetilde{\rho}(s_{3})X_{1},Z_{1})$$
We use Lemma \ref{maxlemma1} to show that the three quadruples given by $\big(Z_{2}, \widetilde{\rho}(s_{1})Y_{2}, \widetilde{\rho}(s_{1})Y_{1},X_{1}\big)$, $\big(X_{2}, \widetilde{\rho}(s_{2})Z_{2}, \widetilde{\rho}(s_{2})Z_{1},Y_{1}\big)$, and $\big(Y_{2}, \widetilde{\rho}(s_{3})X_{2}, \widetilde{\rho}(s_{3})X_{1},Z_{1}\big)$ are maximal . We use Lemma \ref{maxlemma2} to deduce that the 12-tuple is therefore maximal.\\
\underline{\textbf{Assume true for} $\boldsymbol{n}$ \textbf{show true for} $\boldsymbol{n+1}$}: Assuming $\xi_{n}$ monotone we consider the map $\xi_{n+1}: \mathcal{O}^{\H^{2}}_{n+1} \to  \mathcal{O}^{\mathcal{L}}_{n+1}$. We first study the set $\mathcal{O}^{\H^{2}}_{n+1}$. We describe how to obtain it from $\mathcal{O}^{\H^{2}}_{n}$ and how to write its positive order (\emph{Claim 1} and \emph{Claim 2}). The set $\mathcal{O}^{\H^{2}}_{n+1}$ is given by 
$$\mathcal{O}^{\H^{2}}_{n+1}=\mathcal{O}^{\H^{2}}_{n} \cup \big\{ s \cdot H,\ |s|=n+1 \big\}$$
If we fix an element $s \in W_{3}$ such that $|s|=n+1$ and look at the set $\{s \cdot H^{\H^{2}}\}$ we are adding exactly two points inside $\mathcal{O}^{\H^{2}}_{n+1}$ both lying between two points contained in $\mathcal{O}^{\H^{2}}_{n}$. This is made precise in the following two statements:\\

\noindent\underline{\emph{Claim 1}}: For any $s \in W_{3}$ such that $|s|=n+1$ it holds $|\mathcal{O}^{\H^{2}}_{n} \cup \{s \cdot H^{\H^{2}}\} |=|\mathcal{O}^{\H^{2}}_{n}|+2 $.\\

\noindent\underline{\emph{Claim 2}}: Let $s=ws_{i}$ where $|w|=n$ and $s_{i} \in W_{3}$. It holds
\begin{enumerate}
    \item If $s_{i}=s_{1}$ then the two points $ws_{1} \cdot H^{\H^{2}}$ added inside $\mathcal{O}^{\H^{2}}_{n+1}$ are $\{ws_{1}y_{1},ws_{1}y_{2} \}$ where
    \begin{equation*}
    \begin{cases}
(wz_{2},ws_{1}y_{2},ws_{1}y_{1},wx_{1}) \text{ positive if } n \text{ even}\\
(wx_{1},ws_{1}y_{1},ws_{1}y_{2},wz_{2}) \text{ positive if } n \text{ odd}
\end{cases}    
\end{equation*}
\item If $s_{i}=s_{2}$ then the two points $ws_{2} \cdot H^{\H^{2}}$ added inside $\mathcal{O}^{\H^{2}}_{n+1}$ are $\{ws_{2}z_{1},ws_{2}z_{2} \}$ where
\begin{equation*}
    \begin{cases}
(wx_{2},ws_{2}z_{2},ws_{2}z_{1},wy_{1}) \text{ positive if } n \text{ even}\\
(wy_{1},ws_{2}z_{1},ws_{2}z_{2},wx_{2}) \text{ positive if } n \text{ odd}
    \end{cases}
\end{equation*}
\item If $s_{i}=s_{3}$ then the two points $ws_{3} \cdot H^{\H^{2}}$ added inside $\mathcal{O}^{\H^{2}}_{n+1}$ are $\{ws_{3}x_{1},ws_{3}x_{2} \}$ where
\begin{equation*}
    \begin{cases}
(wy_{2},ws_{3}x_{2},ws_{3}x_{1},wz_{1}) \text{ positive if } n \text{ even}\\
(wz_{1},ws_{3}x_{1},ws_{3}x_{2},wy_{2}) \text{ positive if } n \text{ odd}
    \end{cases}
\end{equation*} 
\end{enumerate}

\noindent \underline{\emph{Proof of Claim 1}} Let $s \in W_{3}$ such that $|s|=n+1$ and consider the set $s \cdot H^{\H^{2}}$. Suppose $s$ ends with the element $s_{1}$ i.e. we can write
$s=ws_{1}$ for a $w\in W_{3}$, $|w|=n$. Among the six points $ws_{1}\cdot H^{\H^{2}}=\{ws_{1}x_{1},ws_{1}x_{2},ws_{1}y_{1},ws_{1}y_{2},ws_{1}z_{1},ws_{1}z_{2}\}$ we know  
$$ws_{1}x_{1}=wx_{2}, \ ws_{1}x_{2}=wx_{1}  \text{ and } ws_{1}z_{1}=wz_{2}, \ ws_{1}z_{2}=wz_{1}$$
so that $\{ws_{1}x_{1},ws_{1}x_{2},ws_{1}z_{1},ws_{1}z_{2}\} \subset \mathcal{O}^{\H^{2}}_{n}$. In particular $$\mathcal{O}^{\H^{2}}_{n} \cup \{ws_{1} \cdot H^{\H^{2}}\} = \mathcal{O}^{\H^{2}}_{n} \cup\{ws_{1}y_{1},ws_{2}y_{2}\}$$
A similar proof holds for $s=ws_{2}$ and $s=ws_{3}$. \hfill$\square$
\\

\noindent \underline{\emph{Proof of Claim 2}} Let us show (1) In the proof of \emph{Claim 1} we have already shown that the two points added inside $\mathcal{O}^{\H^{2}}_{n+1}$ are $\{ws_{1}y_{1},ws_{1}y_{2} \}$. We know that we can write $\mathcal{O}_{1}^{\H^{2}}$ as the positive 12-tuple (inductive step $n=1$)
$$\mathcal{O}_{1}^{\H^{2}}=(z_{2},s_{1}y_{2},s_{1}y_{1},x_{1},x_{2},s_{2}z_{2},s_{2}z_{1},y_{1},y_{2},s_{3}x_{2},s_{3}x_{1},z_{1})$$
In particular $(z_{2},s_{1}y_{2},s_{1}y_{1},x_{1})$ is positive. Let $w=w_{n}\cdot...\cdot w_{1}$, where $w_{i} \in \{s_{1},s_{2},s_{3}\}$.  At every step \begin{gather*}
(w_{1}z_{2},w_{1}s_{1}y_{2},w_{1}s_{1}y_{1},w_{1}x_{1})\to(w_{1}w_{2}z_{2},w_{1}w_{2}s_{1}y_{2},w_{1}w_{2}s_{1}y_{1},w_{1}w_{2}x_{1})\to...\\
...\to(w_{1}w_{2}...w_{n}z_{2},w_{1}w_{2}...w_{n}s_{1}y_{2},w_{1}w_{2}...w_{n}s_{1}y_{1},w_{1}w_{2}...w_{n}x_{1})
\end{gather*}
we satisfy the conditions of Proposition \ref{reflinH2}. It follows that the image under $s=ws_{1}$ of the positive quadruple $(z_{2},s_{1}y_{2},s_{1}y_{1},x_{1})$ stays positive if $n$ even and is negative if $n$ odd. Points (2) and (3) are similar. \hfill$\square$
\\

We state similar statements for the set $\mathcal{O}^{\mathcal{L}}_{n+1}$ . 

\noindent \underline{\emph{Claim 3}}: For any $s \in W_{3}$ such that $|s|=n+1$ it holds $|\mathcal{O}^{\mathcal{L}}_{n} \cup \{ \widetilde{\rho}(s) \cdot H^{\mathcal{L}}\} |=|\mathcal{O}^{\mathcal{L}}_{n}|+2$.\\

\noindent\underline{\emph{Claim 4}}: Let $s=ws_{i}$ where $|w|=n$ and $s_{i} \in W_{3}$. It holds
\begin{enumerate}
    \item If $s_{i}=s_{1}$ then the two Lagrangians $ \widetilde{\rho}(ws_{1}) \cdot H^{\mathcal{L}}$ added inside $\mathcal{O}^{\mathcal{L}}_{n+1}$ are $\{ \widetilde{\rho}(ws_{1})Y_{1}, \widetilde{\rho}(ws_{1})Y_{2} \}$ where
    \begin{equation*}
        \begin{cases}
         ( \widetilde{\rho}(w)Z_{2}, \widetilde{\rho}(ws_{1})Y_{2}, \widetilde{\rho}(ws_{1})Y_{1}, \widetilde{\rho}(w)X_{1}) \text{ maximal if } n \text{ even}\\ 
         ( \widetilde{\rho}(w)X_{1}, \widetilde{\rho}(ws_{1})Y_{1}, \widetilde{\rho}(ws_{1})Y_{2}, \widetilde{\rho}(w)Z_{2}) \text{ maximal if } n \text{ odd}
        \end{cases}
    \end{equation*}
\item If $s_{i}=s_{2}$ then the two Lagrangians $\rho(ws_{2}) \cdot H^{\mathcal{L}}$ added inside $\mathcal{O}^{\mathcal{L}}_{n+1}$ are $\{ \widetilde{\rho}(ws_{2})Z_{1}, \widetilde{\rho}(ws_{2})Z_{2} \}$ where
\begin{equation*}
    \begin{cases}
     ( \widetilde{\rho}(w)X_{2}, \widetilde{\rho}(ws_{2})Z_{2}, \widetilde{\rho}(ws_{2})Z_{1}, \widetilde{\rho}(w)Y_{1}) \text{ maximal if } n \text{ even}\\  
     ( \widetilde{\rho}(w)Y_{1}), \widetilde{\rho}(ws_{2})Z_{1}, \widetilde{\rho}(ws_{2})Z_{2}, \widetilde{\rho}(w)X_{2}) \text{ maximal if } n \text{ odd}
    \end{cases}
\end{equation*}
    \item If $s_{i}=s_{3}$ then the two Lagrangians $\rho(ws_{3}) \cdot H^{\mathcal{L}}$ added inside $\mathcal{O}^{\mathcal{L}}_{n+1}$ are $\{ \widetilde{\rho}(ws_{3})X_{1}, \widetilde{\rho}(ws_{3})X_{2} \}$ where
    \begin{equation*}
        \begin{cases}
        ( \widetilde{\rho}(w)Y_{2}, \widetilde{\rho}(ws_{3})X_{2}, \widetilde{\rho}(ws_{3})X_{1}, \widetilde{\rho}(w)Z_{1}) \text{ maximal if } n \text{ even}\\ 
        ( \widetilde{\rho}(w)Z_{1}, \widetilde{\rho}(ws_{3})X_{1}, \widetilde{\rho}(ws_{3})X_{2}, \widetilde{\rho}(w)Y_{2}) \text{ maximal if } n \text{ odd}
        \end{cases}
    \end{equation*}
\end{enumerate}

\noindent \underline{\emph{Proof of Claim 3}} The proof is similar to \emph{Claim 1} where we change $s$ with $ \widetilde{\rho}(s)$ and $x_{i},y_{i},z_{i}$ with $X_{i},Y_{i},Z_{i}$ and follows directly from the definition of $ \widetilde{\rho}$ (Definition \ref{defmaxW3repr}). \hfill$\square$
\\

\noindent \underline{\emph{Proof of Claim 4}} Let us show (1). By definition of $ \widetilde{\rho}$ it is clear that the two Lagrangians added inside $\mathcal{O}^{\mathcal{L}}_{n+1}$ are $\{  \widetilde{\rho}(ws_{1})Y_{1}, \widetilde{\rho}(ws_{1})Y_{2} \}$. We know that $(Z_{2},\rho(s_{1})Y_{2}, \widetilde{\rho}(s_{1})Y_{1},X_{1})$ is maximal (inductive step $n=1$). Let $w=w_{n}\cdot...\cdot w_{1}$, where $w_{i} \in \{s_{1},s_{2},s_{3}\}$.  At every step 
\begin{gather*}
 ( \widetilde{\rho}(w_{1})Z_{2}, \widetilde{\rho}(w_{1}s_{1})Y_{2}, \widetilde{\rho}(w_{1}s_{1})Y_{1}, \widetilde{\rho}(w_{1})X_{1})\to( \widetilde{\rho}(w_{1}w_{2})Z_{2}, \widetilde{\rho}(w_{1}w_{2}s_{1})Y_{2}, \widetilde{\rho}(w_{1}w_{2}s_{1})Y_{1}, \widetilde{\rho}(w_{1}w_{2})X_{1})\to...\\ 
 ...\to( \widetilde{\rho}(w_{1}w_{2}...w_{n})Z_{2}, \widetilde{\rho}(w_{1}w_{2}...w_{n}s_{1})Y_{2}, \widetilde{\rho}(w_{1}w_{2}...w_{n}s_{1})Y_{1}, \widetilde{\rho}(w_{1}w_{2}...w_{n})X_{1})
\end{gather*}
we satisfy the conditions of Proposition \ref{propreflinsiegel}. It follows that the image under $ \widetilde{\rho}(s)= \widetilde{\rho}(ws_{1})$ of the maximal quadruple $(Z_{2}, \widetilde{\rho}(s_{1})Y_{2}, \widetilde{\rho}(s_{1})Y_{1},X_{1})$ is maximal if $n$ is even and is minimal if $n$ is odd. (2) and (3) are similar. \hfill$\square$
\\

\noindent The map $\xi_{n}: \mathcal{O}^{\H^{2}}_{n+1} \to \mathcal{O}^{\mathcal{L}}_{n}$ is monotone by inductive hypothesis. In \emph{Claim 3} and \emph{Claim 4} we have proven that the set $\mathcal{O}^{\mathcal{L}}_{n+1}$ is obtained in the following way: for any $s$ of length $n+1$ we add two Lagrangians $l_{1},l_{2}$ in a way such that $(a,l_{1},l_{2},b)$ maximal for $a,b \in \mathcal{O}^{\mathcal{L}}_{n}$. By Lemma \ref{maxlemma2} it is easy to see that $\xi_{n+1}$ is monotone on the entire set $\mathcal{O}^{\H^{2}}_{n+1}$.\\

\noindent We have proven that $\xi_{n}: \mathcal{O}^{\H^{2}}_{n+1} \to \mathcal{O}^{\mathcal{L}}_{n}$ is monotone for any $n\geq 0$ and it is $ \widetilde{\rho}$-equivariant by definition. In particular when we restrict the map to $\Gamma_{0,3} \leq W_{3}$, that is we consider the map $\rho= \widetilde{\rho}\circ\phi$, we obtain a monotone $\rho$-equivariant map.
Using the same approach of \cite{burger2010surface} it can be shown that $\xi_{n}$ can be extended to a map $\xi$ defined on $S^{1}$ with the same properties. 
\endproof

\begin{definition} \label{defchiS} 
We define $\chi^{\mathcal{S}}:=Im(f)$ where $f$ is the map
$$\begin{aligned}
f: \chi^{\text{max}}(W_{3},\PSp^{\pm}(4,\R))& \to \chi^{\text{max,Shilov}}(\Gamma_{0,3},\PSp(4,\R)) \\
 \bigl[\widetilde{\rho}\bigr]& \mapsto \bigl[\widetilde{\rho}|_{Im(\phi)}\bigr]
\end{aligned}
$$
\end{definition}

\subsection{Parameter space for $\chi^{\text{max}}(W_{3},\PSp^{\pm}(4,\R))$}

Let $\X$ be the symmetric space associated to $\Sp(4,\R)$. Recall that we denote by $R_{st},R_{ex}$ and $\mathcal{K}$ the following matrices:
$$R_{st}=\bpm -\Id&0\\0&\Id\epm, \ R_{ex}=\bpm -r&0\\0&r\epm, \ r=\bpm -1&0\\0&1\epm$$

$$ \mathcal{K}=\Big\{\bpm -K&0\\0&K\epm, \  K \in \PO(2), \ K^{2}=\Id \Big\}$$
Observe that $\{R_{st},R_{ex}\} \subset \mathcal{K}$. Recall also that we denote by $F_{\underline{b},\underline{d},\alpha_{1},\alpha_{2}}$ the malefic map (\ref{defmaleficmap}). 

\begin{thm} \label{paramforchimax}
The set $\chi^{\text{max}}(W_{3},\PSp^{\pm}(4,\R))$  is parametrized by $\mathcal{S}\subset \mathcal{A}\times\mathcal{K}^{3}$ consisting of points $\big( \underline{b},\underline{c},\underline{d},[\alpha_{1},\alpha_{2} ],R_{1},R_{2},R_{3} \big)$ in $\mathcal{A}\times\mathcal{K}^{3}$ such that
\begin{equation*}
    \begin{cases}
        \underline{d} \in \mathfrak{a} \Rightarrow R_{1} \in \{R_{st},R_{ex}\}\\
        \underline{b} \in \mathfrak{a} \Rightarrow R_{2} \in \{R_{st},R_{ex}\}\\
        F_{\underline{b},\underline{d},\alpha_{1},\alpha_{2}}(\underline{c}) \in \mathfrak{a} \Rightarrow R_{3} \in \{R_{st},R_{ex}\}
    \end{cases}
\end{equation*}
\end{thm}

\begin{figure}[!h]
   \centering
\captionsetup{justification=centering,margin=2cm}
   \setlength{\unitlength}{0.1\textwidth}
   \begin{picture}(5.2,5.1)
\put(0,0){\includegraphics[width=8.3cm,height=7.7cm]{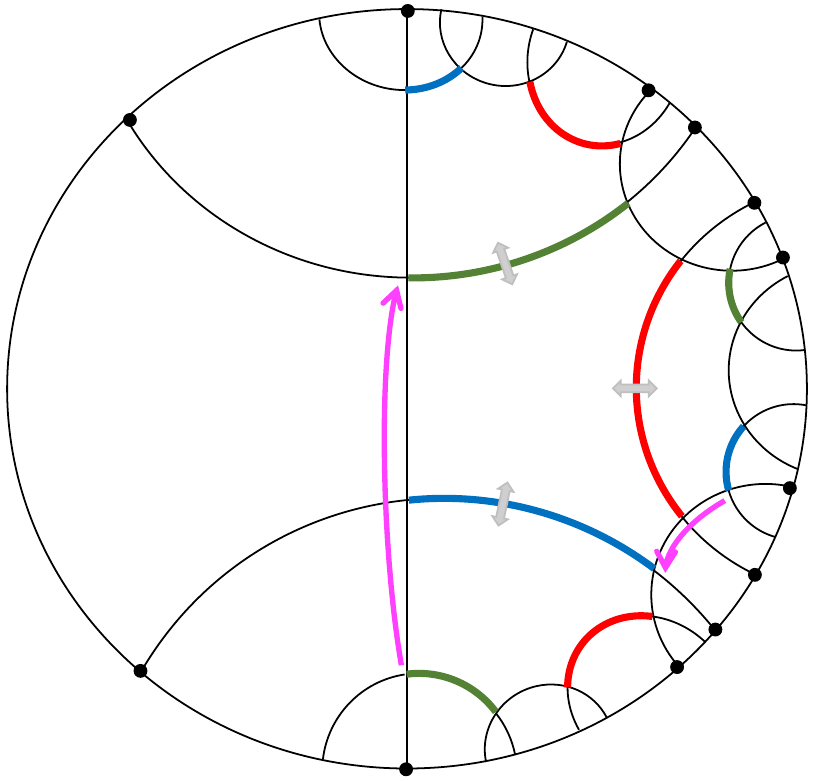}}
\put(2.6,-0.2){$0$}
\put(2.6,5){$\infty$}
\put(4.7,0.8){$A$}
\put(4.3,0.5){$A^{2}$}
\put(0.5,0.5){$-A$}
\put(5.15,3.35){$C$}
\put(5,3.7){$Z_{2}$}
\put(5.25,1.65){$\Id$}
\put(5,1.2){$Z_{1}$}
\put(4.6,4.2){$D$}
\put(4.24,4.45){$DC^{-1}D$}
\put(0.4,4.2){$-D$}
\put(3,1){$\overline{R}_{2}(H)$}
\put(3.2,2.4){$H$}
\put(2.3,2.5){\textcolor{magenta}{$g$}}
\put(4.5,1.3){\textcolor{magenta}{$h$}}
\put(4.25,2.65){$_{\overline{R}_{3}(H)}$}
\put(3,3.8){$\overline{R}_{1}(H)$}  
   \end{picture}
\caption{Configuration of hexagons $H$, $\overline{R}_{1}(H)$, $ \overline{R}_{2}(H)$, $\overline{R}_{3}(H)$}
   \label{maxrepr}
\end{figure}

\proof
\underline{\textbf{From parameters to representations}}: 
Let $\big( \underline{b},\underline{c},\underline{d},[\alpha_{1},\alpha_{2} ],R_{1},R_{2},R_{3} \big)$ be in $\mathcal{S}$. Let $(H,\Yi)=(0,A,\Id,C,D,\infty)$ be a right-angled hexagon with $\mathcal{A}(H,\Yi)=(\underline{b},\underline{c},\underline{d},[\alpha_{1},\alpha_{2}])$. In particular the maximal 12-tuple associated to $(H,\Yi)$ is given by
$$H=(\infty,-D,-A,0,A^{2},A,Z_{1},\Id,C,Z_{2},D,DC^{-1}D)$$
where $Z_{1},Z_{2}$ are uniquely determined by $\Y_{A^{2},\Id} \perp \Y_{Z_{1},Z_{2}} \perp \Y_{C,DC^{-1}D}$. Let $g_{1},g_{2},g_{3}$ be isometries such that
\begin{gather*}
 g_{1}(D,DC^{-1}D,\infty,-D)=(0,\Id,Y_{1},\infty)\\ 
 g_{2}(-A,0,A^{2},A)=(0,\Id,Y_{2},\infty)\\ 
 g_{3}(Z_{1},\Id,C,Z_{2})=(0,\Id,Y_{3},\infty) 
\end{gather*}
for $Y_{1},Y_{2},Y_{3}$ diagonal matrices inside $\Sym^{+}(2,\R)$. 
Put
$$\overline{R}_{1}:=g_{1}^{-1}R_{1}g_{1}$$
$$\overline{R}_{2}:=g_{2}^{-1}R_{2}g_{2}$$
$$\overline{R}_{3}:=g_{3}^{-1}R_{3}g_{3}$$
By Corollary \ref{reflsetwitharcoords} we know that $\overline{R}_{1},\overline{R}_{2},\overline{R}_{3}$ belong to $\mathcal{R}_{\Y_{-D,D}}^{\Yi,\Y_{C,DC^{-1}D}}$, $\mathcal{R}_{\Y_{-A,A}}^{\Yi,\Y_{A^{2},\Id}}$ and $\mathcal{R}_{\Y_{Z_{1},Z_{2}}}^{\Y_{A^{2},\Id},\Y_{C,DC^{-1}D}}$ respectively. Let $W_{3}= \langle s_{1},s_{2},s_{3}| \ s_{1}^{2}=s_{2}^{2}=s_{3}^{2}=1 \rangle$ and define $\widetilde{\rho}$ as the representation
\begin{equation*}
    \begin{aligned}
        \widetilde{\rho}:W_{3}&\to \PSp^{\pm}(4,\R)\\
        s_{1}&\mapsto \overline{R}_{1}\\
        s_{2}&\mapsto \overline{R}_{2}\\
        s_{3}&\mapsto \overline{R}_{3}
    \end{aligned}
\end{equation*}
\noindent The representation $\widetilde{\rho}$ is maximal by construction. The images $\overline{R}_{i}(H)$ for $i \in \{1,2,3\}$ are drawn in Figure \ref{maxrepr}. The maps $g$ and $h$ appearing in the Figure are the generators of the representation restricted to the group $\Gamma_{0,3}$. This will be explained in Theorem \ref{thesis}.\\

\noindent\underline{\textbf{From representations to parameters}}: Let $\widetilde{\rho} \in \chi^{\text{max}}(W_{3},\PSp^{\pm}(4,\R))$. We know that $\widetilde{\rho}$ has the properties described in Definition \ref{defmaxW3repr}: we can determine a right-angled hexagon $(H,\Y_{X_{1},X_{2}})$ where $X_{1},X_{2}$ are as in Figure \ref{maxrefl}. We compute the arc coordinates $\mathcal{A}(H,\Y_{X_{1},X_{2}})$.  By definition we know that $\widetilde{\rho}(s_{1}),\widetilde{\rho}(s_{2})$ and $\widetilde{\rho}(s_{3})$ belong to the reflection sets of three alternating sides of this hexagon. We compute the corresponding elements in $\mathcal{K}^{3}$ (one for every side) using Corollary \ref{reflsetwitharcoords}.
\endproof

\subsection{Parameter space for $\chi^{\mathcal{S}}$}

To parametrize $\chi^{\mathcal{S}}$ we impose an equivalent relation on $\mathcal{S}$ identifying the points that have same image under $f$. Recall the set $\mathcal{D} \subset \mathcal{A}$ corresponds to right-angled hexagons in $\X$ lying inside a maximal polydisk. 
\begin{definition}
We define $\mathcal{S}_{0} \subset \mathcal{A}\times\mathcal{K}^{3}$ as the set
$$\mathcal{S}_{0}:= \Big \{ (\underline{b},\underline{c},\underline{d},[\alpha_{1},\alpha_{2}],R_{1},R_{2},R_{3}) \in \mathcal{A}\times \mathcal{K}^{3} | \ (\underline{b},\underline{c},\underline{d},[\alpha_{1},\alpha_{2}]) \in \mathcal{D} \Big\} $$
We further define the following equivalent relation on $\mathcal{S}_{0}$:
\begin{equation} \label{eqrelS0}
(p,R_{1},R_{2},R_{3})\sim(p,R'_{1},R'_{2},R'_{3}) \iff 
\begin{cases}
    R_{1}R_{2}=R'_{1}R'_{2}\\
    R_{2}R_{3}=R'_{2}R'_{3}
\end{cases}
\end{equation}
\end{definition}

\begin{thm} \label{thesis}
The set $\chi^{\mathcal{S}}$ is parametrized by the parameter space $\mathcal{S}/_{\sim}$ where 
$\mathcal{S}$ consists of points $\big( \underline{b},\underline{c},\underline{d},[\alpha_{1},\alpha_{2} ],R_{1},R_{2},R_{3} \big)$ in $\mathcal{A}\times\mathcal{K}^{3}$ such that
\begin{equation*}
    \begin{cases}
        \underline{d} \in \mathfrak{a} \Rightarrow R_{1} \in \{R_{st},R_{ex}\}\\
        \underline{b} \in \mathfrak{a} \Rightarrow R_{2} \in \{R_{st},R_{ex}\}\\
        F_{\underline{b},\underline{d},\alpha_{1},\alpha_{2}}(\underline{c}) \in \mathfrak{a} \Rightarrow R_{3} \in \{R_{st},R_{ex}\}
    \end{cases}
\end{equation*}
and if $\big( p,R_{1},R_{2},R_{3} \big) \in \mathcal{S}_{0}$ we put
$(p,R_{1},R_{2},R_{3}) \sim (p',R'_{1},R'_{2},R'_{3})$ where $\sim$ is the equivalent relation in (\ref{eqrelS0}).
\end{thm}

\proof 
\underline{\textbf{From parameters to representations}}:
Let $\big( \underline{b},\underline{c},\underline{d},[\alpha_{1},\alpha_{2} ],R_{1},R_{2},R_{3} \big)$ be in $\mathcal{S}$. The construction of a maximal $\widetilde{\rho}:W_{3}\to \PSp^{\pm}(4,\R)$ is identical to the proof of Theorem \ref{paramforchimax}. 
Let $\rho$ be the restriction $\rho=f(\widetilde{\rho})$ where $f$ is the map of Definition \ref{defchiS}. Then $\rho$ is inside $\chi^{\mathcal{S}}$ and we put
$\big( \underline{b},\underline{c},\underline{d},[\alpha_{1},\alpha_{2} ],R_{1},R_{2},R_{3} \big) = \rho \in \chi^{\mathcal{S}}$. In Figure \ref{maxrepr} we have denoted
$g=\widetilde{\rho}(s_{1}s_{2})$, $h=\widetilde{\rho}(s_{2}s_{3})$.

\noindent In the github repository \url{https://github.com/martamagnani/Arc-coord/blob/main/Param_for_chiS.py} we provide a Python program with output the maps $g,h$ uniquely determining the maximal representation constructed above. We only provide the case $R_{i} \in \{R_{st},R_{ex}\}$. \\

\noindent \underline{\textbf{The equivalence relation on $\mathcal{S}_{0}$}}:
Let $\Gamma_{0,3}=\langle \alpha, \beta \rangle$. A representation $\rho: \Gamma_{0,3} \to \PSp(4,\R)$ is uniquely determined by the maps $\rho(\alpha)$, $ \rho(\beta)$, which are exactly the maps $g$ and $h$ of Figure \ref{maxrepr}. The isometry $g$ is sending the hexagon $\overline{R}_{2}(H)$ to $\overline{R}_{1}(H)$ and the isometry $h$ is sending the hexagon $\overline{R}_{3}(H)$ to $\overline{R}_{2}(H)$.
The equivalence relation on $\mathcal{S}_{0}$ identifies the points for which the map $f$ of Definition \ref{defchiS} is not injective. More precisely for two points 
$s=(p,R_{1},R_{2},R_{3})$ and $s'=(p,R'_{1},R'_{2},R'_{3})$ inside $S_{0}$ we denote $\overline{R}_{i},\overline{R'}_{i}$ the reflections constructed from the parameters $s$ and $s'$ respectively as shown in Figure \ref{maxrepr}. In Remark \ref{reflsamehex} we have detected the points for which $\overline{R}_{2}(H)=\overline{R'}_{2}(H)$, that is
$$\overline{R}_{2}=g_{2}^{-1}R_{st}g_{2}, \ \overline{R'}_{2}=g_{2}^{-1}R_{ex}g_{2} \ \ \text{ if } p \in \mathcal{D} \backslash \mathcal{D}_{\H^{2}}$$
where $H$ and $\overline{R}_{2}(H)=\overline{R'}_{2}(H)$ are both contained in the model polydisk, and
$$\overline{R}_{2},\overline{R'}_{2}\in g_{2}^{-1}\mathcal{K}g_{2} \ \ \text{ if } p \in \mathcal{D}_{\H^{2}}$$
where $H$ and $\overline{R}_{2}(H)=\overline{R'}_{2}(H)$ are both contained in the diagonal disc. It is not hard to show a similar result for $\overline{R}_{i},\overline{R'}_{i}$ when $i \in \{1,3\}$. 
Take two points $s,s'$ in $\mathcal{S}_{0}$ such that for the corresponding reflections $\overline{R}_{i},\overline{R'}_{i}$ constructed in proof of Theorem \ref{paramforchimax} it holds
\begin{equation*}
    \begin{cases}
       \overline{R}_{1}(H)=\overline{R'}_{1}(H)\\ 
       \overline{R}_{2}(H)=\overline{R'}_{2}(H)\\
       \overline{R}_{3}(H)=\overline{R'}_{3}(H)
    \end{cases}
\end{equation*}
The hexagons $H$, $\overline{R}_{i}(H)$, $\overline{R'}_{i}(H)$ of Figure \ref{maxrepr} are all contained in the model polydisk. All the points of the polygonal chains drawn in Figure \ref{combpolchains} are aligned.  Then there exists exactly two maps $g,\widehat{g}$ sending $\overline{R}_{2}(H)$ to $\overline{R}_{1}(H)$ and two maps $h,\widehat{h}$ sending $\overline{R}_{3}(H)$ to $\overline{R}_{2}(H)$. This follows directly from Proposition \ref{stabofhex}. We obtain four elements $\langle g,h \rangle, \langle g,\widehat{h} \rangle,\langle \widehat{g},h \rangle, \langle \widehat{g},\widehat{h} \rangle$ inside $\chi^{\mathcal{S}}$. But the parameter space $\mathcal{S}$ produces $|\{R_{st},R_{ex}\}|^{3}=2^{3}$ different maps. If we don't put the equivalence relation we would be over-counting the number of representations, that is we would construct $\widetilde{\rho},\widetilde{\rho}' \in \chi^{\text{max}}(W_{3},\PSp^{\pm}(4,\R))$ that have the same image under $f$.\\

\noindent \underline{\textbf{From representations to parameters}}: This is as Theorem \ref{paramforchimax}.

 \endproof

\begin{cor} \label{mapfnotinjnorsurj}
The map $f$ of Definition \ref{defchiS} is neither injective nor surjective. 
\end{cor}

\proof It is clear that $f$ is not injective.  We show $\chi^{\mathcal{S}}  \subsetneq \chi^{\text{max,Shilov}}(\Gamma_{0,3},\PSp(4,\R))$. The space $\chi^{\text{max,Shilov}}(\Gamma_{0,3},\PSp(4,\R))$ is 10-dimensional (see for example \cite{alessandrini2019noncommutative}). In the parametrization of Theorem \ref{thesis} we see that the set $\mathcal{S}$ is 8-dimensional. To see this observe that the space of right-angled hexagons is 8 dimensional. In the construction of Theorem \ref{thesis} if the hexagon is generic then the $\mathcal{K}$-component is 0-dimensional, and whenever one length-parameter lands in $\mathfrak{d}$ (the hexagon becomes non-generic) then the $\mathcal{K}$-component is one-dimensional (it is $\PO(2)$).
\endproof

\subsection{Connected components}
The parameter space of $\chi^{\mathcal{S}}$ distinguishes 8 connected components. The space $\chi^{\text{max,Shilov}}(\Gamma_{0,3},\PSp(4,\R))$ has $4$ connected components (see \cite{alessandrini2019noncommutative}). Let us see this by studying the geometrical behaviour of the generators of the representations. 

\noindent For a $g$ in $\PSp(4,\R)$ which fixes two transverse Lagrangians $l_{1},l_{2}$ there is a natural way to associate to $g$ a number  $sgn(g)\in \{\pm 1\}$ depending on $g$ being a reflecting or a non-reflecting isometry as in Definition \ref{isomref}. More precisely $sgn(g)=1$ if $\det A>0$ where $g \sim \bpm A&0\\0&A^{-T}\epm$ and $sgn(g)=-1$ if $\det A<0$. The number $sgn(g)$ tells us whether or not $g$ is reversing the orientation of the hyperbolic component of the tube $\Y_{l_{1},l_{2}}$. Given $\Gamma_{0,3}=\langle \alpha,\beta\rangle$, we consider a continuous surjective map $\delta$ which distinghuishes the connected components of $\chi^{\text{max,Shilov}}(\Gamma_{0,3},\PSp(4,\R))$:

\begin{equation}\label{delta}
    \begin{aligned}
\delta: \chi^{\text{max,Shilov}}(\Gamma_{0,3},\PSp(4,\R))& \to \{\pm 1\}\times \{\pm 1\} \\
 \rho& \mapsto \big(sgn(\rho(\alpha)),sgn(\rho(\beta))\big)
\end{aligned}
\end{equation}

\begin{lem}
    Let $\delta$ be the map in (\ref{delta}). Then for a point of parameters $(p,R_{1},R_{2},R_{3})$ inside $ \chi^{\mathcal{S}}$ it holds $\delta(p,R_{1},R_{2},R_{3})=\big(sgn(R_{1}R_{2}),sgn(R_{2}R_{3})\big)$. 
\end{lem}

\proof In the proof of Theorem \ref{thesis} we have explicitly constructed the representation $\rho \in \chi^{\mathcal{S}}$ from the parameters $(p,R_{1},R_{2},R_{3}) \in \mathcal{S}$. The generators of the representation are given by $\overline{R}_{1}\overline{R}_{2}$ and $\overline{R}_{2}\overline{R}_{3}$, where 
$\overline{R}_{1}:=g_{1}^{-1}R_{1}g_{1}$, 
$\overline{R}_{2}:=g_{2}^{-1}R_{2}g_{2}$, and $\overline{R}_{3}:=g_{3}^{-1}R_{3}g_{3}$ for a suitable choice of $g_{1},g_{2},g_{3}$ (see proof of Theorem \ref{paramforchimax}). 
We need to prove $sgn(\overline{R}_{i}\overline{R}_{j})=sgn(R_{i}R_{j})$. Let us concentrate on the case $i=1,j=2$. The geometrical behaviour of $\overline{R}_{1}\overline{R}_{2}$ is shown in Figure \ref{combpolchains} in the case where $p$ determines a generic hexagon. The property of reflecting the hyperbolic component (i.e. being a reflecting isometry as in Definition \ref{isomref}) or not only depends on $sgn(R_{1}R_{2})$. Sending a polygonal chain with internal angles $\alpha_{1},\alpha_{2}$ to one with same internal angles implies being a non-reflecting isometry. This can be generalized for any $i,j$.
\endproof

\begin{prop} \label{propconncomp}
    The set $\chi^{\mathcal{S}}$ hits all connected component inside $\chi^{\text{max,Shilov}}(\Gamma_{0,3},\PSp(4,\R))$.
\end{prop}

\proof This is clear. In particular each of the four connected components of $\chi^{\text{max,Shilov}}(\Gamma_{0,3},\PSp(4,\R))$ contains two connected components of $\chi^{s}$. 
\endproof

\printbibliography

@book{martelli2016introduction,
  title={An Introduction to Geometric Topology},
  author={Martelli, Bruno},
  publisher={Create Space Independent Publishing Platform},
  year={2016},
}

@article{fanoni2020basmajian,
  title={Basmajian-type inequalities for maximal representations},
  author={Fanoni, Federica and Pozzetti, Maria Beatrice},
  journal={Journal of Differential Geometry},
  volume={116},
  number={3},
  pages={405--458},
  year={2020},
  publisher={Lehigh University}
}

@article{burger2017maximal,
  title={Maximal representations, non-Archimedean Siegel spaces, and buildings},
  author={Burger, Marc and Pozzetti, Maria Beatrice},
  journal={Geometry \& Topology},
  volume={21},
  number={6},
  pages={3539--3599},
  year={2017},
  publisher={Mathematical Sciences Publishers}
}

@inproceedings{souriau2005construction,
  title={Construction explicite de l'indice de Maslov. Applications},
  author={Souriau, Jean-Marie},
  booktitle={Group Theoretical Methods in Physics: Fourth International Colloquium, Nijmegen 1975},
  pages={117--148},
  year={2005},
  organization={Springer}
}

@article{burger2010surface,
  title={Surface group representations with maximal Toledo invariant},
  author={Burger, Marc and Iozzi, Alessandra and Wienhard, Anna},
  journal={Annals of Mathematics},
  pages={517--566},
  year={2010},
  publisher={JSTOR}
}

@article{strubel2015fenchel,
  title={Fenchel--Nielsen coordinates for maximal representations},
  author={Strubel, Tobias},
  journal={Geometriae Dedicata},
  volume={176},
  pages={45--86},
  year={2015},
  publisher={Springer}
}

@article{toledo1989representations,
  title={Representations of surface groups in complex hyperbolic space},
  author={Toledo, Domingo},
  journal={Journal of Differential Geometry},
  volume={29},
  number={1},
  pages={125--133},
  year={1989},
  publisher={Lehigh University}
}

@article{parreau2010distance,
  title={La distance vectorielle dans les immeubles affines et les espaces sym{\'e}triques},
  author={Parreau, Anne},
  journal={preprint},
  year={2010}
}

@article{kapovich2017dynamics,
  title={Dynamics on flag manifolds: domains of proper discontinuity and cocompactness},
  author={Kapovich, Michael and Leeb, Bernhard and Porti, Joan},
  journal={Geometry \& Topology},
  volume={22},
  number={1},
  pages={157--234},
  year={2017},
  publisher={Mathematical Sciences Publishers}
}

@article{alessandrini2019noncommutative,
  title={Noncommutative coordinates for symplectic representations},
  author={Alessandrini, Daniele and Guichard, Olivier and Rogozinnikov, Eugen and Wienhard, Anna},
  journal={arXiv preprint arXiv:1911.08014},
  year={2019}
}

@book{goldman1980discontinuous,
  title={Discontinuous groups and the Euler class},
  author={Goldman, William M.},
  year={1980},
  publisher={University of California, Berkeley}
}

@article{fock2006moduli,
  title={Moduli spaces of local systems and higher Teichm{\"u}ller theory},
  author={Fock, Vladimir and Goncharov, Alexander},
  journal={Publications Math{\'e}matiques de l'IH{\'E}S},
  volume={103},
  pages={1--211},
  year={2006}
}

@article{siegel1943symplectic,
  title={Symplectic Geometry},
  author={Siegel, Carl Ludwig},
  journal={American Journal of Mathematics},
  volume={65},
  number={1},
  pages={1--86},
  year={1943},
  publisher={JSTOR}
}

@article{wolf1972fine,
  title={Fine structure of Hermitian symmetric spaces},
  author={Wolf, Joseph A},
  journal={Symmetric spaces (Short Courses, Washington Univ., St. Louis, Mo., 1969--1970)},
  volume={8},
  pages={271--357},
  year={1972},
  publisher={Dekker New York}
}

@Book{wienhard2004bounded,
 Author = {Wienhard, Anna},
 Title = {Bounded cohomology and geometry},
 FSeries = {Bonner Mathematische Schriften},
 Series = {Bonn. Math. Schr.},
 ISSN = {0524-045X},
 Volume = {368},
 Year = {2004},
 Publisher = {Bonn: Univ. Bonn, Mathematisches Institut (Dissertation)},
}

@article{lion1980weil,
  title={The Weil representation},
  author={Lion, Gerard and Vergne, Mich{\`e}le},
  journal={Maslov index and theta series," Birkkhauser, Boston},
  year={1980}
}

@Article{burger2005maximal,
 Author = {Burger, Marc and Iozzi, Alessandra and Labourie, Fran{\c{c}}ois and Wienhard, Anna},
 Title = {Maximal representations of surface groups: symplectic {Anosov} structures},
 FJournal = {Pure and Applied Mathematics Quarterly},
 Journal = {Pure and Applied Mathematics Quarterly},
 Volume = {1},
 Number = {3},
 Pages = {543--590},
 Year = {2005},
}

@Article{guichardwienhard2012,
 Author = {Guichard, Olivier and Wienhard, Anna},
 Title = {Anosov representations: domains of discontinuity and applications.},
 FJournal = {Inventiones Mathematicae},
 Journal = {Inventiones Mathematicae},
 Volume = {190},
 Number = {2},
 Pages = {357--438},
 Year = {2012},
}

@article {Har86,
    AUTHOR = {Harer, John L.},
     TITLE = {The virtual cohomological dimension of the mapping class group
              of an orientable surface},
   JOURNAL = {Inventiones Mathematicae},
  FJOURNAL = {Inventiones Mathematicae},
    VOLUME = {84},
      YEAR = {1986},
    NUMBER = {1},
     PAGES = {157--176},
}

@article{penner1987decorated,
  title={The decorated Teichm{\"u}ller space of punctured surfaces},
  author={Penner, Robert},
  journal={Communications in Mathematical Physics},
  volume={108},
  number={1},
  pages={299--339},
  year={1987}
}

@article{penner2002,
author = {Penner, Robert},
year = {2002},
month = {11},
pages = {},
title = {Decorated Teichm\"uller Theory of Bordered Surfaces},
volume = {12},
journal = {Communications in Analysis and Geometry},
}

@article{ushijima1999canonical,
  title={A canonical cellular decomposition of the Teichm{\"u}ller space of compact surfaces with boundary},
  author={Ushijima, Akira},
  journal={Communications in Mathematical Physics},
  volume={201},
  pages={305--326},
  year={1999},
  publisher={Springer}
}

@article {Mondello09,
    AUTHOR = {Mondello, Gabriele},
     TITLE = {Triangulated {R}iemann surfaces with boundary and the
              {W}eil-{P}etersson {P}oisson structure},
   JOURNAL = {Journal of Differential Geometry},
  FJOURNAL = {Journal of Differential Geometry},
    VOLUME = {81},
      YEAR = {2009},
    NUMBER = {2},
     PAGES = {391--436},
}

@article{luo2007teichmuller,
  title={On Teichm{\"u}ller spaces of surfaces with boundary},
  author={Luo, Feng},
  journal={Duke Mathematical Journal},
  volume={139},
  number={3},
  pages={463},
  year={2007},
  publisher={Duke University Press}
}

@article{guo2009parameterizations,
  title={On parameterizations of Teichm{\"u}ller spaces of surfaces with boundary},
  author={Guo, Ren},
  journal={Journal of Differential Geometry},
  volume={82},
  number={3},
  pages={629--640},
  year={2009},
  publisher={Lehigh University}
}

@article {BIPP21,
    AUTHOR = {Burger, Marc and Iozzi, Alessandra and Parreau, Anne and
              Pozzetti, Maria Beatrice},
     TITLE = {The real spectrum compactification of character varieties:
              characterizations and applications},
   JOURNAL = {Comptes Rendus Math\'{e}matique. Acad\'{e}mie des Sciences.
              Paris},
  FJOURNAL = {Comptes Rendus Math\'{e}matique. Acad\'{e}mie des Sciences.
              Paris},
    VOLUME = {359},
      YEAR = {2021},
     PAGES = {439--463},
}

@article{burger2023real,
  title={The real spectrum compactification of character varieties},
  author={Burger, Marc and Iozzi, Alessandra and Parreau, Anne and Pozzetti, Maria Beatrice},
  journal={arXiv preprint arXiv:2311.01892},
  year={2023}
}

@article{burger2021weyl,
  title={Weyl chamber length compactification of the {$\rm
              PSL(2,\R) \times \rm
              PSL(2,\R)$} maximal character variety},
  author={Burger, Marc and Iozzi, Alessandra and Parreau, Anne and Pozzetti, Maria Beatrice},
  journal={arXiv preprint arXiv:2112.13624},
  year={2021}
}

@article {OuyTamb23,
    AUTHOR = {Ouyang, Charles and Tamburelli, Andrea},
     TITLE = {Length spectrum compactification of the {$\rm
              SO_0(2,3)$}-{H}itchin component},
   JOURNAL = {Advances in Mathematics},
  FJOURNAL = {Advances in Mathematics},
    VOLUME = {420},
      YEAR = {2023},
     PAGES = {Paper No. 108997, 37},
}

@article {OuyangTambGoth,
    AUTHOR = {Ouyang, Charles and Tamburelli, Andrea},
     TITLE = {Boundary of the {G}othen components},
   JOURNAL = {Topology and its Applications},
  FJOURNAL = {Topology and its Applications},
    VOLUME = {326},
      YEAR = {2023},
     PAGES = {Paper No. 108420, 11},
}

\end{document}